\documentclass[11pt,a4paper,reqno]{amsart}
\usepackage{a4wide}
\usepackage{amsaddr}
\usepackage{amsmath}
\usepackage{amsfonts, amssymb, amsthm}
\usepackage{paralist}
\usepackage[colorlinks=true]{hyperref}
\hypersetup{urlcolor=blue, citecolor=red}
\usepackage{graphics,graphicx,color}

\makeatletter
\renewcommand{\email}[2][]{%
	\ifx\emails\@empty\relax\else{\g@addto@macro\emails{,\space}}\fi%
	\@ifnotempty{#1}{\g@addto@macro\emails{\textrm{(#1)}\space}}%
	\g@addto@macro\emails{#2}%
}
\makeatother

\theoremstyle{definition}

\numberwithin{equation}{section}

\newcommand{\R}{{\mathbb R}}
\newcommand{\eps}{{\varepsilon}}

\begin{document}
\title[Nonlinear Choquard  equation with  combined nonlinearities]
{Asymptotic profiles for Choquard equations with combined attractive nonlinearities}
\author{Shiwang Ma}\email{shiwangm@nankai.edu.cn}
\address{School of Mathematical Sciences and LPMC, Nankai University\\ 
	Tianjin 300071, China}
\author{Vitaly Moroz}\email{v.moroz@swansea.ac.uk}
\address{Department of Mathematics, Swansea University\\ 
	Fabian Way, Swansea SA1 8EN, 	Wales, UK}




\date{}

\begin{abstract}
We study asymptotic behaviour of positive ground state solutions of the nonlinear  Choquard equation
$$
-\Delta u+\varepsilon u=(I_\alpha \ast |u|^{p})|u|^{p-2}u+ |u|^{q-2}u 
\quad {\rm in} \ \mathbb R^N,
$$
where $N\ge 3$ is an integer, $p\in [\frac{N+\alpha}{N}, \frac{N+\alpha}{N-2}]$, $q\in (2,\frac{2N}{N-2}]$, $I_\alpha$ is the Riesz potential  and  $\varepsilon>0$ is a parameter. We  show that as $\varepsilon\to 0$ (resp. $\varepsilon\to \infty$), after {\em a suitable rescaling} the ground state solutions of $(P_\varepsilon)$ converge in $H^1(\mathbb R^N)$ to a particular solution of some limit equations. We also  establish a sharp asymptotic characterisation of such a rescaling, and the exact asymptotic behaviours of $u_\varepsilon(0), \|\nabla u_\varepsilon\|_2^2, \|u_\varepsilon\|_2^2, 
\int_{\mathbb R^N}(I_\alpha\ast |u_\varepsilon|^p)|u_\varepsilon|^p$ and $\|u_\varepsilon\|_q^q$, which depend in a non-trivial way on the exponents $p, q$ and the space dimension $N$. We also discuss a connection of our results with an associated mass constrained problem with normalization constraint $\int_{\mathbb R^N}|u|^2=c^2$. As a consequence of the main results, we obtain the existence, multiplicity  and exact asymptotic behaviour of positive normalized solutions of such a problem as  $c\to 0$ and $c\to \infty$.
\end{abstract}

\maketitle

\newpage

\section*{1. Introduction}

We study standing--wave solutions of the nonlinear Schr\"odinger equation with a combined  nonlinearity
\[
i\psi_t=\Delta \psi+(I_\alpha\ast |\psi|^{p})|\psi|^{p-2}\psi+a|\psi|^{q-2}\psi\quad\text{in $\R^N\times\R$,}
\eqno(1.1)
\]
where $N\ge 3$ is an integer, $\psi: \mathbb R^N\times \mathbb R\to \mathbb C$,  $p\in [\frac{N+\alpha}{N}, \frac{N+\alpha}{N-2}]$, $q\in (2,2^*]$ with $2^*=\frac{2N}{N-2}$, and $I_\alpha$ is the Riesz potential defined for every $x\in \mathbb R^N\setminus \{0\}$ by 
$$
I_\alpha (x)=\frac{A_\alpha(N)}{|x|^{N-\alpha}}, \quad A_\alpha(N)=\frac{\Gamma(\frac{N-\alpha}{2})}{\Gamma(\frac{\alpha}{2})\pi^{N/2}2^\alpha},
$$
where $\Gamma$ denotes the Gamma function.

A theory of NLS with local combined power nonlinearities was developed by Tao, Visan and Zhang \cite{Tao} and attracted a lot attention during the past decade (cf. \cite{Akahori-2,Akahori-3, Cazenave-1, Coles, Jeanjean-5,   Li-2, Li-1, Li-3, Li-4, Liu-1, Moroz-1, Sun-1} and further references therein).  
Nonlocal equation (1.1) in the case $p=2$ and $\alpha=2$ was proposed in cosmology (1.1), under the name of the Gross--Pitaevskii--Poisson equation,  as a model to describes the dynamics of the Cold Dark Matter in the form of the Bose--Einstein Condensate \cite{Wang,Bohmer-Harko,Chavanis-11}.
The nonlocal convolution term in (1.1) represents the Newtonian gravitational attraction between bosonic particles. The local term takes into account the short--range self--interaction between bosons. The non-interacting case $a=0$ corresponds to the Schr\"odinger--Newton (Choquard) model of self--gravitating bosons \cite{Ruffini}, which is mathematically well--studied \cite{MvS-survey}. When $a=1$ the quantum self--interaction between bosons is focusing/attractive, while for $a=-1$ the self--interaction is defocusing/repulsive, see surveys  \cite{Chavanis-15,Paredes} for the astrophysical background. Mathematically, the repulsive case $a=-1$ was recently studied in \cite{Liu-1}, see also further references therein. In this work we are concerned with the attractive case $a=1$.

A standing--wave solutions of (1.1) with a frequency $\varepsilon >0$ is a finite energy solution in the form
$$
\psi(t,x)=e^{-i\varepsilon t}u(x).
$$
This ansatz yields   the equation for $u$ in the form
\[
-\Delta u+\varepsilon u=(I_\alpha \ast |u|^{p})|u|^{p-2}u+|u|^{q-2}u
\quad {\rm in} \ \mathbb R^N.
 \eqno(P_\varepsilon)
\]
A solution $u_\varepsilon\in H^1(\mathbb R^N)$ of $(P_\varepsilon)$ is a critical point of the Action functional defined by
$$
I(u)=\frac{1}{2}\int_{\mathbb R^N}|\nabla u|^2+\frac{\varepsilon}{2}\int_{\mathbb R^N}|u|^2-\frac{1}{2p}\int_{\mathbb R^N}(I_\alpha\ast |u|^p)|u|^p-\frac{1}{q}\int_{\mathbb R^N}|u|^q.
\eqno(1.2)
$$
The existence and properties of ground state $u_\varepsilon\in H^1(\mathbb R^N)$ to $(P_\varepsilon)$ have been studied in \cite{Li-2,Li-1}.
We are interested in the limit asymptotic profile of the ground states $u_\varepsilon$ of the  problem $(P_\varepsilon)$, and in the asymptotic behaviour of different norms of $u_\varepsilon$, as $\varepsilon\to 0$ and $\varepsilon\to\infty$.
Of particular importance is the $L^2$--mass of the ground state
\[
M(\varepsilon):=\|u_\varepsilon\|_2^2,
\eqno(1.3)
\]
which plays a key role in the analysis of stability of the corresponding standing--wave solution of the time--dependent NLS (1.1). More precisely, the importance of $M(\varepsilon)$ is for instance seen in the Grillakis-Shatah-Strauss theory \cite{Grillakis-1,Grillakis-2, Shatah-1, Weinstein-1} of stability for these solutions within the time-dependent Schr\"odinger equation. The latter says that the solution $u_\varepsilon$ is orbitally stable when $M'(\varepsilon) > 0$ and that it is unstable when $M'(\varepsilon) < 0$. Therefore the intervals where $M(\varepsilon)$ is increasing furnish stable solutions whereas those where $M(\varepsilon)$ is decreasing correspond to unstable solutions. The Grillakis-Shatah-Strauss theory relies on another conserved quantity, the energy, which is defined below and for which the variations of $M(\varepsilon)$ also play a crucial role.

Alternatively, one can search for solutions to $(P_\varepsilon)$ having prescribed mass, and in this case $\varepsilon\in \mathbb R$ is part of the unknown.  That is, for a fixed $c>0$,  search for $u\in H^1(\mathbb R^N)$ and $\lambda\in \mathbb R$ satisfying
\[
\left\{\begin{array}{rl}
&-\Delta u=\lambda u+(I_\alpha\ast |u|^p)|u|^{p-2}u+|u|^{q-2}u \  \ in  \  \mathbb R^N,\smallskip\\
&u\in H^1(\mathbb R^N),  \   \    \   \int_{\mathbb R^N}|u|^2=c^2.
\end{array}\right.
\eqno(1.4)
\]
The solution of (1.4) is usually denoted by a pair $(u,\lambda)$ and  called a normalized solution.
This approach seems to be particularly meaningful from the physical point of view \cite{Chavanis-15,Paredes}, and often offers a good insight of the dynamical properties of the standing--wave solutions for (1.1), such as stability or instability.

It is standard to verify that the Energy functional 
\[
E(u)=\frac{1}{2}\int_{\mathbb R^N}|\nabla u|^2-\frac{1}{2p}\int_{\mathbb R^N}(I_\alpha\ast |u|^p)|u|^p-\frac{1}{q}\int_{\mathbb R^N}|u|^q
\eqno(1.5)
\]
is of class $C^1$ on $H^1(\R^N)$ and that a critical point of $E$ restricted to the (mass) constraint 
\[
S(c)=\{ u\in H^1(\mathbb R^N): \ \|u\|_2^2=c^2\}
\eqno(1.6)
\]
gives a solution to $(P_\varepsilon)$ and (1.4).  Here   $\lambda:=-\varepsilon$ arises as a Lagrange multiplier. We refer the readers to \cite{Li-3,Li-4, Sun-1} and the references therein. While main results of this paper are obtained for the unconstrained problem $(P_\varepsilon)$, we also study normalised solutions of (1.4).

In what follows we denote 
$$M(0):=\lim_{\varepsilon\to 0}M(\varepsilon),\qquad M(\infty):=\lim_{\varepsilon\to \infty}M(\varepsilon).$$ The following propositions outline the basic properties of solutions of $(P_\eps)$ and are direct consequences of Theorems A, B,  C, Lemma 1.1 and the main results in this paper.
\smallskip

\noindent
{\bf Proposition 1.1.}  {\it  If $p=\frac{N+\alpha}{N}$ and $q\in (2, 2+\frac{4}{N})$,  then the problem  $(P_\varepsilon)$ admits a positive ground state $u_\varepsilon\in H^1(\mathbb R^N)$,  which is radially symmetric and radially nonincreasing. 
Furthermore, the following statements hold true:

$(I)$  As $\varepsilon\to 0$,  then 
$$
u_\varepsilon(0)\sim \left\{\begin{array}{rcl}
\varepsilon^{\frac{1}{q-2}}, \qquad  \qquad &{\rm if}& \    \  q\in (2,2+\frac{4\alpha}{N(2+\alpha)}],\\
  \varepsilon^{\frac{2N}{\alpha[4-N(q-2)]}}, \    \quad  &{\rm if}& \    \ q\in (2+\frac{4\alpha}{N(2+\alpha)},2+\frac{4}{N}),
 \end{array}\right. 
$$
$$
\|u_\varepsilon\|_2^2\sim \left\{\begin{array}{rcl} 
\varepsilon^{\frac{4-N(q-2)}{2(q-2)}}, \quad  \quad &{\rm if}& \    \ q\in (2,2+\frac{4\alpha}{N(2+\alpha)}],\\
  \varepsilon^{\frac{N}{\alpha}}, \   \quad  \quad  \qquad &{\rm if}& \    \  q\in (2+\frac{4\alpha}{N(2+\alpha)},2+\frac{4}{N}),
 \end{array}\right. 
$$
$$
 \|\nabla u_\varepsilon \|_2^2\sim \left\{\begin{array}{rcl}
\varepsilon^{\frac{2N-q(N-2)}{2(q-2)}},   \quad  \   \  &{\rm if}&  q\in (2,2+\frac{4\alpha}{N(2+\alpha)}],\\
 \varepsilon^{\frac{N[2N-q(N-2)]}{\alpha[4-N(q-2)]}},  \   \  &{\rm if}& q\in (2+\frac{4\alpha}{N(2+\alpha)},2+\frac{4}{N}),
 \end{array}\right. 
$$
$$
E(u_\varepsilon)=\left\{\begin{array}{rcl}
\varepsilon^{\frac{2N-q(N-2)}{2(q-2)}}\left[-\frac{4-N(q-2)}{4q}S_q^{\frac{q}{q-2}}+O(\varepsilon^{-\frac{N(2+\alpha)(q-2)-4\alpha}{2N(q-2)}})\right],    &{\rm if}&  q\in (2,2+\frac{4\alpha}{N(2+\alpha)}),\\
\varepsilon^{\frac{N+\alpha}{\alpha}}\left[-\frac{N}{2(N+\alpha)}S_1^{\frac{N+\alpha}{\alpha}}+O(\varepsilon^{\frac{N(2+\alpha)(q-2)-4\alpha}{\alpha[4-N(q-2)]}})\right],  \quad \quad &{\rm if}& q\in (2+\frac{4\alpha}{N(2+\alpha)},2+\frac{4}{N}).
 \end{array}\right. 
$$

$(II)$   As $\varepsilon\to \infty$, then 
$$
u_\varepsilon(0)\sim \left\{\begin{array}{rcl}
\varepsilon^{\frac{2N}{\alpha[4-N(q-2)]}}, \quad  &{\rm if}& \    \  q\in (2,2+\frac{4\alpha}{N(2+\alpha)}],\\
\varepsilon^{\frac{1}{q-2}}, \    \quad  \qquad &{\rm if}& \    \ q\in (2+\frac{4\alpha}{N(2+\alpha)},2+\frac{4}{N}),
 \end{array}\right. 
$$
$$
\|u_\varepsilon\|_2^2\sim \left\{\begin{array}{rcl} 
\varepsilon^{\frac{N}{\alpha}}, \qquad  \qquad &{\rm if}& \    \ q\in (2,2+\frac{4\alpha}{N(2+\alpha)}],\\
  \varepsilon^{\frac{4-N(q-2)}{2(q-2)}}, \  \  \quad   &{\rm if}& \    \  q\in (2+\frac{4\alpha}{N(2+\alpha)},2+\frac{4}{N}),
 \end{array}\right. 
$$
$$
 \|\nabla u_\varepsilon \|_2^2\sim \left\{\begin{array}{rcl}
\varepsilon^{\frac{N[2N-q(N-2)]}{\alpha[4-N(q-2)]}},   \   \   \  &{\rm if}&  q\in (2,2+\frac{4\alpha}{N(2+\alpha)}],\\
\varepsilon^{\frac{2N-q(N-2)}{2(q-2)}},  \quad  \   \  &{\rm if}& q\in (2+\frac{4\alpha}{N(2+\alpha)},2+\frac{4}{N}),
 \end{array}\right. 
$$
$$
E(u_\varepsilon)=\left\{\begin{array}{rcl}
\varepsilon^{\frac{N+\alpha}{\alpha}}\left[-\frac{N}{2(N+\alpha)}S_1^{\frac{N+\alpha}{\alpha}}+O(\varepsilon^{\frac{N(2+\alpha)(q-2)-4\alpha}{\alpha[4-N(q-2)]}})\right],   \qquad &{\rm if}&  q\in (2,2+\frac{4\alpha}{N(2+\alpha)}),\\
\varepsilon^{\frac{2N-q(N-2)}{2(q-2)}}\left[-\frac{4-N(q-2)}{4q}S_q^{\frac{q}{q-2}}+O(\varepsilon^{-\frac{N(2+\alpha)(q-2)-4\alpha}{2N(q-2)}})\right],  &{\rm if}& q\in (2+\frac{4\alpha}{N(2+\alpha)},2+\frac{4}{N}).
 \end{array}\right. 
$$

$(III)$  $M(0)= 0, M(+\infty)=+\infty$. Additionally, if $M(\varepsilon)$ is of class $C^1$ for small $\varepsilon>0$ and large $\varepsilon>0$, then there exist some small  $\varepsilon_0>0$ and some large $\varepsilon_\infty>0$ such that 
$$
M'(\varepsilon)>0, \qquad  for \ all \ \varepsilon \in (0,\varepsilon_0)\cup (\varepsilon_\infty,+\infty).
$$ }

\smallskip
\noindent
{\bf Proposition 1.2.} {\it  If $p=\frac{N+\alpha}{N-2}$, $q\in (2,2^*)$ for $N\ge 4$ and $q\in (4,6)$ for $N=3$, then  the problem  $(P_\varepsilon)$ admits a positive ground state $u_\varepsilon\in H^1(\mathbb R^N)$,  which is  radially symmetric and radially nonincreasing. Furthermore, the following statements hold true:

$(I)$  As $\varepsilon\to \infty$,  then 
$$
u_\varepsilon(0)\sim \left\{\begin{array}{rcl}
\varepsilon^{\frac{1}{q-2}}, \qquad  \qquad &{\rm if}& \    \ N\ge 5, \ q\in (2,2^*),\\
 \varepsilon^{\frac{1}{q-2}}(\ln\varepsilon)^{\frac{2}{q-2}}, \quad &{\rm if}& \   \ N=4, \ q\in (2,4), \\
  \varepsilon^{\frac{1}{2(q-4)}}, \    \quad  \qquad &{\rm if}& \    \  N=3, \  q\in (4,6),
 \end{array}\right. 
$$
$$
\|u_\varepsilon\|_2^2\sim \left\{\begin{array}{rcl} 
\varepsilon^{-\frac{4}{(N-2)(q-2)}}, \quad \   \ \quad &{\rm if}& \    \ N\ge 5, \ q\in (2,2^*),\\
\varepsilon^{-\frac{2}{q-2}}(\ln\varepsilon)^{-\frac{4-q}{q-2}}, \quad &{\rm if}& \   \ N=4, \ q\in (2,4),\\
 \varepsilon^{-\frac{q-2}{2(q-4)}}, \   \quad  \quad  \qquad &{\rm if}& \    \  N=3, \  q\in (4,6),
 \end{array}\right. 
$$
$$
\|\nabla u_\varepsilon\|_2^2= S_\alpha^{\frac{N+\alpha}{2+\alpha}}+\left\{\begin{array}{rcl}
O(\varepsilon^{-\frac{2N-q(N-2)}{(N-2)(q-2)}}),  \  \  \quad {\rm  if} \   \ N\ge 5, \ q\in (2, 2^*),\\
O((\varepsilon\ln\varepsilon)^{-\frac{4-q}{q-2}}),   \  \  \quad {\rm if} \   \ N=4,\ q\in (2,4),\\
O(\varepsilon^{-\frac{6-q}{2(q-4)}}),  \quad \quad \quad {\rm if}  \    \  N=3, \ q\in (4,6),
\end{array}\right.
$$
$$
E(u_\varepsilon)=\frac{2+\alpha}{2(N+\alpha)}S_\alpha^{\frac{N+\alpha}{2+\alpha}}- \left\{\begin{array}{rcl} \Theta(\varepsilon^{-\frac{2N-q(N-2)}{(N-2)(q-2)}}), \quad  \quad  {\rm if} \   \ N\ge 5, \ q\in (2,2^*), \\
\Theta((\varepsilon\ln\varepsilon)^{-\frac{4-q}{q-2}}),  \quad \quad  {\rm if} \   \ N=4,\ q\in (2,4),  \\
\Theta( \varepsilon^{-\frac{6-q}{2(q-4)}}), \quad  \quad  \quad  {\rm if}  \    \  N=3, \ q\in (4, 6).
 \end{array}\right.
$$

$(II)$   As $\varepsilon\to 0$, then 
$$
u_\varepsilon(0)\sim \varepsilon^{\frac{1}{q-2}}, \quad  \quad {\rm if} \   \   \left\{\begin{array}{rcl} 
 N\ge 4, \  q\in (2,2^*),\\
 N=3, \  \   \ q\in (4, 6),
   \end{array}\right. 
$$
$$
\| u_\varepsilon\|_2^2\sim \varepsilon^{\frac{4-N(q-2)}{2(q-2)}},  \quad {\rm if} \   \  \left\{\begin{array}{rcl} 
 N\ge 4, \  q\in (2,2^*),\\
 N=3, \  \   \ q\in (4, 6),
   \end{array}\right. 
$$
$$
\|\nabla u_\varepsilon\|_2^2\sim \varepsilon^{\frac{2N-q(N-2)}{2(q-2)}},  \quad {\rm if} \   \  \left\{\begin{array}{rcl} 
 N\ge 4, \  q\in (2,2^*),\\
 N=3, \  \   \ q\in (4, 6),
   \end{array}\right. 
$$
$$
E(u_\varepsilon)=\varepsilon^{\frac{2N-q(N-2)}{2(q-2)}}\left[\frac{N(q-2)-4}{4q}S_q^{\frac{q}{q-2}}+O(\varepsilon^{\frac{(2+\alpha)(2N-q(N-2))}{2(q-2)}})\right].
$$

$(III)$  If $N\ge 4$, $q\in (2,2+\frac{4}{N})$, then $M(0)=M(+\infty)=0$. Additionally,  if $M(\varepsilon)$ is of class $C^1$ for small $\varepsilon>0$ and large $\varepsilon>0$, then there exist some small  $\varepsilon_0>0$ and  large $\varepsilon_\infty>0$ such that 
$$
M'(\varepsilon)>0, \quad  for  \ \varepsilon \in (0,\varepsilon_0),\qquad 
M'(\varepsilon)<0, \quad for \  \varepsilon\in (\varepsilon_\infty,+\infty).
$$
If $N\ge 4$, $q\in (2+\frac{4}{N}, 2^*)$, or $N=3$, $q\in (4,6)$, then $M(0)=+\infty, M(+\infty)=0$, and  if $M(\varepsilon)$ is of class $C^1$ for small $\varepsilon>0$ and large $\varepsilon>0$, then there exist some small  $\varepsilon_0>0$ and  large $\varepsilon_\infty>0$ such that 
$$
M'(\varepsilon)<0, \qquad  for  \ \varepsilon \in (0,\varepsilon_0)\cup (\varepsilon_\infty,+\infty).
$$ }

\smallskip
\noindent
{\bf Proposition  1.3.} {\it  If $q=2^*$, $p\in (1+\frac{\alpha}{N-2}, \frac{N+\alpha}{N-2})$ for $N\ge 4$ and  $ p\in (2+\alpha, 3+\alpha)$ for $N=3$,  then the problem  $(P_\varepsilon)$ admits a positive ground state $u_\varepsilon\in H^1(\mathbb R^N)$,  which is radially symmetric and radially nonincreasing. Furthermore, the following statements hold true:

$(I)$  As $\varepsilon\to \infty$, then 
$$
u_\varepsilon(0)\sim \left\{\begin{array}{rcl}
\varepsilon^{\frac{N-2}{2(N-2)(p-1)-2\alpha}},  \   &{\rm if}& N\ge 5, \  p\in (1+\frac{\alpha}{N-2}, \frac{N+\alpha}{N-2}),\\
(\varepsilon\ln\varepsilon)^{\frac{1}{2p-2-\alpha}}, \quad   &{\rm if}& N=4, \   p\in (\max\{2, 1+\frac{\alpha}{2}\}, 2+\frac{\alpha}{2}), \\
  \varepsilon^{\frac{1}{4(p-2-\alpha)}}, \    \qquad  &{\rm if}&  N=3, \ p\in (2+\alpha, 3+\alpha),
   \end{array}\right. 
$$
$$
\|u_\varepsilon\|_2^2\sim \left\{\begin{array}{rcl} 
 \varepsilon^{-\frac{2}{(N-2)(p-1)-\alpha}}, \qquad  \quad   &{\rm if}&  N\ge 5, \  p\in (1+\frac{\alpha}{N-2}, \frac{N+\alpha}{N-2}),\\
\varepsilon^{-\frac{2}{2p-2-\alpha}}(\ln\varepsilon)^{-\frac{4+\alpha-2p}{2p-2-\alpha}},   &{\rm if}&  N=4, \ p\in (\max\{2, 1+\frac{\alpha}{2}\}, 2+\frac{\alpha}{2}),\\
  \varepsilon^{-\frac{p-1-\alpha}{2(p-2-\alpha)}}, \   \quad  \quad  \qquad    &{\rm if}&  N=3, \  \ p\in (2+\alpha, 3+\alpha),
   \end{array}\right. 
$$
$$
\|\nabla u_\varepsilon\|_2^2= S^{\frac{N}{2}}+\left\{\begin{array}{rcl} 
O(\varepsilon^{-\frac{N+\alpha-p(N-2)}{(N-2)(p-1)-\alpha}}),  \   \quad &{\rm  if}&  N\ge 5,\  p\in (1+\frac{\alpha}{N-2}, \frac{N+\alpha}{N-2}),\\
O((\varepsilon\ln\varepsilon)^{-\frac{4+\alpha-2p}{2p-2-\alpha}}), \  \  &{\rm if}&  N=4,\ p\in (\max\{2, 1+\frac{\alpha}{2}\}, 2+\frac{\alpha}{2}),\\
 O( \varepsilon^{-\frac{3+\alpha-p}{2(p-2-\alpha)}}), \     \  \quad \quad  &{\rm if}&   N=3,\ p\in (2+\alpha, 3+\alpha), \end{array}\right. 
$$
$$
E(u_\varepsilon)=\frac{1}{N}S^{\frac{N}{2}}-\left\{\begin{array}{rcl} \Theta(\varepsilon^{-\frac{N+\alpha-p(N-2)}{(N-2)(p-1)-\alpha}}), \     \  &{\rm if}&  N\ge 5,\  p\in (1+\frac{\alpha}{N-2}, \frac{N+\alpha}{N-2}),\\
\Theta((\varepsilon\ln\varepsilon)^{-\frac{4+\alpha-2p}{2p-2-\alpha}}),    \   \  &{\rm if}&  N=4,\ p\in (\max\{2, 1+\frac{\alpha}{2}\}, 2+\frac{\alpha}{2}),\\
 \Theta( \varepsilon^{-\frac{3+\alpha-p}{2(p-2-\alpha)}}), \  \quad \quad    &{\rm if}&  N=3, \ p\in (2+\alpha, 3+\alpha), 
 \end{array}\right. 
$$

(II) As $\varepsilon\to 0$, then 
$$
u_\varepsilon(0)\sim \varepsilon^{\frac{2+\alpha}{4(p-1)}}, \quad {\rm if} \   \  \left\{\begin{array}{rcl} 
 N\ge 4, \  p\in (1+\frac{\alpha}{N-2}, \frac{N+\alpha}{N-2}),\\
 N=3, \  \   \ p\in (2+\alpha, 3+\alpha),
   \end{array}\right. 
$$
$$
\|u_\varepsilon\|_2^2\sim    \varepsilon^{\frac{2+\alpha-N(p-1)}{2(p-1)}}, \quad {\rm if} \    \ \left\{\begin{array}{rcl} 
 N\ge 4, \  p\in (1+\frac{\alpha}{N-2}, \frac{N+\alpha}{N-2}),\\
 N=3, \  \    \  p\in (2+\alpha, 3+\alpha),
   \end{array}\right. 
$$
$$
\|\nabla u_\varepsilon\|_2^2\sim \varepsilon^{\frac{N+\alpha-p(N-2)}{2(p-1)}}, \quad {\rm if} \    \ \left\{\begin{array}{rcl} 
 N\ge 4, \  p\in (1+\frac{\alpha}{N-2}, \frac{N+\alpha}{N-2}),\\
 N=3, \  \    \  p\in (2+\alpha, 3+\alpha),
   \end{array}\right. 
 $$
$$
E(u_\varepsilon)=\varepsilon^{\frac{N+\alpha-p(N-2)}{2(p-1)}}\left[\frac{N(p-1)-2-\alpha}{4p}S_p^{\frac{p}{p-1}}+O(\varepsilon^{\frac{N+\alpha-p(N-2)}{(N-2)(p-1)}})\right].
$$

$(III)$  If $N\ge 4$, $\alpha<N-2$,  $p\in (1+\frac{\alpha}{N-2}, 1+\frac{2+\alpha}{N})$, then $M(0)=M(+\infty)=0$. Additionally, if $M(\varepsilon)$ is of class $C^1$ for small $\varepsilon>0$ and large $\varepsilon>0$, then there exist some small  $\varepsilon_0>0$ and  large $\varepsilon_\infty>0$ such that 
$$
M'(\varepsilon)>0, \quad  for \ \varepsilon \in (0,\varepsilon_0), \qquad 
M'(\varepsilon)<0, \quad for  \  \varepsilon\in (\varepsilon_\infty,+\infty).
$$
If $N\ge 4$,  $p\in (\max\{1+\frac{\alpha}{N-2},1+\frac{2+\alpha}{N}\}, \frac{N+\alpha}{N-2})$, or $N=3$, $p\in (2+\alpha, 3+\alpha)$, then $M(0)=+\infty, M(+\infty)=0$. Additionally, if $M(\varepsilon)$ is of class $C^1$ for small $\varepsilon>0$ and large $\varepsilon>0$, then there exist some small  $\varepsilon_0>0$ and  large $\varepsilon_\infty>0$ such that 
$$
M'(\varepsilon)<0, \qquad  for  \ \varepsilon \in (0,\varepsilon_0)\cup (\varepsilon_\infty,+\infty).
$$}

\noindent
\smallskip
{\bf Proposition 1.4.} {\it If $p\in (\frac{N+\alpha}{N},  \frac{N+\alpha}{N-2})$ and $q\in (2,2^*)$, then the problem  $(P_\varepsilon)$ admits a positive ground state $u_\varepsilon\in H^1(\mathbb R^N)$,  which is radially symmetric and radially nonincreasing. Furthermore, let $S_p$ and $S_q$ be the constants given in Theorems 2.4 and 2.5, respectively, then the following statements hold true:

(I) As $\varepsilon\to 0$, then
$$
u_\varepsilon(0)\sim \left\{\begin{array}{rcl} 
\varepsilon^{\frac{1}{q-2}}, \quad \quad  &{\rm if}&   q\le \frac{2(2p+\alpha)}{2+\alpha},\\
\varepsilon^{\frac{2+\alpha}{4(p-1)}},  \quad  &{\rm if}&   q>\frac{2(2p+\alpha)}{2+\alpha},
 \end{array}\right. 
$$
$$
\|u_\varepsilon\|_2^2\sim \left\{\begin{array}{rcl} 
\varepsilon^{\frac{4-N(q-2)}{2(q-2)}},
 \quad \quad  &{\rm if}&  q\le \frac{2(2p+\alpha)}{2+\alpha},\\
\varepsilon^{\frac{2+\alpha-N(p-1)}{2(p-1)}},  \quad   &{\rm if}&  q>\frac{2(2p+\alpha)}{2+\alpha},
 \end{array}\right. 
$$
and 
$$
\|\nabla u_\varepsilon\|_2^2\simeq
\frac{N(q-2)}{2q}S_q^{\frac{q}{q-2}}\varepsilon^{\frac{2N-q(N-2)}{2(q-2)}}, \quad {\rm if}  \  \  q\not=\frac{2(2p+\alpha)}{2+\alpha}.
$$
If $q\not=\frac{2(2p+\alpha)}{2+\alpha}$, then as $\varepsilon\to 0$
$$
M(\varepsilon)= \left\{\begin{array}{rcl} 
\varepsilon^{\frac{4-N(q-2)}{2(q-2)}}\left(\frac{2N-q(N-2)}{2q}S_q^{\frac{q}{q-2}}-\Theta(\varepsilon^{\frac{2(2p+\alpha)-q(2+\alpha)}{2(q-2)}})\right), \   \   \  \    &{\rm if}&  q<\frac{2(2p+\alpha)}{2+\alpha},\\
\varepsilon^{\frac{2+\alpha-N(p-1)}{2(p-1)}}\left(\frac{N+\alpha-p(N-2)}{2p}S_p^{\frac{p}{p-1}}-\Theta(\varepsilon^{\frac{q(2+\alpha)-2(2p+\alpha)}{4(p-1)}})\right),   &{\rm if}& q> \frac{2(2p+\alpha)}{2+\alpha},
\end{array}\right.
$$
$$
E(u_\varepsilon)=\left\{\begin{array}{rcl} 
\varepsilon^{\frac{2N-q(N-2)}{2(q-2)}}\left(\frac{N(q-2)-4}{4q}S_q^{\frac{q}{q-2}}+O(\varepsilon^{\frac{2(2p+\alpha)-q(2+\alpha)}{2(q-2)}})\right),\   \  \  \   \ &{\rm if}& q< \frac{2(2p+\alpha)}{2+\alpha},\\
\varepsilon^{\frac{N+\alpha-p(N-2)}{2(p-1)}}\left(\frac{N(p-1)-2-\alpha}{4p}S_p^{\frac{p}{p-1}}+O(\varepsilon^{\frac{q(2+\alpha)-2(2p+\alpha)}{4(p-1)}})\right),  &{\rm if}& q> \frac{2(2p+\alpha)}{2+\alpha}.
\end{array}\right.
$$

(II) As $\varepsilon\to \infty$, then
$$
u_\varepsilon(0)\sim \left\{\begin{array}{rcl} 
\varepsilon^{\frac{2+\alpha}{4(p-1)}},  \quad  &{\rm if}& q\le \frac{2(2p+\alpha)}{2+\alpha},\\
\varepsilon^{\frac{1}{q-2}}, \quad \quad  &{\rm if}&  q>\frac{2(2p+\alpha)}{2+\alpha},
 \end{array}\right. 
$$
$$
\|u_\varepsilon\|_2^2\sim \left\{\begin{array}{rcl} 
\varepsilon^{\frac{2+\alpha-N(p-1)}{2(p-1)}}, \quad  &{\rm if}&  q\le \frac{2(2p+\alpha)}{2+\alpha},\\
\varepsilon^{\frac{4-N(q-2)}{2(q-2)}},  \quad \quad  &{\rm if}&  q>\frac{2(2p+\alpha)}{2+\alpha},
 \end{array}\right. 
$$
and 
$$
\|\nabla u_\varepsilon\|_2^2\simeq 
\frac{N(p-1)-\alpha}{2p}S_p^{\frac{p}{p-1}}\varepsilon^{\frac{N+\alpha-p(N-2)}{2(p-1)}},   \quad {\rm if} \  \  q\not=\frac{2(2p+\alpha)}{2+\alpha}.
$$
If $q\not=\frac{2(2p+\alpha)}{2+\alpha}$, then as $\varepsilon\to \infty$
$$
M(\varepsilon)= \left\{\begin{array}{rcl} 
\varepsilon^{\frac{2+\alpha-N(p-1)}{2(p-1)}}\left(\frac{N+\alpha-p(N-2)}{2p}S_p^{\frac{p}{p-1}}+O(\varepsilon^{-\frac{2(2p+\alpha)-q(2+\alpha)}{4(p-1)}})\right),
   &{\rm if}& q<\frac{2(2p+\alpha)}{2+\alpha},\\
\varepsilon^{\frac{4-N(q-2)}{2(q-2)}}\left(\frac{2N-q(N-2)}{2q}S_q^{\frac{q}{q-2}}+O(\varepsilon^{-\frac{q(2+\alpha)-2(2p+\alpha)}{2(q-2)}})\right),  \   \  \   \  &{\rm if}&  q> \frac{2(2p+\alpha)}{2+\alpha},
\end{array}\right.
$$
$$
E(u_\varepsilon)=\left\{\begin{array}{rcl} 
\varepsilon^{\frac{N+\alpha-p(N-2)}{2(p-1)}}\left(\frac{N(p-1)-2-\alpha}{4p}S_p^{\frac{p}{p-1}}+O(\varepsilon^{-\frac{2(2p+\alpha)-q(2+\alpha)}{4(p-1)}})\right),  &{\rm if}& q< \frac{2(2p+\alpha)}{2+\alpha},\\
\varepsilon^{\frac{2N-q(N-2)}{2(q-2)}}\left(\frac{N(q-2)-4}{4q}S_q^{\frac{q}{q-2}}+O(\varepsilon^{-\frac{q(2+\alpha)-2(2p+\alpha)}{2(q-2)}})\right), \    \    \quad  &{\rm if}& q> \frac{2(2p+\alpha)}{2+\alpha}.
\end{array}\right.
$$

$(III)$   Let $p_0:=1+\frac{2+\alpha}{N}$ and $q_0:=2+\frac{4}{N}$, then
$$
M(0)= \left\{\begin{array}{rcl} 
0, \quad \quad \qquad \quad  &{\rm if}& \   \ q<q_0 \   \   {\rm or} \    \  p<p_0,\\
\frac{2}{N+2}S_{q_0}^{\frac{N+2}{2}},  \quad \quad  &{\rm if}& \   \ q=q_0  \  \  {\rm and} \   \ p>p_0,\\
\frac{2+\alpha}{N+2+\alpha}S_{p_0}^{\frac{N+2+\alpha}{2+\alpha}},   &{\rm if}& \   \ q>q_0  \  \  {\rm and} \   \ p=p_0,\\
\infty, \quad  \qquad  \  \quad \  &{\rm if}& \   \ q>q_0 \    \  {\rm and}  \   \ p>p_0,
 \end{array}
\right.
$$
and 
$$
M(\infty)= \left\{\begin{array}{rcl} 
0, \quad \quad \qquad  \   &{\rm if}& \   \ q>q_0 \   \   {\rm or} \    \  p>p_0,\\
\frac{2}{N+2}S_{q_0}^{\frac{N+2}{2}},  \quad \quad &{\rm if}& \   \ q=q_0  \  \  {\rm and} \   \ p<p_0,\\
\frac{2+\alpha}{N+2+\alpha}S_{p_0}^{\frac{N+2+\alpha}{2+\alpha}},  &{\rm if}& \   \ q<q_0  \  \  {\rm and} \   \ p=p_0,\\
\infty, \quad  \qquad  \   \   \    &{\rm if}& \   \ q<q_0 \    \  {\rm and}  \   \ p<p_0.
 \end{array}
\right.
$$
Moreover, if  $q\not=\frac{2(2p+\alpha)}{2+\alpha}$ and additionally,  $M(\varepsilon)$ is of class $C^1$ for small $\varepsilon>0$ and large $\varepsilon>0$, then there exists a small $\varepsilon_0>0$  such that
for any $\varepsilon\in (0,\varepsilon_0)$,
$$
 \left\{\begin{array}{rcl} 
M'(\varepsilon)>0, \quad   &{\rm if}& \   \ q<q_0 \   \   {\rm or} \    \  p<p_0,\\
M'(\varepsilon)<0,  \quad   &{\rm if}& \  \ q=q_0  \    \  {\rm and} \   \  p>p_0\\
M'(\varepsilon)<0,  \quad    &{\rm if}& \  \ q>q_0  \    \  {\rm and} \   \  p=p_0\\
M'(\varepsilon)<0, \quad   &{\rm if}& \   \ q>q_0 \    \  {\rm and}  \   \ p>p_0,
 \end{array}\right.
 \eqno(1.7)
$$
and there exists a large $\varepsilon_\infty>0$ such that for any $\varepsilon\in (\varepsilon_\infty, +\infty)$, 
$$
 \left\{\begin{array}{rcl} 
M'(\varepsilon)<0, \quad   &{\rm if}& \   \ q>q_0 \   \   {\rm or} \    \  p>p_0,\\
M'(\varepsilon)>0, \quad   &{\rm if}& \   \ q<q_0 \    \  {\rm and}  \   \ p<p_0.
 \end{array}
\right.
\eqno(1.8)
$$}

\smallskip
Figure 1 outlines  the limits of $M(\varepsilon)$ as $\varepsilon\to 0$ and $\varepsilon\to\infty$ and reveals the variation of $M(\lambda)$ for small $\varepsilon>0$ and large $\varepsilon>0$ when $(p,q)$ 
belongs to different regions in the $(p,q)$ plane, as described in Propositions 1.1-1.4.


\begin{figure}
	\centering
	\includegraphics[width=0.7\linewidth]{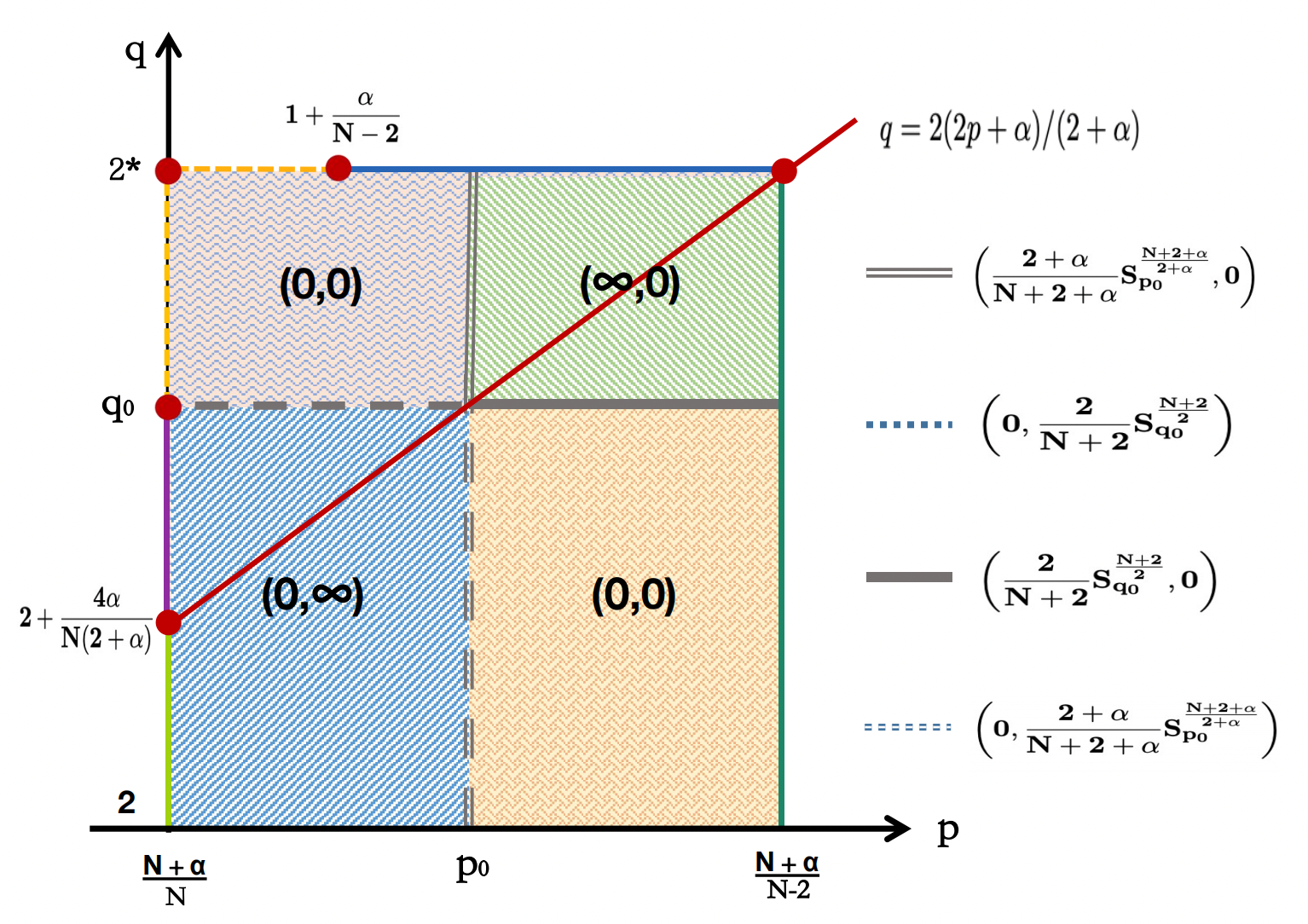}
	\caption{The variation of $M(\varepsilon)$ for small and large $\varepsilon$, here $(\cdot,\cdot)=(M(0), M(\infty))$.}
	\label{fig:figure1}
\end{figure}





\vskip 5mm
\noindent
{\bf Remark 1.1.}   Clearly, if $u_\varepsilon\in H^1(\mathbb R^N)$ is a ground state of $(P_\varepsilon)$, and for some $c>0$ there holds
\[
M(\varepsilon)=\|u_\varepsilon\|_2^2=c^2,
\eqno(1.9)
\] 
then $u_\varepsilon$ is a normalized solution of (1.4) with $\lambda=-\varepsilon$.  We  denote this normalized solution by  a pair $(u_c, \lambda_c)$ with $\lambda_c=-\varepsilon$, or $u_c$ for simplicity. 
As direct consequences of  Propositions 1.1-1.4, we have the following results.

\smallskip
\noindent
{\bf Corollary 1.1.}   {\it Let $p=\frac{N+\alpha}{N}$, $q\in (2,2+\frac{4}{N})$, then for any  $c>0$ the problem (1.4) has at least one positive  normalized solution $u_c\in H^1(\mathbb R^N)$, which is  radially symmetric and radially nonincreasing. 
Moreover,  as $c\to 0$, 
$$
 \|\nabla u_c \|_2^2\sim 
c^{\frac{2(2N-q(N-2))}{4-N(q-2)}}\to 0,
$$
$$
E(u_c)\sim \left\{\begin{array}{rcl}
-c^{\frac{2(2N-q(N-2))}{4-N(q-2)}}\to 0^-,       &{\rm if}&   q\in (2,2+\frac{4\alpha}{N(2+\alpha)}),\\
-c^{\frac{2(N+\alpha)}{N}}\to 0^-,  \qquad &{\rm if} & q\in (2+\frac{4\alpha}{N(2+\alpha)},2+\frac{4}{N}).
 \end{array}\right. 
$$
As $c\to \infty$, 
$$
 \|\nabla u_c \|_2^2\sim 
c^{\frac{2(2N-q(N-2))}{4-N(q-2)}}\to +\infty,
$$
$$
E(u_c)\sim \left\{\begin{array}{rcl}
-c^{\frac{2(N+\alpha)}{N}}\to -\infty,  \qquad      &{\rm if}&   q\in (2,2+\frac{4\alpha}{N(2+\alpha)}),\\
-c^{\frac{2(2N-q(N-2))}{4-N(q-2)}}\to -\infty,   &{\rm if} & q\in (2+\frac{4\alpha}{N(2+\alpha)},2+\frac{4}{N}).
 \end{array}\right. 
$$

}

\smallskip
\noindent
{\bf Corollary 1.2.}  {\it If $p=\frac{N+\alpha}{N-2}$, $q\in (2,2^*)$ for $N\ge 4$ and $q\in (4,6)$ for $N=3$, then  the following statements hold true:

(I)  If $q<q_0$, then there  exists  a constant $c_0>0$ such that for any $c\in (0,c_0)$,  the problem  (1.4)  has at least  two  positive normalized solutions  $u_c^1, u_c^2\in H^1(\mathbb R^N)$,  which are  radially symmetric and radially nonincreasing. Moreover, 
$$
\|\nabla u_c^1\|_2^2\sim c^{\frac{2(2N-q(N-2))}{4-N(q-2)}}\to 0, \  \ as \  c\to 0,
$$
$$
E(u_c^1)\simeq \frac{4q}{N(q-2)}S_q^{\frac{q}{q-2}}c^{\frac{2(2N-q(N-2))}{4-N(q-2)}}\to 0^-, \  \  as \ c\to 0,
$$
$$
\|\nabla u_c^2\|_2^2\simeq S_\alpha^{\frac{N+\alpha}{2+\alpha}},  \quad 
E(u_c^2)\simeq \frac{2+\alpha}{2(N+\alpha)}S_\alpha^{\frac{N+\alpha}{2+\alpha}}, \   \   as \ c\to 0.
$$

(II)  If $q>q_0$, then for any $c>0$ the problem  (1.4)  has at least  one  positive normalized solution  $u_c\in H^1(\mathbb R^N)$,  which is   radially symmetric and radially nonincreasing. Moreover, 
$$
\|\nabla u_c\|_2^2\simeq S_\alpha^{\frac{N+\alpha}{2+\alpha}},  \quad 
E(u_c)\simeq \frac{2+\alpha}{2(N+\alpha)}S_\alpha^{\frac{N+\alpha}{2+\alpha}}, \   \   as \ c\to 0,
$$
$$
\|\nabla u_c\|_2^2\sim c^{\frac{2(2N-q(N-2))}{4-N(q-2)}}\to 0, \  \ as \  c\to \infty,
$$
$$
E(u_c)\simeq \frac{4q}{N(q-2)}S_q^{\frac{q}{q-2}}c^{\frac{2(2N-q(N-2))}{4-N(q-2)}}\to 0^+, \  \  as \ c\to \infty.
$$
}

\smallskip
\noindent
{\bf Corollary 1.3.}  {\it  If $q=2^*$, $p\in (1+\frac{\alpha}{N-2}, \frac{N+\alpha}{N-2})$ for $N\ge 5$, $p\in (\max\{2, 1+\frac{\alpha}{2}\}, 2+\frac{\alpha}{2})$ for $N=4$ and  $ p\in (2+\alpha, 3+\alpha)$ for $N=3$, then  the following statements hold true:

(I)  If $p<p_0$, then there  exists  a constant $c_0>0$ such that for any $c\in (0,c_0)$,  the problem  (1.4)  has at least  two  positive normalized solutions  $u_c^1, u_c^2\in H^1(\mathbb R^N)$,  which are  radially symmetric and radially nonincreasing. Moreover, 
$$
\|\nabla u_c^1\|_2^2\sim c^{\frac{2(N+\alpha-p(N-2))}{2+\alpha-N(p-1)}}\to 0, \  \ as \  c\to 0,
$$
$$
E(u_c^1)\simeq \frac{N(p-1)-2-\alpha}{4p}S_p^{\frac{p}{p-1}}c^{\frac{2(N+\alpha-p(N-2))}{2+\alpha-N(p-1)}}\to 0^-, \  \  as \ c\to 0,
$$
$$
\|\nabla u_c^2\|_2^2\simeq S^{\frac{N}{2}},  \quad
E(u_c^2)\simeq \frac{1}{N}S^{\frac{N}{2}}, \   \   as \ c\to 0.
$$

(II)  If $p>p_0$, then for any $c>0$ the problem  (1.4)  has at least  one  positive normalized solution  $u_c\in H^1(\mathbb R^N)$,  which is   radially symmetric and radially nonincreasing. Moreover, 
$$
\|\nabla u_c\|_2^2\simeq S^{\frac{N}{2}},  \quad 
E(u_c)\simeq \frac{1}{N}S^{\frac{N}{2}}, \   \   as \ c\to 0,
$$
$$
\|\nabla u_c\|_2^2\sim c^{\frac{2(N+\alpha-p(N-2))}{2+\alpha-N(p-1)}}\to 0, \  \ as \  c\to \infty,
$$
$$
E(u_c)\simeq \frac{N(p-1)-2-\alpha}{4p}S_p^{\frac{p}{p-1}}c^{\frac{2(N+\alpha-p(N-2))}{2+\alpha-N(p-1)}}\to 0^+, \  \  as \ c\to \infty.
$$
}

\smallskip
\noindent
{\bf Corollary 1.4.} {\it Let $p\in (\frac{N+\alpha}{N}, \frac{N+\alpha}{N-2})$, $q\in (2,2^*)$, then the following statements hold true:

(I) If $p<p_0, q<q_0$  or  $p>p_0, q>q_0$,  then for any  $c>0$ the problem (1.4)  has at least one positive  normalized solution $u_c\in H^1(\mathbb R^N)$.
Moreover, if  $p<p_0, q<q_0$ and $ q< \frac{2(2p+\alpha)}{2+\alpha}$, then
$$
\|\nabla u_c\|_2^2\simeq \left\{\begin{array}{lcl}
 \frac{N(q-2)}{2q}S_q^{\frac{q}{q-2}} \left(\frac{2q}{2N-q(N-2)}S_q^{-\frac{q}{q-2}}c^2\right)^{\frac{2N-q(N-2)}{4-N(q-2)}}\to 0, &as& c\to 0,\\
\frac{N(p-1)-\alpha}{2p}S_p^{\frac{p}{p-1}}  \left(\frac{2p}{N+\alpha-p(N-2)}S_p^{-\frac{p}{p-1}}c^2\right)^{\frac{N+\alpha-p(N-2)}{2+\alpha-N(p-1)}}\to +\infty, &as& c\to \infty,
  \end{array}\right.
 $$
 $$
E(u_c)\simeq \left\{\begin{array}{lcl}
 \frac{N(q-2)-4}{4q}S_q^{\frac{q}{q-2}}\left(\frac{2q}{2N-q(N-2)}S_q^{-\frac{q}{q-2}}c^2\right)^{\frac{2N-q(N-2)}{4-N(q-2)}}\to 0^-,  &as&  c\to 0,\\
\frac{N(p-1)-2-\alpha}{4p}S_p^{\frac{p}{p-1}}\left(\frac{2p}{N+\alpha-p(N-2)}S_p^{-\frac{p}{p-1}}c^2\right)^{\frac{N+\alpha-p(N-2)}{2+\alpha-N(p-1)}}\to-\infty,  &as&  c\to \infty. 
 \end{array}\right.
 $$
 If $p<p_0, q<q_0$ and  $ q>\frac{2(2p+\alpha)}{2+\alpha}$, then 
 $$
\|\nabla u_c\|_2^2\simeq \left\{\begin{array}{lcl}
 \frac{N(p-1)-\alpha}{2p}S_p^{\frac{p}{p-1}}  \left(\frac{2q}{2N-q(N-2)}S_q^{-\frac{q}{q-2}}c^2\right)^{\frac{(q-2)(N+\alpha-p(N-2))}{(p-1)(4-N(q-2))}}\to +\infty, &as& c\to \infty,\\
 \frac{N(q-2)}{2q}S_q^{\frac{q}{q-2}}   \left(\frac{2p}{N+\alpha-p(N-2)}S_p^{-\frac{p}{p-1}}c^2\right)^{\frac{(p-1)(2N-q(N-2))}{(q-2)(2+\alpha-N(p-1))}}\to 0, &as& c\to 0,
  \end{array}\right.
 $$
  $$
E(u_c)\simeq \left\{\begin{array}{lcl}
 \frac{N(q-2)-4}{4q}S_q^{\frac{q}{q-2}}\left(\frac{2q}{2N-q(N-2)}S_q^{-\frac{q}{q-2}}c^2\right)^{\frac{2N-q(N-2)}{4-N(q-2)}}\to -\infty,  &as&  c\to \infty,\\
\frac{N(p-1)-2-\alpha}{4p}S_p^{\frac{p}{p-1}}\left(\frac{2p}{N+\alpha-p(N-2)}S_p^{-\frac{p}{p-1}}c^2\right)^{\frac{N+\alpha-p(N-2)}{2+\alpha-N(p-1)}}\to 0^-,  &as&  c\to 0. 
 \end{array}\right.
 $$
 If $p>p_0, q>q_0$ and $ q< \frac{2(2p+\alpha)}{2+\alpha}$, then
$$
\|\nabla u_c\|_2^2\simeq \left\{\begin{array}{lcl}
  \frac{N(q-2)}{2q}S_q^{\frac{q}{q-2}}\left(\frac{2q}{2N-q(N-2)}S_q^{-\frac{q}{q-2}}c^2\right)^{\frac{2N-q(N-2)}{4-N(q-2)}}\to 0, &as& c\to \infty,\\
 \frac{N(p-1)-\alpha}{2p}S_p^{\frac{p}{p-1}}  \left(\frac{2p}{N+\alpha-p(N-2)}S_p^{-\frac{p}{p-1}}c^2\right)^{\frac{N+\alpha-p(N-2)}{2+\alpha-N(p-1)}}\to +\infty, &as& c\to 0,
  \end{array}\right.
 $$
 $$
E(u_c)\simeq \left\{\begin{array}{lcl}
 \frac{N(q-2)-4}{4q}S_q^{\frac{q}{q-2}}\left(\frac{2q}{2N-q(N-2)}S_q^{-\frac{q}{q-2}}c^2\right)^{\frac{2N-q(N-2)}{4-N(q-2)}}\to 0^+,  &as&  c\to \infty,\\
\frac{N(p-1)-2-\alpha}{4p}S_p^{\frac{p}{p-1}}\left(\frac{2p}{N+\alpha-p(N-2)}S_p^{-\frac{p}{p-1}}c^2\right)^{\frac{N+\alpha-p(N-2)}{2+\alpha-N(p-1)}}\to +\infty,  &as&  c\to 0. 
 \end{array}\right.
 $$
 If $p>p_0, q>q_0$ and $ q>\frac{2(2p+\alpha)}{2+\alpha}$, then 
$$
\|\nabla u_c\|_2^2\simeq \left\{\begin{array}{lcl}
  \frac{N(q-2)}{2q}S_q^{\frac{q}{q-2}}  \left(\frac{2p}{N+\alpha-p(N-2)}S_p^{-\frac{p}{p-1}}c^2\right)^{\frac{(p-1)(2N-q(N-2))}{(q-2)(2+\alpha-N(p-1))}}\to 0, &as& c\to \infty,\\
   \frac{N(p-1)-\alpha}{2p}S_p^{\frac{p}{p-1}}  \left(\frac{2q}{2N-q(N-2)}S_q^{-\frac{q}{q-2}}c^2\right)^{\frac{(q-2)(N+\alpha-p(N-2))}{(p-1)(4-N(q-2))}}\to +\infty, &as& c\to 0,
     \end{array}\right.
 $$
 $$
E(u_c)\simeq \left\{\begin{array}{lcl}
\frac{N(p-1)-2-\alpha}{4p}S_p^{\frac{p}{p-1}}\left(\frac{2p}{N+\alpha-p(N-2)}S_p^{-\frac{p}{p-1}}c^2\right)^{\frac{N+\alpha-p(N-2)}{2+\alpha-N(p-1)}}\to 0^+,  &as&  c\to \infty,\\
 \frac{N(q-2)-4}{4q}S_q^{\frac{q}{q-2}}\left(\frac{2q}{2N-q(N-2)}S_q^{-\frac{q}{q-2}}c^2\right)^{\frac{2N-q(N-2)}{4-N(q-2)}}\to +\infty,  &as&  c\to 0.
  \end{array}\right.
 $$

(II)  If $p=p_0$ and $q\not=q_0$ \ $( \ resp. \  q=q_0$ and $p\not=p_0\ )$,  then  for any $c\in (0,\sqrt{\frac{2+\alpha}{N+2+\alpha}}S_{p_0}^{\frac{N+2+\alpha}{2(2+\alpha)}})$  \ $(resp. \  c\in (0, \sqrt{\frac{2}{N+2}}S_{q_0}^{\frac{N+2}{4}})\ )$,  the problem (1.4) has at least one positive  normalized solution $u_c\in H^1(\mathbb R^N)$.  Moreover,
if $p=p_0, q<q_0$, then as $c\to 0$,
$$
\|\nabla u_c\|_2^2\simeq     \frac{N(q-2)}{2q}S_q^{\frac{q}{q-2}} \left(\frac{2q}{2N-q(N-2)}S_q^{-\frac{q}{q-2}}c^2\right)^{\frac{2N-q(N-2)}{4-N(q-2)}}\to 0,
$$
$$
E(u_c)\simeq  \frac{N(q-2)-4}{4q}S_q^{\frac{q}{q-2}}\left(\frac{2q}{2N-q(N-2)}S_q^{-\frac{q}{q-2}}c^2\right)^{\frac{2N-q(N-2)}{4-N(q-2)}}\to 0^-.
$$
If $p=p_0, q>q_0$, then as $c\to 0$,
$$
\|\nabla u_c\|_2^2\simeq \frac{N}{N+2+\alpha}S_{p_0}^{\frac{N+2+\alpha}{2+\alpha}} \left(\frac{2q}{2N-q(N-2)}S_{q}^{-\frac{q}{q-2}}c^2\right)^{\frac{2(q-2)}{4-N(q-2)}}\to +\infty,
$$
$$
E(u_c)\simeq  \frac{N(q-2)-4}{4q}S_q^{\frac{q}{q-2}}\left(\frac{2q}{2N-q(N-2)}S_q^{-\frac{q}{q-2}}c^2\right)^{\frac{2N-q(N-2)}{4-N(q-2)}}\to +\infty,
$$
and as $c\to \sqrt{\frac{2+\alpha}{N+2+\alpha}}S_{p_0}^{\frac{N+2+\alpha}{2(2+\alpha)}}$,
$$
\|\nabla u_c\|_2^2\sim \frac{N(q-2)}{2q}S_{q}^{\frac{q}{q-2}}\left(\frac{2+\alpha}{N+2+\alpha}S_{p_0}^{\frac{N+2+\alpha}{2+\alpha}}-c^2\right)^{\frac{2(2N-q(N-2))}{(q-2)(N(q-2)-4)}}\to 0, 
$$
$$
E(u_c)=O\left((\frac{2+\alpha}{N+2+\alpha}S_{p_0}^{\frac{N+2+\alpha}{2+\alpha}}-c^2)^{\frac{N(q-2)}{N(q-2)-4}}\right)\to 0.
$$
If $p<p_0, q=q_0$, then as $c\to 0$,
$$
\|\nabla u_c\|_2^2\simeq \frac{N}{N+2}S_{q_0}^{\frac{N+2}{2}}  \left(\frac{2p}{N+\alpha-p(N-2)}S_p^{-\frac{p}{p-1}}c^2\right)^{\frac{2(p-1)}{2+\alpha-N(p-1)}}\to 0, 
$$
$$
E(u_c)\simeq \frac{N(p-1)-2-\alpha}{4p}S_p^{\frac{p}{p-1}}\left(\frac{2p}{N+\alpha-p(N-2)}S_p^{-\frac{p}{p-1}}c^2\right)^{\frac{N+\alpha-p(N-2)}{2+\alpha-N(p-1)}}\to 0^-.
$$
If $p>p_0, q=q_0$, then as $c\to 0$, 
$$
\|\nabla u_c\|_2^2\simeq  \frac{N(p-1)-\alpha}{2p}S_{p}^{\frac{p}{p-1}} \left(\frac{2p}{N+\alpha-p(N-2)}S_p^{-\frac{p}{p-1}}c^2\right)^{\frac{N+\alpha-p(N-2)}{2+\alpha-N(p-1)}}\to +\infty,
$$
$$
E(u_c)\simeq  \frac{N(p-1)-2-\alpha}{4p}S_p^{\frac{p}{p-1}}\left(\frac{2p}{N+\alpha-p(N-2)}S_p^{-\frac{p}{p-1}}c^2\right)^{\frac{N+\alpha-p(N-2)}{2+\alpha-N(p-1)}}\to +\infty, 
$$
and as $c\to \sqrt{\frac{2}{N+2}}S_{q_0}^{\frac{N+2}{4}}$, 
$$
\|\nabla u_c\|_2^2\sim \frac{N}{N+2}S_{q_0}^{\frac{N+2}{2}}\left(\frac{2}{N+2}S_{q_0}^{\frac{N+2}{2}}-c^2\right)^{\frac{2}{N(p-1)-2-\alpha}}\to 0,
 $$
$$
E(u_c)=O\left((\frac{2}{N+2}S_{q_0}^{\frac{N+2}{2}}-c^2)^{\frac{N(p-1)-\alpha}{N(p-1)-2-\alpha}}\right)\to 0.
$$

(III) If $p<p_0, q>q_0$, then for any $c\in (0,\sup_{\varepsilon>0}\sqrt{M(\varepsilon)})$ the problem (1.4) has two positive normalized solutions $u_c^1$ and $u_c^2$ satisfying 
$$
\|\nabla u_c^1\|_2^2\simeq \frac{N(q-2)}{2q}S_q^{\frac{q}{q-2}} \left(\frac{2p}{N+\alpha-p(N-2)}S_p^{-\frac{p}{p-1}}c^2\right)^{\frac{(p-1)(2N-q(N-2))}{(q-2)(2+\alpha-N(p-1))}}\to 0, \ as \ c\to 0,
$$
$$
E(u_c^1)\simeq \frac{N(p-1)-2-\alpha}{4p}S_p^{\frac{p}{p-1}}\left(\frac{2p}{N+\alpha-p(N-2)}S_p^{-\frac{p}{p-1}}c^2\right)^{\frac{N+\alpha-p(N-2)}{2+\alpha-N(p-1)}}\to 0^-, \ as \ c\to 0,
$$
$$
\|\nabla u_c^2\|_2^2\simeq   \frac{N(p-1)-\alpha}{2p}S_p^{\frac{p}{p-1}} \left(\frac{2q}{2N-q(N-2)}S_q^{-\frac{q}{q-2}}c^2\right)^{\frac{(q-2)(N+\alpha-p(N-2))}{(p-1)(4-N(q-2))}}\to +\infty, \ as \ c\to 0,
$$
$$
E(u_c^2)\simeq \frac{N(q-2)-4}{4q}S_q^{\frac{q}{q-2}}\left(\frac{2q}{2N-q(N-2)}S_q^{-\frac{q}{q-2}}c^2\right)^{\frac{2N-q(N-2)}{4-N(q-2)}}\to +\infty, \ as \ c\to 0.
$$
If $p>p_0, q<q_0$, then for any $c\in (0,\sup_{\varepsilon>0}\sqrt{M(\varepsilon)})$, the problem (1.4) has two positive normalized solutions $u_c^1$ and $u_c^2$ satisfying 
$$
\|\nabla u_c^1\|_2^2\simeq \frac{N(q-2)}{2q}S_q^{\frac{q}{q-2}}\left(\frac{2q}{2N-q(N-2)}S_q^{-\frac{q}{q-2}}c^2\right)^{\frac{2N-q(N-2)}{4-N(q-2)}}\to 0, \ as \ c\to 0,
$$
$$
E(u_c^1)\simeq \frac{N(q-2)-4}{4q}S_q^{\frac{q}{q-2}}\left(\frac{2q}{2N-q(N-2)}S_q^{-\frac{q}{q-2}}c^2\right)^{\frac{2N-q(N-2)}{4-N(q-2)}}\to 0^-, \ as \ c\to 0,
$$
$$
\|\nabla u_c^2\|_2^2\simeq   \frac{N(p-1)-\alpha}{2p}S_p^{\frac{p}{p-1}} \left(\frac{2p}{N+\alpha-p(N-2)}S_p^{-\frac{p}{p-1}}c^2\right)^{\frac{N+\alpha-p(N-2)}{2+\alpha-N(p-1)}}\to +\infty, \ as \ c\to 0,
$$
$$
E(u_c^2)\simeq \frac{N(p-1)-2-\alpha}{4p}S_p^{\frac{p}{p-1}}\left(\frac{2p}{N+\alpha-p(N-2)}S_p^{-\frac{p}{p-1}}c^2\right)^{\frac{N+\alpha-p(N-2)}{2+\alpha-N(p-1)}}\to +\infty, \ as \ c\to 0.
$$}

Similar existence results on normalized solutions  are already obtained in \cite{ Li-3, Li-4, Li-5, Li-6, Sun-1, Sun-2} in the whole possible range of  parameters, but the precise asymptotic behaviour of the normalized solutions was not addressed there and is a new result of this paper.

\smallskip
\noindent
{\bf Remark 1.2.}    In the case $N=3$, $p=2$, $q=4$ and $\alpha\in(0,3)$ the problem $(P_\varepsilon)$  is known as in astrophysics as the Gross-Pitaevskii-Poisson equation (see the Introduction). By Proposition 1.4, we see that $(P_\varepsilon)$ admits a positive ground state $u_\varepsilon\in H^1(\mathbb R^N)$, which is radially symmetric and radially nonincreasing. Moreover,  as $\varepsilon\to 0$, we have 
$$
u_\varepsilon(0)\sim \varepsilon^{\frac{2+\alpha}{4}},\quad 
\|\nabla u_\varepsilon\|_2^2\simeq\frac{3}{4}S_q^2 \varepsilon^{\frac{1}{2}},
$$
$$
M(\varepsilon):=\|u_\varepsilon\|_2^2=\varepsilon^{\frac{\alpha-1}{2}}\left(\frac{1+\alpha}{4}S_p^2-\Theta(\varepsilon^{\frac{\alpha}{2}})\right),
\quad 
E(u_\varepsilon)=\varepsilon^{\frac{1+\alpha}{2}}\left(\frac{1-\alpha}{8}S_p^2+O(\varepsilon^{\frac{\alpha}{2}})\right).
$$
As $\varepsilon\to \infty$, we have 
$$
u_\varepsilon(0)\sim \varepsilon^{\frac{1}{2}},
\quad 
\|\nabla u_\varepsilon\|_2^2\simeq \frac{3-\alpha}{4}S_p^2\varepsilon^{\frac{1+\alpha}{2}},
$$
$$
M(\varepsilon):=\|u_\varepsilon\|_2^2=\varepsilon^{-\frac{1}{2}}\left(\frac{1}{4}S_q^2+O(\varepsilon^{-\frac{\alpha}{2}})\right),
\quad 
E(u_\varepsilon)=\varepsilon^{\frac{1}{2}}\left(\frac{1}{8}S_q^2+O(\varepsilon^{-\frac{\alpha}{2}})\right).
$$
Therefore, we conclude that 

(1)  when $\alpha\in (1,3)$ and $c\in (0,\sup_{\varepsilon>0}\sqrt{M(\varepsilon)})$,  (1.4) admits two positive normalized solutions  $u_c^1$ and $u_c^2$  satisfying 
$$
\|u_c^1\|_\infty\sim \left(\frac{4}{1+\alpha}S_p^{-2}c^2\right)^{\frac{2+\alpha}{2(\alpha-1)}}\to 0,  \  \   as \  c\to 0,
$$
$$
\|\nabla u_c^1\|_2^2\simeq \frac{3}{4}S_q^2\left(\frac{4}{1+\alpha}S_p^{-2}c^2\right)^{\frac{1}{\alpha-1}}\to 0, \  \   as \  c\to 0,
$$
$$
E(u_c^1)\simeq -\frac{\alpha-1}{8}(\frac{4}{1+\alpha})^{\frac{1+\alpha}{\alpha-1}}S_p^{-\frac{4}{\alpha-1}}c^{\frac{2(1+\alpha)}{\alpha-1}}\to 0^-, \  \   as \  c\to 0,
$$
$$
\|u_c^2\|_\infty\sim \frac{1}{4}S_q^2c^{-2}\to +\infty,  \  \   as \  c\to 0,
$$
$$
\|\nabla u_c^2\|_2^2\simeq \frac{3-\alpha}{4}S_p^2\left(\frac{1}{16}S_q^4c^{-4}\right)^{\frac{1+\alpha}{2}}\to +\infty,  \  \   as \  c\to 0,
$$
$$
E(u_c^2)\simeq \frac{1}{32}S_q^4c^{-2}\to +\infty,  \  \   as \  c\to 0.
$$
Moreover, as $c\to 0$, the rescaled family
$$
w_c^1(x):=\varepsilon_{1c}^{-\frac{2+\alpha}{4}}u_c^1(\varepsilon_{1c}^{-\frac{1}{2}}x), \quad \varepsilon_{1c}\simeq (\frac{4}{1+\alpha})^{\frac{2}{\alpha-1}}S_p^{-\frac{4}{\alpha-1}}c^{\frac{4}{\alpha-1}},
$$
converges in $H^1(\mathbb R^3)$, up to a subsequence, to a positive solution $w_0$ of the Choquard equation
\[
-\Delta w+w=(I_\alpha\ast |w|^2)w, \quad x\in\mathbb R^3.
\eqno(1.10)
\]
If $\alpha=2$, then positive solution $w_0$ of (1.10)  is unique  \cite{Lieb-1} and $w_{c}^1\to w_0$ in $H^1(\mathbb R^3)$. As $c\to 0$, the rescaled family 
$$
w_c^2(x):=\varepsilon_{2c}^{-\frac{1}{2}}u_{c}^2(\varepsilon_{2c}^{-\frac{1}{2}}x), \quad  \varepsilon_{2c}\simeq \frac{1}{16}S_q^4c^{-4},
$$
converges in $H^1(\mathbb R^3)$ to the unique positive solution $w_\infty$ of the Gross-Pitaevskii equation
$$
-\Delta w+w=w^3, \quad x\in\mathbb R^3.\eqno(1.11)
$$

(2)  when $\alpha=1$ and $c\in (0,\frac{\sqrt 2}{2}S_p)$, (1.4) admits a positive normalized solution $u_c$ satisfying  
$$
\|u_c\|_\infty\sim \frac{1}{4}S_q^2c^{-2}\to +\infty,  \  \   as \  c\to 0,
$$
$$
\|\nabla u_c\|_2^2\simeq \frac{1}{32}S_p^2S_q^4c^{-4}\to +\infty,  \  \   as \  c\to 0,
$$
$$
E(u_c)\simeq \frac{1}{32}S_q^4c^{-2}\to +\infty,  \  \   as \  c\to 0,
$$
$$
\|u_c\|_\infty\sim (\frac{1}{2}S_p^2-c^2)^{\frac{3}{2}}\to 0, \  \  as \  c\to\frac{\sqrt 2}{2}S_p,
$$
$$
\|\nabla u_c\|_2^2\sim \frac{3}{4}S_q^2(\frac{1}{2}S_p^2-c^2)\to 0,   \  \  as \  c\to\frac{\sqrt 2}{2}S_p,
$$
$$
E(u_c)= O\left((\frac{1}{2}S_p^2-c^2)^3\right)\to 0,   \  \  as \  c\to\frac{\sqrt 2}{2}S_p.
$$
Moreover, as $c\to \frac{\sqrt 2}{2}S_p$, the rescaled family
$$
w_c(x):=\varepsilon_c^{-\frac{3}{4}}u_c(\varepsilon_c^{-\frac{1}{2}}x), \quad \varepsilon_c\sim (\frac{1}{2}S_p^2-c^2)^2,
$$
converges in $H^1(\mathbb R^3)$, up to a subsequence, to a positive solution $w_0$ of  (1.10) with $\alpha=1$,
and as $c\to 0$, while the rescaled family 
$$
w_c(x):=\varepsilon_c^{-\frac{1}{2}}u_c(\varepsilon_c^{-\frac{1}{2}}x), \quad  \varepsilon_c\simeq \frac{1}{16}S_q^4c^{-4},
$$
converges in $H^1(\mathbb R^3)$ to the unique positive solution  $w_\infty$ of (1.11).

(3) when $\alpha\in (0,1)$ and any $c>0$,  (1.4) admits a positive normalized solution $u_c$ satisfying 
$$
\|u_c\|_\infty\sim \frac{1}{4}S_q^2c^{-2}\to +\infty,  \  \   as \  c\to 0,
$$
$$
\|\nabla u_c\|_2^2\simeq \frac{3-\alpha}{4}S_p^2\left(\frac{1}{4}S_q^2c^{-2}\right)^{1+\alpha}\to +\infty,  \  \   as \  c\to 0,
$$
$$
E(u_c)\simeq \frac{1}{32}S_q^4c^{-2}\to +\infty,  \  \   as \  c\to 0,
$$
$$
\|u_c\|_\infty\sim \left(\frac{4}{1+\alpha}S_p^{-2}c^2\right)^{\frac{2+\alpha}{2(\alpha-1)}}\to 0,  \  \   as \  c\to \infty,
$$
$$
\|\nabla u_c\|_2^2\simeq \frac{3}{4}S_q^2\left(\frac{4}{1+\alpha}S_p^{-2}c^2\right)^{\frac{1}{\alpha-1}}\to 0, \  \   as \  c\to \infty,
$$
$$
E(u_c)\simeq \frac{1-\alpha}{8}(\frac{4}{1+\alpha})^{\frac{1+\alpha}{\alpha-1}}S_p^{-\frac{4}{\alpha-1}}c^{\frac{2(1+\alpha)}{\alpha-1}}\to 0^+, \  \   as \  c\to \infty.
$$
Moreover, as $c\to \infty$, the rescaled family
$$
w_c(x):=\varepsilon_c^{-\frac{2+\alpha}{4}}u_c(\varepsilon_c^{-\frac{1}{2}}x), \quad \varepsilon_c\simeq (\frac{4}{1+\alpha})^{\frac{2}{\alpha-1}}S_p^{-\frac{4}{\alpha-1}}c^{\frac{4}{\alpha-1}},
$$
converges in $H^1(\mathbb R^3)$, up to a subsequence,  to a positive solution $w_0$ of  (1.10),
and as $c\to 0$, the rescaled family 
$$
w_c(x):=\varepsilon_c^{-\frac{1}{2}}u_c(\varepsilon_c^{-\frac{1}{2}}x), \quad  \varepsilon_c\simeq \frac{1}{16}S_q^4c^{-4},
$$
converges in $H^1(\mathbb R^3)$ to the unique positive solution $w_\infty$ of (1.11).

\smallskip
\noindent
{\bf Remark 1.3.} In \cite{Moroz-1}, the second author and  Muratov first  study the asymptotic properties and limit profiles of ground states for  a class of local scalar field equations  
$$
-\Delta u+\varepsilon u=|u|^{p-2}u-|u|^{q-2}u \quad {\rm in} \ \mathbb R^N,
$$
where $N\ge 3$ and $q>p>2$.  Later, in \cite{Lewin-1},  M. Lewin and S. Nodari prove a general result about the uniqueness and non-degeneracy of positive radial solutions to the above equation. The non-degeneracy of the unique solution $u_\varepsilon$ allows to derive its behaviour in the two limits $\varepsilon\to 0$ and $\varepsilon\to \varepsilon_*$, where $\varepsilon_*$ is a  threshold for the existence.  Amongst other things, a precise asymptotic  expression of $M(\varepsilon)=\|u_\varepsilon\|_2^2$ is obtained in \cite{Lewin-1}, which gives the uniqueness of energy minimizers at fixed mass in certain regimes. 
 
In \cite{Liu-1}, Zeng Liu and the second author extended the results in \cite{Moroz-1} to a class of Choquard type equation with attractive nonlocal and repulsive local interaction terms,
\[
-\Delta u+\varepsilon u=(I_\alpha \ast |u|^{p})|u|^{p-2}u- |u|^{q-2}u 
\quad {\rm in} \ \mathbb R^N,
\eqno(1.12)
 \]
where $N\ge 3$, $p>\frac{N+\alpha}{N}$ and $q>2$. Under some assumptions on the exponents $p$ and $q$, the limit profiles of the ground states are discussed in the two cases $\varepsilon\to 0$ and $\varepsilon\to \infty$. 

In \cite{Akahori-2, Akahori-3}, T. Akahori et al. consider the following Schr\"odinger equation with two focusing  exponents $p$ and $q$:
\[
-\Delta u+\varepsilon u=|u|^{p-2}u+\mu|u|^{q-2}u, \quad {\rm in} \ \mathbb R^N,
\eqno(1.13)
\]
where $N\ge 3$, $2<q<p\le 2^*$ and $\mu>0$ is a parameter.  When  $\mu=1$, $p=2^*$ and $q\in (2,2^*)$,  the authors in \cite{Akahori-3}   proved that for small $\varepsilon>0$ the ground state is unique and as $\varepsilon\to 0$, the unique ground state $u_\varepsilon$ tends to the unique positive solution of the equation $-\Delta u+u=u^{q-1}$. After a suitable rescaling, authors in \cite{Akahori-2}  establish a uniform decay estimate for the ground state $u_\varepsilon$, and then   prove  the uniqueness and nondegeneracy of ground states $u_\varepsilon$  for $N\ge 5$ and large $\varepsilon>0$,  and show that for $N\ge 3$, as $\varepsilon\to\infty$, $u_\varepsilon$ tends to a particular solution of the critical Emden--Fowler equation.
More recently, Jeanjean, Zhang and Zhong \cite{Jeanjean-4} also study the asymptotic behaviour of solutions as $\varepsilon\to 0$ and $\varepsilon\to \infty$  for the equation  (1.13) with a general subcritical nonlinearity and discuss the cornection  to the existence, non-existence and multiplicity of  prescribed mass positive solutions to (1.13) with  the  associated $L^2$ constraint condition $\int_{\mathbb R^N}|u|^2=c^2$. For  other papers relevant to (1.13), we refer the reader to \cite{Ma-1} and the references therein.

 For quite a long time the paper \cite{Jeanjean-2} was the only one dealing with existence of normalized solutions in cases when the energy is unbounded from below on the $L^2$-constraint. More recently, however, problems of this type received much attention. We refer the readers to  \cite{ Soave-1,  Soave-2, Wei-1, Wei-2} and references therein  for the existence and multiplicity of  normalized solutions to the  equations (1.13).
 In \cite{Wei-1}, Wei and Wu study the existence and asymptotic behavior of normalized solutions  for (1.13) with $p=2^*$, and  obtain in some cases precise asymptotic behavior of ground states and mountain-pass solutions as $\mu\to 0$ and $\mu$ goes to its upper bound,  here  $\lambda:=-\varepsilon$ arises as a Lagrange multiplier.  We refer the readers to 
 \cite{ Jeanjean-3, Jeanjean-4, Soave-1,  Soave-2, Wei-2} for  the  asymptotic behaviour of normalized solutions when the parameter $\mu$ varies in its range.  Roughly speaking, the parameter $\mu$ changes some thresholds for the existence but it does not change the qualitative properties of solutions. In a sense,  changing the parameter $\mu$ is equivalent to changing the mass $c>0$.  More precisely,  it follows from the reduction given in \cite{Wei-2} that finding a normalized solution of (1.13) with $p=2^*$ and mass constrained  condition $\int_{\mathbb R^N}|u|^2=c^2$
is equivalent to finding a normalized solution of the problem
   $$
\left\{\begin{array}{rl}
&-\Delta u+\varepsilon u=|u|^{2^*-2}u+|u|^{q-2}u, \  \ in  \  \mathbb R^N,\\
&u\in H^1(\mathbb R^N),  \   \    \   \int_{\mathbb R^N}|u|^2=c^2\mu^{\frac{2}{q-q\gamma_q}},
\end{array}\right.
\eqno(1.14)
$$ 
where $\gamma_q=\frac{N(q-2)}{2q}$. In particular,  sending $\mu\to 0$ (resp. $\mu\to \infty$) is equavelent to sending the mass $c\to 0$ (resp. $c\to \infty$).

  In  \cite{Jeanjean-3},  taking the mass $c>0$ as a parameter,  Jeanjean and  Le also discuss the asymptotic behaviour of normalized solutions of (1.13) with  $p=2^*$. In  the case $N\ge 4$, $2 < q < 2 + \frac{4}{N}$,  
   amongst other things,  Jeanjean and Le obtained a normalized  solution $u_c$ of mountain pass type for small $c>0$ and proved  that 
\[
\lim_{c\to 0}\|\nabla u_c\|_2^2=S^{\frac{N}{2}}, \quad  \lim_{c\to 0}E(u_c)=\frac{1}{N}S^{\frac{N}{2}}.
\eqno(1.15)
\] 
 The relationship between action ground state and energy ground state is discussed in \cite{Dovetta-1,  Ilyasov-1, Jeanjean-3, Jeanjean-1, Wei-2}. In particular, it is shown in those papers that the energy ground state is 
necessarily an action ground state. So  the results in the present paper imply the precise asymptotic behaviour of normalized ground states, and we shall address this problem in a forthcoming paper. 

\smallskip
\noindent
{\bf Remark 1.4.} 
By the precise asymptotic expression of $M(\varepsilon)$ obtained in Propositions 1.1-1.4 (also see Corollaries 2.1-2.4), 
 if $M(\varepsilon)$ is of class $C^1$ for small $\varepsilon>0$ and large $\varepsilon>0$,
then the sign of $M'(\varepsilon)$  follows from the following Lemma 1.1.  Note that  by Proposition 1.4 (I),  as $\varepsilon\to 0$, we have 
$$
M(\varepsilon)=\frac{2}{N+2}S_{q_0}^{\frac{N+2}{2}}-\Theta(\varepsilon^{\frac{N(p-1)-2-\alpha}{2}}), \quad {\rm if} \ q=q_0 \ {\rm and} \ p>p_0,
$$
$$
M(\varepsilon)=\frac{2+\alpha}{N+2+\alpha}S_{p_0}^{\frac{N+2+\alpha}{2+\alpha}}-\Theta(\varepsilon^{\frac{N(q-2)-4}{4}}), \quad {\rm if} \ q>q_0 \ {\rm and} \ p=p_0.
$$
Therefore, to prove the second and third inequalities in (1.7), we replace $M(\varepsilon)$ by 
$$
\frac{2}{N+2}S_{q_0}^{\frac{N+2}{2}}-M(\varepsilon) \  \   \   {\rm and}  \   \   \  \frac{2+\alpha}{N+2+\alpha}S_{p_0}^{\frac{N+2+\alpha}{2+\alpha}}-M(\varepsilon),
$$
in the following Lemma 1.1,  respectively.

\smallskip
\noindent
{\bf Lemma 1.1.} {\it Let $\eta>0$ is a constant and  $M(\varepsilon)$ is of class $C^1$ on $ (0,\infty)$.  Then the following statements hold true:

(1)  If $M(\varepsilon)\sim \varepsilon^\eta$ as $\varepsilon\to 0$, then there is $\varepsilon_0>0$ such that  $M'(\varepsilon)>0$ 
for $\varepsilon\in (0,\varepsilon_0)$.

(2)  If $M(\varepsilon)\sim \varepsilon^{-\eta}$ as $\varepsilon\to 0$, then there is $\varepsilon_0>0$ such that  $M'(\varepsilon)<0$ 
for $\varepsilon\in (0,\varepsilon_0)$. 

(3)  If $M(\varepsilon)\sim \varepsilon^{-\eta}$ as $\varepsilon\to\infty$, then there is $\varepsilon_\infty>0$ such that  $M'(\varepsilon)<0$ 
for $\varepsilon\in (\varepsilon_\infty, \infty)$.

(4)  If $M(\varepsilon)\sim \varepsilon^{\eta}$ as $\varepsilon\to\infty$, then there is  $\varepsilon_\infty>0$ such that  $M'(\varepsilon)>0$ 
for $\varepsilon\in (\varepsilon_\infty, \infty)$.

}

\begin{proof}
Since $\eta>0$, we can choose $\delta>0$ such that
$$
\max\{1,\eta\}<\delta<\eta+1.
$$
To prove (1), we let $\tilde M(\varepsilon)=\varepsilon^{1-\delta}M(\varepsilon)$ for $\varepsilon\in (0,\infty)$ and $\tilde M(0)=0$. Since 
$$
0\le \tilde M(\varepsilon)\le C\varepsilon^{1-\delta+\eta}\to 0, \quad {\rm as} \  \varepsilon\to 0,
$$
it follows that  $\tilde M(\varepsilon)$ is of class $C^1$ in $(0,\infty)$ and continuous on $[0,\infty)$. Moreover,
$$
\tilde M'(0):=\lim_{\varepsilon\to 0^+}\frac{\tilde M(\varepsilon)-\tilde M(0)}{\varepsilon}=\lim_{\varepsilon\to 0^+}\frac{M(\varepsilon)}{\varepsilon^\delta}
\ge \lim_{\varepsilon\to 0^+}c\varepsilon^{\eta-\delta}=+\infty.
$$
Therefore, there exists $\varepsilon_0>0$ such that 
$\tilde M'(\varepsilon)>0$ for any $\varepsilon\in (0,\varepsilon_0)$.
That is, 
$$
\tilde M'(\varepsilon)=(1-\delta)\varepsilon^{-\delta}M(\varepsilon)+\varepsilon^{1-\delta}M'(\varepsilon)>0,
$$
which together with $\delta>1$ implies that
$$
M'(\varepsilon)>(\delta-1)\varepsilon^{-1}M(\varepsilon)>0.
$$
The proof  of (1) is complete.

To prove (2), we let $\tilde M(\varepsilon)=-\frac{\varepsilon^{1-\delta}}{M(\varepsilon)}$ for $\varepsilon\in (0,\infty)$ and $\tilde M(0)=0$. Since 
$$
0\ge \tilde M(\varepsilon)\ge -\frac{\varepsilon^{1-\delta}}{c\varepsilon^{-\eta}}\to 0, \quad {\rm as} \  \varepsilon\to 0,
$$
it follows that  $\tilde M(\varepsilon)$ is of class $C^1$ in $(0,\infty)$ and continuous on $[0,\infty)$. Moreover,
$$
\tilde M'(0):=\lim_{\varepsilon\to 0^+}\frac{\tilde M(\varepsilon)-\tilde M(0)}{\varepsilon}= -\lim_{\varepsilon\to 0^+}\frac{\varepsilon^{-\delta}}{M(\varepsilon)}
\le -\lim_{\varepsilon\to 0^+}\frac{\varepsilon^{-\delta}}{C\varepsilon^{-\eta}}=-\infty.
$$
Therefore, there exists $\varepsilon_0>0$ such that 
$\tilde M'(\varepsilon)<0$ for any $\varepsilon\in (0,\varepsilon_0)$.
That is, 
$$
\tilde M'(\varepsilon)=-\frac{(1-\delta)\varepsilon^{-\delta}M(\varepsilon)-\varepsilon^{1-\delta}M'(\varepsilon)}{M(\varepsilon)^2}<0,
$$
which together with $\delta>1$ implies that
$$
M'(\varepsilon)<(1-\delta)\varepsilon^{-1}M(\varepsilon)<0.
$$
The proof of (2) is complete. 

To prove (3), we let $\tilde M(\varepsilon)=-\varepsilon^{1-\delta}M(\frac{1}{\varepsilon})$ for $\varepsilon\in (0,\infty)$ and $\tilde M(0)=0$. Since 
$$
0\ge \tilde M(\varepsilon)\ge -C\varepsilon^{1-\delta+\eta}\to 0, \quad {\rm as} \  \varepsilon\to 0,
$$
it follows that  $\tilde M(\varepsilon)$ is of class $C^1$ in $(0,\infty)$ and continuous on $[0,\infty)$. Moreover,
$$
\tilde M'(0):=\lim_{\varepsilon\to 0^+}\frac{\tilde M(\varepsilon)-\tilde M(0)}{\varepsilon}\le -c\lim_{\varepsilon\to 0^+}\varepsilon^{\mu-\delta}=-\infty.
$$
Therefore, there exists $\varepsilon_0>0$ such that 
$\tilde M'(\varepsilon)<0$ for any $\varepsilon\in (0,\varepsilon_0)$.
That is, 
$$
\tilde M'(\varepsilon)=-(1-\delta)\varepsilon^{-\delta}M(\frac{1}{\varepsilon})+\varepsilon^{-1-\delta}M'(\frac{1}{\varepsilon})<0,
$$
which together with $\delta>1$ implies that
$$
M'(\frac{1}{\varepsilon})<(1-\delta)\varepsilon M(\frac{1}{\varepsilon})<0.
$$
Take $\varepsilon_\infty=\frac{1}{\varepsilon_0}$, then $M'(\varepsilon)<0$ for $\varepsilon\in (\varepsilon_\infty,\infty)$. The proof of (3) is complete.

To prove (4), we let $\tilde M(\varepsilon)=\frac{\varepsilon^{1-\delta}}{M(\frac{1}{\varepsilon})}$ for $\varepsilon\in (0,\infty)$ and $\tilde M(0)=0$. Since 
$$
0\le \tilde M(\varepsilon)\le \frac{\varepsilon^{1-\delta}}{c\varepsilon^{-\eta}}\to 0, \quad {\rm as} \  \varepsilon\to 0,
$$
it follows that  $\tilde M(\varepsilon)$ is of class $C^1$ in $(0,\infty)$ and continuous on $[0,\infty)$. Moreover,
$$
\tilde M'(0):=\lim_{\varepsilon\to 0^+}\frac{\tilde M(\varepsilon)-\tilde M(0)}{\varepsilon}\ge \lim_{\varepsilon\to 0^+}\frac{\varepsilon^{-\delta}}{C\varepsilon^{-\eta}}=+\infty.
$$
Therefore, there exists $\varepsilon_0>0$ such that 
$\tilde M'(\varepsilon)>0$ for any $\varepsilon\in (0,\varepsilon_0)$.
That is, 
$$
\tilde M'(\varepsilon)=\frac{(1-\delta)\varepsilon^{-\delta}M(\frac{1}{\varepsilon})+\varepsilon^{-1-\delta}M'(\frac{1}{\varepsilon})}{M(\frac{1}{\varepsilon})^2}>0,
$$
which together with $\delta>1$ implies that
$$
M'(\frac{1}{\varepsilon})>(\delta-1)\varepsilon M(\frac{1}{\varepsilon})>0.
$$
Take $\varepsilon_\infty=\frac{1}{\varepsilon_0}$, then $M'(\varepsilon)>0$ for $\varepsilon\in (\varepsilon_\infty,\infty)$. The proof of (4) is complete. 
\end{proof}

\textbf{Organization of the paper}. In Section 2, we state  the main results in this paper. In Section 3, we give  some preliminary results which are needed in the proof of our main results. Sections 4--6 are devoted to the proofs of Theorems 2.1--2.3, respectively.  Finally, in the last section, we prove Theorem 2.4--2.5, and present some further results.

\smallskip

\textbf{Basic notations}. Throughout this paper, we assume $N\geq
3$. $ C_c^{\infty}(\mathbb{R}^N)$ denotes the space of the functions
infinitely differentiable with compact support in $\mathbb{R}^N$.
$L^p(\mathbb{R}^N)$ with $1\leq p<\infty$ denotes the Lebesgue space
with the norms
$\|u\|_p=\left(\int_{\mathbb{R}^N}|u|^p\right)^{1/p}$.
 $ H^1(\mathbb{R}^N)$ is the usual Sobolev space with norm
$\|u\|_{H^1(\mathbb{R}^N)}=\left(\int_{\mathbb{R}^N}|\nabla
u|^2+|u|^2\right)^{1/2}$. $ D^{1,2}(\mathbb{R}^N)=\{u\in
L^{2^*}(\mathbb{R}^N): |\nabla u|\in
L^2(\mathbb{R}^N)\}$. $H_r^1(\mathbb{R}^N)=\{u\in H^1(\mathbb{R}^N):
u\  \mathrm{is\ radially \ symmetric}\}$.  $B_r$ denotes the ball in $\mathbb R^N$ with radius $r>0$ and centered at the origin,  $|B_r|$ and $B_r^c$ denote its Lebesgue measure  and its complement in $\mathbb R^N$, respectively.  As usual, $C$, $c$, etc., denote generic positive constants.
 For any  small $\epsilon>0$ and two nonnegative functions $f(\epsilon)$ and  $g(\epsilon)$, we write

(1)  $f(\epsilon)\lesssim g(\epsilon)$ or $g(\epsilon)\gtrsim f(\epsilon)$ if there exists a positive constant $C$ independent of $\epsilon$ such that $f(\epsilon)\le Cg(\epsilon)$.

(2) $f(\epsilon)\sim g(\epsilon)$ if $f(\epsilon)\lesssim g(\epsilon)$ and $f(\epsilon)\gtrsim g(\epsilon)$.

\smallskip
\noindent
If $|f(\epsilon)|\lesssim |g(\epsilon)|$, we write $f(\epsilon)=O((g(\epsilon))$. We also denote by $\Theta=\Theta(\epsilon)$ a  generic  positive function satisfying
$
C_1\epsilon \le \Theta(\epsilon)\le C_2\epsilon
$
for some positive numbers $C_1,C_2>0$, which are independent of $\epsilon$. Finally, if $\lim f(\epsilon)/g(\epsilon)=1$ as $\epsilon\to \epsilon_0$, then we write $f(\epsilon)\simeq g(\epsilon)$ as $\epsilon\to \epsilon_0$.

\section*{2. Main Results}

Let 
$$
v(x)=\varepsilon^{s}u(\varepsilon^{t}x).
$$
Then it is easy to see that if $s=-\frac{2+\alpha}{4(p-1)}, t=-\frac{1}{2}$, $(P_\varepsilon)$ is reduced to 
$$
-\Delta v+v=(I_\alpha\ast |v|^p)|v|^{p-2}v+\varepsilon^{\Lambda_1}|v|^{q-2}v,
$$
where 
$$
\Lambda_1=\frac{q(2+\alpha)-2(2p+\alpha)}{4(p-1)}=\left\{\begin{array}{rcl}
>0,   \quad {\rm  if}  \   \    q>\frac{2(2p+\alpha)}{2+\alpha},\\
=0,  \quad {\rm if}   \    \   q=\frac{2(2p+\alpha)}{2+\alpha},\\
<0,   \quad  {\rm if}  \    \  q<\frac{2(2p+\alpha)}{2+\alpha}.
\end{array}\right.
$$

If  $s=-\frac{1}{q-2}, t=-\frac{1}{2}$, $(P_\varepsilon)$ is reduced to 
$$
-\Delta v+v=\varepsilon^{\Lambda_2}(I_\alpha\ast |v|^p)|v|^{p-2}v+|v|^{q-2}v,
$$
where 
$$
\Lambda_2=\frac{2(2p+\alpha)-q(2+\alpha)}{2(q-2)}=\left\{\begin{array}{rcl}
<0,   \quad {\rm if}  \   \    q>\frac{2(2p+\alpha)}{2+\alpha},\\
=0,  \quad {\rm if}   \    \   q=\frac{2(2p+\alpha)}{2+\alpha},\\
>0,   \quad {\rm if}  \    \  q<\frac{2(2p+\alpha)}{2+\alpha}.
\end{array}\right.
$$

Motivated by this, in the present paper, we consider  the following  equations
$$
-\Delta v+v=(I_\alpha \ast |v|^{p})|v|^{p-2}v+\lambda |v|^{q-2}v, 
\quad {\rm in} \ \mathbb R^N,
 \eqno(Q_\lambda)
$$
and 
$$
-\Delta v+v=\mu(I_\alpha \ast |v|^{p})|v|^{p-2}v+|v|^{q-2}v, 
\quad {\rm in} \ \mathbb R^N,
 \eqno(Q_\mu)
$$
where $p\in [\frac{N+\alpha}{N}, \frac{N+\alpha}{N-2}]$, $q\in (2,2^*]$ and $\lambda, \mu>0$ are  parameters. 

The corresponding  functionals are defined by 
$$
I_\lambda(v):=\frac{1}{2}\int_{\mathbb R^N}|\nabla v|^2+|v|^2-\frac{1}{2p}\int_{\mathbb R^N}(I_\alpha\ast |v|^p)|v|^p-\frac{\lambda}{q}\int_{\mathbb R^N}|v|^q
$$
and 
$$
I_\mu(v):=\frac{1}{2}\int_{\mathbb R^N}|\nabla v|^2+|v|^2-\frac{\mu}{2p}\int_{\mathbb R^N}(I_\alpha\ast |v|^p)|v|^p-\frac{1}{q}\int_{\mathbb R^N}|v|^q,
$$
respectively. The energy of the ground states given by
$$
m_\lambda:=\inf_{v\in \mathcal M_\lambda}I_\lambda(v)
\quad {\rm and } \quad 
m_\mu:=\inf_{v\in \mathcal M_\mu}I_\mu(v)
$$
are well-defined and positive, where $\mathcal M_\lambda$ and $\mathcal M_\mu$ denote the correspoding Nehari manifolds 
$$
\mathcal M_\lambda:=\left\{ v\in H^1(\mathbb R^N)\setminus\{0\}  \ \left | \ \int_{\mathbb R^N}|\nabla v|^2+|v|^2=\int_{\mathbb R^N}(I_\alpha\ast |v|^p)|v|^p+\lambda\int_{\mathbb R^N}|v|^q \right. \right\},
$$
$$
\mathcal M_\mu:=\left\{ v\in H^1(\mathbb R^N)\setminus\{0\}  \ \left | \ \int_{\mathbb R^N}|\nabla v|^2+|v|^2=\mu\int_{\mathbb R^N}(I_\alpha\ast |v|^p)|v|^p+\int_{\mathbb R^N}|v|^q \right. \right\}.
$$

\vskip 5mm 

The ground state solutions of $(Q_\lambda)$ and $(Q_\mu)$ are denoted by $v_\lambda$ and $v_\mu$ respectively. The existence of these kind solutions are proved in \cite{Li-2, Li-1}.  More precisely,  the following theorems are proved in \cite{Li-1}.

{\bf Theorem A.} {\it 
Let $N\geq 3$,\ $\alpha\in(0,N)$, $p=\frac{N+\alpha}{N}$ and
$\lambda>0$. Then there is a constant $\lambda_0>0$ such that $(Q_\lambda)$ admits a positive ground state
$v_\lambda\in H^1(\mathbb{R}^N)$ which is radially symmetric and radially
nonincreasing if one of the following conditions holds:

(1) $q\in(2,2+\frac{4}{N})$;

(2) $q\in[2+\frac{4}{N},\frac{2N}{N-2})$ and $\lambda>\lambda_0$.
}

{\bf Theorem B.} {\it 
 Let $N\geq 3$,\ $\alpha\in(0,N)$, $p=\frac{N+\alpha}{N-2}$ and
$\lambda>0$. Then there is a constant $\lambda_1>0$ such that $(Q_\lambda)$ admits a positive ground state
$v_\lambda\in H^1(\mathbb{R}^N)$ which is radially symmetric and radially
nonincreasing if one of the following conditions holds:

(1) $N\geq 4$ and $q\in (2,\frac{2N}{N-2})$;

(2) $N=3$ and  $q\in (4,\frac{2N}{N-2})$;

(3) $N=3$,  $q\in (2,4]$ and $\lambda>\lambda_1$.
}

{\bf Theorem C. } {\it 
Let $N\geq 3$, $\alpha\in(0,N)$, $q=2^*$ and $\mu>0$. Then there are two  constants $\mu_0, \mu_1>0$ such that $(Q_\mu)$
admits a positive ground state  $v_\mu\in H^1(\mathbb{R}^N)$ which is
radially symmetric and radially nonincreasing if one of the
following conditions holds:

(1) $N\geq 4$ and $p\in
(1+\frac{\alpha}{N-2},\frac{N+\alpha}{N-2})$;

(2) $N\geq 4$, $p\in (\frac{N+\alpha}{N},1+\frac{\alpha}{N-2}]$ and
$\mu>\mu_0$;

(3) $N=3$ and $p\in (2+\alpha,\frac{N+\alpha}{N-2})$;

(4) $N=3$, $p\in (\frac{N+\alpha}{N},2+\alpha]$ and $\mu>\mu_1$.
}

\vskip 5mm 

When $p=\frac{N+\alpha}{N}$ and $\lambda>0$ is small, we are going to show that after a suitable rescaling the limit equation for $(Q_\lambda)$ is given by the critical equation
$$
 U=(I_\alpha\ast |U|^{\frac{N+\alpha}{N}})U^{\frac{\alpha}{N}} \quad {\rm in } \   \ \mathbb R^N.
\eqno(2.1)
$$
It is well-known that the radial ground states of $(2.1)$ are given by the function 
$$
U_1(x):=\left(\frac{A}{1+|x|^2}\right)^{\frac{N}{2}}
\eqno(2.2)
$$
and the family of its rescalings 
$$
U_\rho(x):=\rho^{-\frac{N}{2}}U_1(x/\rho),  \quad \rho>0.
\eqno(2.3)
$$
\vskip 5mm 

{\bf Theorem 2.1.} {\it If  $p=\frac{N+\alpha}{N}, q\in (2, 2+\frac{4}{N})$, and  $\{v_\lambda\}$ is  a family of radial ground states of $(Q_\lambda)$,
 then for small $\lambda>0$
 $$
v_\lambda(0)\sim \lambda^{\frac{N}{4-N(q-2)}},
$$
$$
\|\nabla v_\lambda\|^2_2\sim \lambda^{\frac{4}{4-N(q-2)}}, \quad \|v_\lambda\|^q_q\sim \lambda^{\frac{N(q-2)}{4-N(q-2)}}, 
$$
$$
\int_{\mathbb R^N}(I_\alpha\ast |v_\lambda|^{\frac{N+\alpha}{N}})|v_\lambda|^{\frac{N+\alpha}{N}}=S_1^{\frac{N+\alpha}{\alpha}}+ O(\lambda^{\frac{4}{4-N(q-2)}}),
$$ 
$$
\|v_\lambda\|^2_2=S_1^{\frac{N+\alpha}{\alpha}}+ O(\lambda^{\frac{4}{4-N(q-2)}}),
$$
Moreover, there exists $\zeta_\lambda\in (0,+\infty)$
verifying  
$$
\zeta_\lambda \sim \lambda^{-\frac{2}{4-N(q-2)}}
$$
such that for small $\lambda>0$, the rescaled family of ground states
$$
 w_\lambda(x)=\zeta_\lambda^{\frac{N}{2}}v_\lambda(\zeta_\lambda x)
$$
satisfies 
$$
\|\nabla  w_\lambda\|^2_2\sim \| w_\lambda\|_{q}^{q} \sim \int_{\mathbb R^N}(I_\alpha \ast |w_\lambda|^{\frac{N+\alpha}{N}})| w_\lambda|^{\frac{N+\alpha}{N}}\sim \| w_\lambda\|_2^2\sim 1,
$$
and as $\lambda\to 0$, $w_\lambda$ converges in $H^1(\mathbb R^N)$ to the extremal function  $U_{\rho_0}$, where 
$$
\rho_0=\left(\frac{2q\int_{\mathbb R^N}|\nabla U_1|^2}{N(q-2)\int_{\mathbb R^N}|U_1|^q}\right)^{\frac{2}{4-N(q-2)}}.
\eqno(2.4)
$$
Furthermore,  the least energy $m_\lambda$ of the ground sate  satisfies 
$$
\frac{\alpha}{2(N+\alpha)}S_1^{\frac{N+\alpha}{\alpha}}-m_\lambda\sim \lambda^{\frac{4}{4-N(q-2)}},  
$$
as $\lambda\to 0$,  where 
$$
S_1=\inf_{u\in H^1(\mathbb{R}^N)\setminus\{0\}}\frac{\int_{\mathbb{R}^N}| u|^2}{(\int_{\mathbb{R}^N}(I_\alpha\ast |u|^{\frac{N+\alpha}{N}})|u|^{\frac{N+\alpha}{N}})^{\frac{N}{N+\alpha}}}.
\eqno(2.5)
$$}

\vskip 5mm

When $p=\frac{N+\alpha}{N-2}$ and $\lambda>0$ is small, we are going to show that after a suitable rescaling the limit equation for $(Q_\lambda)$ is given by the critical 
Choquard equation
$$
-\Delta V=(I_\alpha\ast |V|^{\frac{N+\alpha}{N-2}})V^{\frac{2+\alpha}{N-2}} \quad {\rm in } \   \ \mathbb R^N.
\eqno(2.6)
$$
It is well-known that the radial ground states of $(2.6)$ are given by the function 
$$
V_1(x):=[N(N-2)]^{\frac{N-2}{4}}\left(\frac{1}{1+|x|^2}\right)^{\frac{N-2}{2}}
\eqno(2.7)
$$
and the family of its rescalings 
$$
V_\rho(x):=\rho^{-\frac{N-2}{2}}V_1(x/\rho),  \quad \rho>0.
\eqno(2.8)
$$
\vskip 5mm

{\bf Theorem 2.2.} {\it Let $p=\frac{N+\alpha}{N-2}$ and $q\in (2, 2^*)$, and   $\{v_\lambda\}$ be a family of radial ground states of $(Q_\lambda)$.
 If $N\ge 5$, then for small $\lambda>0$
 $$
v_\lambda(0)\sim \lambda^{-\frac{1}{q-2}},
$$
$$
 \|v_\lambda\|_q^q\sim \lambda^{\frac{2N-q(N-2)}{(N-2)(q-2)}}, \quad \|v_\lambda\|^2_2\sim \lambda^{\frac{4}{(N-2)(q-2)}}, 
$$
$$
\|\nabla v_\lambda\|^2_2=S_\alpha^{\frac{N+\alpha}{2+\alpha}}+ O(\lambda^{\frac{4}{(N-2)(q-2)}}),
$$
$$
\int_{\mathbb R^N}(I_\alpha\ast |v_\lambda|^{\frac{N+\alpha}{N-2}})|v_\lambda|^{\frac{N+\alpha}{N-2}}=S_\alpha^{\frac{N+\alpha}{2+\alpha}}+ O(\lambda^{\frac{4}{(N-2)(q-2)}}).
$$
Moreover,  there exists $\zeta_\lambda\in (0,+\infty)$
verifying  
$$
\zeta_\lambda\sim \lambda^{\frac{2}{(N-2)(q-2)}},  
$$
such that for small $\lambda>0$, the rescaled family of ground states
$$
w_\lambda(x)=\zeta_\lambda^{\frac{N-2}{2}}v_\lambda(\zeta_\lambda x)
$$
satisfies 
$$
\|\nabla w_\lambda\|^2_2\sim \|w_\lambda\|_{q}^{q} \sim \int_{\mathbb R^N}(I_\alpha \ast |w_\lambda|^{\frac{N+\alpha}{N-2}})|w_\lambda|^{\frac{N+\alpha}{N-2}}\sim \|w_\lambda\|_2^2\sim 1,
$$
and as $\lambda\to 0$, $w_\lambda$ converges in $H^1(\mathbb R^N)$ to $V_{\rho_0}$ with 
 $$
 \rho_0=\left(\frac{2(2^*-q)\int_{\mathbb R^N}|V_1|^q}{q(2^*-2)\int_{\mathbb R^N}|V_1|^2}\right)^ {\frac{2}{(N-2)(q-2)}}.
 \eqno(2.9)
 $$
 In the lower dimension cases, we assume that $q\in (2, 4)$ if $N=4$, and $q\in (4, 6)$ if $N=3$, then for small $\lambda>0$
 $$
v_\lambda(0)  \sim\left\{\begin{array}{rcl}
\lambda^{-\frac{1}{q-2}}(\ln\frac{1}{\lambda})^{\frac{1}{q-2}},   \  \quad if \   \ N=4,\\
\lambda^{-\frac{1}{q-4}},  \quad \qquad \qquad if  \    \  N=3,
\end{array}\right.
$$
$$
 \|v_\lambda\|_q^q\sim\left\{\begin{array}{rcl}
\lambda^{\frac{4-q}{q-2}}(\ln\frac{1}{\lambda})^{-\frac{4-q}{q-2}},   \  \quad if \   \ N=4,\\
\lambda^{\frac{6-q}{q-4}},  \qquad \quad \qquad if  \    \  N=3,
\end{array}\right.
$$
$$
\|v_\lambda\|_2^2\sim\left\{\begin{array}{rcl}
\lambda^{\frac{2}{q-2}}(\ln\frac{1}{\lambda})^{-\frac{4-q}{q-2}},   \  \quad if \   \ N=4,\\
\lambda^{\frac{2}{q-4}},  \quad \qquad \qquad if  \    \  N=3,
\end{array}\right.
$$
$$
\|\nabla v_\lambda\|^2_2=S_\alpha^{\frac{N+\alpha}{2+\alpha}}+ \left\{\begin{array}{rcl} 
O(\lambda^{\frac{2}{q-2}}(\ln\frac{1}{\lambda})^{-\frac{4-q}{q-2}}),   \quad  if  \   \ N=4,\\
O( \lambda^{\frac{2}{q-4}}), \quad  \qquad  \qquad  if \    \  N=3,
 \end{array}\right.
$$
$$ 
\int_{\mathbb R^N}(I_\alpha\ast |v_\lambda|^{\frac{N+\alpha}{N-2}})|v_\lambda|^{\frac{N+\alpha}{N-2}}=S_\alpha^{\frac{N+\alpha}{2+\alpha}}+ \left\{\begin{array}{rcl} 
O(\lambda^{\frac{2}{q-2}}(\ln\frac{1}{\lambda})^{-\frac{4-q}{q-2}}),   \quad  if  \   \ N=4,\\
O( \lambda^{\frac{2}{q-4}}), \quad  \qquad  \qquad  if \    \  N=3,
 \end{array}\right.
$$
and there exists $\zeta_\lambda\in (0,+\infty)$
verifying  
$$
\zeta_\lambda \sim\left\{\begin{array}{rcl}
\lambda^{\frac{1}{q-2}}(\ln\frac{1}{\lambda})^{-\frac{1}{q-2}},   \  \quad if \   \ N=4,\\
\lambda^{\frac{2}{q-4}},  \qquad \quad \qquad if  \    \  N=3,
\end{array}\right.
$$
such that for small $\lambda>0$, the rescaled family of ground states
$$
w_\lambda(x)=\zeta_\lambda^{\frac{N-2}{2}}v_\lambda(\zeta_\lambda x)
$$
satisfies 
$$
\|\nabla w\|_2^2\sim \|w_\lambda\|_{q}^{q}\sim \int_{\R^N}(I_\alpha\ast |w_\lambda|^{\frac{N+\alpha}{N-2}})|w_\lambda|^{\frac{N+\alpha}{N-2}}\sim 1, 
$$
$$
  \|w_\lambda\|_2^2\sim \left\{\begin{array}{rcl}
\ln\frac{1}{\lambda},   \quad  \quad if \   \ N=4,\\
\lambda^{-\frac{2}{q-4}},  \quad if  \    \  N=3,
\end{array}\right.
$$
and as $\lambda\to 0$, $w_\lambda$ converges in $D^{1,2}(\mathbb R^N)$ and $L^{q}(\mathbb R^N)$  to $V_1$.
Furthermore, the least energy $m_\lambda$ of the ground state satisfies 
$$
\frac{2+\alpha}{2(N+\alpha)}S_\alpha^{\frac{N+\alpha}{2+\alpha}}-m_\lambda\sim \left\{\begin{array}{rcl} \lambda^{\frac{4}{(N-2)(q-2)}}, \qquad  \quad  if \   \ N\ge 5,\\
\lambda^{\frac{2}{q-2}}(\ln\frac{1}{\lambda})^{-\frac{4-q}{q-2}},   \quad  if \   \ N=4,\\
 \lambda^{\frac{2}{q-4}}, \quad  \qquad  \qquad  if  \    \  N=3,
 \end{array}\right.
$$
as $\lambda\to 0$, where 
$$
S_\alpha:=\inf_{v\in D^{1,2}(\mathbb R^N)\setminus\{0\}}\frac{\int_{\mathbb R^N}|\nabla v|^2}{\left(\int_{\mathbb R^N}(I_\alpha\ast |v|^{\frac{N+\alpha}{N-2}})|v|^{\frac{N+\alpha}{N-2}}\right)^{\frac{N-2}{N+\alpha}}}.
\eqno(2.10)
$$}
\vskip 5mm 

{\bf Remark 2.1.}  If $N\ge 5$ and $\alpha>N-4$, we can choose 
$
\zeta_\lambda= \lambda^{\frac{2}{(N-2)(q-2)}}
$
in Theorem 2.2.
\vskip 5mm 

When $q=2^*$ and $\mu>0$ is small, we are going to show that after a suitable rescaling the limit equation for $(Q_\mu)$ is given by the critical Emden-Fowler equation
$$
-\Delta W=W^{2^*-1} \quad {\rm in } \   \ \mathbb R^N.
\eqno(2.11)
$$
It is well-known that the radial ground states of $(2.11)$ are given by the function 
$$
W_1(x):=[N(N-2)]^{\frac{N-2}{4}}\left(\frac{1}{1+|x|^2}\right)^{\frac{N-2}{2}}
\eqno(2.12)
$$
and the family of its rescalings 
$$
W_\rho(x):=\rho^{-\frac{N-2}{2}}W_1(x/\rho),  \quad \rho>0.
\eqno(2.13)
$$
\vskip 5mm 

{\bf Theorem 2.3.} {\it Let $q=2^*, \ p\in (1+\frac{\alpha}{N-2}, \frac{N+\alpha}{N-2})$, and   $\{v_\mu\}$ be a family of radial ground states of $(Q_\mu)$.
 If $N\ge 5$, then for small $\mu>0$
 $$
v_\mu(0)\sim \mu^{-\frac{N-2}{2(N-2)(p-1)-2\alpha}},
$$
$$
\int_{\mathbb R^N}(I_\alpha\ast |v_\mu|^p)|v_\mu|^p\sim \mu^{\frac{N+\alpha-p(N-2)}{(N-2)(p-1)-\alpha}}, \quad  \|v_\mu\|_2^2\sim \mu^{\frac{2}{(N-2)(p-1)-\alpha}}, 
$$
$$
\|\nabla v_\mu\|_2^2=S^{\frac{N}{2}}+O(\mu^{\frac{2}{(N-2)(p-1)-\alpha}}), 
\quad
\| v_\mu\|_{2^*}^{2^*}=S^{\frac{N}{2}}+O(\mu^{\frac{2}{(N-2)(p-1)-\alpha}}).
$$
Moreover,  there exists $\zeta_\mu\in (0,+\infty)$
verifying  
$$
\zeta_\mu \sim \mu^{\frac{1}{(N-2)(p-1)-\alpha}},
$$
such that for small $\mu>0$, the rescaled family of ground states
$$
w_\mu(x)=\zeta_\mu^{\frac{N-2}{2}}v_\mu(\zeta_\mu x)
$$
satisfies 
$$
\|\nabla w_\mu\|^2_2\sim \|w_\mu\|_{2^*}^{2^*} \sim \int_{\mathbb R^N}(I_\alpha \ast |w_\mu|^p)|w_\mu|^p\sim \|w_\mu\|_2^2\sim 1,
$$
and as $\mu\to 0$, $w_\mu$ converges in $H^1(\mathbb R^N)$ to $W_{\rho_0}$ with 
 $$
 \rho_0=\left(\frac{[(N+\alpha)-p(N-2)]\int_{\mathbb R^N}(I_\alpha\ast |W_1|^p)|W_1|^p}{2p\int_{\mathbb R^N}|W_1|^2}\right)^ {\frac{1}{(N-2)(p-1)-\alpha}}.
 \eqno(2.14)
 $$
 In the lower dimension cases, we assume that $p\in (\max\{2, 1+\frac{\alpha}{2}\}, 2+\frac{\alpha}{2})$ if $N=4$, and $p\in (2+\alpha, 3+\alpha)$ if $N=3$, 
  then for small $\mu>0$
 $$
v_\mu(0)  \sim\left\{\begin{array}{rcl}
\mu^{-\frac{1}{2p-2-\alpha}}(\ln\frac{1}{\mu})^{\frac{1}{2p-2-\alpha}},   \  \quad if \   \ N=4,\\
\mu^{-\frac{1}{2(p-2-\alpha)}},  \qquad \qquad \qquad if  \    \  N=3,
\end{array}\right.
$$
$$
\int_{\mathbb R^N}(I_\alpha\ast |v_\mu|^p)|v_\mu|^p\sim\left\{\begin{array}{rcl}
\mu^{\frac{4+\alpha-2p}{2p-2-\alpha}}(\ln\frac{1}{\mu})^{-\frac{4+\alpha-2p}{2p-2-\alpha}},   \  \quad if \   \ N=4,\\
\mu^{\frac{3+\alpha-p}{p-2-\alpha}},  \qquad \qquad \qquad if  \    \  N=3,
\end{array}\right.
$$
$$
 \|v_\mu\|_2^2\sim\left\{\begin{array}{rcl}
\mu^{\frac{2}{2p-2-\alpha}}(\ln\frac{1}{\mu})^{-\frac{4+\alpha-2p}{2p-2-\alpha}},   \  \quad if \   \ N=4,\\
\mu^{\frac{1}{p-2-\alpha}},  \qquad \qquad \qquad if  \    \  N=3,
\end{array}\right.
$$
$$
\|\nabla v_\mu\|_2^2=S^{\frac{N}{2}}+\left\{\begin{array}{rcl} 
O(\mu^{\frac{2}{2p-2-\alpha}}(\ln\frac{1}{\mu})^{-\frac{4+\alpha-2p}{2p-2-\alpha}}),   \quad  if \   \ N=4,\\
 O( \mu^{\frac{1}{p-2-\alpha}}), \  \quad \quad  \qquad  \qquad   if  \    \  N=3,
 \end{array}\right. 
$$
$$
\| v_\mu\|_{2^*}^{2^*}=S^{\frac{N}{2}}+\left\{\begin{array}{rcl} 
O(\mu^{\frac{2}{2p-2-\alpha}}(\ln\frac{1}{\mu})^{-\frac{4+\alpha-2p}{2p-2-\alpha}}),   \quad if \   \ N=4,\\
 O( \mu^{\frac{1}{p-2-\alpha}}), \  \quad \quad  \qquad  \qquad  if  \    \  N=3.
 \end{array}\right. 
 $$
and  there exists $\zeta_\mu\in (0,+\infty)$
verifying  
$$
\zeta_\mu \sim\left\{\begin{array}{rcl}
\mu^{\frac{1}{2p-2-\alpha}}(\ln\frac{1}{\mu})^{-\frac{1}{2p-2-\alpha}},   \  \quad if \   \ N=4,\\
\mu^{\frac{1}{p-2-\alpha}},  \qquad \qquad \qquad if  \    \  N=3,
\end{array}\right.
$$
such that for small $\mu>0$, the rescaled family of ground states
$$
w_\mu(x)=\zeta_\mu^{\frac{N-2}{2}}v_\mu(\zeta_\mu x)
$$
satisfies 
$$
\|\nabla w_\mu\|^2_2\sim \|w_\mu\|_{2^*}^{2^*}\sim \int_{\mathbb R^N}(I_\alpha\ast |w_\mu|^p)|w_\mu|^p\sim 1, 
$$ 
$$
\|w_\mu\|_2^2\sim \left\{\begin{array}{rcl}
\ln\frac{1}{\mu},   \quad  \quad if \   \ N=4,\\
\mu^{-\frac{1}{p-2-\alpha}},  \quad if  \    \  N=3,
\end{array}\right.
$$
and as $\mu\to 0$, $w_\mu$ converges in $D^{1,2}(\mathbb R^N)$ and $L^{\frac{2Np}{N+\alpha}}(\mathbb R^N)$  to $W_1$.
Furthermore, the least energy $m_\mu$ of the ground state satisfies 
$$
\frac{1}{N}S^{\frac{N}{2}}-m_\mu \sim \left\{\begin{array}{rcl} \mu^{\frac{2}{(N-2)(p-1)-\alpha}}, \   \quad  \qquad  \quad  if  \   \ N\ge 5,\\
\mu^{\frac{2}{2p-2-\alpha}}(\ln\frac{1}{\mu})^{-\frac{4+\alpha-2p}{2p-2-\alpha}},   \quad  if  \   \ N=4,\\
  \mu^{\frac{1}{p-2-\alpha}}, \  \quad \quad  \qquad  \qquad  if \    \  N=3,
 \end{array}\right. 
$$
as $\mu\to 0$, where 
$$
S:=\inf_{w\in D^{1,2}(\mathbb R^N)\setminus \{0\}}\frac{\int_{\mathbb R^N}|\nabla w|^2}{\left(\int_{\mathbb R^N}|w|^{2^*}\right)^{\frac{2}{2^*}}}.
\eqno(2.15)
$$}
\vskip 5mm 

{\bf Remark 2.2.}  If $N\ge 5$ and $p>\max\{2, 1+\frac{\alpha}{N-2}\}$, we can choose $\zeta_\mu = \mu^{\frac{1}{(N-2)(p-1)-\alpha}}$ in Theorem 2.3.  For $N=3$ and $p\in (2+\alpha, 3+\alpha)$, we obtain the exact asymptotic behaviour of ground states as $\mu\to 0$.

\begin{figure}
	\centering
	\includegraphics[width=0.6\linewidth]{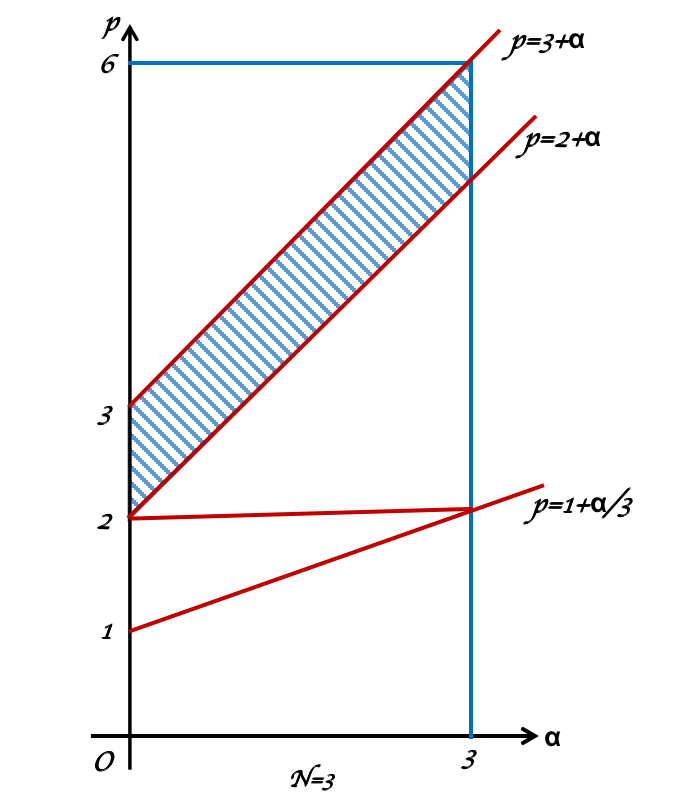}
	\caption{Admissible regions in  $(\alpha,p)$ plane for $N=3$.}
	\label{fig:figure2}
\end{figure}


For $N=4$ and $p\in (\max\{2, 1+\frac{\alpha}{2}\}, 2+\frac{\alpha}{2})$, we obtain the exact asymptotic behaviour of ground states as $\mu\to 0$. If $N= 4$, $\alpha<2$ and  $p\in (1+\frac{\alpha}{2}, 2]$,  the exact asymptotic behaviour 
of ground state solution is still open.

\begin{figure}
	\centering
	\includegraphics[width=0.6\linewidth]{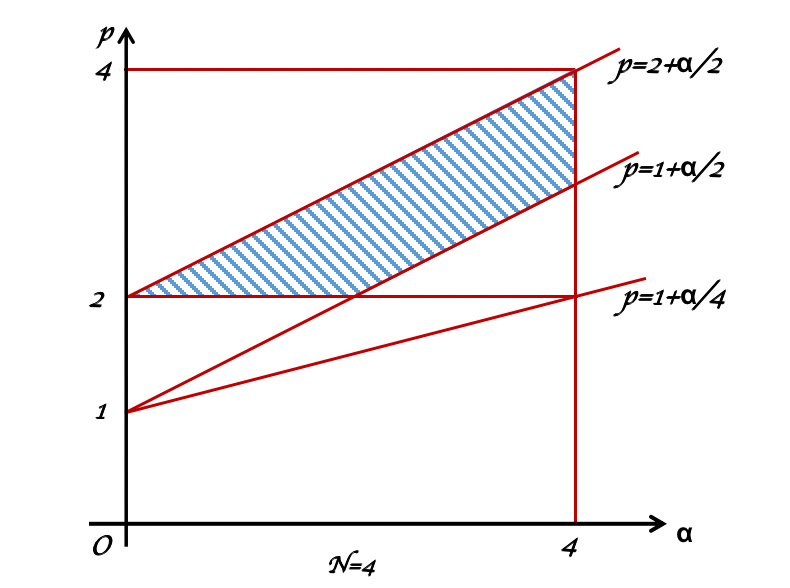}
	\caption{Admissible regions in  $(\alpha,p)$ plane for $N=4$.}
	\label{fig:figure3}
\end{figure}


For $N\ge 5$ and $p\in (1+\frac{\alpha}{N-2}, \frac{N+\alpha}{N-2})$, the exact asymptotic behaviour 
of ground states as $\mu\to 0$ are obtained.

\begin{figure}
	\centering
	\includegraphics[width=0.6\linewidth]{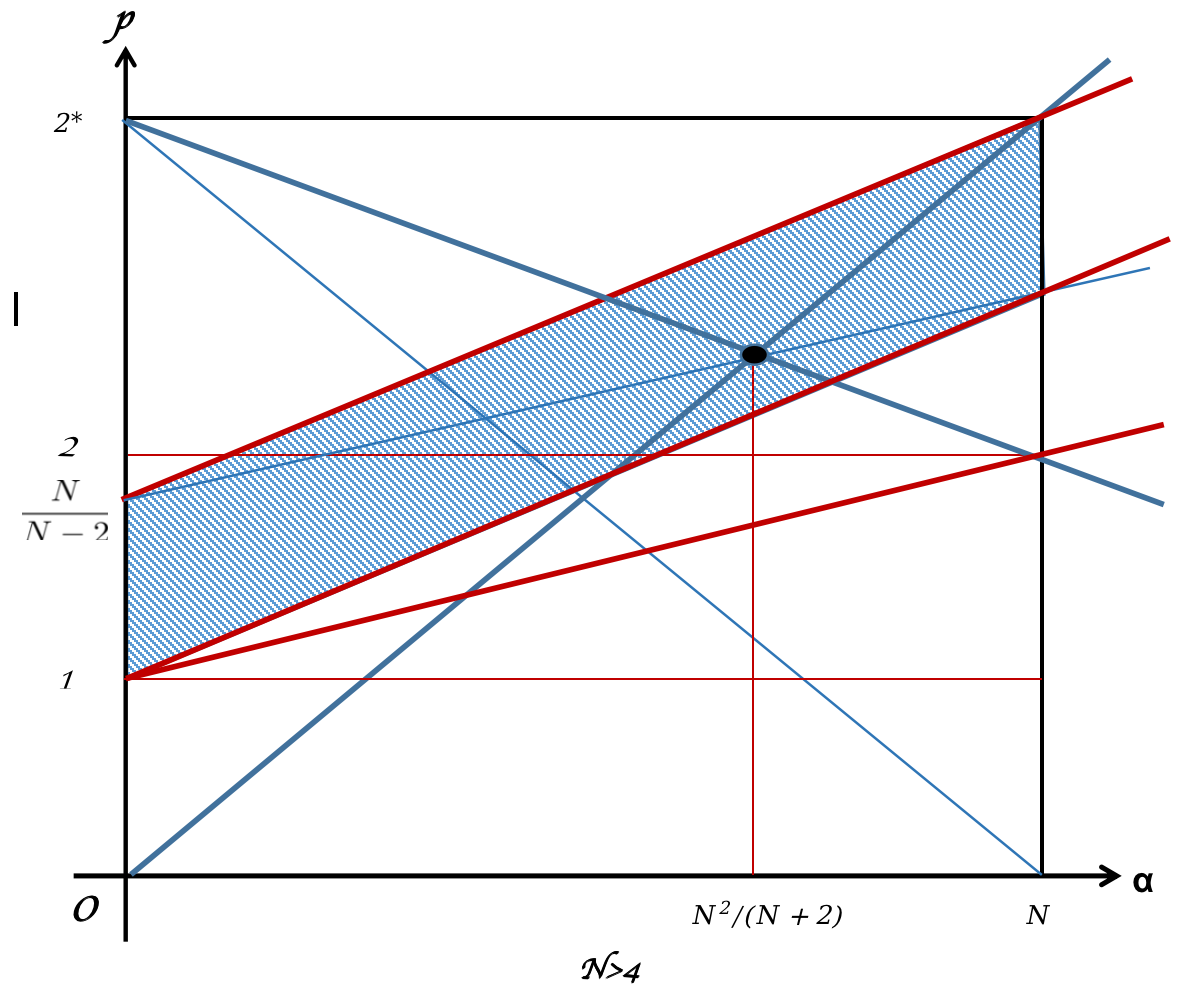}
	\caption{Admissible regions in  $(\alpha,p)$ plane for $N\ge 5$.}
	\label{fig:figure4}
\end{figure}


{\bf Theorem 2.4.} {\it  If $p\in (\frac{N+\alpha}{N}, \frac{N+\alpha}{N-2})$ and $q\in (2,2^*]$. Let $v_\lambda$ be  the radial ground state of $(Q_\lambda)$, then for any sequence $\lambda_n\to 0$, there exists a subsequence, still denoted by $\lambda_n$, such that $v_{\lambda_n}$ 
converges in $H^1(\mathbb R^N)$ to a positive solution $v_0\in H^1(\mathbb R^N)$ of the equation
$$
-\Delta v+v=(I_\alpha \ast |v|^p)v^{p-1}.
\eqno(2.16)
$$
Moreover, as $\lambda\to 0$, there holds
$$
\|v_\lambda\|_2^2=\frac{N+\alpha-p(N-2)}{2p}S_p^{\frac{p}{p-1}}+O(\lambda),  \quad if \  \ q< \frac{2(2p+\alpha)}{2+\alpha},
$$
$$
\|v_\lambda\|_2^2=\frac{N+\alpha-p(N-2)}{2p}S_p^{\frac{p}{p-1}}-\Theta(\lambda),   \quad if \  \ q\ge \frac{2(2p+\alpha)}{2+\alpha},
$$
$$
\|\nabla v_\lambda\|_2^2=\frac{N(p-1)-\alpha}{2p}S_p^{\frac{p}{p-1}}+O(\lambda),
$$
and  the least energy $m_\lambda$ of the ground state satisfies 
$$
\frac{p-1}{2p}S_p^{\frac{p}{p-1}}-m_\lambda\sim \lambda,
$$
as $\lambda\to 0$, where 
$$
S_p=\inf_{v\in H^1(\mathbb R^N)\setminus \{0\}}\frac{\int_{\mathbb R^N}|\nabla v|^2+|v|^2}{(\int_{\mathbb R^N}(I_\alpha\ast |v|^p)|v|^p)^{\frac{1}{p}}}.
\eqno(2.17)
$$
}
\vskip 5mm 

{\bf Theorem 2.5.} {\it  If $p\in [\frac{N+\alpha}{N}, \frac{N+\alpha}{N-2}]$ and $q\in (2,2^*)$. Let $v_\mu$ be  the radial ground state of $(Q_\mu)$, then as $\mu\to 0$, 
$v_\mu$ converges in $H^1(\mathbb R^N)$ to the unique positive solution $v_0\in H^1(\mathbb R^N)$ of the equation
$$
-\Delta v+v=v^{q-1}.
\eqno(2.18)
$$
Moreover, as $\mu\to 0$, there holds
$$
\|v_\mu\|_2^2=\frac{2N-q(N-2)}{2q}S_q^{\frac{q}{q-2}}+O(\mu), \quad if \  \ q> \frac{2(2p+\alpha)}{2+\alpha},
$$
$$
\|v_\mu\|_2^2=\frac{2N-q(N-2)}{2q}S_q^{\frac{q}{q-2}}-\Theta(\mu), \quad if \  \ q\le \frac{2(2p+\alpha)}{2+\alpha},
$$
$$
\|\nabla v_\mu\|_2^2=\frac{N(q-2)}{2q}S_q^{\frac{q}{q-2}}+O(\mu),
$$
$$
\|u_\mu\|_q^q\sim \int_{\mathbb R^N}(I_\alpha\ast |u_\mu|^p)|u_\mu|^p\sim 1,
$$
and the least energy $m_\mu$ of the ground state satisfies 
$$
\frac{q-2}{2q}S_q^{\frac{q}{q-2}}-m_\mu\sim\mu,
$$
as $\mu\to 0$, where 
$$
S_q=\inf_{v\in H^1(\mathbb R^N)\setminus \{0\}}\frac{\int_{\mathbb R^N}|\nabla v|^2+|v|^2}{(\int_{\mathbb R^N}|v|^q)^{\frac{2}{q}}}.
\eqno(2.19)
$$}
\vskip 5mm

{\bf Remark 2.3.} The contributions of the present paper are mainly as follows. 

(I) Firstly, we obtain optimal explicit rescaling in a sense that it is unique up to a multiplicative constant such that the rescaled family of ground states converges in $H^1(\mathbb R^N)$ to a particular solution of the limit equation. 
Such as the lower critical case in Theorem 2.1 for all $N\ge 3$, the upper critical case in Theorem 2.2  for $N\ge 5$ and the critical case in Theorem 2.3 for $N\ge 5$. 

(II) The present paper is inspired by \cite{Moroz-1},  where the second
author and C. Muratov studied the asymptotic properties of ground states for a combined
powers Schr\"odinger equations with a focusing exponent $p>2$ and a defocusing exponent $q>p$.   In \cite{Moroz-1}, by considering a constrained minimization problem, the authors obtain a precise estimate of least energy which implies the uniform boundedness of the rescaled family of ground states in $L^q(\mathbb R^N)$. However, in the present paper, the technique is totally different from that used in \cite{Moroz-1}. For example, in the proof of Theorem 2.2, we first use the Nehari  manifold and Poho\v zaev identity to obtain the uniform boundedness  of the rescaled family of ground states in $L^q(\mathbb R^N)$  and then give a precise estimate of least energy. 

(III)  In the proof of Theorem 2.3, we need an uniform optimal decay estimate of the rescaled family of ground states in the lower dimension cases $N=3$ and $N=4$.  This is a  delicate and difficult  task. 
\vskip 5mm

According to our main results, we draw the following figures which depict the limit equations of  $(P_\varepsilon)$  as $\varepsilon\to \infty$ and $\varepsilon\to 0$ when $(p,q)$ 
belongs to different regions in the the $(p,q)$ plane,  and in a sense
reveals the asymptotic behaviour of rescaled family of ground states to 
$(P_\varepsilon)$ as $\varepsilon\to \infty$ and $\varepsilon\to 0$, respectively.

\begin{figure}
	\centering
	\includegraphics[width=0.7\linewidth]{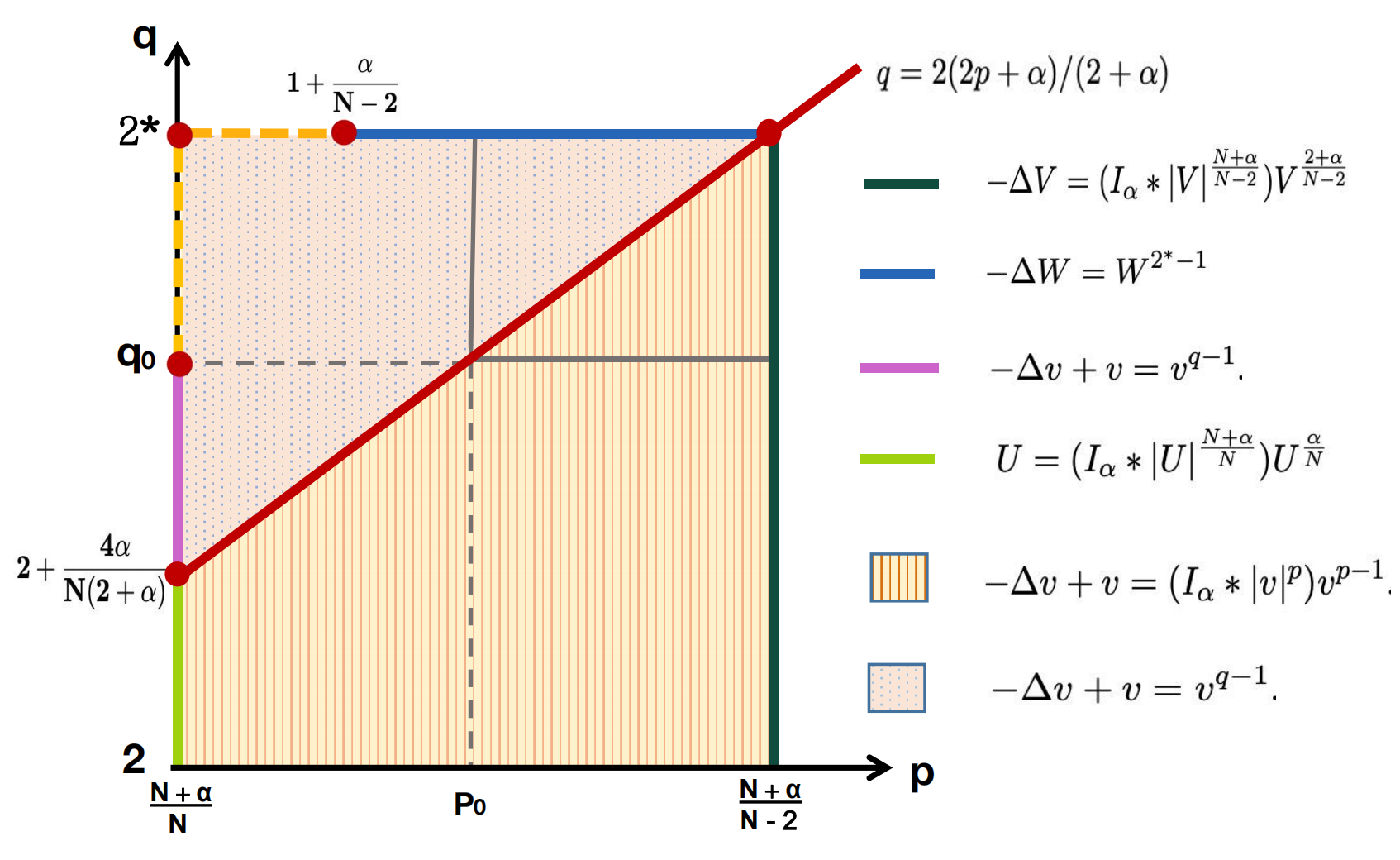}
	\caption{The limit equations of $(P_\varepsilon)$ as $\varepsilon\to \infty$ when $(p,q)$ belongs to different regions.}
	\label{fig:figure5}
\end{figure}

\begin{figure}
	\centering
	\includegraphics[width=0.7\linewidth]{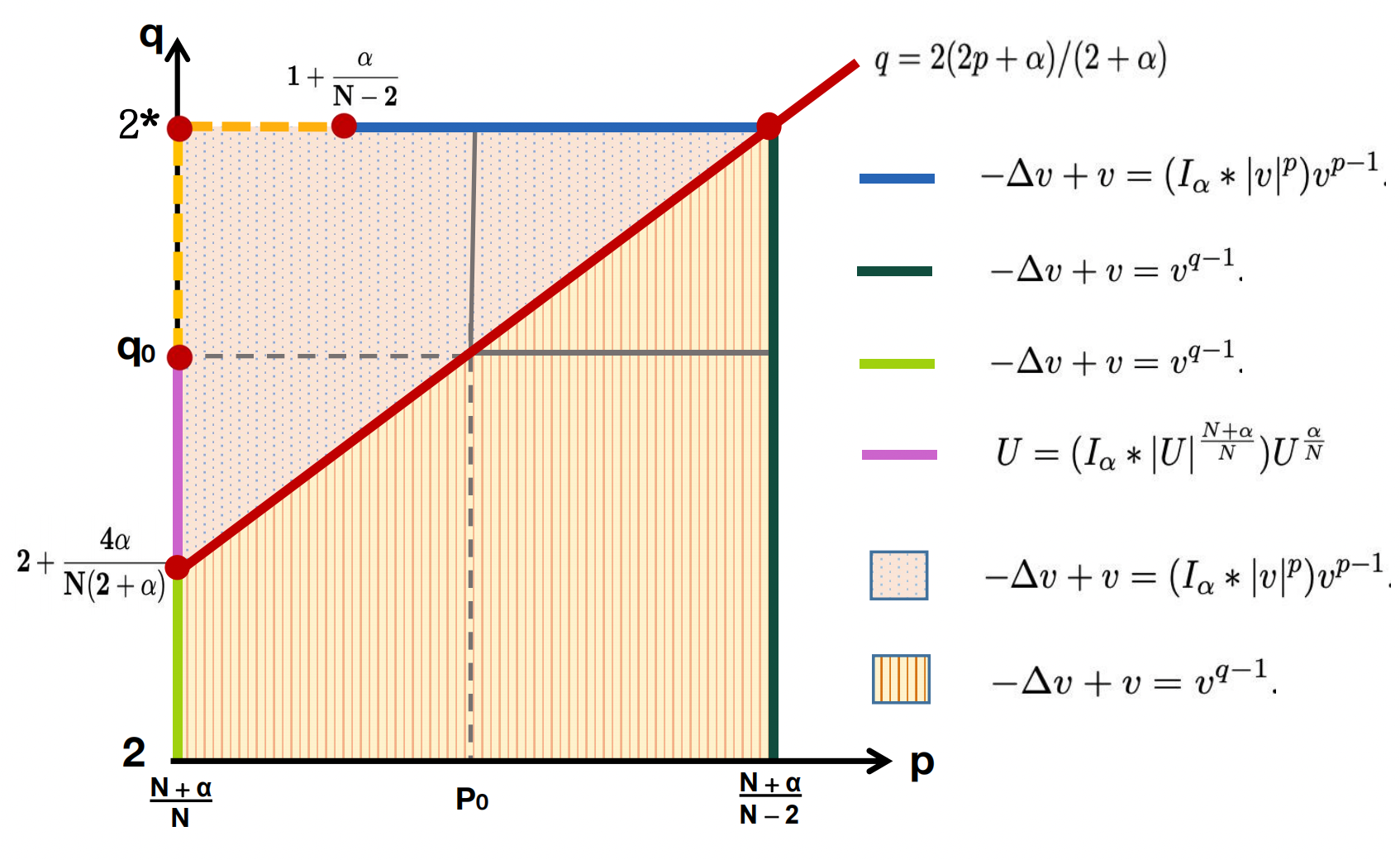}
	\caption{The limit equations of $(P_\varepsilon)$ as $\varepsilon\to 0$ when $(p,q)$ belongs to different regions.}
	\label{fig:figure6}
\end{figure}


The following results are direct consequences of Theorems 2.1--2.5.

{\bf Corollary 2.1.} {\it  If $p=\frac{N+\alpha}{N}$, then the problem  $(P_\varepsilon)$ admits a positive ground states $u_\varepsilon\in H^1(\mathbb R^N)$,  which are radially symmetric and radially nonincreasing. Furthermore, the following statements hold true:

$(I)$  If $q\in (2, 2+\frac{4\alpha}{N(2+\alpha)})$ and $\varepsilon\to 0 \ ( resp. \ q\in (2+\frac{4\alpha}{N(2+\alpha)}, 2+\frac{4}{N})$ and $\varepsilon\to \infty)$,  then 
the ground state solutions satisfy 
$$
u_\varepsilon(0)\sim \varepsilon^{\frac{1}{q-2}},
$$
$$
\|u_\varepsilon\|_q^q\sim \varepsilon^{\frac{2N-q(N-2)}{2(q-2)}},\quad   \int_{\mathbb R^N}(I_\alpha\ast |u_\varepsilon|^{\frac{N+\alpha}{N}})|u_\varepsilon|^{\frac{N+\alpha}{N}}\sim \varepsilon^{\frac{(N+\alpha)[4-N(q-2)]}{2N(q-2)}},
$$
$$
\| u_\varepsilon\|_2^2= \varepsilon^{\frac{4-N(q-2)}{2(q-2)}}\left[\frac{2N-q(N-2)}{2q}S_q^{\frac{q}{q-2}}+O(\varepsilon^{-\frac{N(2+\alpha)(q-2)-4\alpha}{2N(q-2)}})\right],
$$
$$
 \|\nabla u_\varepsilon \|_2^2=\varepsilon^{\frac{2N-q(N-2)}{2(q-2)}}\left[\frac{N(q-2)}{2q}S_q^{\frac{q}{q-2}}+O(\varepsilon^{-\frac{N(2+\alpha)(q-2)-4\alpha}{2N(q-2)}})\right].
$$
Moreover, as $\varepsilon\to 0  \  ( resp.   \  \varepsilon\to \infty )$, the rescaled family of ground states 
$$
w_\varepsilon(x)=\varepsilon^{-\frac{1}{q-2}}u_\varepsilon(\varepsilon^{-\frac{1}{2}}x)
$$
converges in $H^1(\mathbb R^N)$ to the unique positive solution of the equation
$$
-\Delta w+w=w^{q-1}.
$$
Furthermore, if $q\in (2, 2+\frac{4\alpha}{N(2+\alpha)})$ and $\varepsilon\to 0 \ ( resp. \ q\in (2+\frac{4\alpha}{N(2+\alpha)}, 2+\frac{4}{N})$ and $\varepsilon\to \infty)$, the least energy $m_\varepsilon$ satisfies 
$$
m_\varepsilon=\varepsilon^{\frac{2N-q(N-2)}{2(q-2)}}\left[\frac{q-2}{2q}S_q^{\frac{q}{q-2}}-\Theta(\varepsilon^{-\frac{N(2+\alpha)(q-2)-4\alpha}{2N(q-2)}})\right].
$$

(II) If $q\in (2, 2+\frac{4\alpha}{N(2+\alpha)})$ and $\varepsilon\to \infty  \ ( resp. \  q\in (2+\frac{4\alpha}{N(2+\alpha)}, 2+\frac{4}{N})$ and $\varepsilon\to 0) $, then the ground state solutions satisfy 
$$
u_\varepsilon(0)\sim \varepsilon^{\frac{2N}{\alpha[4-N(q-2)]}},
$$
$$
\|u_\varepsilon\|_q^q\sim\|\nabla u_\varepsilon \|_2^2\sim \varepsilon^{\frac{N[2N-q(N-2)]}{\alpha[4-N(q-2)]}},
$$
$$
\| u_\varepsilon\|_2^2= \varepsilon^{\frac{N}{\alpha}}\left[S_1^{\frac{N+\alpha}{\alpha}}+O(\varepsilon^{\frac{N(2+\alpha)(q-2)-4\alpha}{\alpha[4-N(q-2)]}})\right],
$$
$$
 \int_{\mathbb R^N}(I_\alpha\ast |u_\varepsilon|^{\frac{N+\alpha}{N}})|u_\varepsilon|^{\frac{N+\alpha}{N}}=\varepsilon^{\frac{N+\alpha}{\alpha}}\left[S_1^{\frac{N+\alpha}{\alpha}}+O(\varepsilon^{\frac{N(2+\alpha)(q-2)-4\alpha}{\alpha[4-N(q-2)]}})\right].
 $$
Moreover,  there exists $\xi_\varepsilon\in (0,+\infty)$
verifying  
$$
\xi_\varepsilon \sim \varepsilon^{-\frac{N(q-2)}{\alpha[4-N(q-2)]}}
$$
such that as $\varepsilon\to \infty  \  ( resp.   \  \varepsilon\to 0 )$,  the rescaled family of ground states
$$
 w_\varepsilon(x)=\varepsilon^{-\frac{N}{2\alpha}}\xi_\varepsilon^{\frac{N}{2}}u_\varepsilon(\xi_\varepsilon x)
$$
 converges in $H^1(\mathbb R^N)$  to the extremal function $U_{\rho_0}$ of $S_1$ with $\rho_0$ being given in (2.4).
 
 Furthermore, if $q\in (2, 2+\frac{4\alpha}{N(2+\alpha)})$ and $\varepsilon\to \infty  \ ( resp. \  q\in (2+\frac{4\alpha}{N(2+\alpha)}, 2+\frac{4}{N})$ and $\varepsilon\to 0) $, the least energy $m_\varepsilon$ satisfies 
 $$
m_\varepsilon=\varepsilon^{\frac{N+\alpha}{\alpha}}\left[\frac{\alpha}{2(N+\alpha)}S_1^{\frac{N+\alpha}{\alpha}}-\Theta(\varepsilon^{\frac{N(2+\alpha)(q-2)-4\alpha}{\alpha[4-N(q-2)]}})\right].
$$ 
 }
 
 \vskip 5mm 
  
{\bf Corollary 2.2.} {\it  If $p=\frac{N+\alpha}{N-2}$, then the problem  $(P_\varepsilon)$ admits  positive ground states $u_\varepsilon\in H^1(\mathbb R^N)$,  which are radially symmetric and radially nonincreasing. Furthermore, the following statements hold true:

$(I)$  If $N\ge 5$, $q\in (2, 2^*)$ and $\varepsilon\to \infty$, then 
the ground states satisfy 
$$
u_\varepsilon(0)\sim \varepsilon^{\frac{1}{q-2}},
$$
$$
\|u_\varepsilon\|_q^q\sim \varepsilon^{-\frac{2N-q(N-2)}{(N-2)(q-2)}}, \quad \| u_\varepsilon\|_2^2\sim \varepsilon^{-\frac{4}{(N-2)(q-2)}}, 
$$
$$
\|\nabla u_\varepsilon \|_2^2= S_\alpha^{\frac{N+\alpha}{2+\alpha}}+O(\varepsilon^{-\frac{2N-q(N-2)}{(N-2)(q-2)}}),
$$
$$ 
\int_{\mathbb R^N}(I_\alpha\ast |u_\varepsilon|^{\frac{N+\alpha}{N-2}})|u_\varepsilon|^{\frac{N+\alpha}{N-2}}=S_\alpha^{\frac{N+\alpha}{2+\alpha}}
+O(\varepsilon^{-\frac{2N-q(N-2)}{(N-2)(q-2)}}).
$$
Moreover, there exists $\xi_\varepsilon\in (0,+\infty)$ verifying 
$$
\xi_\varepsilon \sim \varepsilon^{-\frac{2}{(N-2)(q-2)}}
$$
such that, as $\varepsilon\to \infty$, the rescaled family of ground states 
$$
w_\varepsilon(x)=\xi_\varepsilon^{\frac{N-2}{2}}u_\varepsilon(\xi_\varepsilon x), 
$$
converges   in $H^1(\mathbb R^N)$ to $V_{\rho_0}$ with  $\rho_0$ being given in (2.9).

If $N=4, q\in (2,4)$ and $N=3, q\in (4,6)$, then  as $\varepsilon\to \infty$, the ground states satisfy 
$$
u_\varepsilon(0)\sim \left\{\begin{array}{rcl} \varepsilon^{\frac{1}{q-2}}(\ln\varepsilon)^{\frac{2}{q-2}}, \quad &{\rm if}& \   \ N=4,\\
  \varepsilon^{\frac{1}{2(q-4)}}, \    \quad  \qquad &{\rm if}& \    \  N=3,
 \end{array}\right. 
$$
$$
\|u_\varepsilon\|_2^2\sim \left\{\begin{array}{rcl}\varepsilon^{-\frac{2}{q-2}}(\ln\varepsilon)^{-\frac{4-q}{q-2}}, \quad &{\rm if}& \   \ N=4,\\
\varepsilon^{-\frac{q-2}{2(q-4)}}, \   \quad  \qquad  \quad &{\rm if}& \    \  N=3,
 \end{array}\right. 
$$
$$
\|u_\varepsilon\|_q^q\sim \left\{\begin{array}{rcl} \varepsilon^{-\frac{4-q}{q-2}}(\ln\varepsilon)^{-\frac{4-q}{q-2}}, \quad &{\rm if}& \   \ N=4,\\
  \varepsilon^{-\frac{6-q}{2(q-4)}}, \     \qquad  \qquad &{\rm if}& \    \  N=3,
 \end{array}\right. 
$$
$$
\|\nabla u_\varepsilon\|_2^2= S_\alpha^{\frac{N+\alpha}{2+\alpha}}+\left\{\begin{array}{rcl}
O((\varepsilon\ln\varepsilon)^{-\frac{4-q}{q-2}}),   \  \  \quad if \   \ N=4,\\
O(\varepsilon^{-\frac{6-q}{2(q-4)}}),  \quad \quad \quad if  \    \  N=3,
\end{array}\right.
$$
$$
 \int_{\mathbb R^N}(I_\alpha\ast |u_\varepsilon|^{\frac{N+\alpha}{N-2}})|u_\varepsilon|^{\frac{N+\alpha}{N-2}}=S_\alpha^{\frac{N+\alpha}{2+\alpha}}
+\left\{\begin{array}{rcl}
O((\varepsilon\ln\varepsilon)^{-\frac{4-q}{q-2}}),   \  \  \quad if \   \ N=4,\\
O(\varepsilon^{-\frac{6-q}{2(q-4)}}),  \quad \quad \quad if  \    \  N=3,
\end{array}\right.
$$
Moreover, and there exists $\xi_\varepsilon\in (0,+\infty)$
verifying  
$$
\xi_\varepsilon \sim\left\{\begin{array}{rcl}
\varepsilon^{-\frac{4-q}{2(q-2)}}(\ln\varepsilon)^{-\frac{1}{q-2}},   \  \quad if \   \ N=4,\\
\varepsilon^{-\frac{6-q}{2(q-4)}},  \qquad \quad \qquad if  \    \  N=3,
\end{array}\right.
$$
such that, as $\varepsilon\to \infty$, the rescaled family of ground states
$$
w_\varepsilon(x)=\xi_\varepsilon^{\frac{N-2}{2}}u_\varepsilon(\xi_\varepsilon x)
$$
converges in $D^{1,2}(\mathbb R^N)$ and $L^{q}(\mathbb R^N)$  to $V_1$.

Furthermore, for large $\varepsilon>0$, the least energy $m_\varepsilon$ satisfies 
$$
\frac{2+\alpha}{2(N+\alpha)}S_\alpha^{\frac{N+\alpha}{2+\alpha}}-m_\varepsilon\sim \left\{\begin{array}{rcl} \varepsilon^{-\frac{2N-q(N-2)}{(N-2)(q-2)}}, \qquad  \quad  if \   \ N\ge 5,\\
(\varepsilon\ln\varepsilon)^{-\frac{4-q}{q-2}},  \qquad \quad  if \   \ N=4,\\
 \varepsilon^{-\frac{6-q}{2(q-4)}}, \quad  \quad  \qquad  if  \    \  N=3.
 \end{array}\right.
$$

$(II)$   If $N\ge 3$, $q\in (2, 2^*)$ and $\varepsilon\to 0$, then 
the ground state solutions satisfy 
$$
u_\varepsilon(0)\sim \varepsilon^{\frac{1}{q-2}},
$$
$$
\|u_\varepsilon\|_q^q\sim \varepsilon^{\frac{2N-q(N-2)}{2(q-2)}},\quad 
\int_{\mathbb R^N}(I_\alpha\ast |u_\varepsilon|^{\frac{N+\alpha}{N-2}})|u_\varepsilon|^{\frac{N+\alpha}{N-2}}\sim 
\varepsilon^{\frac{(N+\alpha)(2N-q(N-2))}{2(N-2)(q-2)}},
$$
$$
\| u_\varepsilon\|_2^2= \varepsilon^{\frac{4-N(q-2)}{2(q-2)}}\left[\frac{2N-q(N-2)}{2q}S_q^{\frac{q}{q-2}}+O(\varepsilon^{\frac{(2+\alpha)(2N-q(N-2))}{2(N-2)(q-2)}})\right],                
$$
$$
 \|\nabla u_\varepsilon \|_2^2= \varepsilon^{\frac{2N-q(N-2)}{2(q-2)}}\left[\frac{N(q-2)}{2q}S_q^{\frac{q}{q-2}}+O(\varepsilon^{\frac{(2+\alpha)(2N-q(N-2))}{2(N-2)(q-2)}})\right].
$$
Moreover, as $\varepsilon \to 0$, the rescaled family of ground states 
$$
w_\varepsilon(x)=\varepsilon^{-\frac{1}{q-2}}u_\varepsilon(\varepsilon^{-\frac{1}{2}}x)
$$
converges  in $H^1(\mathbb R^N)$ to the unique positive solution of the equation
$$
-\Delta w+w=w^{q-1}.
$$
Furthermore, for small $\varepsilon>0$, the least energy $m_\varepsilon$ satisfies 
$$
m_\varepsilon=\varepsilon^{\frac{2N-q(N-2)}{2(q-2)}}\left[\frac{q-2}{2q}S_q^{\frac{q}{q-2}}-\Theta(\varepsilon^{\frac{(2+\alpha)(2N-q(N-2))}{2(q-2)}})\right].
$$
}

\vskip 5mm 

{\bf Corollary 2.3.} {\it  If $q=2^*$, then the problem  $(P_\varepsilon)$ admits a positive ground states $u_\varepsilon\in H^1(\mathbb R^N)$,  which are radially symmetric and radially nonincreasing. Furthermore, the following statements hold true:

$(I)$  If $N\ge 5, \ p\in (1+\frac{\alpha}{N-2}, \frac{N+\alpha}{N-2})$ and $\varepsilon\to \infty$, then 
the ground state solutions satisfy 
$$
u_\varepsilon(0)\sim \varepsilon^{\frac{N-2}{2(N-2)(p-1)-2\alpha}},
$$
$$
\int_{\mathbb R^N}(I_\alpha\ast |u_\varepsilon|^p)|u_\varepsilon|^p\sim   \varepsilon^{-\frac{N+\alpha-p(N-2)}{(N-2)(p-1)-\alpha}},
\quad \|u_\varepsilon\|_2^2\sim \varepsilon^{-\frac{2}{(N-2)(p-1)-\alpha}},
$$
$$
\|\nabla  u_\varepsilon\|_2^2=S^{\frac{N}{2}}+O(\varepsilon^{-\frac{N+\alpha-p(N-2)}{(N-2)(p-1)-\alpha}}),
$$
$$
\|u_\varepsilon\|_{2^*}^{2^*}=S^{\frac{N}{2}}+O(\varepsilon^{-\frac{N+\alpha-p(N-2)}{(N-2)(p-1)-\alpha}}).
$$
Moreover, there exists $\xi_\varepsilon\in (0,+\infty)$ verifying 
$$
\xi_\varepsilon\sim \varepsilon^{-\frac{1}{(N-2)(p-1)-\alpha}},
$$
such that, as $\varepsilon\to \infty$, the rescaled family of ground states 
$$
w_\varepsilon(x)=\xi_\varepsilon^{\frac{N-2}{2}}u_\varepsilon(\xi_\varepsilon x), 
$$
converges in $H^1(\mathbb R^N)$ to $W_{\rho_0}$  with $\rho_0$ being given in (2.14).
 
If $N=4, \ p\in (\max\{2, 1+\frac{\alpha}{2}\}, 2+\frac{\alpha}{2})$, and $N=3, \ p\in (2+\alpha, 3+\alpha)$, then as $\varepsilon\to\infty$,
$$
u_\varepsilon(0)\sim \left\{\begin{array}{rcl} (\varepsilon\ln\varepsilon)^{\frac{1}{2p-2-\alpha}}, \quad  if  \   \ N=4,\\
  \varepsilon^{\frac{1}{4(p-2-\alpha)}}, \   \quad  \quad   if  \    \  N=3,
 \end{array}\right. 
$$
$$
\|u_\varepsilon\|_2^2\sim \left\{\begin{array}{rcl}\varepsilon^{-\frac{2}{2p-2-\alpha}}(\ln\varepsilon)^{-\frac{4+\alpha-2p}{2p-2-\alpha}}, \quad  if  \   \ N=4,\\
 \varepsilon^{-\frac{p-1-\alpha}{2(p-2-\alpha)}}, \   \qquad \quad  \quad  if  \    \  N=3,
 \end{array}\right. 
$$
$$
\int_{\mathbb R^N}(I_\alpha\ast |u_\varepsilon|^p)|u_\varepsilon|^p\sim \left\{\begin{array}{rcl} (\varepsilon\ln\varepsilon)^{-\frac{4+\alpha-2p}{2p-2-\alpha}}, \quad \   \  if  \   \ N=4,\\
 \varepsilon^{-\frac{3+\alpha-p}{2(p-2-\alpha)}}, \    \quad  \quad   if  \    \  N=3,
 \end{array}\right. 
$$
$$
\|\nabla u_\varepsilon\|_2^2= S^{\frac{N}{2}}+\left\{\begin{array}{rcl} O((\varepsilon\ln\varepsilon)^{-\frac{4+\alpha-2p}{2p-2-\alpha}}), \quad  if  \   \ N=4,\\
 O( \varepsilon^{-\frac{3+\alpha-p}{2(p-2-\alpha)}}), \       \quad \quad  if  \    \  N=3,
 \end{array}\right. 
$$
$$
 \|u_\varepsilon\|_{2^*}^{2^*}=S^{\frac{N}{2}}+ \left\{\begin{array}{rcl} O((\varepsilon\ln\varepsilon)^{-\frac{4+\alpha-2p}{2p-2-\alpha}}), \quad  if  \   \ N=4,\\
 O( \varepsilon^{-\frac{3+\alpha-p}{2(p-2-\alpha)}}), \     \quad \quad  if  \    \  N=3.
 \end{array}\right. 
$$
Moreover, there exists $\xi_\varepsilon\in (0,+\infty)$ verifying
$$
\xi_\varepsilon\sim \left\{\begin{array}{rcl}(\varepsilon\ln\varepsilon)^{-\frac{1}{2p-2-\alpha}}, \quad  if  \   \ N=4,\\
  \varepsilon^{-\frac{1}{2(p-2-\alpha)}}, \     \quad  \quad   if  \    \  N=3,
 \end{array}\right. 
$$
such that, as $\varepsilon\to \infty$,  the rescaled family of ground states
$$
w_\varepsilon(x)=\xi_\varepsilon^{\frac{N-2}{2}}u_\varepsilon(\xi_\varepsilon x)
$$
converges in $D^{1,2}(\mathbb R^N)$ and $L^{\frac{2Np}{N+\alpha}}(\mathbb R^N)$  to $W_1$.  

Furthermore, for large $\varepsilon>0$, the least energy $m_\varepsilon$ satisfies 
$$
\frac{1}{N}S^{\frac{N}{2}}-m_\varepsilon \sim \left\{\begin{array}{rcl} \varepsilon^{-\frac{N+\alpha-p(N-2)}{(N-2)(p-1)-\alpha}}, \   \quad  \qquad  \quad  if  \   \ N\ge 5,\\
(\varepsilon\ln\varepsilon)^{-\frac{4+\alpha-2p}{2p-2-\alpha}},   \qquad \qquad  if  \   \ N=4,\\
  \varepsilon^{-\frac{3+\alpha-p}{2(p-2-\alpha)}}, \  \quad \quad  \quad  \qquad  if \    \  N=3.
 \end{array}\right. 
$$

(II) If $p\in (\frac{N+\alpha}{N}, \frac{N+\alpha}{N-2})$ and $\varepsilon\to 0$, then the ground state solutions satisfy 
$$
u_\varepsilon(0)\sim \varepsilon^{\frac{2+\alpha}{4(p-1)}},
$$
$$
\|u_\varepsilon\|_2^2=  \varepsilon^{\frac{2+\alpha-N(p-1)}{2(p-1)}}\left[\frac{N+\alpha-p(N-2)}{2p}S_p^{\frac{p}{p-1}}+O(\varepsilon^{\frac{N+\alpha-p(N-2)}{(N-2)(p-1)}})\right],
$$
$$ 
\|\nabla u_\varepsilon \|_2^2= \varepsilon^{\frac{N+\alpha-p(N-2)}{2(p-1)}}\left[\frac{N(p-1)-\alpha}{2p}S_p^{\frac{p}{p-1}}+O(\varepsilon^{\frac{N+\alpha-p(N-2)}{(N-2)(p-1)}})\right].
$$
Moreover,  for any sequence $\varepsilon_n\to 0$,  the rescaled family of ground states
$$
w_{\varepsilon_n}(x)=\varepsilon_n^{-\frac{2+\alpha}{4(p-1)}}u_{\varepsilon_n}(\varepsilon_n^{-\frac{1}{2}}x)
$$
converges up to a subsequence in $H^1(\mathbb R^N)$ to a positive solution of the equation
$$
-\Delta w+w=(I_\alpha\ast |w|^p)w^{p-1}.
$$
Furthermore, for small $\varepsilon>0$, the least energy $m_\varepsilon$ satisfies 
$$
m_\varepsilon=\varepsilon^{\frac{N+\alpha-p(N-2)}{2(p-1)}}\left[\frac{p-1}{2p}S_p^{\frac{p}{p-1}}-\Theta(\varepsilon^{\frac{N+\alpha-p(N-2)}{(N-2)(p-1)}})\right].
$$
}
\vskip 5mm

{\bf Corollary 2.4.} {\it If $p\in (\frac{N+\alpha}{N},  \frac{N+\alpha}{N-2})$ and $q\in (2,2^*)$, then the problem  $(P_\varepsilon)$ admits a positive ground states $u_\varepsilon\in H^1(\mathbb R^N)$,  which are radially symmetric and radially nonincreasing. Furthermore, the following statements hold true:

(I) If $q<\frac{2(2p+\alpha)}{2+\alpha}$, then for small $\varepsilon>0$,
$$
\|u_\varepsilon\|_2^2=\varepsilon^{\frac{4-N(q-2)}{2(q-2)}}\left[\frac{2N-q(N-2)}{2q}S_q^{\frac{q}{q-2}}-\Theta(\varepsilon^{\frac{2(2p+\alpha)-q(2+\alpha)}{2(q-2)}})\right],
$$
$$
\|\nabla u_\varepsilon\|_2^2=\varepsilon^{\frac{2N-q(N-2)}{2(q-2)}}
\left[\frac{N(q-2)}{2q}S_q^{\frac{q}{q-2}}+O(\varepsilon^{\frac{2(2p+\alpha)-q(2+\alpha)}{2(q-2)}})\right],
$$
and  as $\varepsilon\to 0$, the rescaled family of ground states 
$$
w_\varepsilon(x)=\varepsilon^{-\frac{1}{q-2}}u_\varepsilon(\varepsilon^{-\frac{1}{2}}x),
$$
converges in $H^1(\mathbb R^N)$ to the unique positive solution  $w_0$ of the equation
$$
-\Delta w+w=w^{q-1}.
$$
For large $\varepsilon>0$, 
$$
\|u_\varepsilon\|_2^2=\varepsilon^{\frac{2+\alpha-N(p-1)}{2(p-1)}}\left[\frac{N+\alpha-p(N-2)}{2p}S_p^{\frac{p}{p-1}}+O(\varepsilon^{-\frac{2(2p+\alpha)-q(2+\alpha)}{4(p-1)}})\right],
$$
$$
\|\nabla u_\varepsilon\|_2^2=\varepsilon^{\frac{N+\alpha-p(N-2)}{2(p-1)}}
\left[\frac{N(p-1)-\alpha}{2p}S_p^{\frac{p}{p-1}}+O(\varepsilon^{-\frac{2(2p+\alpha)-q(2+\alpha)}{4(p-1)}})\right],
$$
and as 
$\varepsilon\to \infty$, the rescaled family of ground states
$$
w_{\varepsilon}(x)=\varepsilon^{-\frac{2+\alpha}{4(p-1)}}u_{\varepsilon}(\varepsilon^{-\frac{1}{2}}x), 
$$
converges up to a subsequence in   $H^1(\mathbb R^N)$ to a positive solution $w_\infty$ of the equation
$$
-\Delta w+w=(I_\alpha\ast |w|^p)w^{p-1}.
$$
Furthermore, for small $\varepsilon>0$, the least energy $m_\varepsilon$ satisfies 
$$
m_\varepsilon=\varepsilon^{\frac{2N-q(N-2)}{2(q-2)}}\left[\frac{q-2}{2q}S_q^{\frac{q}{q-2}}-\Theta(\varepsilon^{\frac{2(2p+\alpha)-q(2+\alpha)}{2(q-2)}})\right],
$$
and for large $\varepsilon>0$, the least energy $m_\varepsilon$ satisfies 
$$
m_\varepsilon=\varepsilon^{\frac{N+\alpha-p(N-2)}{2(p-1)}}\left[\frac{p-1}{2p}S_p^{\frac{p}{p-1}}-\Theta(\varepsilon^{-\frac{2(2p+\alpha)-q(2+\alpha)}{4(p-1)}})\right].
$$

(II) If $q>\frac{2(2p+\alpha)}{2+\alpha}$, then 
for small  $\varepsilon> 0$, 
$$
\|u_\varepsilon\|_2^2=\varepsilon^{\frac{2+\alpha-N(p-1)}{2(p-1)}}\left[\frac{N+\alpha-p(N-2)}{2p}S_p^{\frac{p}{p-1}}-\Theta(\varepsilon^{\frac{q(2+\alpha)-2(2p+\alpha)}{4(p-1)}})\right],
$$
$$
\|\nabla u_\varepsilon\|_2^2=\varepsilon^{\frac{2N-q(N-2)}{2(q-2)}}
\left[\frac{N(q-2)}{2q}S_q^{\frac{q}{q-2}}+O(\varepsilon^{\frac{q(2+\alpha)-2(2p+\alpha)}{4(p-1)}})\right],
$$
and as $\varepsilon\to 0$,  the rescaled family of ground states
$$
w_{\varepsilon}(x)=\varepsilon^{-\frac{2+\alpha}{4(p-1)}}u_{\varepsilon}(\varepsilon^{-\frac{1}{2}}x), 
$$
converges up to a subsequence in   $H^1(\mathbb R^N)$ to a positive solution $w_0$ of the equation
$$
-\Delta w+w=(I_\alpha\ast |w|^p)w^{p-1}.
$$
For large $\varepsilon>0$, 
$$
\|u_\varepsilon\|_2^2=\varepsilon^{\frac{4-N(q-2)}{2(q-2)}}\left[\frac{2N-q(N-2)}{2q}S_q^{\frac{q}{q-2}}+O(\varepsilon^{-\frac{q(2+\alpha)-2(2p+\alpha)}{2(q-2)}})\right],
$$
$$
\|\nabla u_\varepsilon\|_2^2=\varepsilon^{\frac{N+\alpha-p(N-2)}{2(p-1)}}
\left[\frac{N(p-1)-\alpha}{2p}S_p^{\frac{p}{p-1}}+O(\varepsilon^{-\frac{q(2+\alpha)-2(2p+\alpha)}{2(q-2)}})\right],
$$
and  as $\varepsilon\to \infty$, the rescaled family of ground states 
$$
w_\varepsilon(x)=\varepsilon^{-\frac{1}{q-2}}u_\varepsilon(\varepsilon^{-\frac{1}{2}}x),
$$
converges in $H^1(\mathbb R^N)$ to the unique positive solution $w_\infty$ of the equation
$$
-\Delta w+w=w^{q-1}.
$$
Furthermore, for small $\varepsilon>0$, the least energy $m_\varepsilon$ satisfies 
$$
m_\varepsilon=\varepsilon^{\frac{N+\alpha-p(N-2)}{2(p-1)}}\left[\frac{p-1}{2p}S_p^{\frac{p}{p-1}}-\Theta(\varepsilon^{\frac{q(2+\alpha)-2(2p+\alpha)}{4(p-1)}})\right],
$$
and for large $\varepsilon>0$, the least energy $m_\varepsilon$ satisfies 
$$
m_\varepsilon=\varepsilon^{\frac{2N-q(N-2)}{2(q-2)}}\left[\frac{q-2}{2q}S_q^{\frac{q}{q-2}}-\Theta(\varepsilon^{-\frac{q(2+\alpha)-2(2p+\alpha)}{2(q-2)}})\right].
$$
}
\vskip 5mm

\section*{3. Preliminaries } 

In this section, we present some preliminary results which are needed in the proof of our main results. Firstly,  we consider the following Choquard type  equation with combined nonlinearites:
$$
-\Delta u+u=\mu(I_{\alpha}\ast |u|^p)|u|^p+\lambda |u|^{q-2}u,  \quad \mathrm{in}\
\mathbb{R}^N,
\eqno(Q_{\mu,\lambda})
$$
where $N\geq 3,\ \alpha\in(0,N)$, $p\in [\frac{N+\alpha}{N}, \frac{N+\alpha}{N-2}], \ q\in (2, 2^*]$,
 $\mu>0$ and $\lambda>0$ are two parameters. 

It has been proved in \cite{Li-1} that any weak solution of $(Q_{\mu,\lambda})$ in
$H^1(\mathbb{R}^N)$ has additional regularity properties, which
allows us to establish the Poho\v{z}aev identity for all finite
energy solutions.

{\bf Lemma 3.1.}
{\it  If
$u\in H^1(\mathbb{R}^N)$ is a solution of $(Q_{\mu,\lambda})$, then $u\in
W_{\mathrm{loc}}^{2,r}(\mathbb{R}^N)$ for every $r>1$. Moreover, $u$
satisfies the Poho\v{z}aev identity
$$
P_{\mu, \lambda}(u):=\frac{N-2}{2}\int_{\mathbb{R}^N}|\nabla
u|^2+\frac{N}{2}\int_{\mathbb{R}^N}|u|^2-\frac{N+\alpha}{2p}\mu\int_{\mathbb{R}^N}(I_\alpha\ast
|u|^p)|u|^p-\frac{N}{q}\lambda\int_{\R^N}|u|^q=0.
$$}

It is well known that any weak solution of  $(Q_{\mu,\lambda})$ corresponds to a critical point of the action  functionals $I_{\mu,\lambda}$ defined by 
$$
I_{\mu,\lambda}(u):=\frac{1}{2}\int_{\mathbb R^N}|\nabla u|^2+|u|^2-\frac{\mu}{2p}\int_{\mathbb R^N}(I_\alpha\ast |u|^p)|u|^p-\frac{\lambda}{q}\int_{\mathbb R^N}|u|^q, 
\eqno(3.1)
$$
which is well defined and is of $C^1$ in $H^1(\R^N)$.  A nontrivial solution $u_{\mu,\lambda}\in H^1(\R^N)$  is called a ground state if 
$$
I_{\mu,\lambda}(u_{\mu,\lambda})=m_{\mu, \lambda}:=\inf\{ I_{\mu,\lambda}(u): \ u\in H^1(\R^N) \ {\rm and}  \  I'_{\mu,\lambda}(u)=0\}.
\eqno(3.2)
$$

In \cite{Li-2, Li-1} (see also the proof of the main results in \cite{Li-1}),  it has been shown that 
$$
m_{\mu, \lambda}=\inf_{u\in \mathcal M_{\mu, \lambda}}I_{\mu, \lambda}(u)=\inf_{u\in \mathcal P_{\mu, \lambda}}I_{\mu, \lambda}(u),
\eqno(3.3)
$$
where $\mathcal M_{\mu, \lambda}$ and $\mathcal P_{\mu, \lambda}$ are the correspoding Nehari  and Poho\v{z}aev manifolds defined by
$$
\mathcal M_{\mu, \lambda}:=\left\{ u\in H^1(\mathbb R^N)\setminus\{0\}  \ \left | \ \int_{\mathbb R^N}|\nabla u|^2+|u|^2=\mu \int_{\mathbb R^N}(I_\alpha\ast |u|^p)|u|^p+\lambda\int_{\mathbb R^N}|u|^q \right. \right\}
$$
and 
$$
\mathcal P_{\mu, \lambda}:=\left\{ u\in H^1(\mathbb R^N)\setminus\{0\}  \ \left | \ P_{\mu,\lambda}(u)=0 \right. \right\},
$$
respectively.
Moreover, the following min-max descriptions are valid:

{\bf Lemma 3.2.}  {\it  Let
 $$
u_t(x)=\left\{\begin{array}{rcl} u(\frac{x}{t}), \quad if \  t>0,\\
0, \quad \  if \  \ t=0,
\end{array}\right.
$$
then 
$$
m_{\mu, \lambda}=\inf_{u\in H^1(\mathbb R^N)\setminus\{0\}}\sup_{t\ge 0}I_{\mu, \lambda}(tu)=\inf_{u\in H^1(\mathbb R^N)\setminus\{0\}}\sup_{t\ge 0}I_{\mu, \lambda}(u_t).
\eqno(3.4)
$$
In particular, we have $m_{\mu, \lambda}=I_{\mu,\lambda}(u_{\mu, \lambda})=\sup_{t>0}I_{\mu, \lambda}(tu_{\mu,\lambda})=\sup_{t>0}I_{\mu, \lambda}((u_{\mu,\lambda})_t)$.}

When $\mu=1$ and $\lambda=0$, then the equation  $(Q_{\mu,\lambda})$  reduces to  
$$
-\Delta u+u=(I_\alpha \ast |u|^{p})|u|^{p-2}u, 
\quad {\rm in} \ \  \mathbb R^N,
 \eqno(Q_{1,0})
$$
when $\mu=0$ and $\lambda=1$, then the equation  $(Q_{\mu,\lambda})$ reduces to  
$$
\Delta u+u=|u|^{q-2}u, 
\quad {\rm in} \ \  \mathbb R^N,
 \eqno(Q_{0,1})
$$
Then the correspoding Nehari manifolds are as follows.
$$
\mathcal M_{1,0}=\left\{ u\in H^1(\mathbb R^N)\setminus\{0\}  \ \left | \ \int_{\mathbb R^N}|\nabla u|^2+|u|^2=\int_{\mathbb R^N}(I_\alpha\ast |u|^p)|u|^p\right. \right\}.
$$
$$
\mathcal M_{0,1}=\left\{ u\in H^1(\mathbb R^N)\setminus\{0\}  \ \left | \ \int_{\mathbb R^N}|\nabla u|^2+|u|^2=\int_{\mathbb R^N}|u|^q\right. \right\}.
$$
It is known  that
$$
m_{1,\lambda}=\inf_{u\in \mathcal M_{1,\lambda}}I_{1,\lambda}(u), \qquad 
m_{1,0}:=\inf_{u\in \mathcal M_{1,0}} \ I_{1,0}(u),
\eqno(3.5)
$$
and 
$$
m_{\mu,1}:=\inf_{u\in \mathcal M_{\mu,1}} \ I_{\mu,1}(u), \qquad 
m_{0,1}:=\inf_{u\in \mathcal M_{0,1}} \ I_{0,1}(u)
\eqno(3.6)
$$
are well-defined and positive.

Let $u_{\mu,\lambda}$ be the ground state for $(Q_{\mu, \lambda})$, then we have the following

{\bf Lemma 3.3.} {\it The solution sequences  $\{u_{1,\lambda}\}$ and $\{u_{\mu,1}\}$ are  bounded in $H^1(\mathbb R^N)$.}
\begin{proof}  It is not hard to see  that $m_{1,\lambda}\le m_{1,0}\le C<+\infty$. If $q\ge 2p$, then 
$$
\begin{array}{rcl}
m_{1,\lambda}&=&I_{1,\lambda}(u_{1,\lambda})=I_{1,\lambda}(u_{1,\lambda})-\frac{1}{2p}I'_{1,\lambda}(u_{1,\lambda})u_{1,\lambda}\\
&=&\left(\frac{1}{2}-\frac{1}{2p}\right)\int_{\mathbb R^N}|\nabla u_{1,\lambda}|^2+|u_{1,\lambda}|^2+\left(\frac{1}{2p}-\frac{1}{q}\right)\lambda\int_{\mathbb R^N}|u_{1,\lambda}|^q\\
&\ge& \left(\frac{1}{2}-\frac{1}{2p}\right)\int_{\mathbb R^N}|\nabla u_{1,\lambda}|^2+|u_{1,\lambda}|^2,
\end{array}
$$
and if $q<2p$, then
$$
\begin{array}{rcl}
m_{1,\lambda}&=&I_{1,\lambda}(u_{1,\lambda})=I_{1,\lambda}(u_{1,\lambda})-\frac{1}{q}I'_{1,\lambda}(u_{1,\lambda})u_{1,\lambda}\\
&=&\left(\frac{1}{2}-\frac{1}{q}\right)\int_{\mathbb R^N}|\nabla u_{1,\lambda}|^2+|u_{1,\lambda}|^2+\left(\frac{1}{q}-\frac{1}{2p}\right)\int_{\mathbb R^N}(I_\lambda\ast |u_{1,\lambda}|^p)|u_{1,\lambda}|^p\\
&\ge& \left(\frac{1}{2}-\frac{1}{q}\right)\int_{\mathbb R^N}|\nabla u_{1,\lambda}|^2+|u_{1,\lambda}|^2.
\end{array}
$$
Therefore, we conclude that $\{u_{1,\lambda}\}$ is bounded in $H^1(\mathbb R^N)$. 

Arguing as above, we show that $\{u_{\mu,1}\}$ is bounded in $H^1(\mathbb R^N)$. The proof is completed. 
\end{proof}

The following well
known Hardy-Littlewood-Sobolev inequality can be found in
\cite{Lieb-Loss 2001}.

{\bf Lemma 3.4.} {\it 
Let $p, r>1$ and $0<\alpha<N$ with $1/p+(N-\alpha)/N+1/r=2$. Let
$u\in L^p(\mathbb{R}^N)$ and $v\in L^r(\mathbb{R}^N)$. Then there
exists a sharp constant $C(N,\alpha,p)$, independent of $u$ and $v$,
such that
$$
\left|\int_{\mathbb{R}^N}\int_{\mathbb{R}^N}\frac{u(x)v(y)}{|x-y|^{N-\alpha}}\right|\leq
C(N,\alpha,p)\|u\|_p\|v\|_r.
$$
If $p=r=\frac{2N}{N+\alpha}$, then
$$
C(N,\alpha,p)=C_\alpha(N)=\pi^{\frac{N-\alpha}{2}}\frac{\Gamma(\frac{\alpha}{2})}{\Gamma(\frac{N+\alpha}{2})}\left\{\frac{\Gamma(\frac{N}{2})}{\Gamma(N)}\right\}^{-\frac{\alpha}{N}}.
$$}

{\bf Remark 3.1. }  By the Hardy-Littlewood-Sobolev inequality, for any $v\in L^s(\mathbb R^N)$ with $s\in (1,\frac{N}{\alpha})$, $I_\alpha\ast v\in L^{\frac{Ns}{N-\alpha s}}(\mathbb R^N)$ and 
$$
\|I_\alpha\ast v\|_{\frac{Ns}{N-\alpha s}}\le A_\alpha(N)C(N,\alpha, s)\|v\|_{s}.
\eqno(3.7)
$$

 {\bf Lemma 3.5.} ( P. L. Lions \cite{Lions-1} )
{\it Let $r>0$ and $2\leq q\leq 2^{*}$. If $(u_{n})$ is bounded in $H^{1}(\mathbb{R}^N)$ and if
$$\sup_{y\in\mathbb{R}^N}\int_{B_{r}(y)}|u_{n}|^{q}\to0,\,\,\textrm{as\ }n\to\infty,$$
then $u_{n}\to0$ in $L^{s}(\mathbb{R}^N)$ for $2<s<2^*$. Moreover, if $q=2^*$, then $u_{n}\to0$ in $L^{2^{*}}(\mathbb{R}^N)$.
}
 
{\bf  Lemma 3.6.} {\it  Let $r>0$, $N\ge 3$, $\alpha\in (0,N)$ and $\frac{N+\alpha}{N}\le p\le \frac{N+\alpha}{N-2}$. If  $(u_n)$ be bounded in $H^1(\mathbb R^N)$ and if 
$$
\lim_{n\to \infty}\sup_{z\in \mathbb R^N} \int_{B_r(z)}\int_{B_r(z)}\frac{|u_n(x)|^p|u_n(y)|^p}{|x-y|^{N-\alpha}} dx dy = 0, 
$$
then 
$$
\lim_{n\to \infty}\int_{\mathbb R^N}|u_n|^s dx = \lim_{n\to\infty} \int_{\mathbb R^N}(I_\alpha\ast |u_n|^t)|u_n|^t dx = 0,
$$
for any $2<s<2^*$ and $\frac{N+\alpha}{N}<t<\frac{N+\alpha}{N-2}$. Moreover, if $p=\frac{N+\alpha}{N-2}$, then 
$$
\lim_{n\to \infty}\int_{\mathbb R^N}|u_n|^{2^*} dx = \lim_{n\to\infty}\int_{\mathbb R^N} (I_\alpha\ast |u_n|^{\frac{N+\alpha}{N-2}})|u_n|^{\frac{N+\alpha}{N-2}} dx = 0.
$$} 
\begin{proof}  Similar to the proof of   \cite[Lemma 3.8]{Cassani-1}, for any $p\in [\frac{N+\alpha}{N}, \frac{N+\alpha}{N-2}]$,  it is easy to show that  
$$
\lim_{n\to \infty}\sup_{z\in \mathbb R^N} \int_{B_r(z)}\int_{B_r(z)}\frac{|u_n(x)|^p|u_n(y)|^p}{|x-y|^{N-\alpha}} dx dy = 0
$$
is equivalent to the following condition
$$
\lim_{n\to \infty}\sup_{z\in \mathbb R^N}\int_{B_r(z)}|u_n|^{\frac{2Np}{N+\alpha}} dx =0.
$$
Then the conclusion follows from Lemma 3.5. The proof is complete. 
\end{proof}

{\bf Lemma 3.7.}  (Radial Lemma A.II, H. Berestycki and P. L. Lions \cite{Berestycki-1})
{\it  Let $N\ge 2$, then every radial function $u\in H^1(\mathbb R^N)$ is almost everywhere equal to a function $\tilde u(x)$, continuous for $x\not=0$, such that
$$
 |\tilde u(x)|\le  C_N |x|^{(1-N)/2}\|u\|_{H^1(\mathbb R^N)} \qquad for  \ |x|\ge \alpha_N,
\eqno(3.8)
$$
 where $C_N$ and $\alpha_N$ depend only on the dimension $N$.}

 {\bf Lemma 3.8.} ( Radial Lemma A.III. H. Berestycki and P. L. Lions \cite{Berestycki-1} ) 
 {\it Let $N\ge 3$, then every radial function $u$ in $D^{1,2}(\mathbb R^N)$ is almost everywhere equal to a function $\tilde u(x)$, continuous for $x\not=0$, such that
 $$
  |\tilde u(x)|\le C_N |x|^{(2-N)/2}\|u\|_{D^{1,2}(\mathbb R^N)} \qquad for  \ |x|\ge 1, 
\eqno(3.9)
  $$
  where $C_N$ only depends on $N$.} 
  \vskip 5mm

{\bf Lemma 3.9.}  {\it Let $0<\alpha<N$ and $0\le f\in L^1(\R^N)$. Assume that
$$
\lim_{|x|\to \infty}\frac{\int_{|y|\le |x|}f(y)|y|dy}{|x|}=0,
\eqno(3.10)
$$
$$
\lim_{|x|\to \infty}\int_{|y-x|\le |x|/2}\frac{f(y)dy}{|x-y|^{N-\alpha}}=0.
\eqno(3.11)
$$
Then as $|x|\to\infty$, 
$$
\int_{\R^N}\frac{f(y)dy}{|x-y|^{N-\alpha}}=\frac{\|f\|_{L^1}}{|x|^{N-\alpha}}+o\left(\frac{1}{|x|^{N-\alpha}}\right).
\eqno(3.12)
$$}

Note that $f\in L^1(\R^N)$ alone is not sufficient to obtain (3.12) even if $f$ is radially symmetric, see \cite{Sie-1}.

\begin{proof}  Fix $0\not=x\in \R^N$, we  decompose $\R^N$ as the union of three sets
 $B=\{y: \ |y-x|<|x|/2\}$, 
$A=\{y\notin B: \ |y|\le |x| \}$ and $C=\{y\not\in B: \ |y|>|x| \}$. 

We want to estimate the quantity
$$
\left|\int_{A\cup C}f(y)\left(\frac{1}{|x-y|^{N-\alpha}}-\frac{1}{|x|^{N-\alpha}}\right)dy\right|\le 
\int_{A\cup C}f(y)\left|\frac{1}{|x-y|^{N-\alpha}}-\frac{1}{|x|^{N-\alpha}}\right|dy.
$$
Since $|x|/2\le |x-y|\le 2|x|$ for all $y\in A$, by the  Mean Value Theorem we have
$$
\left|\frac{1}{|x-y|^{N-\alpha}}-\frac{1}{|x|^{N-\alpha}}\right|\le \frac{c_1|y|}{|x|^{N-\alpha+1}},  \  (y\in A), 
$$
where $c_1=(N-\alpha)2^{N-\alpha+1}.$  Thus 
$$
\left|\int_{A}f(y)\left(\frac{1}{|x-y|^{N-\alpha}}-\frac{1}{|x|^{N-\alpha}}\right)dy\right|\le \frac{c_1}{|x|^{N-\alpha+1}}\int_Af(y)|y|dy.
$$
On the other hand, since $|x-y|>|x|/2$ for all $y\in C$, then
$$
\left|\frac{1}{|x-y|^{N-\alpha}}-\frac{1}{|x|^{N-\alpha}}\right|\le \frac{1}{|x|^{N-\alpha}}, \  (y\in C),
$$
from which we compute that
$$
\left|\int_{C}f(y)\left(\frac{1}{|x-y|^{N-\alpha}}-\frac{1}{|x|^{N-\alpha}}\right)dy\right|\le \frac{1}{|x|^{N-\alpha}}\int_{A\cup C}f(y)dy.
$$
Then
$$
\begin{array}{cl}
&\left|\int_{\R^N}\frac{f(y)}{|x-y|^{N-\alpha}}dy-\frac{\|f\|_{L^1}}{|x|^{N-\alpha}}\right|\\
&\le 
\frac{c_1}{|x|^{N-\alpha+1}}\int_Af(y)|y|dy +\int_B\frac{f(y)}{|x-y|^{N-\alpha}}dy+\frac{1}{|x|^{N-\alpha}}\int_{B\cup C}f(y)dy.\\
\end{array}
$$
The conclusion follows from (3.10), (3.11) and since $f\in L^1(\R^N)$.
\end{proof}

{\bf Lemma 3.10.} {\it Let $0<\alpha<N$, $0\le f(x)\in L^1(\mathbb R^N)$ be a radially symmetric function such that
$$
\lim_{|x|\to +\infty}f(|x|)|x|^N=0.
\eqno(3.13)
$$
If $\alpha\le 1$, we additionally assume that $f$ is monotone non-increasing. Then as $|x|\to +\infty$, we have 
$$
\int_{\mathbb R^N}\frac{f(y)}{|x-y|^{N-\alpha}}dy=\frac{\|f\|_{L^1}}{|x|^{N-\alpha}}+o\left(\frac{1}{|x|^{N-\alpha}}\right).
\eqno(3.14)
$$}
\begin{proof}
Using (3.13) by I'Hospital rule we conclude that
$$
\int_{|y|\le |x|}f(y)|y|dy=\int_0^{|x|}f(r)r^Ndr=o(|x|),  \   (|x|\to \infty), 
$$
so (3.10)  holds.

For $|x|\gg 1$, using radial estimates on the Riesz kernels in \cite[Lemma 2.2]{Duo-1} and (3.13) we obtain
for $\alpha >1$: 
$$
\int_{|y-x|\le |x|/2}\frac{f(y)dy}{|x-y|^{N-\alpha}}\lesssim |x|^{\alpha-1}\int_{|x|/2}^{3|x|/2}f(r)dr=o(|x|^{-(N-\alpha)});
$$
for $\alpha=1$, additionally using monotonocity of $f$:
$$
\begin{array}{lcl}
\int_{|y-x|\le |x|/2}\frac{f(y)dy}{|x-y|^{N-\alpha}}&\lesssim &\int_{|x|/2}^{3|x|/2}f(r)\log\frac{1}{1-r/|x|}dr\\
&\le & f(|x|/2)\int_{|x|/2}^{3|x|/2}\log\frac{1}{1-r/|x|}dr=o(|x|^{-(N-1)});
\end{array}
$$
for $\alpha<1$, additionally using monotonicity of $f$:
$$
\begin{array}{lcl}
\int_{|y-x|\le |x|/2}\frac{f(y)dy}{|x-y|^{N-\alpha}}&\lesssim & \int_{|x|/2}^{3|x|/2}\frac{f(r)}{|r-|x||^{1-\alpha}}dr\\
&\le & f(|x|/2)\int_{|x|/2}^{3|x|/2}\frac{1}{|1-|x||^{1-\alpha}}dr= o(|x|^{-(N-\alpha)});
\end{array}
$$
so (3.12) holds. This completes the proof.
\end{proof}

{\bf Lemma 3.11.} {\it Assume that $m'>\frac{N}{2}, m\in (0,m']$ and $x\in \mathbb{R}^N$ tends to infinity. Then 
$$
\int_{\mathbb{R}^N}\frac{1}{(1+|y-x|^2)^m(1+|y|^2)^{m'}}dy\sim \frac{1}{|x|^{2m}}.
\eqno(3.15)
$$ }
\begin{proof}  We suppose, by using a convenient rotation, that $x=(x_0,0,\cdots,0)$.  Without loss of generality, we also assume that $x_0>0$ is large.
The first conclusion follows from the following two claims.

Claim 1:  For some positive constant $c>0$ and large $x$, we have
$$
\int_{\mathbb{R}^N}\frac{1}{(1+|y-x|^2)^m(1+|y|^2)^{m'}}dx\ge  \frac{c}{|x|^{2m}}.
\eqno(3.16)
$$
In fact, we have 
$$
\begin{array}{rl}
 \displaystyle{}\int_{\mathbb{R}^N}\frac{1}{(1+|y-x|^2)^m(1+|y|^2)^{m'}}dy
 &\ge  \displaystyle{}\int_{B_1(0)}\frac{1}{(1+|y-x|^2)^m(1+|y|^2)^{m'}}dy \\ 
&\ge  \displaystyle{}\frac{1}{2^{m'}}\int_{B_1(0)}\frac{1}{(1+|y-x|^2)^m}dy\\ 
&\ge  \displaystyle{}\frac {1}{2^{m'}}\int_{B_1(0)}\frac{1}{(5|x|^2)^m}dy,
\end{array}
$$
from which  Claim  1 follows.

Claim 2:  For some positive constant $c>0$ and large $x$, we also have
$$
\int_{\mathbb{R}^N}\frac{1}{(1+|y-x|^2)^m(1+|y|^2)^{m'}}dy\le   \frac{c}{|x|^{2m}}.
\eqno(3.17)
$$
To prove this claim, we use the notation: $y=(r,y')$ with $r\in \mathbb{R}$ and $y'\in\mathbb{R}^{N-1}$.  Then
$$
\begin{array}{rl}
 \displaystyle{}\int_{\mathbb{R}^N}\frac{1}{(1+|y-x|^2)^m(1+|y|^2)^{m'}}dy
 &=  \displaystyle{}\int_{-\infty}^{+\infty}\int_{\mathbb{R}^{N-1}}\frac{1}{(1+|y-x|^2)^m(1+|y|^2)^{m'}}dy'dr\\ 
&=   \displaystyle{}\int_{-\infty}^{x_0/2}\int_{\mathbb{R}^{N-1}}\frac{1}{(1+|y-x|^2)^m(1+|y|^2)^{m'}}dy'dr\\ 
&\mbox{} \   \displaystyle{}+ \int_{x_0/2}^{+\infty}\int_{\mathbb{R}^{N-1}}\frac{1}{(1+|y-x|^2)^m(1+|y|^2)^{m'}}dy'dr.
\end{array}
$$
If $r\le x_0/2$, then we have 
$$
|y-x|=\sqrt{(x_0-r)^2+|y'|^2}\ge x_0-r\ge x_0/2=|x|/2.
$$
Therefore, noting that $m'>N/2$, we obtain
$$
\begin{array}{lcl}
J_1(x):&=&\displaystyle{}\int_{-\infty}^{x_0/2}\int_{\mathbb{R}^{N-1}}\frac{1}{(1+|y-x|^2)^m(1+|y|^2)^{m'}}dy'dr\\
&\le& \displaystyle{}\frac{4^m}{(4+|x|^2)^m}\int_{-\infty}^{x_0/2}\int_{\mathbb{R}^{N-1}}\frac{1}{(1+|y|^2)^{m'}}dy'dr\\ 
&\le& \displaystyle{}\frac{4^m}{(4+|x|^2)^m}\int_{\mathbb{R}^{N}}\frac{1}{(1+|y|^2)^{m'}}dy\le \frac{C}{|x|^{2m}}.
\end{array}
$$
Since
$$
\frac{1}{1+|y|^2}\le \frac{1}{1+|y-x|^2},  \qquad {\rm for \ all} \ y\in \{y=(r,y'):  \ r\ge x_0/2, \ y'\in \mathbb{R}^{N-1}\},
$$
a similar  argument shows that
$$
\begin{array}{lcl}
J_2(x):&=&\displaystyle{}\int_{x_0/2}^{+\infty}\int_{\mathbb{R}^{N-1}}\frac{1}{(1+|y-x|^2)^m(1+|y|^2)^{m'}}dy'dr\\
& \displaystyle{}\le&  \int_{x_0/2}^{+\infty}\int_{\mathbb{R}^{N-1}}\frac{1}{(1+|y-x|^2)^{m'}(1+|y|^2)^{m}}dy'dr\\ 
&= &\displaystyle{}\int_{-x_0/2}^{+\infty}\int_{\mathbb{R}^{N-1}}\frac{1}{(1+|y+x|^2)^{m}(1+|y|^2)^{m'}}dy'dr\le  \frac{C}{|x|^{2m}}.
 \end{array}
 $$
This establishes  Claim 2. The proof is complete.
\end{proof}

{\bf Remark 3.2.} In Lemma 3.11, if $m'>\frac{N}{2}$ can be chosen in such a way that $C_1\ge m'\ge C_2>\frac{N}{2}$ with $C_1,C_2$ being independent of $m'$ and $m$, then there exists  constants $C_0$ and $R_0$ dependent of $C_1$ and $C_2$ such that  
$$
\int_{\mathbb{R}^N}\frac{1}{(1+|y-x|^2)^m(1+|y|^2)^{m'}}dy\le \frac{C_0}{|x|^{2m}},\quad |x|\ge R_0.
$$

\vskip 5mm 

The following Moser iteration lemma  is given in \cite[Proposition B.1]{Akahori-2}. See also  \cite{LiuXQ}  and \cite{GT}.

{\bf Lemma 3.12.} {\it Assume $N\ge 3$. Let $a(x)$ and $b(x)$ be functions on $B_4$, and let $u\in H^1(B_4)$ be a weak solution to 
$$
-\Delta u+a(x)u=b(x)u \qquad  in \  \ B_4.
\eqno(3.18)
$$
Suppose that $a(x)$ and $u$ satisfy that 
$$
a(x)\ge 0 \quad for \ a. e. \ x\in B_4, 
\eqno(3.19)
$$
and 
$$
 \int_{B_4}a(x)|u(x)v(x)|dx<\infty \quad for \ each \ v\in H_0^1(B_4).
\eqno(3.20)
$$
(i) Assume that for any $\varepsilon\in (0,1)$, there exists $t_\varepsilon>0$ such  that
$$
\|\chi_{[|b|>t_\varepsilon]}b\|_{L^{N/2}(B_4)}\le \varepsilon,
$$
where $[|b|>t]:=\{x\in B_4: \ |b(x)|>t\},$ and $\chi_A(x)$ denotes the characteristic function of $A\subset \mathbb R^N$. Then for any $r\in (0,\infty)$, there exists a constant $C(N,r,t_\varepsilon)$ such that 
$$
\||u|^{r+1}\|_{H^1(B_1)}\le C(N, r,t_\varepsilon)\|u\|_{L^{2^*}(B_4)}.
$$
(ii) \ Let $s>N/2$ and assume that $b\in L^s(B_4)$. Then there exits a constant $C(N,s,\|b\|_{L^s(B_4)})$ such that
$$
\|u\|_{L^\infty(B_1)}\le C(N,s,\|b\|_{L^s(B_4)})\|u\|_{L^{2^*}(B_4)}.
$$
Here, the constants $C(N,r, t_\varepsilon)$ and $C(N,s,\|b\|_{L^s(B_4)})$ in (i) and (ii) remain bounded as long as $r, t_\varepsilon$ and $\|b\|_{L^s(B_4)}$ are bounded.}

\vskip 5mm

\section*{4. Proof of Theorem 2.1}

In this section, we always assume that $p=\frac{N+\alpha}{N}$, $q\in (2, 2+\frac{4}{N})$ and $\lambda>0$ is a small parameter.

It is easy to see that under the rescaling 
$$
w(x)=\lambda^{-\frac{N}{4-N(q-2)}}v(\lambda^{-\frac{2}{4-N(q-2)}}x), 
\eqno(4.1)
$$
the equation $(Q_\lambda)$ is reduced to 
$$
-\lambda^\sigma\Delta w+ w=(I_\alpha\ast |w|^p)|w|^{p-2}w+\lambda^\sigma|w|^{q-2}w, 
\eqno(\bar Q_\lambda)
$$
where $ \sigma:=\frac{4}{4-N(q-2)}>1$.       The corresponding functional is given by 
$$
J_\lambda(w):=\frac{1}{2}\int_{\mathbb R^N}\lambda^\sigma|\nabla w|^2+|w|^2-\frac{1}{2p}\int_{\mathbb R^N}(I_\alpha\ast |w|^p)|w|^p-\frac{1}{q}\lambda^\sigma\int_{\mathbb R^N}|w|^q.
$$

{\bf Lemma 4.1.}  {\it  Let $\lambda>0, v\in H^1(\mathbb R^N)$ and $w$ is the rescaling (4.1) of $v$. Then

(1)  \ $ \ \|w\|_2^2=\|v\|_2^2, \quad \int_{\mathbb R^N}(I_\alpha\ast |w|^p)|w|^p= \int_{\mathbb R^N}(I_\alpha\ast |v|^p)|v|^p,$

(2) \  $\lambda^{\sigma}\|\nabla w\|_2^2=\|\nabla v\|_2^2, \quad   \lambda^\sigma\|w\|_q^q=\lambda\|v\|_q^q$,

(3) \   $I_\lambda(v)=J_\lambda(w)$. }

We define the Nehari manifolds as follows.
$$
\mathcal{N}_\lambda=
\left\{w\in H^1(\mathbb R^N)\setminus\{0\} \ \left | \ \lambda^\sigma\int_{\mathbb R^N}|\nabla w|^2+\int_{\mathbb R^N}|w|^2=\int_{\mathbb R^N}(I_\alpha\ast |w|^p)|w|^p+\lambda^\sigma
\int_{\mathbb R^N}|w|^q\  \right. \right\}
$$
and 
$$
\mathcal{N}_0=
\left\{w\in H^1(\mathbb R^N)\setminus\{0\} \ \left | \ \int_{\mathbb R^N}|w|^2=\int_{\mathbb R^N}(I_\alpha\ast |w|^p)|w|^p\  \right. \right\}. 
$$
Then 
$$
m_\lambda:=\inf_{w\in \mathcal {N}_\lambda}J_\lambda(w), \quad {\rm and} \quad 
m_0:=\inf_{u\in \mathcal {N}_0}J_0(u)
$$
are well-defined and positive. Moreover, $J_0$ is attained on $\mathcal N_0$ and 
$m_0:=\inf_{w\in \mathcal {N}_0}J_0(w)=\frac{\alpha}{2(N+\alpha)}S_1^{\frac{N+\alpha}{\alpha}}$.

 For $w\in H^1(\mathbb R^N)\setminus \{0\}$, we set
$$
\tau_1(w)=\frac{\int_{\mathbb R^N}|w|^2}{\int_{\mathbb R^N}(I_\alpha\ast |w|^p)|w|^p}.
\eqno(4.2)
$$
Then $(\tau_1(w))^{\frac{N}{2\alpha}}w\in \mathcal N_0$ for any $w\in H^1(\mathbb R^N)\setminus\{0\}$,  and $w\in \mathcal N_0$  if and only if $\tau_1(w)=1$.

Define the Poho\v{z}aev manifold as follows
$$
\mathcal P_\lambda:=\{w\in H^1(\mathbb R^N)\setminus\{0\} \   | \  P_\lambda(w)=0 \},
$$
where 
$$
\begin{array}{rcl}
P_\lambda(w):&=&\frac{\lambda^\sigma(N-2)}{2}\int_{\mathbb R^N}|\nabla w|^2+\frac{ N}{2}\int_{\mathbb R^N}|w|^2\\ \\
&\quad &-\frac{N+\alpha}{2p}\int_{\mathbb R^N}(I_\alpha\ast |w|^p)|w|^p-\frac{\lambda^\sigma N}{q}\int_{\mathbb R^N}|w|^q.
\end{array}
\eqno(4.3)
$$
Let $v_\lambda\in H^1(\mathbb R^N)$ be the ground state for $(Q_\lambda)$ and 
$$w_\lambda(x)=\lambda^{-\frac{N}{4-N(q-2)}}v_\lambda(\lambda^{-\frac{2}{4-N(q-2)}}x).$$ Then by Lemma 3.1, $w_\lambda\in \mathcal P_\lambda$. Moreover,  we have the following minimax characterizations for the least energy level $m_\lambda$.
$$
m_\lambda=\inf_{w\in H^1(\mathbb R^N)\setminus\{0\}}\sup_{t\ge 0}J_\lambda(tw)=\inf_{w\in H^1(\mathbb R^N)\setminus\{0\}}\sup_{t\ge 0}J_\lambda(w_t).
\eqno(4.4)
$$
In particular, we have $m_\lambda=J_\lambda(w_\lambda)=\sup_{t>0}J_\lambda(tw_\lambda)=\sup_{t>0}J_\lambda((w_\lambda)_t)$. A similar result also holds for $m_0$ and $J_0$.

{\bf Lemma 4.2.} {\it The rescaled family of solutions $\{w_\lambda\}$ is bounded in $H^1(\mathbb R^N)$.}
\begin{proof}   Since $\{w_\lambda\}$ is bounded in $L^2(\mathbb R^N)$, it suffices to show that it is also bounded 
in $D^{1,2}(\mathbb R^N)$.  By $w_\lambda\in \mathcal N_\lambda\cap\mathcal P_\lambda$,  we obtain
$$
\lambda^\sigma\int_{\mathbb R^N}|\nabla w_\lambda|^2+\int_{\mathbb R^N}|w_\lambda|^2=\int_{\mathbb R^N}(I_\alpha\ast |w_\lambda|^p)|w_\lambda|^p+\lambda^\sigma\int_{\mathbb R^N}|w_\lambda|^q,
$$
$$
\frac{\lambda^\sigma(N-2)}{2}\int_{\mathbb R^N}|\nabla w_\lambda|^2+\frac{N}{2}\int_{\mathbb R^N}|w_\lambda|^2=\frac{N+\alpha}{2p}\int_{\mathbb R^N}(I_\alpha\ast |w_\lambda|^p)|w_\lambda|^p+\frac{\lambda^\sigma N}{q}\int_{\mathbb R^N}|w_\lambda|^q.
$$
Therefore, we have 
$$
\|\nabla v_\lambda\|_2^2=\lambda^\sigma\|\nabla w_\lambda\|_2^2=\frac{N(q-2)}{2q}\lambda^\sigma\int_{\mathbb R^N}|w_\lambda|^q=\frac{N(q-2)}{2q}\lambda\int_{\mathbb R_N}|v_\lambda|^q.
$$
Particularly, we have 
$$
\|\nabla w_\lambda \|_2^2=\frac{N(q-2)}{2q}\|w_\lambda\|_q^q.
\eqno(4.5)
$$

By the Gagliardo-Nirenberg Inequality, we obtain 
$$
\|v_\lambda\|_q^q\le C\|\nabla v_\lambda\|_2^{\frac{N(q-2)}{2}}\|v_\lambda\|_2^{\frac{2N-q(N-2)}{2}}.
$$
Therefore, we get
$$
\|\nabla v_\lambda\|_2^{\frac{4-N(q-2)}{2}}\le C\frac{N(q-2)}{2q}\lambda\|v_\lambda\|_2^{\frac{2N-q(N-2)}{2}}.
$$
Hence, 
$$
\lambda^\sigma\|\nabla w_\lambda\|_2^2=\|\nabla v_\lambda\|_2^2\le \tilde C\lambda^\sigma \|v_\lambda\|_2^{\frac{2[2N-q(N-2)]}{4-N(q-2)}},
$$
which together with the boundedness of $\|v_\lambda\|_2$ implies that $w_\lambda$ is bounded in $D^{1,2}(\mathbb R^N)$. 
\end{proof}

Now, we give  the following estimation on $\tau_1(w_\lambda)$ and  the least energy $m_\lambda$:

{\bf Lemma 4.3.} {\it $1<\tau_1(w_\lambda)\le 1+O(\lambda)$ and $m_\lambda=m_0+O(\lambda)$ as $\lambda \to 0$}.
\begin{proof} First, since $w_\lambda\in \mathcal{N}_\lambda$,  by (4.5), it follows that 
$$
\tau_1(w_\lambda)=\frac{\int_{\mathbb R^N}|w_\lambda|^2}{\int_{\mathbb R^N}(I_\alpha\ast |w_\lambda|^p)|w_{\lambda}|^p}
>  \frac{\int_{\mathbb R^N}\lambda^\sigma|\nabla w_\lambda|^2+|w_\lambda|^2-\lambda^\sigma\int_{\mathbb R^N}|w_\lambda|^q}{\int_{\mathbb R^N}(I_\alpha\ast |w_\lambda|^p)|w_{\lambda}|^p}=1
\eqno(4.6)
$$
and   by Lemma 4.1 and the Sobolev  inequality, we have
$$
\begin{array}{rcl}
\tau_1(w_\lambda)&=&\frac{\int_{\mathbb R^N}|w_\lambda|^2}{\int_{\mathbb R^N}(I_\alpha\ast |w_\lambda|^p)|w_{\lambda}|^p}\\ &\le& \frac{\int_{\mathbb R^N}\lambda^\sigma|\nabla w_\lambda|^2+|w_\lambda|^2}{\int_{\mathbb R^N}(I_\alpha\ast |w_\lambda|^p)|w_{\lambda}|^p+\lambda^\sigma\int_{\mathbb R^N}|w_\lambda|^q-\lambda^\sigma\int_{\mathbb R^N}|w_\lambda|^q}\\
&=&\frac{\|v_\lambda\|^2}{\|v_\lambda\|^2-\lambda\|v_\lambda\|_q^q}\le \frac{1}{1-\lambda C\|v_\lambda\|^{q-2}}.
\end{array}
\eqno(4.7)
$$
Since $\|v_\lambda\|$ is bounded, it follows that $1<\tau_1(w_\lambda)\le 1+O(\lambda)$ as $\lambda\to 0$.

For $w\in H^1(\mathbb R^N)$, let 
$$
w_t(x)=\left\{\begin{array}{cl}
w(x/t) \     & \ t>0,\\
0, \     & t=0.
\end{array}\right.
$$
Then by Lemma 3.2 and Poho\v{z}aev's identity, it is easy to show that $\sup_{t\ge 0}J_\lambda((w_\lambda)_t)=J_\lambda(w_\lambda)=m_\lambda$. Therefore, we get 
$$
\begin{array}{rcl}
m_0&\le& \sup_{t\ge 0}J_0((w_\lambda)_t)=J_0((w_\lambda)_t)|_{t=\tau_1((w_\lambda))^{1/\alpha}}\\
&\le & \sup_{t\ge 0}J_\lambda((w_\lambda)_t)+\lambda^\sigma(\tau_1(w_\lambda))^{N/\alpha}\int_{\mathbb R^N}|w_\lambda|^q\\
&\le& m_\lambda+\lambda(1+O(\lambda))^{N/\alpha}\|v_\lambda\|_q^q\\
&=& m_\lambda+O(\lambda).
\end{array}
$$

On the other hand,  let $U\in \mathcal N\subset H^1(\mathbb R^N)$ be such that
$$
S_1=\frac{\int_{\mathbb R^N}|U|^2}{\left(\int_{\mathbb R^N}(I_\alpha\ast |U|^p)|U|^p\right)^{1/p}}.
$$
Then $\int_{\mathbb R^N}|U|^2=\int_{\R^N}(I_\alpha\ast |U|^{\frac{N+\alpha}{N}})|U|^{\frac{N+\alpha}{N}}=S_1^{\frac{N+\alpha}{\alpha}}$ and $m_0=J_0(U)=\frac{\alpha}{2(N+\alpha)}S_1^{\frac{N+\alpha}{\alpha}}$. 
By Lemma 3.2 again, we obtain
$$
\begin{array}{rcl}
m_\lambda&\le& \sup_{t\ge 0}J_\lambda(tU)\\
&=&
\sup_{t\ge 0}\left\{\frac{t^2}{2}\int_{\mathbb R^N}\lambda^\sigma|\nabla U|^2+|U|^2-\frac{t^{2p}}{2p}\int_{\mathbb R^N}(I_\alpha\ast |U|^p)|U|^p
-\frac{\lambda^\sigma t^q}{q}\int_{\mathbb R^N}|U|^q\right\}\\
&\le & \sup_{t\ge 0}\left\{\frac{t^2}{2}\int_{\mathbb R^N}|U|^2-\frac{t^{2p}}{2p}\int_{\mathbb R^N}(I_\alpha\ast |U|^p)|U|^p\right\}
+\lambda^\sigma\sup_{t\ge 0}\left\{\frac{t^2}{2}\int_{\mathbb R^N}|\nabla U|^2-\frac{t^q}{q}\int_{\mathbb R^N}|U|^q\right\}\\
&=& J_0(U)+O(\lambda^\sigma)\\
&=&m_0+o(\lambda).
\end{array}
$$
The proof is completed. 
\end{proof}

{\bf Lemma 4.4.} {\it  $\|w_\lambda\|_2^2\sim1$ as $\lambda\to 0$.}
\begin{proof}  
By the definition of $\tau_1(w_\lambda)$, Lemma 4.3  and the  Hardy-Sobolev-Littlewood inequality, for small $\lambda>0$, we have 
$$
\|w_\lambda\|_2^2=\tau_1(w_\lambda)\int_{\mathbb R^N}(I_\alpha\ast |w_\lambda|^p)|w_\lambda|^p\le 2S_1^{-p}\|w_\lambda\|_2^{2p},
$$
and thus it follows that 
$$
\|w_\lambda\|_2^2\ge 2^{-\frac{N}{\alpha}}S_1^{\frac{N+\alpha}{\alpha}},
$$
which together with the boundedness of $w_\lambda$ implies that $\|w_\lambda\|_2^2\sim1$ as $\lambda\to 0$. The proof is completed. 
\end{proof}

Now, we give  the following estimation on the least energy:

{\bf Lemma 4.5.} {\it  Let $N\ge 3$ and $q\in (2,2+\frac{4}{N})$, then
  $$
  m_0-m_\lambda \sim \lambda^\sigma\qquad  as \quad \lambda\to 0.
  $$}
\begin{proof}  By Lemma 3.2, Lemma 4.4 and the boundedness of $\{w_\lambda\}$,  we find
$$
\begin{array}{rcl}
m_0&\le& \sup_{t\ge 0}J_0((w_\lambda)_t)=J_0((w_\lambda)_{t_\lambda})\\
&\le &\sup_{t\ge 0} J_\lambda((w_\lambda)_t)+\lambda^\sigma  \left(\frac{t_\lambda^N}{q}\int_{\mathbb R^N}|w_\lambda|^q-\frac{t_\lambda^{N-2}}{2}\int_{\mathbb R^N}|\nabla w_\lambda|^{2}\right)\\
&\le &m_\lambda +C\lambda^\sigma,
\end{array}
\eqno(4.8)
$$
where 
$$
t_\lambda=\left(\frac{\int_{\mathbb R^N}| w_\lambda|^2}{\int_{\mathbb R^N}(I_\alpha\ast |w_\lambda|^p)|w_\lambda|^p}\right)^{\frac{1}{\alpha}}=(\tau_1(w_\lambda))^{\frac{1}{\alpha}}.
$$

For each $\rho>0$, the family  $U_\rho(x):=\rho^{-\frac{N}{2}}U_1(x/\rho)$ are radial ground states of  $ v=(I_\alpha\ast |v|^p)v^{p-1}$, and verify that
 $$
 \|\nabla U_\rho\|_2^2=\rho^{-2}\|\nabla U_1\|_2^2, \qquad \int_{\mathbb R^N}|U_\rho|^q=\rho^{N-\frac{N}{2}q}\int_{\mathbb R^N}|U_1|^q.
 \eqno(4.9)
 $$ 
 Let $g_0(\rho)=\frac{1}{q}\int_{\mathbb R^N}|U_\rho|^q-\frac{1}{2}\int_{\mathbb R^N}|\nabla U_\rho|^2$. Then there exists $\rho_0=\rho(q)\in (0,+\infty)$ with
 $$
 \rho_0=\left(\frac{2q\int_{\mathbb R^N}|\nabla U_1|^2}{N(q-2)\int_{\mathbb R^N}|U_1|^q}\right)^{\frac{2}{4-N(q-2)}}
 $$
such that 
 $$
 g_0(\rho_0)=\sup_{\rho>0}g_0(\rho)=\frac{4-N(q-2)}{2N(q-2)}\left(\frac{N(q-2)\int_{\mathbb R^N}|U_1|^q}{2q\int_{\mathbb R^N}|\nabla U_1|^2}\right)^{\frac{4}{4-N(q-2)}}\int_{\mathbb R^N}|\nabla U_1|^2.
 $$ 
 Let $U_0=U_{\rho_0}$, then there exists $t_\lambda\in (0,+\infty)$ such that
 $$
\begin{array}{rcl}
m_\lambda&\le &\sup_{t\ge 0}J_\lambda(tU_0)=J_\lambda(t_\lambda U_0)\\
&=&\frac{t^2_\lambda}{2}\int_{\mathbb R^N}| U_0|^2-\frac{t^{2p}_\lambda}{2p}\int_{\R^N}(I_\alpha\ast |U_0|^p)|U_0|^p-\lambda^\sigma\{\frac{t^{q}_\lambda}{q}\int_{\mathbb R^N}|U_0|^{q}-\frac{t^2_\lambda}{2}\int_{\mathbb R^N}|\nabla U_0|^2\}\\
&\le &\sup_{t\ge 0}\left(\frac{t^2}{2}-\frac{t^{2p}}{2p}\right)\int_{\mathbb R^N}|U_0|^2 -\lambda^\sigma\{\frac{t^{q}_\lambda}{q}\int_{\mathbb R^N}|U_0|^{q}-\frac{t^2_\lambda}{2}\int_{\mathbb R^N}|\nabla U_0|^2\}\\&=& m_0-\lambda^\sigma\{\frac{t^{q}_\lambda}{q}\int_{\mathbb R^N}|U_0|^{q}-\frac{t^2_\lambda}{2}\int_{\mathbb R^N}|\nabla U_0|^2\}.
\end{array}
\eqno(4.10)
$$
If $t_\lambda\ge 1$, then
$$
\int_{\mathbb R^N}|U_0|^2+\lambda^\sigma\int_{\mathbb R^N}|\nabla U_0|^2\ge t_\lambda^{\min\{\frac{2\alpha}{N}, q-2\}}\left\{\int_{\mathbb R^N}(I_\alpha\ast |U_0|^p)|U_0|^p+\lambda^\sigma\int_{\mathbb R^N}|U_0|^q\right\}.
$$
Hence
$$
t_\lambda\le \left(\frac{\int_{\mathbb R^N}| U_0|^2+\lambda^\sigma\int_{\mathbb R^N}|\nabla U_0|^2}{\int_{\mathbb R^N}(I_\alpha\ast |U_0|^p)|U_0|^p+\lambda^\sigma\int_{\mathbb R^N} |U_0|^q}\right)^{\frac{1}{\min\{\frac{2\alpha}{N}, q-2\} }}.
$$
If $t_\lambda\le  1$, then
$$
\int_{\mathbb R^N}|U_0|^2+\lambda^\sigma\int_{\mathbb R^N}|\nabla U_0|^2\le t_\lambda^{\min\{\frac{2\alpha}{N}, q-2\}}\left\{\int_{\mathbb R^N}(I_\alpha\ast |U_0|^p)|U_0|^p+\lambda^\sigma\int_{\mathbb R^N}|U_0|^q\right\}.
$$
Hence
$$
t_\lambda\ge \left(\frac{\int_{\mathbb R^N}|U_0|^2+\lambda^\sigma\int_{\mathbb R^N}|\nabla U_0|^2}{\int_{\mathbb R^N}(I_\alpha\ast |U_0|^p)|U_0|^p+\lambda^\sigma\int_{\mathbb R^N}|U_0|^q}\right)^{\frac{1}{\min\{\frac{2\alpha}{N}, q-2\}}}.
$$
Since
 $$
 \int_{\mathbb R^N}(I_\alpha\ast |U_0|^p)|U_0|^p=\int_{\mathbb R^N}|U_0|^2   \qquad {\rm and} \qquad \int_{\mathbb R^N}|U_0|^q>\int_{\mathbb R^N}|\nabla U_0|^2,
 $$
we conclude that 
$$
\left(\frac{\int_{\mathbb R^N}|U_0|^2+\lambda^\sigma\int_{\mathbb R^N}|\nabla U_0|^2}{\int_{\mathbb R^N}(I_\alpha\ast |U_0|^p)|U_0|^p+\lambda^\sigma\int_{\mathbb R^N}|U_0|^q}\right)^{\frac{1}{\min\{\frac{2\alpha}{N}, q-2\} }}\le t_\lambda
\le 1.
\eqno(4.11)
$$
Therefore, $\lim_{\lambda\to 0}t_\lambda=1$ and hence there exists a constant $C>0$ such that
$$
m_\lambda\le m_0-C\lambda^\sigma,
$$
for small $\lambda>0$. The proof is complete. 
\end{proof}

Let 
$$
\mathbb D(u):=\int_{\mathbb R^N}(I_\alpha\ast |u|^{\frac{N+\alpha}{N}})|u|^{\frac{N+\alpha}{N}}.
$$
The following result is  a special case of the classical Brezis-Lieb lemma \cite{Brezis-1} for Riesz potentials, for a proof, we refer the reader to \cite[Lemma 2.4]{Moroz-2}.

{\bf Lemma 4.6.} {\it Let $N \in\mathbb  N, \alpha\in  (0,N)$, and $(w_n)_{n\in \mathbb N}$ be a bounded sequence in $L^2(\mathbb R^N)$. If $w_n \to  w$ almost everywhere on $\mathbb R^N$ as $n\to \infty$, then
$$
\lim_{n\to\infty}\mathbb{D}(w_n)-\mathbb D(w_n-w_0) = \mathbb D(w_0).
$$}

{\bf Lemma 4.7.} {\it $\|w_\lambda\|_2^2\sim \mathbb D(w_\lambda)\sim \|\nabla w_\lambda \|_2^2\sim \|w_\lambda\|_q^q\sim 1$ as $\lambda\to 0$.}
\begin{proof}  
It   follows from Lemma 4.3  that
$$
\begin{array}{rcl}
m_0&\le &J_0((\tau_1(w_\lambda))^{\frac{N}{2\alpha}}w_\lambda)\\
&=&\frac{1}{2}(\tau_1(w_\lambda))^{\frac{N}{\alpha}}\|w_\lambda\|_2^2-\frac{1}{2p}(\tau_1(w_\lambda))^{\frac{N+\alpha}{\alpha}}\int_{\mathbb R^N}(I_\alpha\ast |w_\lambda|^p)|w_\lambda|^p\\
&\le& (\tau_1(w_\lambda))^{\frac{N}{\alpha}}\left[\frac{1}{2}\|w_\lambda\|_2^2-\frac{1}{2p}\int_{\mathbb R^N}(I_\alpha\ast |w_\lambda|^p)|w_\lambda|^p\right].
\end{array}
\eqno(4.12)
$$
Since $w_\lambda\in \mathcal N_\lambda$,  by Lemma 4.4, we obtain
$$
\begin{array}{rcl}
\tau_1(w_\lambda)&=&\frac{\int_{\mathbb R^N}|w_\lambda|^2}{\int_{\mathbb R^N}(I_\alpha\ast |w_\lambda|^p)|w_\lambda|^p}\\
&=&\frac{\int_{\mathbb R^N}|w_\lambda|^2}{\int_{\mathbb R^N}|w_\lambda|^2+\lambda^\sigma(\int_{\mathbb R^N}|\nabla w_\lambda|^2-\int_{\mathbb R^N}|w_\lambda|^q)}\\
&=&\frac{\int_{\mathbb R^N}|w_\lambda|^2}{\int_{\mathbb R^N}|w_\lambda|^2-\lambda^\sigma\frac{2N-q(N-2)}{2q}\int_{\mathbb R^N}|w_\lambda|^q}\\
&\le & 1+C_1\lambda^\sigma \|w_\lambda\|_q^q.
\end{array}
\eqno(4.13)
$$
Therefore, by (4.5), (4.12) and  (4.13),  we obtain 
$$
\begin{array}{rcl}
m_\lambda&=&\frac{1}{2}\int_{\mathbb R^N}\lambda^{\sigma}|\nabla w_\lambda|_2^2+|w_\lambda|_2^2-\frac{1}{2p}\int_{\mathbb R^N}(I_\alpha\ast |w_\lambda|^p)|w_\lambda|^p-\frac{\lambda^\sigma}{q}\int_{\mathbb R^N}|w_\lambda|^q\\
&\ge& \lambda^\sigma(\frac{1}{2}\|\nabla w_\lambda\|_2^2-\frac{1}{q}\|w_\lambda\|_q^q)+\frac{m_0}{(\tau_1(w_\lambda))^{\frac{N}{\alpha}}}\\
&\ge &\lambda^\sigma\frac{N(q-2)-4}{4q}\|w_\lambda\|_q^q+m_0-C_1\lambda^\sigma\|w_\lambda\|_q^q.
\end{array}
$$
Recall that by Lemma 4.5, we have $m_0-m_\lambda\ge C_2\lambda^\sigma$. Hence,  we get
$$
\frac{4-N(q-2)}{4q}\|w_\lambda\|_q^q\ge C_2-C_1\|w_\lambda\|_q^q, 
$$
which yields 
$$
\|w_\lambda\|_q^q\ge \frac{4qC_2}{4-N(q-2)+4qC_1}>0.
$$
Since $w_\lambda$ is bounded in $H^1(\mathbb R^N)$, it follows that $\|w_\lambda\|_q^q\sim 1$ as $\lambda\to 0$.

Since $\|w_\lambda\|_2^2\sim 1$ as $\lambda\to 0$, the Gagliardo-Nirenberg inequality implies 
$$
\|w_\lambda\|_q^q\le C\|\nabla w_\lambda\|_2^{\frac{N(q-2)}{2}},
$$
which together with the boundedness of $w_\lambda$ in $H^1(\mathbb R^N)$  yields $\|\nabla w_\lambda \|_2^2\sim 1$ as $\lambda\to 0$.

Finally, by the definition of $\tau_1(w_\lambda)$, Lemma 4.3 and Lemma 4,4,  it follows  that 
$$
\mathbb D(w_\lambda)=\int_{\mathbb R^N}(I_\alpha\ast |w_\lambda|^{\frac{N+\alpha}{N}})|w_\lambda|^{\frac{N+\alpha}{N}}=(\tau_1(w_\lambda))^{-1} \|w_\lambda\|_2^2\sim 1,
$$
as $\lambda\to 0$. The proof is complete. 
\end{proof}

{\bf Lemma 4.8.} {\it Let $N\ge 3$ and $q\in (2, 2+\frac{4}{N})$,  then for any $\lambda_n\to 0$, there exists $\rho\in [\rho_0, +\infty)$, such that, up to a subsequence, $w_{\lambda_n}\to U_{\rho}$ in $L^2(\mathbb R^N)$, where 
$$
  \rho_0=\rho_0(q):= \left(\frac{2q\int_{\mathbb R^N}|\nabla U_1|^2}{N(q-2)\int_{\mathbb R^N}|U_1|^q}\right)^{\frac{2}{4-N(q-2)}}.
 $$
Moreover, $w_{\lambda_n}\to U_\rho$ in $D^{1,2}(\mathbb R^N)$ if and only if $\rho=\rho_0$. }
\begin{proof} Note that $w_{\lambda_n}$ is a positive radially symmetric function, and by Lemma 4.2, $\{w_{\lambda_n}\}$ is bounded in $H^1(\mathbb R^N)$. Then there exists $w_0\in H^1(\mathbb R^N)$ such that 
$$
w_{\lambda_n} \rightharpoonup w_0   \quad {\rm weakly \ in} \  H^1(\mathbb R^N), \quad w_{\lambda_n}\to w_0 \quad {\rm in} \ L^p(\mathbb R^N) \quad {\rm for \ any} \ p\in (2,2^*),
\eqno(4.14)
$$
and 
$$
w_{\lambda_n}(x)\to w_0(x) \quad a. \ e. \  {\rm on} \ \mathbb R^N,  \qquad w_{\lambda_n}\to w_0 \quad {\rm in} \   L^2_{loc}(\mathbb R^N).
\eqno(4.15)
$$
Observe that
$$
J_0(w_{\lambda_n})=J_{\lambda_n}(w_{\lambda_n})+\frac{\lambda_n^\sigma}{q}\int_{\mathbb R^N}|w_{\lambda_n}|^q-\frac{\lambda_n^\sigma}{2}\int_{\mathbb R^N}|\nabla w_{\lambda_n}|^2=m_{\lambda_n}+o_n(1)=m_0+o_n(1),
$$
and 
$$
J'_0(w_{\lambda_n})w=J'_{\lambda_n}(w_{\lambda_n})w+\lambda_n^\sigma\int_{\mathbb R^N}|w_{\lambda_n}|^{q-2}w_{\lambda_n} w-\lambda_n^\sigma\int_{\mathbb R^N}\nabla w_{\lambda_n} \nabla w=o_n(1).
$$
Therefore, $\{w_{\lambda_n}\}$ is a PS sequence of $J_0$ at level $m_0=\frac{\alpha}{2(N+\alpha)}S_1^{\frac{N+\alpha}{\alpha}}$.

By Lemma 4.7, we have $w_0\not=0$, and hence $m_0\le J_0(w_0)$.
By Lemma 4.5 and Lemma 4.6, we have
$$
\begin{array}{rcl}
o_n(1)&=&m_{\lambda_n}-m_0\\
&\ge &\frac{\lambda_n^\sigma}{2}\|\nabla w_{\lambda_n}\|_2^2+\frac{1}{2}[\|w_{\lambda_n}\|_2^2-\|w_0\|_2^2]-\frac{N}{2(N+\alpha)}[\mathbb D(w_{\lambda_n})-\mathbb D(w_0)]-\frac{\lambda_n^\sigma}{q}\|w_{\lambda_n}\|_q^q\\
&=&\frac{1}{2}\|w_{\lambda_n}-w_0\|_2^2-\frac{N}{2(N+\alpha)}\mathbb D(w_{\lambda_n}-w_0)+o_n(1),
\end{array}
$$
$$
\begin{array}{rcl}
0&=&\langle J_{\lambda_n}'(w_{\lambda_n})-J_0'(w_0), w_{\lambda_n}-w_0\rangle\\
&=& \lambda_n^\sigma\int_{\mathbb R^N}\nabla w_{\lambda_n}(\nabla w_{\lambda_n}-\nabla w_0)+\int_{\mathbb R^N}|w_{\lambda_n}-w_0|^2\\
&\mbox{}&\quad 
-\int_{\mathbb R^N}(I_\alpha\ast |w_{\lambda_n}|^{\frac{N+\alpha}{N}})w_{\lambda_n}^{\frac{\alpha}{N}}(w_{\lambda_n}-w_0)
+\int_{\mathbb R^N}(I_\alpha\ast |w_0|^{\frac{N+\alpha}{N}})w_0^{\frac{\alpha}{N}}(w_{\lambda_n}-w_0)\\
&\mbox{}&\quad
-\lambda_n^\sigma\int_{\mathbb R^N}|w_{\lambda_n}|^{q-2}w_{\lambda_n}(w_{\lambda_n}-w_0)\\
&=& \|w_{\lambda_n}-w_0\|_2^2-\mathbb D(w_{\lambda_n}-w_0)+o_n(1).
\end{array}
$$
Hence, it follows that
$$
\|w_{\lambda_n}-w_0\|_2^2\le \frac{N}{N+\alpha}\mathbb D(w_{\lambda_n}-w_0)+o_n(1)=\frac{N}{N+\alpha}\|w_{\lambda_n}-w_0\|_2^2+o_n(1),
$$
and hence
$$
\|w_{\lambda_n}-w_0\|_2\to 0, \qquad  {\rm as} \  \  \lambda_n\to 0.
$$
By the Hardy-Littlewood-Sobolev inequality and Lemma 4.6, it follows that 
$$
\lim_{\lambda_n\to 0}\mathbb D(w_{\lambda_n})=\mathbb D(w_0).
$$
Since $\tau_1(w_{\lambda_n})\to 1$ as $\lambda_n\to 0$ by Lemma 4.3,  it follows that $w_0\in \mathcal N_0$.

On the other hand, by Lemma 4.1 and the boundedness of $v_{\lambda_n}$ in $H^1(\mathbb R^N)$, we have 
$$
\begin{array}{rcl}
m_{\lambda_n}&=&J_{\lambda_n}(w_{\lambda_n})\\
&=&\frac{1}{2}\int_{\mathbb R^N}\lambda_n^\sigma|\nabla w_{\lambda_n}|^2+|w_{\lambda_n}|^2-\frac{1}{2p}\int_{\mathbb R^N}
(I_\alpha\ast |w_{\lambda_n}|^{\frac{N+\alpha}{N}})|w_{\lambda_n}|^{\frac{N+\alpha}{N}}-\frac{\lambda_n^\sigma}{q}\int_{\mathbb R^N}|w_{\lambda_n}|^q\\
&=&  \frac{1}{2}\int_{\mathbb R^N}|\nabla v_{\lambda_n}|^2+|w_{\lambda_n}|^2-\frac{1}{2p}\int_{\mathbb R^N}
(I_\alpha\ast |w_{\lambda_n}|^{\frac{N+\alpha}{N}})|w_{\lambda_n}|^{\frac{N+\alpha}{N}}-\frac{\lambda_n}{q}\int_{\mathbb R^N}|v_{\lambda_n}|^q\\
&\ge &  \frac{1}{2}\int_{\mathbb R^N}|w_{\lambda_n}|^2-\frac{1}{2p}\int_{\mathbb R^N}
(I_\alpha\ast |w_{\lambda_n}|^{\frac{N+\alpha}{N}})|w_{\lambda_n}|^{\frac{N+\alpha}{N}}-C\lambda_n.
\end{array}
$$
Sending $\lambda_n\to 0$, it then follows from Lemma 4.5 that
$$
m_0\ge \frac{1}{2}\int_{\mathbb R^N}|w_0|^2-\frac{1}{2p}\int_{\mathbb R^N}
(I_\alpha\ast |w_0|^{\frac{N+\alpha}{N}})|w_0|^{\frac{N+\alpha}{N}}
=J_0(w_0).
$$
Therefore, note that $w_0\in \mathcal N_0$, we obtain $J_0(w_0)=m_0$. 
Thus, $w_0=U_{\rho}$ for some $\rho\in (0,+\infty)$.

 Moreover, by (4.5), we obtain
$$
\|\nabla w_0\|_2^2\le \lim_{\lambda_n\to 0}\|\nabla w_{\lambda_n}\|_2^2=\frac{N(q-2)}{2q}\int_{\mathbb R^N}|w_0|^q,
$$
from which it follows that  
 $$
 \rho\ge\left(\frac{2q\int_{\mathbb R^N}|\nabla U_1|^2}{N(q-2)\int_{\mathbb R^N}|U_1|^q}\right)^{\frac{2}{4-N(q-2)}}.
 $$
If $\rho=\rho_0$, then (4.5) implies that $\lim_{n\to\infty}\|\nabla w_{\lambda_n}\|_2^2=\|\nabla U_{\rho_0}\|_2^2$, and hence $w_{\lambda_n}\to U_{\rho_0}$ in $D^{1,2}(\mathbb R^N)$.
\end{proof}

\begin{proof}[Proof of Theorem 2.1.]   
Let 
$$
M_\lambda=w_\lambda(0),  \quad z_\lambda=M_\lambda [U_{\rho_0}(0)]^{-1},  
$$
where $\rho_0$ is given in Lemma 4.8. We further perform a scaling 
$$
\tilde w_\lambda(x)=z_\lambda^{-1}w_\lambda(z_\lambda^{-\frac{2}{N}}x),
$$
then 
$$
\tilde w_\lambda(0)=z_\lambda^{-1}w_\lambda(0)=U_{\rho_0}(0)M_\lambda^{-1}w_\lambda(0)=U_{\rho_0}(0),
$$
and $\tilde w_\lambda$ satisfies  the rescaled equation 
$$
-\lambda^\sigma z_\lambda^{\frac{4}{N}}\Delta \tilde w+ \tilde w=(I_\alpha\ast |\tilde w|^p)|\tilde w|^{p-2}\tilde w+\lambda^\sigma z_\lambda^{q-2}|\tilde w|^{q-2}\tilde w, 
$$
The corresponding functional is given by 
$$
 J_\lambda( w_\lambda)=\frac{1}{2}\int_{\mathbb R^N}\lambda^\sigma z_\lambda^{\frac{4}{N}}|\nabla \tilde w_\lambda|^2+|\tilde w_\lambda|^2-\frac{1}{2p}\int_{\mathbb R^N}(I_\alpha\ast |\tilde w_\lambda|^p)|\tilde w_\lambda|^p-\frac{1}{q}\lambda^\sigma z_\lambda^{q-2}\int_{\mathbb R^N}|\tilde w_\lambda|^q.
$$
Moreover, we have 

(1)  \ $ \ \|\tilde w\|_2^2=\|w\|_2^2, \quad \int_{\mathbb R^N}(I_\alpha\ast |\tilde w|^p)|\tilde w|^p= \int_{\mathbb R^N}(I_\alpha\ast |w|^p)|w|^p,$

(2) \  $z_\lambda^{\frac{4}{N}}\|\nabla \tilde w\|_2^2=\|\nabla w\|_2^2, \quad   z_\lambda^{q-2}\|\tilde w\|_q^q=\|w\|_q^q$.

By Lemma 4.8, for any $\lambda_n\to 0$, there exists $\rho\ge \rho_0$ such that
$$
M_{\lambda_n}=w_{\lambda_n}(0)\to U_{\rho}(0)=\rho^{-\frac{N}{2}}U_1(0)\le \rho_0^{-\frac{N}{2}}U_1(0)<+\infty,
$$
which yields that $M_\lambda\le C$ for some $C>0$ and any small $\lambda>0$. 

Suppose that there exists a sequence $\lambda_n\to 0$ such that $M_{\lambda_n}\to 0$.
By Lemma 4.8,  up to a subsequence, $M_{\lambda_n}=w_{\lambda_n}(0)\to U_{\rho}(0)\not=0$ for some $\rho\in (0,+\infty)$.
This leads to a contradiction. Therefore, there exists some $c>0$ such that
$M_\lambda\ge c>0$.

Let 
$$
\zeta_\lambda =z_\lambda^{-\frac{2}{N}} \lambda^{-\frac{2}{4-N(q-2)}}.
$$
Then 
$$
\zeta_\lambda \sim \lambda^{-\frac{2}{4-N(q-2)}}
$$
and  for small $\lambda>0$, the rescaled family of ground states
$$
\tilde w_\lambda(x)=\zeta_\lambda^{\frac{N}{2}}v_\lambda(\zeta_\lambda x)
$$
satisfies 
$$
\|\nabla \tilde w_\lambda\|^2_2\sim \|\tilde w_\lambda\|_{q}^{q} \sim \int_{\mathbb R^N}(I_\alpha \ast |\tilde w_\lambda|^{\frac{N+\alpha}{N}})|\tilde w_\lambda|^{\frac{N+\alpha}{N}}\sim \|\tilde w_\lambda\|_2^2\sim 1,
$$
and as $\lambda\to 0$, $\tilde w_\lambda$ converges in $L^2(\mathbb R^N)$ to the extremal function  $U_{\rho_0}$. 
Then by Lemma 4.8, we also have $\tilde w_\lambda\to U_{\rho_0}$ in $D^{1,2}(\mathbb R^N)$. Thus we conclude that $\tilde w_\lambda\to U_{\rho_0}$ in $H^1(\mathbb R^N)$. 
 
 Since $w_\lambda\in \mathcal N_\lambda$, it follows that
$$
\begin{array}{lcl}
m_\lambda &=&(\frac{1}{2}-\frac{1}{2p})\lambda^\sigma\int_{\R^N}|\nabla  w_\lambda|^2+(\frac{1}{2}-\frac{1}{2p})\int_{\R^N}|w_\lambda|^2-(\frac{1}{q}
-\frac{1}{2p})\lambda^\sigma\int_{\R^N}|w_\lambda|^q\\
&=&\frac{\alpha}{2(N+\alpha)}\lambda^\sigma\int_{\R^N}|\nabla  w_\lambda|^2+\frac{\alpha}{2(N+\alpha)}\int_{\R^N}|w_\lambda|^2-\frac{2p-q}{2pq}\lambda^\sigma\int_{\R^N}|w_\lambda|^q\\
&=&\frac{\alpha}{2(N+\alpha)}\int_{\R^N}|w_\lambda|^2+O(\lambda^\sigma).
\end{array}
$$
Similarly, we also have 
$$
m_0=\frac{\alpha}{2(N+\alpha)}\int_{\R^N}|U_1|^2.
$$
Then it follows  from Lemma 4.5 that 
$$
\int_{\R^N}| U_1|^2-\int_{\R^N}|w_\lambda|^2=\frac{2(N+\alpha)}{\alpha}(m_0-m_\lambda)+O(\lambda^\sigma)=O(\lambda^\sigma).
$$
Since 
$
\| U_1\|_2^2=\int_{\R^N}(I_\alpha\ast |U_1|^p)|U_1|^p=S_1^{\frac{N+\alpha}{\alpha}},
$
we conclude  that 
$$
\|\tilde w_\lambda\|_2^2=\| w_\lambda\|_2^2=S_1^{\frac{N+\alpha}{\alpha}}+O(\lambda^\sigma).
$$
Finally, by $w_\lambda\in \mathcal N_\lambda$, we also have
$$
\int_{\R^N}(I_\alpha\ast |\tilde w_\lambda|^{\frac{N+\alpha}{N}})|\tilde w_\lambda|^{\frac{N+\alpha}{N}}=\int_{\R^N}(I_\alpha\ast |w_\lambda|^{\frac{N+\alpha}{N}})|w_\lambda|^{\frac{N+\alpha}{N}}=\|w_\lambda\|_2^2+O(\lambda^\sigma)=  S_1^{\frac{N+\alpha}{\alpha}}+O(\lambda^\sigma).
$$
The statements on $v_\lambda$ follow from the corresponding results on $w_\lambda$ and $\tilde w_\lambda$, and 
the proof is complete.  \end{proof}

\vskip 5mm

\section*{5.  Proof of Theorem 2.2}

 In this section, we always assume that $p=\frac{N+\alpha}{N-2}$ and  $q\in (2, 2^*)$ if $N\ge 4$, and $q\in (4,6)$ if $N=3$. It is easy to see that under the rescaling 
$$
w(x)=\lambda^{\frac{1}{q-2}}v(\lambda^{\frac{2^*-2}{2(q-2)}}x), 
\eqno(5.1)
$$
the equation $(Q_\lambda)$ is reduced to 
$$
-\Delta w+ \lambda^\sigma w=(I_\alpha\ast |w|^p)|w|^{p-2}w+\lambda^\sigma |w|^{q-2}w, 
\eqno(\bar Q_\lambda)
$$
where $\sigma:=\frac{2^*-2}{q-2}>1.$ 

The associated functional is defined by 
$$
J_\lambda(w):=\frac{1}{2}\int_{\mathbb R^N}|\nabla w|^2+\lambda^\sigma|w|^2-\frac{1}{2p}\int_{\mathbb R^N}(I_\alpha\ast |w|^p)|w|^p-\frac{1}{q}\lambda^\sigma\int_{\mathbb R^N}|w|^q.
$$

{\bf Lemma  5.1.}  {\it Let $\lambda>0, v\in H^1(\mathbb R^N)$ and $w$ is the rescaling (5.1) of $v$.  Then

(1) $ \    \   \|\nabla w\|_2^2=\|\nabla v\|_2^2, \quad \int_{\mathbb R^N}(I_\alpha\ast |w|^p)|w|^p= \int_{\mathbb R^N}(I_\alpha\ast |v|^p)|v|^p,$

(2)  $ \   \   \lambda^\sigma\|w\|_2^2=\|v\|_2^2,  \quad   \lambda^\sigma\|w\|_q^q=\lambda\|v\|_q^q$,

(3)  $\    \  I_\lambda(u)=J_\lambda(w)$.}

We define the Nehari manifolds as follows.
$$
\mathcal{N}_\lambda=
\left\{w\in H^1(\mathbb R^N)\setminus\{0\} \ \left | \ \int_{\mathbb R^N}|\nabla w|^2+\lambda^\sigma\int_{\mathbb R^N}|w|^2=\int_{\mathbb R^N}(I_\alpha\ast |u|^p)|w|^p+\lambda^\sigma
\int_{\mathbb R^N}|w|^q\  \right. \right\}
$$
and 
$$
\mathcal{N}_0=
\left\{w\in H^1(\mathbb R^N)\setminus\{0\} \ \left | \ \int_{\mathbb R^N}|\nabla w|^2=\int_{\mathbb R^N}(I_\alpha\ast |w|^p)|w|^p\  \right. \right\}. 
$$
Then 
$$
m_\lambda:=\inf_{w\in \mathcal {N}_\lambda}J_\lambda(w), \quad {\rm and} \quad 
m_0:=\inf_{u\in \mathcal {N}_0}J_0(u)
$$
are well-defined and positive. Moreover, $J_0$ is attained on $\mathcal N_0$.

 For $w\in H^1(\mathbb R^N)\setminus \{0\}$, we set
$$
\tau_2(w)=\frac{\int_{\mathbb R^N}|\nabla w|^2}{\int_{\mathbb R^N}(I_\alpha\ast |w|^p)|w|^p}.
\eqno(5.2)
$$
Then $(\tau_2(w))^{\frac{N-2}{2(2+\alpha)}}w\in \mathcal N_0$ for any $w\in H^1(\mathbb R^N)\setminus\{0\}$,  and $w\in \mathcal N_0$  if and only if $\tau_2(w)=1$.

Define the Pohozaev manifold as follows
$$
\mathcal P_\lambda:=\{w\in H^1(\mathbb R^N)\setminus\{0\} \   | \  P_\lambda(w)=0 \},
$$
where 
$$
\begin{array}{rcl}
P_\lambda(w):&=&\frac{N-2}{2}\int_{\mathbb R^N}|\nabla w|^2+\frac{ \lambda^\sigma N}{2}\int_{\mathbb R^N}|w|^2\\ \\
&\quad &-\frac{N+\alpha}{2p}\int_{\mathbb R^N}(I_\alpha\ast |w|^p)|w|^p-\frac{\lambda^\sigma N}{q}\int_{\mathbb R^N}|w|^q.
\end{array}
\eqno(5.3)
$$
Then by Lemma 3.1, $w_\lambda\in \mathcal P_\lambda$. Moreover,  we have a similar minimax characterizations for the least energy level $m_\lambda$ as in Lemma 3.2.

{\bf Lemma 5.2.} {\it The rescaled family of solutions $\{w_\lambda\}$ is bounded in $H^1(\mathbb R^N)$.}

{\bf Proof.}  First, we show that $\{w_\lambda\}$ is bounded in $H^1(\mathbb R^N)$.  Since $\{w_\lambda\}$ is bounded in $D^{1,2}(\mathbb R^N)$, it suffices to show that it is also bounded 
in $L^2(\mathbb R^N)$.  By $w_\lambda\in \mathcal N_\lambda\cap \mathcal P_\lambda$, we obtain
$$
\int_{\mathbb R^N}|\nabla w_\lambda|^2+\lambda^\sigma \int_{\mathbb R^N}|w_\lambda|^2=\int_{\mathbb R^N}(I_\alpha\ast |w_\lambda|^p)|w_\lambda|^p+\lambda^\sigma\int_{\mathbb R^N}|w_\lambda|^q,
$$
$$
\frac{N-2}{2}\int_{\mathbb R^N}|\nabla w_\lambda|^2+\frac{\lambda^\sigma N}{2}\int_{\mathbb R^N}|w_\lambda|^2=\frac{N+\alpha}{2p}\int_{\mathbb R^N}(I_\alpha\ast |w_\lambda|^p)|w_\lambda|^p+\frac{\lambda^\sigma N}{q}\int_{\mathbb R^N}|w_\lambda|^q.
$$
Thus, we obtain 
$$
\int_{\mathbb R^N}|w_\lambda|^2=\frac{2(2^*-q)}{q(2^*-2)}\int_{\mathbb R^N}|w_\lambda|^q.
\eqno(5.4)
$$
By the Sobolev embedding theorem  and the  interpolation inequality, we obtain 
$$
\int_{\mathbb R^N}|w_\lambda|^q\le \left(\int_{\mathbb R^N}|w_\lambda|^2\right)^{\frac{2^*-q}{2^*-2}}\left(\int_{\mathbb R^N}|w_\lambda|^{2^*}\right)^{\frac{q-2}{2^*-2}}
\le \left(\int_{\mathbb R^N}|w_\lambda|^2\right)^{\frac{2^*-q}{2^*-2}}\left(\frac{1}{S}\int_{\mathbb R^N}|\nabla w_\lambda|^2\right)^{\frac{2^*(q-2)}{2(2^*-2)}},
$$
where $S$ is the best Sobolev constant.  Therefore,  we have 
$$
\left(\int_{\mathbb R^N}|w_\lambda|^2\right)^{\frac{q-2}{2^*-2}}\le \frac{2(2^*-q)}{q(2^*-2)}\left(\frac{1}{S}\int_{\mathbb R^N}|\nabla w_\lambda|^2\right)^{\frac{2^*(q-2)}{2(2^*-2)}}.
$$
It then follows from Lemma 5.1 that
$$
\int_{\mathbb R^N}|w_\lambda|^2\le \left(\frac{2(2^*-q)}{q(2^*-2)}\right)^{\frac{2^*-2}{q-2}}\left(\frac{1}{S}\int_{\mathbb R^N}|\nabla v_\lambda|^2\right)^{2^*/2},
\eqno(5.5)
$$
which together with the boundedness of $v_\lambda$ in $H^1(\mathbb R^N)$ implies that $w_\lambda$ is bounded in $L^2(\mathbb R^N)$.

Now, we give  the following estimation on the least energy:

{\bf Lemma 5.3.} {\it  If $N\ge 5$ and $q\in (2,2^*)$, then
  $$
  m_0-m_\lambda \sim \lambda^\sigma\qquad  as \quad \lambda\to 0.
  $$
  If $N=4$ and $q\in (2,4)$, or $N=3$ and $q\in (4,6)$, then
  $$
  m_0-m_\lambda\lesssim \lambda^\sigma,  \qquad as \quad \lambda\to 0.
  $$
  }
\begin{proof}  First, we claim that there exists a comstant $C>0$ such that
$$
1<\tau_2(w_\lambda)\le 1+C\lambda^\sigma.
\eqno(5.6)
$$
In fact, since  $w_\lambda\in \mathcal{N}_\lambda$,  we see that
$$
\tau_2(w_\lambda)=\frac{\int_{\mathbb R^N}|\nabla w_\lambda|^2}{\int_{\mathbb R^N}(I_\alpha \ast |w_\lambda|^p)|w_\lambda|^p}
=1+\lambda^\sigma \frac{\int_{\mathbb R^N}|w_\lambda|^q-\int_{\mathbb R^N}|w_\lambda|^2}{\int_{\mathbb R^N}(I_\alpha\ast |w_\lambda|^p)|w_\lambda|^p}.
$$
Since 
$$
\int_{\mathbb R^N}|w_\lambda|^q\le \left(\int_{\mathbb R^N}|w_\lambda|^2\right)^{\frac{2^*-q}{2^*-2}}\left(\int_{\mathbb R^N}|w_\lambda|^{2^*}\right)^{\frac{q-2}{2^*-2}}, 
$$
we see that 
$$
\frac{\int_{\mathbb R^N}|w_\lambda|^q-\int_{\mathbb R^N}|w_\lambda|^2}{\int_{\mathbb R^N}|w_\lambda|^{2^*}}\le \zeta_\lambda^{\theta_q}(1-\zeta_\lambda^{1-\theta_q})\le \theta_q^{\frac{\theta_q}{1-\theta_q}}(1-\theta_q):=G(q), 
$$
where 
$$
\theta_q=\frac{2^*-q}{2^*-2}, \qquad \zeta_\lambda=\frac{\int_{\mathbb R^N}|v_\lambda|^2}{\int_{\mathbb R^N}|v_\lambda|^{2^*}}.
$$
Therefore, by the boundedness of $w_\lambda$ in $D^{1,2}(\R^N)$, we get
$$
\begin{array}{rcl}
\tau_2(w_\lambda)&\le& 1+\lambda^\sigma G(q)\frac{\int_{\mathbb R^N}|w_\lambda|^{2^*}}{\int_{\mathbb R^N}(I_\alpha\ast |w_\lambda|^p)|w_\lambda|^p}\\
&\le & 1+\lambda^\sigma G(q)S^{-\frac{N}{N-2}}\frac{(\int_{\mathbb R^N}|\nabla w_\lambda|^{2})^{\frac{N}{N-2}}}{\int_{\mathbb R^N}(I_\alpha\ast |w_\lambda|^p)|w_\lambda|^p}\\
&=&1+\lambda^\sigma G(q)S^{-\frac{N}{N-2}}\tau_2(w_\lambda)(\int_{\mathbb R^N}|\nabla w_\lambda|^{2})^{\frac{2}{N-2}}\\
&\le & 1+\lambda^\sigma C\tau_2(w_\lambda),
\end{array}
$$
and hence for small $\lambda>0$, there holds
$$
\tau_2(w_\lambda)\le \frac{1}{1-\lambda^\sigma C}=1+\lambda^\sigma\frac{C}{1-\lambda^\sigma C}\le 1+\frac{1}{2}C\lambda^\sigma.
$$
On the other hand,  by (5.4), we have that $\int_{\mathbb R^N}|w_\lambda|^2=\frac{2(2^*-q)}{q(2^*-2)}\int_{\mathbb R^N}|w_\lambda|^q<\int_{\mathbb R^N}|w_\lambda|^q$, therefore, we get 
$\tau_2(w_\lambda)>1$. This proved the claim.

 If $N\ge 3$, by Lemma 3.2 and the boundedness of $\{w_\lambda\}$, we find
$$
\begin{array}{rcl}
m_0&\le& \sup_{t\ge 0} J_\lambda((w_\lambda)_t)+\lambda^\sigma t^N_\lambda \left(\frac{1}{q}\int_{\mathbb R^N}|w_\lambda|^q-\frac{1}{2}\int_{\mathbb R^N}|w_\lambda|^{2}\right)\\
&\le &m_\lambda +C\lambda^\sigma,
\end{array}
\eqno(5.7)
$$
where 
$$
t_\lambda=\left(\frac{\int_{\mathbb R^N}|\nabla w_\lambda|^2}{\int_{\mathbb R^N}(I_\alpha\ast |w_\lambda|^p)|w_\lambda|^p}\right)^{\frac{1}{2+\alpha}}=(\tau_2(w_\lambda))^{\frac{1}{2+\alpha}}.
$$

For each $\rho>0$, the family  $V_\rho(x):=\rho^{-\frac{N-2}{2}}V_1(x/\rho)$ are radial ground states of  $-\Delta v=(I_\alpha\ast |v|^p)v^{p-1}$, and verify that
 $$
 \|V_\rho\|_2^2=\rho^2\|V_1\|_2^2, \qquad \int_{\mathbb R^N}|V_\rho|^q=\rho^{N-\frac{N-2}{2}q}\int_{\mathbb R^N}|V_1|^q.
 \eqno(5.8)
 $$ 
 Let $g_0(\rho)=\frac{1}{q}\int_{\mathbb R^N}|V_\rho|^q-\frac{1}{2}\int_{\mathbb R^N}|V_\rho|^2$. Then there exists $\rho_0=\rho(q)\in (0,+\infty)$ with
 $$
 \rho_0=\left(\frac{[2N-q(N-2)]\int_{\mathbb R^N}|V_1|^q}{2q\int_{\mathbb R^N}|V_1|^2}\right)^{\frac{2}{(N-2)(q-2)}}
 $$
such that 
 $$
 g_0(\rho_0)=\sup_{\rho>0}g_0(\rho)=\frac{(N-2)(q-2)}{4q}\left(\frac{[2N-q(N-2)]\int_{\mathbb R^N}|V_1|^q}{2q\int_{\mathbb R^N}|V_1|^2}\right)^{\frac{2^*-q}{q-2}}\int_{\mathbb R^N}|V_1|^q.
 $$ 
 Let $V_0=V_{\rho_0}$, then there exists $t_\lambda\in (0,+\infty)$ such that
 $$
\begin{array}{rcl}
m_\lambda&\le &\sup_{t\ge 0}J_\lambda(tV_0)=J_\lambda(t_\lambda V_0)\\
&=&\frac{t^2_\lambda}{2}\int_{\mathbb R^N}|\nabla V_0|^2-\frac{t^{2p}_\lambda}{2p}(I_\alpha\ast |V_0|^p)|V_0|^p-\lambda^\sigma\{\frac{t^{q}_\lambda}{q}\int_{\mathbb R^N}|V_0|^{q}-\frac{t^2_\lambda}{2}\int_{\mathbb R^N}|V_0|^2\}\\
&\le &\sup_{t\ge 0}\left(\frac{t^2}{2}-\frac{t^{2p}}{2p}\right)\int_{\mathbb R^N}|\nabla V_0|^2 -\lambda^\sigma\{\frac{t^{q}_\lambda}{q}\int_{\mathbb R^N}|V_0|^{q}-\frac{t^2_\lambda}{2}\int_{\mathbb R^N}|V_0|^2\}\\&=& m_0-\lambda^\sigma\{\frac{t^{q}_\lambda}{q}\int_{\mathbb R^N}|V_0|^{q}-\frac{t^2_\lambda}{2}\int_{\mathbb R^N}|V_0|^2\}.
\end{array}
\eqno(5.9)
$$
If $t_\lambda\ge 1$, then
$$
\int_{\mathbb R^N}|\nabla V_0|^2+\lambda^\sigma\int_{\mathbb R^N}|V_0|^2\ge t_\lambda^{q-2}\left\{\int_{\mathbb R^N}(I_\alpha\ast |V_0|^p)|V_0|^p+\lambda^\sigma\int_{\mathbb R^N}|V_0|^q\right\}.
$$
Hence
$$
t_\lambda\le \max\left\{1, \left(\frac{\int_{\mathbb R^N}|\nabla V_0|^2+\lambda^\sigma\int_{\mathbb R^N}|V_0|^2}{\int_{\mathbb R^N}(I_\alpha\ast |V_0|^p)|V_0|^p+\lambda^\sigma \int_{\mathbb R^N}|V_0|^q}\right)^{\frac{1}{q-2}}\right\}.
$$
If $t_\lambda\le  1$, then
$$
\int_{\mathbb R^N}|\nabla V_0|^2+\lambda^\sigma\int_{\mathbb R^N}|V_0|^2\le t_\lambda^{q-2}\left\{\int_{\mathbb R^N}(I_\alpha\ast |V_0|^p)|V_0|^p+\lambda^\sigma\int_{\mathbb R^N}|V_0|^q\right\}.
$$
Hence
$$
t_\lambda\ge \min\left\{1, \left(\frac{\int_{\mathbb R^N}|\nabla V_0|^2+\lambda^\sigma\int_{\mathbb R^N}|V_0|^2}{\int_{\mathbb R^N}(I_\alpha\ast |V_0|^p)|V_0|^p+\lambda^\sigma\int_{\mathbb R^N}|V_0|^q}\right)^{\frac{1}{q-2}}\right\}.
$$
Since
 $$
 \int_{\mathbb R^N}(I_\alpha\ast |V_0|^p)|V_0|^p=\int_{\mathbb R^N}|\nabla V_0|^2   \qquad {\rm and} \qquad \int_{\mathbb R^N}|V_0|^q>\int_{\mathbb R^N}|V_0|^2,
 $$
we conclude that 
$$
\left(\frac{\int_{\mathbb R^N}|\nabla V_0|^2+\lambda^\sigma\int_{\mathbb R^N}|V_0|^2}{\int_{\mathbb R^N}(I_\alpha\ast |V_0|^p)|V_0|^p+\lambda^\sigma\int_{\mathbb R^N}|V_0|^q}\right)^{\frac{1}{q-2}}\le t_\lambda
\le 1.
\eqno(5.10)
$$
Therefore, $\lim_{\lambda\to 0}t_\lambda=1$ and hence there exists a constant $C>0$ such that
$$
m_\lambda\le m_0-C\lambda^\sigma,
$$
for small $\lambda>0$. The proof is complete. 
\end{proof}

{\bf Lemma 5.4.}  {\it In the lower dimension cases, there exists a constant $\varpi=\varpi(q)>0$ such that  for $\lambda>0$ small, 
$$
m_\lambda\le \left\{\begin{array}{rcl}
m_0-\lambda^\sigma \left(\ln\frac{1}{\lambda}\right)^{-\frac{4-q}{q-2}}\varpi=m_0-\lambda^{\frac{2}{q-2}}\left(\ln\frac{1}{\lambda}\right)^{-\frac{4-q}{q-2}}\varpi, \    
&if&    N=4, \\
m_0-\lambda^{\sigma+\frac{2(6-q)}{(q-4)(q-2)}}\varpi=m_0-\lambda^{\frac{2}{q-4}}\varpi, \quad \quad \quad  \quad \quad  &if&   N=3 \  and \  q\in (4,6).
\end{array} \right.
$$}
\begin{proof}  
Let  $\rho>0$ and $R$ be a large parameter, and $\eta_R\in C_0^\infty(\mathbb R)$ is a cut-off function such that $\eta_R(r)=1$ for $|r|<R$, $0<\eta_R(r)<1$ for $R<|r|<2R$, $\eta_R(r)=0$ for $|r|>2R$ and $|\eta'_R(r)|\le 2/R$.  

For $\ell\gg 1$, a straightforward computation shows that 
$$
\int_{\mathbb R^N}|\nabla (\eta_\ell V_1)|^2=\left\{ \begin{array}{rcl} \frac{2(N+\alpha)}{2+\alpha}m_0+O(\ell^{-2}), \quad \qquad &{\rm if} & \  \ N=4,\\
\frac{2(N+\alpha)}{2+\alpha}m_0+O(\ell^{-1}),  \qquad \  \    \ &{\rm if} & \   \ N=3.
\end{array}\right.
$$ 
$$
\int_{\mathbb R^N} (I_\alpha\ast |\eta_\ell V_1|^p)|\eta_\ell V_1|^p=\frac{2(N+\alpha)}{2+\alpha}m_0+O(\ell^{-N}),
$$
$$
\int_{\mathbb R^N}|\eta_\ell V_1|^2=\left\{ \begin{array}{rcl} \ln \ell(1+o(1)), \quad &{\rm if} & \  \ N=4,\\
\ell(1+o(1)), \quad  \quad &{\rm if} & \   \ N=3.
\end{array}\right.
$$ 
By Lemma 3.2, we find
$$
\begin{array}{rcl}
m_\lambda&\le &\sup_{t\ge 0}J_\lambda((\eta_RV_\rho)_t)=J_\lambda((\eta_RV_\rho)_{t_\lambda})\\
&\le &\sup_{t\ge 0}\left(\frac{t^{N-2}}{2}\int_{\mathbb R^N}|\nabla(\eta_RV_\rho)|^2-\frac{t^{N+\alpha}}{2p}\int_{\mathbb R^N}(I_\alpha\ast |\eta_RV_\rho|^{p})|\eta_RV_\rho|^p\right)\\
&\quad & -\lambda^\sigma
t_\lambda^N\left[\int_{\mathbb R^N}\frac{1}{q}|\eta_RV_\rho|^q-\frac{1}{2}|\eta_RV_\rho|^2\right]\\
&= & (I)-\lambda^\sigma (II),
\end{array}
\eqno(5.11)
$$
where $t_\lambda\in (0, +\infty)$ is the unique critical point of the function $g(t)$ defined by 
$$
\begin{array}{rcl}
g(t)&=&\frac{t^{N-2}}{2}\int_{\mathbb R^N}|\nabla(\eta_RV_\rho)|^2+ \frac{t^N}{2}\lambda^\sigma\int_{\mathbb R^N}|\eta_RV_\rho|^2\\
&\quad &-\frac{t^{N+\alpha}}{2p}\int_{\mathbb R^N}(I_\alpha\ast |\eta_RV_\rho|^p)|\eta_RV_\rho|^p-\frac{t^N}{q}\lambda^\sigma\int_{\mathbb R^N}|\eta_RV_\rho|^{q}.
\end{array}
$$
That is, $t=t_\lambda$ solves the equation $\ell_1(t)=\ell_2(t)$, where
$$
\ell_1(t):=\frac{N-2}{2t^2}\int_{\mathbb R^N}|\nabla(\eta_RV_\rho)|^2
$$
and 
$$
\ell_2(t):=\frac{N+\alpha}{2p}t^\alpha\int_{\mathbb R^N}(I_\alpha\ast |\eta_RV_\rho|^p)|\eta_RV_\rho|^p
+\frac{N}{q}\lambda^\sigma\int_{\mathbb R^N}|\eta_RV_\rho|^{q}-\frac{N}{2}\lambda^\sigma\int_{\mathbb R^N}|\eta_RV_\rho|^2.
$$

If $t_\lambda\ge 1$, then
$$
\begin{array}{rcl}
\frac{N-2}{2t_\lambda^{2}}\int_{\mathbb R^N}|\nabla(\eta_RV_\rho)|^2&\ge& \frac{N+\alpha}{2p}\int_{\mathbb R^N}(I_\alpha\ast |\eta_RV_\rho|^p)|\eta_RV_\rho|^p\\
&\quad &
+\frac{N}{q}\lambda^\sigma\int_{\mathbb R^N}|\eta_RV_\rho|^{q}-\frac{N}{2}\lambda^\sigma\int_{\mathbb R^N}|\eta_RV_\rho|^2,
\end{array}
$$
and hence 
$$
\begin{array}{rcl}
t_\lambda
&\le &\left(\frac{\int_{\mathbb R^N}|\nabla(\eta_RV_\rho)|^2}{\int_{\mathbb R^N}(I_\alpha\ast |\eta_RV_\rho|^p)|\eta_RV_\rho|^p+2^*\lambda^\sigma\{\frac{1}{q}\int_{\mathbb R^N}|\eta_RV_\rho|^q-\frac{1}{2}\int_{\mathbb R^N}|\eta_RV_\rho|^2\}}\right)^{\frac{1}{2}}\\
&\le & \left(\frac{\int_{\mathbb R^N}|\nabla(\eta_RV_\rho)|^2}{\int_{\mathbb R^N}(I_\alpha\ast |\eta_RV_\rho|^p)|\eta_RV_\rho|^p}\right)^{\frac{1}{2}}.
\end{array}
\eqno(5.12)
$$
If $t_\lambda\le 1$, then
$$
\begin{array}{rcl}
\frac{N-2}{2t_\lambda^2}\int_{\mathbb R^N}|\nabla(\eta_RV_\rho)|^2&\le& \frac{N+\alpha}{2p}\int_{\mathbb R^N}(I_\alpha\ast |\eta_RV_\rho|^p)|\eta_RV_\rho|^p\\
&\quad &
+\frac{N}{q}\lambda^\sigma\int_{\mathbb R^N}|\eta_RV_\rho|^{q}-\frac{N}{2}\lambda^\sigma\int_{\mathbb R^N}|\eta_RV_\rho|^2,
\end{array}
$$
and hence 
$$
\begin{array}{rcl}
t_\lambda
&\ge&  \left(\frac{\int_{\mathbb R^N}|\nabla(\eta_RV_\rho)|^2}{\int_{\mathbb R^N}(I_\alpha\ast |\eta_RV_\rho|^p)|\eta_RV_\rho|^p+2^*\lambda^\sigma\{\frac{1}{q}\int_{\mathbb R^N}|\eta_RV_\rho|^q-\frac{1}{2}\int_{\mathbb R^N}|\eta_RV_\rho|^2\}}\right)^{\frac{1}{2}}\\
&\ge & \left(\frac{\int_{\mathbb R^N}|\nabla(\eta_RV_\rho)|^2}{\int_{\mathbb R^N}(I_\alpha\ast |\eta_RV_\rho|^p)|\eta_RV_\rho|^p}\right)^{\frac{1}{2}}
\left\{1
-2^*\lambda^\sigma\frac{\frac{1}{q}\int_{\mathbb R^N}|\eta_RV_\rho|^q-\frac{1}{2}\int_{\mathbb R^N}|\eta_RV_\rho|^2}
{\int_{\mathbb R^N}(I_\alpha\ast |\eta_RV_\rho|^p)|\eta_RV_\rho|^p}\right\}.
\end{array}
\eqno(5.13)
$$
Therefore, we obtain
$$
\begin{array}{rcl}
|t_\lambda-1|&\le &\left|\left(\frac{\int_{\mathbb R^N}|\nabla(\eta_RV_\rho)|^2}{\int_{\mathbb R^N}(I_\alpha\ast |\eta_RV_\rho|^p)|\eta_RV_\rho|^p}\right)^{\frac{1}{2}}-1\right|\\
&\quad &
+\lambda^\sigma \left(\frac{\int_{\mathbb R^N}|\nabla(\eta_RV_\rho)|^2}{\int_{\mathbb R^N}(I_\alpha\ast |\eta_RV_\rho|^p)|\eta_RV_\rho|^p}\right)^{\frac{1}{2}}
\frac{2^*\phi(\rho)}{\int_{\mathbb R^N}(I_\alpha\ast |\eta_RV_\rho|^p)|\eta_RV_\rho|^p},
\end{array}
$$
where $\phi(\rho):=\frac{1}{q}\int_{\mathbb R^N}|\eta_RV_\rho|^q-\frac{1}{2}\int_{\mathbb R^N}|\eta_RV_\rho|^2$.

Set $\ell=R/\rho$, then
$$
\begin{array}{rcl}
(I)&=&\frac{2+\alpha}{2(N+\alpha)}\frac{(\int_{\mathbb R^N}|\nabla(\eta_\ell V_1)|^2)^{\frac{N+\alpha}{2+\alpha}}}{(\int_{\mathbb R^N}(I_\alpha\ast |\eta_\ell V_1|^p)|\eta_\ell V_1|^p)^{\frac{N-2}{2+\alpha}}}\\
&=&\left\{\begin{array}{rcl} \frac{2+\alpha}{2(N+\alpha)}\frac{(\frac{2(N+\alpha)}{2+\alpha}m_0+O(\ell^{-2}))^{\frac{N+\alpha}{2+\alpha}}}{(\frac{2(N+\alpha)}{2+\alpha}m_0+O(\ell^{-4}))^{\frac{N-2}{2+\alpha}}}, \  \  {\rm if} \   \  N=4,\\
\frac{2+\alpha}{2(N+\alpha)}\frac{(\frac{2(N+\alpha)}{2+\alpha}m_0+O(\ell^{-1}))^{\frac{N+\alpha}{2+\alpha}}}{(\frac{2(N+\alpha)}{2+\alpha}m_0+O(\ell^{-3}))^{\frac{N-2}{2+\alpha}}}, \  \  {\rm if} \   \  N=3,
\end{array}\right.
\\
&=&\left\{\begin{array}{rcl} m_0+O(\ell^{-2}), \  \  {\rm if} \   \  N=4,\\
m_0+O(\ell^{-1}),   \   \   {\rm if}  \   \  N=3.\end{array}\right.
\end{array}
\eqno(5.14)
$$
Since 
$$
\begin{array}{rcl}
\phi(\rho)&=&\frac{1}{q}\int_{\mathbb R^N}|\eta_RV_\rho|^q-\frac{1}{2}\int_{\mathbb R^N}|\eta_RV_\rho|^2\\
&=&\frac{1}{q}\rho^{N-\frac{N-2}{2}q}\int_{\mathbb R^N}|\eta_\ell V_1|^q-\frac{1}{2}\rho^2\int_{\mathbb R^N}|\eta_\ell V_1|^2
\end{array}
$$
take its maximum value $\varphi(\rho_\ell)$ at the unique point 
$$
\begin{array}{rcl}
\rho_\ell:&=&\left(\frac{[2N-q(N-2)]\int_{\mathbb R^N}|\eta_\ell V_1|^q}{2q\int_{\mathbb R^N}|\eta_\ell V_1|^2}\right)^{\frac{2}{(N-2)(q-2)}}\\
&\sim &\left\{\begin{array}{rcl}
(\ln \ell)^{-\frac{2}{(N-2)(q-2)} }\quad &{\rm if}&  \   N=4,\\
\ell^{-\frac{2}{(N-2)(q-2)}} \quad &{\rm if}&  \   N=3,
\end{array}\right.
\end{array}
$$
we obtain 
$$
\begin{array}{rcl}
\phi(\rho_\ell)&=&\sup_{\rho\ge 0}\phi(\rho)\\
&=&\frac{4+q(N-2)-2N}{4q}\rho_\ell^{N-\frac{N-2}{2}q}\int_{\mathbb R^N}|\eta_\ell V_1|^q\\
&=&\frac{4+q(N-2)-2N}{4q}\left(\frac{2N-q(N-2)}{2q}\right)^{\frac{2N-q(N-2)}{(N-2)(q-2)}}\frac{(\int_{\mathbb R^N}|\eta_\ell V_1|^q)^{\frac{4}{(N-2)(q-2)}}}{(\int_{\mathbb R^N}|\eta_\ell V_1|^2)^{\frac{2N-q(N-2)}{(N-2)(q-2)}}}\\
&\le &\frac{4+q(N-2)-2N}{4q}\left(\frac{2N-q(N-2)}{2q}\right)^{\frac{2N-q(N-2)}{(N-2)(q-2)}}\int_{\mathbb R^N}|\eta_\ell V_1|^{2^*}\\
&\to &\frac{4+q(N-2)-2N}{4q}\left(\frac{2N-q(N-2)}{2q}\right)^{\frac{2N-q(N-2)}{(N-2)(q-2)}}\int_{\mathbb R^N}| V_1|^{2^*},
\end{array}
$$
as $\ell \to +\infty$, where we have used the interpolation inequality
$$
\int_{\mathbb R^N}|\eta_\ell V_1|^q\le \left(\int_{\mathbb R^N}|\eta_\ell V_1|^2\right)^{\frac{2^*-q}{2^*-2}}\left(\int_{\mathbb R^N}|\eta_\ell V_1|^{2^*}\right)^{\frac{q-2}{2^*-2}}.
$$
Since 
$$
\int_{\mathbb R^N}|\eta_\ell V_1|^q\to \int_{\mathbb R^N}|V_1|^q,
$$ 
as $\ell \to +\infty$, it follows that
$$
\begin{array}{rcl}
\phi(\rho_\ell)&=&\frac{4+q(N-2)-2N}{4q}\left(\frac{[2N-q(N-2)]\int_{\mathbb R^N}|\eta_\ell V_1|^q}{2q\int_{\mathbb R^N}|\eta_\ell V_1|^2}\right)^{\frac{2N-q(N-2)}{(N-2)(q-2)}}\int_{\mathbb R^N}|\eta_\ell V_1|^q\\
&=& \left\{\begin{array}{rcl} C(\ln \ell(1+o(1))^{-\frac{2N-q(N-2)}{(N-2)(q-2)}}  \qquad &{\rm if}&  \   \  N=4,\\
C(\ell(1+o(1))^{-\frac{2N-q(N-2)}{(N-2)(q-2)}} \qquad &{\rm if }& \   \ N=3.
\end{array}\right.
\end{array}
$$
Since  $\phi(\rho)$ is bounded, we find
$$
\begin{array}{rcl}
|t_\lambda-1|&\le& \left|\left(\frac{\int_{\mathbb R^N}|\nabla(\eta_\ell V_1)|^2}{\int_{\mathbb R^N}(I_\alpha\ast |\eta_\ell V_1|^p)|\eta_\ell V_1|^p}\right)^{\frac{1}{2}}-1\right|\\
&\quad &
+\lambda^\sigma \left(\frac{\int_{\mathbb R^N}|\nabla(\eta_\ell V_1)|^2}{\int_{\mathbb R^N}(I_\alpha\ast |\eta_\ell V_1|^p)|\eta_\ell V_1|^p}\right)^{\frac{1}{2}}
\frac{2^*C}{\int_{\mathbb R^N}(I_\alpha\ast |\eta_\ell V_1|^p)|\eta_\ell V_1|^p}\\
&\to &\frac{2^*C\lambda^\sigma}{\int_{\mathbb R^N}(I_\alpha\ast | V_1|^p)|V_1|^p},
\end{array}
$$
as $\ell \to +\infty$. Thus, for small $\lambda>0$, we have 
$$
\begin{array}{rcl}
(II)&=&\phi(\rho_\ell)+(t_\lambda^N-1)\phi(\rho_\ell )\\
&\sim &\left\{\begin{array}{rcl} (\ln \ell)^{-\frac{2N-q(N-2)}{(N-2)(q-2)}},  \qquad &{\rm if}&  \   \  N=4,\\
\ell^{-\frac{2N-q(N-2)}{(N-2)(q-2)}}, \qquad &{\rm if }& \   \ N=3.
\end{array}\right.\end{array}
$$
It follows that if $N=4$, then 
$$
\begin{array}{rcl}
m_\lambda &\le &(I)-\lambda^\sigma (II)\\
&\le &m_0+O(\ell^{-2})-C\lambda^\sigma (\ln \ell)^{-\frac{2N-q(N-2)}{(N-2)(q-2)} }.
\end{array}
\eqno(5.15)
$$
Take $\ell=(1/\lambda)^M$. Then 
$$
m_\lambda\le m_0+O(\lambda^{2M})-C\lambda^\sigma M^{-\frac{2N-q(N-2)}{(N-2)(q-2)}}(\ln\frac{1}{\lambda})^{-\frac{2N-q(N-2)}{(N-2)(q-2) }}.
$$
If $M>\frac{1}{q-2}$, then $2M>\sigma$, and hence
$$
m_\lambda\le m_0-\lambda^\sigma(\ln\frac{1}{\lambda})^{-\frac{2N-q(N-2)}{(N-2)(q-2)} }\varpi=m_0-\lambda^{\frac{2}{q-2}}
(\ln\frac{1}{\lambda})^{-\frac{4-q}{q-2}}\varpi,
\eqno(5.16)
$$
for small $\mu>0$, where 
$$
\varpi=\frac{1}{2}CM^{-\frac{2N-q(N-2)}{(N-2)(q-2)}}.
$$

If  $N=3$, then 
$$
\begin{array}{rcl}
m_\lambda&\le &(I)-\lambda^\sigma (II)\\
&\le &m_0+O(\ell^{-1})-C\lambda^\sigma \ell^{-\frac{2N-q(N-2)}{(N-2)(q-2)}}.
\end{array}
\eqno(5.17)
$$
Take $\ell=\delta^{-1}\lambda^{-\tau}$. Then 
$$
m_\lambda\le m_0+\delta O(\lambda^\tau)-C\lambda^\sigma \delta^{\frac{2N-q(N-2)}{(N-2)(q-2)}}\lambda^{\tau\frac{2N-q(N-2)}{(N-2)(q-2)} }
$$
If $q\in (4,6)$ and 
$$
\tau=\frac{2(N-2)}{2+q(N-2)-2N}=\frac{2}{q-4},
$$
then
$$
m_\lambda\le m_0+(\delta O(1)-C\delta^{\frac{2N-q(N-2)}{(N-2)(q-2)}})\lambda^{\frac{2}{q-4}}.
$$
Since 
$$
1>\frac{2N-q(N-2)}{(N-2)(q-2)},
$$
it follows that for some small $\delta>0$, there exists $\varpi>0$ such that 
$$
m_\lambda\le m_0-\lambda^{\frac{2}{q-4}}\varpi.
$$
This completes the proof. \end{proof}

Combining Lemma 5.3 and Lemma 5.4, we get the following

{\bf Lemma 5.5.} {\it Let $\delta_\lambda:=m_0-m_\lambda$, then 
$$
\lambda^\sigma\gtrsim \delta_\lambda\gtrsim \left\{\begin{array}{rcl} 
 \lambda^\sigma, \  \  \   \quad  \quad  \qquad &{ if}&  \    N\ge 5,\\
 \lambda^{\frac{2}{q-2}}(\ln\frac{1}{\lambda})^{-\frac{4-q}{q-2}},   \ &{ if}& \   N=4,\\
 \lambda^{\frac{2}{q-4}},  \qquad \qquad  &{ if}& \     N=3 \ and \ q\in (4, 6).
 \end{array}\right.
 $$   }

{\bf Lemma 5.6.} {\it 
Assume $N\ge 5$. Then $\|w_\lambda\|_q\sim 1$ as $\lambda\to 0$.}

\begin{proof} 
By (5.7),  we have
$$
m_0\le m_\lambda+\lambda^\sigma (\tau_2(w_\lambda))^{\frac{N}{N+\alpha}}\frac{q-2}{q(2^*-2)}\int_{\mathbb R^N}|w_\lambda|^q.
$$
Therefore, it follows from (5.6) and Corollary 5.5 that 
$$
\|w_\lambda\|_q^q\ge \frac{m_0-m_\lambda}{(\tau_2(w_\lambda))^{\frac{N}{2+\alpha}}}\cdot \frac{q(2^*-2)}{q-2}\lambda^{-\sigma}\ge \frac{Cq(2^*-2)}{(q-2)(\tau_2(w_\lambda))^{\frac{N}{2+\alpha}}}\ge C>0,
$$
which together with the boundedness of $\{w_\lambda\}$ implies the desired  conclusion. 
\end{proof}

{\bf Lemma 5.7.} {\it Let $N\ge 5$, $\alpha>N-4$ and $q\in (2,2^*)$,  then $w_\lambda\to V_{\rho_0}$ in $H^1(\mathbb R^N)$ as $\lambda\to 0$, where $V_{\rho_0}$ is  a positive ground sate of the equation $-\Delta V=(I_\alpha\ast |V|^p)V^{p-1}$ with
$$
 \rho_0=\left(\frac{2(2^*-q)\int_{\mathbb R^N}|V_1|^q}{q(2^*-2)\int_{\mathbb R^N}|V_1|^2}\right)^ {\frac{2}{(N-2)(q-2)}}.
 \eqno(5.18)
 $$
In the lower dimension cases $N=4$ and $N=3$, there exists $\xi_\lambda\in (0,+\infty)$ with $\xi_\lambda\to 0$ such that 
$$
w_\lambda-\xi_\lambda^{-\frac{N-2}{2}}V_1(\xi^{-1}_\lambda\cdot)\to 0
$$
as $\lambda\to 0$ in $D^{1,2}(\mathbb R^N)$ and $L^{2^*}(\mathbb R^N)$. }

\begin{proof} Note that $w_\lambda$ is a positive radially symmetric function, and by Lemma 5.2, $\{w_\lambda\}$ is bounded in $H^1(\mathbb R^N)$. Then there exists $w_0\in H^1(\mathbb R^N)$ verifying $-\Delta w=(I_\alpha \ast w^p)w^{p-1}$ such that 
$$
w_\lambda \rightharpoonup w_0   \quad {\rm weakly \ in} \  H^1(\mathbb R^N), \quad w_\lambda\to w_0 \quad {\rm in} \ L^p(\mathbb R^N) \quad {\rm for \ any} \ p\in (2,2^*),
\eqno(5.19)
$$
and 
$$
w_\lambda(x)\to w_0(x) \quad a. \ e. \  {\rm on} \ R^N,  \qquad w_\lambda\to w_0 \quad {\rm in} \   L^2_{loc}(\mathbb R^N).
\eqno(5.20)
$$
Observe that
$$
J_0(w_\lambda)=J_\lambda(w_\lambda)+\frac{\lambda^\sigma}{q}\int_{\mathbb R^N}|w_\lambda|^q-\frac{\lambda^\sigma}{2}\int_{\mathbb R^N}|w_\lambda|^2=m_\lambda+o(1)=m_0+o(1),
$$
and 
$$
J'_0(w_\lambda)w=J'_\lambda(w_\lambda)w+\lambda^\sigma\int_{\mathbb R^N}|w_\lambda|^{q-2}w_\lambda w-\lambda^\sigma\int_{\mathbb R^N}w_\lambda w=o(1).
$$
Therefore, $\{w_\lambda\}$ is a PS sequence of $J_0$ at level  $m_0=\frac{2+\alpha}{2(N+\alpha)}S_\alpha^{\frac{N+\alpha}{2+\alpha}}$.

By Lemma 3.6,  it is standard to show that there exists $\zeta^{(j)}_\lambda\in (0,+\infty)$, $w^{(j)}\in D^{1,2}(\mathbb R^N)$ with $j=1,2,\cdots, k$, $k$  a non-negative integer, such that
$$
w_\lambda=w_0+\sum_{j=1}^k(\zeta^{(j)}_\lambda)^{-\frac{N-2}{2}}w^{(j)}((\zeta^{(j)}_\lambda)^{-1} x)+\tilde w_\lambda,
\eqno(5.21)
$$
where $\tilde w_\lambda \to 0$ in $L^{2^*}(\mathbb R^N)$ and $w^{(j)}$ are nontrivial solutions of the limit equation $-\Delta v=(I_\alpha\ast v^p)v^{p-1}$. Moreover, we have
$$
\limsup_{\lambda\to 0}\|w_\lambda\|^2_{D^1(\mathbb R^N)}\ge \|w_0\|^2_{D^1(\mathbb R^N)}+\sum_{j=1}^k\|w^{(j)}\|^2_{D^1(\mathbb R^N)}
\eqno(5.22)
$$
and 
$$
m_0=J_0(w_0)+\sum_{j=1}^kJ_0(w^{(j)}).
\eqno(5.23)
$$
Moreover, $J_0(w_0)\ge 0$ and $J_0(w^{(j)})\ge m_0$ for all $j=1,2,\cdots, k.$

If $N\ge 5$, then by Lemma 5.6, we have $w_0\not=0$ and hence $J_0(w_0)=m_0$ and  $k=0$.  Thus $w_\lambda\to w_0$ in $L^{2^*}(\mathbb R^N)$. Since $J_0'(w_\lambda)\to 0$, it follows that
$w_\lambda\to w_0$ in $D^{1,2}(\mathbb R^N)$.  

Since $w_\lambda(x)$ is radial and radially decreasing, for every $x\in\mathbb R^N\setminus \{0\}$, we have 
$$
w^2_\lambda(x)\le \frac{1}{|B_{|x|}|}\int_{B_{|x|}}|w_\lambda|^2\le \frac{1}{|x|^N}\int_{\mathbb R^N}|w_\lambda|^2\le \frac{C}{|x|^N},
$$
then 
$$
w_\lambda(x) \le C|x|^{-\frac{N}{2}}, \qquad |x|\ge 1.
\eqno(5.24)
$$

If  $\alpha>N-4$, then we have $p=\frac{N+\alpha}{N-2}>2$ and hence
$$
|w_\lambda|^p|x|^N\le C|x|^{-\frac{N}{2}p+N}=C|x|^{-\frac{N}{2}(p-2)}\to 0, \qquad {\rm as} \  |x|\to +\infty.
$$
By virtue of  Lemma 3.10, we obtain 
$$
(I_\alpha\ast |w_\lambda|^p)(x)\le C|x|^{-N+\alpha}, \qquad |x|\ge 1.
$$
and then
$$
(I_\alpha\ast |w_\lambda|^p)(x)|w_\lambda|^{p-2}(x)\le C|x|^{-\frac{N^2-N\alpha+4\alpha}{2(N-2)}}, \qquad |x|\ge \tilde R.
\eqno(5.25)
$$
Since
$$
\left(-\Delta-C|x|^{-\frac{N^2-N\alpha+4\alpha}{2(N-2)}}\right)w_\lambda\le \left(-\Delta +\lambda^\sigma-(I_\alpha\ast |w_\lambda|^p)w_\lambda^{p-2}-\lambda^\sigma w_\lambda^{q-2}\right)w_\lambda=0,
$$
for large $|x|$.  We also have
$$
\left(-\Delta-C|x|^{-\frac{N^2-N\alpha+4\alpha}{2(N-2)}}\right)\frac{1}{|x|^{N-2-\varepsilon_0}}=\left(\varepsilon_0(N-2-\varepsilon_0)-C|x|^{-\frac{(N-4)(N-\alpha)+8}{2(N-2)}}\right)\frac{1}{|x|^{N-\varepsilon_0}},
$$
which is positive for $|x|$ large enough. By (5.24) and  the maximum principle on $\mathbb R^N\setminus B_R$, we deduce that 
$$
w_\lambda(x)\le \frac{w_\lambda(R)R^{N-2-\varepsilon_0}}{|x|^{N-2-\varepsilon_0}}\le \frac{CR^{N/2-2-\varepsilon_0}}{|x|^{N-2-\varepsilon_0}},  \qquad {\rm for} \  \ |x|\ge R.
\eqno(5.26)
$$
When $\varepsilon_0>0$ is small enough, the domination is in $L^2(\mathbb R^N)$ for $N\ge 5$, and this shows, by the dominated convergence theorem, that 
$w_\lambda\to w_0$ in $L^2(\mathbb R^N)$.  Thus, we conclude that $w_\lambda\to w_0$ in  $H^1(\mathbb R^N)$. Moreover, by (5.4), we obtain
$$
\|w_0\|_2^2=\frac{2(2^*-q)}{q(2^*-2)}\int_{\mathbb R^N}|w_0|^q,
$$
from which it follows that $w_0=V_{\rho_0}$ with  
$$
 \rho_0=\left(\frac{2(2^*-q)\int_{\mathbb R^N}|V_1|^q}{q(2^*-2)\int_{\mathbb R^N}|V_1|^2}\right)^ {\frac{2}{(N-2)(q-2)}}.
 $$

If $N=4$ or $3$. By Fatou's lemma, we have $\|w_0\|^2_2\le \liminf_{\lambda\to 0}\|w_\lambda\|_2^2<\infty$, therefore,  $w_0=0$ and hence $k=1$. Thus, we obtain  
$J_0(w^{(1)})=m_0$ and hence $w^{(1)}=V_\rho$ for some $\rho\in (0,+\infty)$. Therefore, we conclude that 
$$
w_\lambda-\xi_\lambda^{-\frac{N-2}{2}}V_1(\xi_\lambda^{-1}\cdot )\to 0
$$ 
in $L^{2^*}(\mathbb R^N)$ as $\lambda\to 0$, where 
$\xi_\lambda:=\rho\zeta_\lambda^{(1)}\in (0,+\infty)$ satisfying $\xi_\lambda\to 0$ as $\lambda\to 0$.
Since 
$$
J_0'(w_\lambda-\xi_\lambda^{-\frac{N-2}{2}}V_1(\xi_\lambda^{-1}\cdot ))=J'_0(w_\lambda)+J'_0(V_1)+o(1)=o(1)
$$ 
as $\lambda\to 0$, it follows that $w_\lambda-\xi_\lambda^{-\frac{N-2}{2}}V_1(\xi_\lambda^{-1}\cdot )\to 0$ in $D^{1,2}(\mathbb R^N)$. 
\end{proof}

{\bf Lemma 5.8.} {\it Let $N\ge 5$ and $q\in (2,2^*)$,  then there exists a $\zeta_\lambda\in (0,\infty)$ verifying 
$$
\zeta_\lambda\sim  \lambda^{\frac{2^*-2}{2(q-2)}}
$$
such that the rescaled ground states 
$$
w_\lambda(x)=\zeta_\lambda^{\frac{N-2}{2}}v_\lambda(\zeta_\lambda x)
$$
converges to  $V_{\rho_0}$ in $H^1(\mathbb R^N)$ as $\lambda\to 0$, where $V_{\rho_0}$ is  given in Lemma 5.7.}

\begin{proof}
If $\alpha>N-4$, then the statement is valid with $\zeta_\lambda=\lambda^{\frac{2^*-2}{2(q-2)}}$. 
If $\alpha\le N-4$, then for any $\lambda_n\to 0$, up to a subsequence, we can assume that
$w_{\lambda_n}\to V_\rho$ in $D^{1,2}(\mathbb R^N)$ with $\rho\in (0,\rho_0]$. Moreover, $w_{\lambda_n}\to V_\rho$ in $L^2(\mathbb R^N)$ if and only if
$\rho=\rho_0$. Arguing as in the proof of Theorem 2.1, we can show that there exists a $\zeta_\lambda\sim \lambda^{\frac{2^*-2}{2(q-2)}}$ such that
$w_\lambda(x)=\zeta_\lambda^{\frac{N-2}{2}}v_\lambda(\zeta_\lambda x )$ converges to $V_{\rho_0}$ in $L^2(\mathbb R^N)$, and hence in $H^1(\mathbb R^N)$.
This completes the proof. \end{proof}

\vskip 5mm

In the lower dimension cases $N=4$ and $N=3$, we further perform a  scaling
$$
\tilde w(x)=\xi_\lambda^{\frac{N-2}{2}} w(\xi_\lambda x), 
\eqno(5.27)
$$
where $\xi_\lambda\in (0,+\infty)$ is given in Lemma 5.7.  Then the rescaled equation is as follows
$$
-\Delta \tilde w+\lambda^\sigma\xi_\lambda^{2}\tilde w=(I_\alpha\ast |\tilde w|^{\frac{N+\alpha}{N-2}})\tilde w^{\frac{2+\alpha}{N-2}}+\lambda^\sigma\xi_\lambda^{N-\frac{N-2}{2}q}\tilde w^{q-1}.
\eqno(R_\lambda)
$$
The  corresponding energy functional is given by
$$
\tilde J_\lambda(\tilde w):=\frac{1}{2}\int_{\mathbb R^N}|\nabla\tilde  w|^2+\lambda^\sigma\xi_\lambda^2|\tilde w|^2-\frac{1}{2p}\int_{\mathbb R^N}(I_\alpha\ast |\tilde w|^{p})|\tilde w|^p-\frac{1}{q}\lambda^\sigma\xi_\lambda^{N-\frac{N-2}{2}q}\int_{\mathbb R^N}|\tilde w|^q.
\eqno(5.28)
$$
Clearly, we have $\tilde J_\lambda(\tilde w)=J_\lambda(w)=I_\lambda(v)$. 

Furthermore, we  have the following lemma.

{\bf Lemma 5.9.}  {\it  Let $v,w,\tilde w\in H^1(\mathbb R^N)$ satisfy  (5.1) and (5.27), then the following statements hold true

(1) $ \ \|\nabla \tilde w\|_2^2= \|\nabla w\|_{2}^{2}=\|\nabla v\|_{2}^{2}, \  \int_{\mathbb R^N}(I_\alpha\ast |\tilde w|^p)|\tilde w|^p=\int_{\mathbb R^N}(I_\alpha\ast |w|^p)|w|^p=\int_{\mathbb R^N}(I_\alpha\ast |v|^p)|v|^p,$

(2)  $\xi_\lambda^{2}\|\tilde w\|^2_2=\|w\|_2^2=\lambda^{-\sigma}\| v\|_2^2, \   \   \xi_\lambda^{N-\frac{N-2}{2}q}\|\tilde w\|^q_q=\|w\|_q^q=\lambda^{1-\sigma} \|v\|_q^q$.
}

Set  $\tilde w_\lambda(x)=\xi_\lambda^{\frac{N-2}{2}} w_\lambda(\xi_\lambda x)$, then by Lemma 5.7, we have 
$$
\|\nabla(\tilde w_\lambda-V_1)\|_2\to 0, \qquad \|\tilde w_\lambda-V_1\|_{2^*}\to 0,   \qquad {\rm as} \  \ \lambda\to 0.
\eqno(5.29)
$$

Note that the corresponding Nehari and Pohozaev's identities are as follows
$$
\int_{\mathbb R^N}|\nabla \tilde w_\lambda|^2+\lambda^\sigma\xi_\lambda^{2}\int_{\mathbb R^N}|
\tilde w_\lambda|^2=\int_{\mathbb R^N}(I_\alpha\ast |\tilde w_\lambda|^{p})|\tilde w_\lambda|^p+\lambda^\sigma\xi_\lambda^{N-\frac{N-2}{2}q}\int_{\mathbb R^N}|\tilde w_\lambda|^q
\eqno(5.30)
$$
and 
$$
\frac{1}{2^*}\int_{\mathbb R^N}|\nabla \tilde w_\lambda|^2+\frac{1}{2}
\lambda^\sigma\xi_\lambda^{2}\int_{\mathbb R^N}|\tilde w_\lambda|^2=\frac{1}{2^*}\int_{\mathbb R^N}(I_\alpha\ast |\tilde w_\lambda|^{p})|\tilde w_\lambda|^p+\frac{1}{q}\lambda^\sigma\xi_\lambda^{N-\frac{N-2}{2}q}\int_{\mathbb R^N}|\tilde w_\lambda|^q,
\eqno(5.31)
$$
it follows that 
$$
\left(\frac{1}{2}-\frac{1}{2^*}\right)\lambda^\sigma\xi_\lambda^{2}\int_{\mathbb R^N}|\tilde w_\lambda|^2=\left(\frac{1}{q}-\frac{1}{2^*}\right)\lambda^\sigma\xi_\lambda^{N-\frac{N-2}{2}q}\int_{\mathbb R^N}|\tilde w_\lambda|^q.
$$
Thus, we obtain
$$
\xi_\lambda^{\frac{(N-2)(q-2)}{2}}\int_{\mathbb R^N}|\tilde w_\lambda|^2=\frac{2(2^*-q)}{q(2^*-2)}\int_{\mathbb R^N}|\tilde w_\lambda|^q.
\eqno(5.32)
$$

To control the norm $\|\tilde w_\lambda\|_2$, we note that  for any $\lambda>0$, $\tilde w_\lambda>0$ satisfies the linear inequality
$$
-\Delta \tilde w_\lambda+\lambda^\sigma\xi_\lambda^{2}\tilde w_\lambda=(I_\alpha\ast |\tilde w_\lambda|^{p})\tilde w_\lambda^{p-1}+\lambda^\sigma\xi_\lambda^{N-\frac{N-2}{2}q}\tilde w_\lambda^{q-1}>0,  \qquad x\in \mathbb R^N.
\eqno(5.33)
$$

{\bf Lemma 5.10.}  {\it There exists a constant $c>0$ such that 
$$
\tilde w_\lambda(x)\ge c|x|^{-(N-2)}\exp({-\lambda^{\frac{\sigma}{2}}\xi_\lambda}|x|),  \quad |x|\ge 1.
\eqno(5.34)
$$}
The proof of the above lemma  is similar to that of   \cite[Lemma 4.8]{Moroz-1}.  As consequences, we have the following lemma.

{\bf Lemma 5.11.}  {\it If $N=3$, then $\|\tilde w_\lambda\|_2^2\gtrsim \lambda^{-\frac{\sigma}{2}}\xi_\lambda^{-1}$.}

{\bf Lemma 5.12.}  {\it If $N=4$, then $\|\tilde w_\lambda\|_2^2\gtrsim  - \ln(\lambda^{\sigma}\xi_\lambda^2).$}

We remark that $\tilde w_\lambda$ is only defined for $N=4$ and $N=3$. But the following discussion also apply to the case $N\ge 5$.
To prove our main result, the key point is to show the boundedness of $\|\tilde w_\lambda\|_q$. 

{\bf Lemma 5.13.}  {\it  Assume $N\ge 3$, $\alpha>N-4$ and $ 2<q<2^*$. Then there exist constants $L_0>0$ and $C_0>0$ such that for any small $\lambda>0$ and $|x|\ge L_0\lambda^{-\sigma/2}\xi_\lambda^{-1}$, 
$$
\tilde w_\lambda(x)\le C_0\lambda^{\sigma(N-2)/4}\xi_\lambda^{(N-2)/2}\exp(-\frac{1}{2}\lambda^{\sigma/2}\xi_\lambda |x|).
$$}
\begin{proof} By (5.25) and (5.26),  if $|x|\ge L_0\lambda^{-\sigma/2}\xi_\lambda^{-1}$ with  $L_0>0$ being large enough, we have
$$
\begin{array}{rcl}
(I_\alpha\ast |\tilde w_\lambda|^p)(x)|\tilde w_\lambda|^{p-2}(x)&=&\xi_\lambda^{(N-2)(p-1)-\alpha}(I_\alpha\ast |w_\lambda|^p)(\xi_\lambda x)|w_\lambda|^{p-2}(\xi_\lambda x)\\
&\le &C\xi_\lambda^2L_0^{-\frac{N^2-N\alpha+4\alpha}{2(N-2)}}\lambda^{\sigma\cdot \frac{N^2-N\alpha+4\alpha}{4(N-2)}}\\
&\le &\frac{1}{4}\lambda^\sigma\xi_\lambda^2,
\end{array}
$$
here we have used the fact that 
$$
\frac{N^2-N\alpha+4\alpha}{4(N-2)}>1,
$$
which follows from the inequality $N<N+2<\frac{4(\alpha+2)}{\alpha-N+4}, \ \forall \alpha\in (N-4,N)$. 

By (5.24) and (5.26), for $|x|\ge L_0\lambda^{-\sigma/2}\xi_\lambda^{-1}$, we get
$$
\begin{array}{rcl}
\lambda^\sigma\xi_\lambda^{N-\frac{N-2}{2}q}\tilde w_\lambda^{q-2}(x)&=&\lambda^\sigma\xi_\lambda^{N-\frac{N-2}{2}q}\xi_\lambda^{\frac{N-2}{2}(q-2)}w_\lambda^{q-2}(\xi_\lambda x)
\\
&\le& \lambda^{\sigma}\xi_\lambda^{2}\cdot C|\xi_\lambda x|^{-\frac{N}{2}(q-2)}\\
&\le &CL_0^{-N(q-2)/2}\lambda^{\sigma+\frac{\sigma N}{4}(q-2)}\xi_\lambda^2\\
&\le & \frac{1}{4}\lambda^\sigma\xi_\lambda^2.
\end{array}
$$
Therefore, we obtain 
$$
-\Delta\tilde w_\lambda(x)+\frac{1}{2}\lambda^{\sigma}\xi_\lambda^2\tilde w_\lambda(x)\le 0, \quad 
{\rm for \ all} \  |x|\ge L_0\lambda^{-\sigma/2}\xi_\lambda^{-1}.
$$ 

We adopt an argument as used in  \cite[Lemma 3.2]{Akahori-2}.
 Let $R>L_0\lambda^{-\sigma/2}\xi_\lambda^{-1}$, and introduce a positive function 
$$
\psi_R(r):=\exp(-\frac{1}{2}\lambda^{\sigma/2}\xi_\lambda(r-L_0\lambda^{-\sigma/2}\xi_\lambda^{-1}))+\exp(\frac{1}{2}\lambda^{\sigma/2}\xi_\lambda(r-R)).
$$
It is easy to see that 
$$
|\psi'_R(r)|\le \frac{1}{2}\lambda^{\sigma/2}\xi_\lambda\psi_R(r),  \quad 
\psi''_R(r)=\frac{1}{4}\lambda^{\sigma}\xi_\lambda^2\psi_R(r).
$$
We use the same symbol $\psi_R$ to denote the radial function $\psi_R(|x|)$ on $\mathbb R^N$. Then for $L_0\lambda^{-\sigma/2}\xi_\lambda^{-1}<r<R$,  if $L_0\ge 2(N-1)$, then we have 
$$
\begin{array}{rcl}
-\Delta \psi_R+\frac{1}{2}\lambda^{\sigma}\xi_\lambda^2\psi_R&=&-\psi''_R-\frac{N-1}{r}\psi_R'+\frac{1}{2}\lambda^\sigma\xi_\lambda^2\psi_R\\
&\ge& -\frac{1}{4}\lambda^{\sigma}\xi_\lambda^2\psi_R-\frac{N-1}{L_0}\lambda^{\sigma/2}\xi_\lambda\cdot \frac{1}{2}\lambda^{\sigma/2}\xi_\lambda\psi_R+\frac{1}{2}\lambda^\sigma\xi_\lambda^2\psi_R\\
&\ge &0.
\end{array}
$$
Furthemore, $\psi_R(L_0\lambda^{-\sigma/2}\xi_\lambda^{-1})\ge 1$ and $\psi_R(R)\ge 1$, thus we have 
$$
\tilde w_\lambda(R)\le \tilde w_\lambda(L_0\lambda^{-\sigma/2}\xi_\lambda^{-1})\le CL_0^{-\frac{N-2}{2}}\lambda^{\frac{\sigma(N-2)}{4}}\xi_\lambda^{\frac{N-2}{2}}\min \{\psi_R(L_0\lambda^{-\sigma/2}\xi_\lambda^{-1}), \psi_R(R)\}.
$$
Hence, the comparison principle implies that if $L_0\lambda^{-\sigma/2}\xi_\lambda^{-1}\le |x|\le R$, then
$$
\tilde w_\lambda(x)\le CL_0^{-\frac{N-2}{2}}\lambda^{\frac{\sigma(N-2)}{4}}\xi_\lambda^{\frac{N-2}{2}}\psi_R(|x|).
$$
Since $R>L_0\lambda^{-\sigma/2}\xi_\lambda^{-1}$ is arbitrary,  taking $R\to\infty$, we find that 
$$
\tilde w_\lambda(x)\le CL_0^{-\frac{N-2}{2}}e^{L_0/2}\lambda^{\frac{\sigma(N-2)}{4}}\xi_\lambda^{\frac{N-2}{2}}e^{-\frac{1}{2}\lambda^{\sigma/2}\xi_\lambda |x|},
$$
for all $|x|\ge L_0\lambda^{-\sigma/2}\xi_\lambda^{-1}$. The proof is complete.
\end{proof}

In the following proposition, we establish the optimal decay estimate of $\tilde w_\lambda$ at infinity.

{\bf Proposition 5.14.}  {\it Assume $N\ge 3$, $\alpha>N-4$ and $ 2<q<2^*$. Then there exists a constant $C>0$ such that for small $\lambda>0$, there holds 
$$
\tilde w_\lambda(x)\le C(1+|x|)^{-(N-2)}, \qquad x\in \mathbb R^N.
$$}

We consider the Kelvin transform of $\tilde w_\lambda$. For any $w\in H^1(\mathbb R^N)$, we denote by $K[w]$ the  Kelvin transform of $w$, that is,
$$
K[w](x):=|x|^{-(N-2)}w\left(\frac{x}{|x|^2}\right). 
$$
It is easy to see that $\|K[\tilde w_\lambda]\|_{L^\infty(B_1)}\lesssim 1$ implies that 
$$
\tilde w_\lambda(x)\lesssim |x|^{-(N-2)},  \qquad |x|\ge 1,
$$
uniformly for small $\lambda>0$.

Thus, to prove Proposition  5.14, it suffices to show that there exists $\lambda_0>$ such that 
$$
\sup_{\lambda\in (0,\lambda_0)}\|K[\tilde w_\lambda]\|_{L^{\infty}(B_1)}<\infty.
\eqno(5.35)
$$

It is easy to verify that $K[\tilde w_\lambda]$ satisfies 
$$
-\Delta K[\tilde w_\lambda]+\frac{\lambda^\sigma\xi_\lambda^{2}}{|x|^4}K[\tilde w_\lambda]=\frac{1}{|x|^{4}}(I_\alpha\ast |\tilde w_\lambda|^{p})(\frac{x}{|x|^2})\tilde w_\lambda^{p-2}(\frac{x}{|x|^2})K[\tilde w_\lambda]+\frac{\lambda^\sigma\xi_\lambda^{\gamma/2}}{|x|^{\gamma}}K[\tilde w_\lambda]^{q-1},
\eqno(5.36)
$$
here and in what follows, we set
$$
\gamma:=2N-(N-2)q>0.
$$
 
We also see from Lemma 5.13 that if $|x|\le \lambda^{\sigma/2}\xi_\lambda/L_0$, then 
$$
K[\tilde w_\lambda](x)\lesssim \frac{1}{|x|^{N-2}}\lambda^{\frac{\sigma(N-2)}{4}}\xi_\lambda^{\frac{N-2}{2}}e^{-\frac{1}{2}\lambda^{\sigma/2}\xi_\lambda |x|^{-1}}.
\eqno(5.37)
$$

   Let 
$$
a(x)=\frac{\lambda^\sigma\xi_\lambda^2}{|x|^4},  \quad b(x)=\frac{1}{|x|^{4}}(I_\alpha\ast |\tilde w_\lambda|^p)(\frac{x}{|x|^2})\tilde w_\lambda^{p-2}(\frac{x}{|x|^2})
+\frac{\lambda^\sigma\xi^{\gamma/2}}{|x|^\gamma}K[\tilde w_\lambda]^{q-2}(x).
$$
Then (5.36) reads as 
$$
-\Delta K[\tilde w_\lambda]+a(x)K[\tilde w_\lambda]=b(x)K[\tilde w_\lambda].
$$
We shall apply the Moser iteration to prove  (5.35).

We note that it follows from (5.37) that for any $v\in H^1_0(B_4)$, 
$$
\int_{B_4}\frac{\lambda^\sigma\xi_\lambda^2}{|x|^4}K[\tilde w_\lambda](x)|v(x)|dx<\infty.
\eqno(5.38)
$$
Since $\tilde w_\lambda\to V_1$  in $L^{2^*}(\mathbb R^N)$   as $\lambda\to 0$,  and  for any $s>1$, the Lebesgue space $L^{s}(\mathbb R^N)$ has the Kadets-Klee property,   it is easy to see that
$$
\lim_{\lambda\to 0}\int_{\mathbb R^N}|\tilde w_\lambda^p-V_1^p|^{\frac{2N}{N+\alpha}}dx=0
$$
and 
$$
\lim_{\lambda\to 0}\int_{\mathbb R^N}|\tilde w_\lambda^{p-2}-V_1^{p-2}|^{\frac{2N}{\alpha-N+4}}dx=0.
$$
Therefore, by the H\"older inequality and the Hardy-Littlewood-Sobolev inequality, we find 
$$
\begin{array}{rl}
&\left\|\frac{1}{|x|^4}(I_\alpha\ast (|\tilde w_\lambda|^p-|V_1|^p))(\frac{x}{|x|^2})\tilde w_\lambda^{p-2}(\frac{x}{|x|^2})\right\|_{L^{\frac{N}{2}}(\mathbb R^N)}^{N/2}\\
&\quad =\int_{\mathbb R^N}\frac{1}{|x|^{2N}}|(I_\alpha\ast (|\tilde w_\lambda|^p-|V_1|^p))(\frac{x}{|x|^2})|^{N/2}|\tilde w_\lambda^{p-2}(\frac{x}{|x|^2})|^{N/2}dx\\
&\quad= \int_{\mathbb R^N}|(I_\alpha\ast (|\tilde w_\lambda|^p-|V_1|^p))(z)|^{N/2}|\tilde w_\lambda^{p-2}(z)|^{N/2}dz\\
&\quad \lesssim \left(\int_{\mathbb R^N}|\tilde w_\lambda^p-V_1^p|^{\frac{2N}{N+\alpha}}\right)^{\frac{N^2-\alpha^2}{4(N-2)}}\left(\int_{\mathbb R^N}|\tilde w_\lambda|^{2^*}\right)^{\frac{\alpha-N+4}{4}}\\
&\quad \to 0, \quad {\rm as} \ \ \lambda\to 0
\end{array}
\eqno(5.39)
$$
and 
$$
\begin{array}{rl}
&\left\|\frac{1}{|x|^4}(I_\alpha\ast |V_1|^p)(\frac{x}{|x|^2})[\tilde w_\lambda^{p-2}(\frac{x}{|x|^2})-V_1^{p-2}(\frac{x}{|x|^2})]\right\|_{L^{\frac{N}{2}}(\mathbb R^N)}^{N/2}\\
&\quad =\int_{\mathbb R^N}\frac{1}{|x|^{2N}}|(I_\alpha\ast |V_1|^p))(\frac{x}{|x|^2})|^{N/2}|\tilde w_\lambda^{p-2}(\frac{x}{|x|^2})-V_1^{p-2}(\frac{x}{|x|^2})|^{N/2}dx\\
&\quad= \int_{\mathbb R^N}|(I_\alpha\ast |V_1|^p)(z)|^{N/2}|\tilde w_\lambda^{p-2}(z)-V_1^{p-2}(z)|^{N/2}dz\\
&\quad \lesssim \left(\int_{\mathbb R^N}|V_1|^{2^*}\right)^{\frac{N^2-\alpha^2}{4(N-2)}}\left(\int_{\mathbb R^N}|\tilde w_\lambda^{p-2}(z)-V_1^{p-2}(z)|^{\frac{2N}{\alpha-N+4}}\right)^{\frac{\alpha-N+4}{4}}\\
&\quad \to 0, \quad {\rm as} \ \ \lambda\to 0.
\end{array}
\eqno(5.40)
$$
It  follows from (5.39) and (5.40) that
$$
\begin{array}{rl}
&\left\|\frac{1}{|x|^{4}}(I_\alpha\ast |\tilde w_\lambda|^p)(\frac{x}{|x|^2})\tilde w_\lambda^{p-2}(\frac{x}{|x|^2})-\frac{1}{|x|^{4}}(I_\alpha\ast |V_1|^p)(\frac{x}{|x|^2})V_1^{p-2}(\frac{x}{|x|^2})\right\|_{L^{\frac{N}{2}}(\mathbb R^N)}\\
&\le \left\| \frac{1}{|x|^{4}}(I_\alpha\ast |V_1|^p)(\frac{x}{|x|^2})\left[\tilde w_\lambda^{p-2}(\frac{x}{|x|^2})-V_1^{p-2}(\frac{x}{|x|^2})\right]\right\|_{L^{\frac{N}{2}}(\mathbb R^N)}\\
&\quad +\left\|\frac{1}{|x|^{4}}(I_\alpha\ast ( |\tilde w_\lambda|^p-|V_1|^p))(\frac{x}{|x|^2})\tilde w_\lambda^{p-2}(\frac{x}{|x|^2})\right\|_{L^{\frac{N}{2}}(\mathbb R^N)}\\
&\quad \to 0, \quad {\rm as} \ \ \lambda\to 0.
\end{array}
\eqno(5.41)
$$

{\bf Lemma 5.15.}  {\it Assume $N\ge 3$ and $2<q<2^*$. Then it holds that 
$$
\lim_{\lambda\to 0}\int_{|x|\le 4}\left|\frac{\lambda^\sigma\xi_\lambda^{\gamma/2}}{|x|^\gamma}K[\tilde w_\lambda]^{q-2}(x)\right|^{N/2}dx=0.
$$}

\begin{proof}  We divide the integral into two parts:
$$
I_\lambda^{(1)}(\frac{N}{2}):=\int_{|x|\le \lambda^{\sigma/2}\xi_\lambda/L_0}\left |\frac{\lambda^\sigma\xi_\lambda^{\gamma/2}}{|x|^\gamma}K[\tilde w_\lambda]^{q-2}(x)\right |^{N/2}dx,
$$
$$
I_\lambda^{(2)}(\frac{N}{2}):=\int_{\lambda^{\sigma/2}\xi_\lambda/L_0\le |x|\le 4}\left |\frac{\lambda^\sigma\xi_\lambda^{\gamma/2}}{|x|^\gamma}K[\tilde w_\lambda]^{q-2}(x)\right |^{N/2}dx.
$$

It follows from (5.37) and the H\"older inequality  that 
$$
\begin{array}{rcl}
I_\lambda^{(1)}(\frac{N}{2})&\le& \lambda^{\frac{\sigma N}{8}(4+(N-2)(q-2))}\xi_\lambda^{\frac{N}{4}(\gamma+(N-2)(q-2))}\int_{|x|\le \lambda^{\sigma/2}\xi_\lambda/L_0}|x|^{-\frac{\gamma N}{2}}e^{-\frac{N(q-2)}{4}\lambda^{\sigma/2}\xi_\lambda |x|^{-1}}dx\\
&=&\lambda^{\frac{\sigma N}{8}(8+2\gamma+(N-2)(q-2))}\xi_\lambda^{\frac{N}{4}(4+3\gamma+(N-2)(q-2))}\int_{L_0}^{+\infty}s^{\frac{\gamma N}{2}-N-1}e^{-\frac{N(q-2)}{4}s}ds\\
&\lesssim & \lambda^{\frac{\sigma N}{8}(12+2N-q(N-2))}\xi_\lambda^{\frac{N}{2}(8+4N-2q(N-2))},
\end{array}
$$
and 
$$
\begin{array}{rcl}
I_\lambda^{(2)}(\frac{N}{2})&=&\lambda^{\frac{\sigma N}{2}}\xi_\lambda^{\frac{\gamma N}{4}}\int_{\lambda^{\sigma/2}\xi_\lambda/L_0\le |x|\le 4}|x|^{-\frac{\gamma N}{2}}K[\tilde w_\lambda]^{\frac{N}{2}(q-2)}(x)dx\\
&\le &\lambda^{\frac{\sigma N}{2}}\xi_\lambda^{\frac{\gamma N}{4}}\left(\int_{\lambda^{\sigma/2}\xi_\lambda/L_0\le |x|\le 4}K[\tilde w_\lambda]^{2^*}dx\right)^{1-\frac{\gamma}{4}}
\left(\int_{\lambda^{\sigma/2}\xi_\lambda/L_0\le |x|\le 4} |x|^{-2N}dx\right)^{\frac{\gamma}{4}}\\
&\lesssim &\lambda^{\frac{\sigma N}{2}}\xi_\lambda^{\frac{\gamma N}{4}}\left(\int_{\lambda^{\sigma/2}\xi_\lambda/L_0}^4r^{-N-1}dr\right)^{\frac{\gamma}{4}}\\
&\lesssim &\lambda^{\frac{\sigma N}{8}(N-2)(q-2)}.
\end{array}
$$
From which the conclusion follows. 
\end{proof}

\begin{proof}[Proof of Proposition 5.14.] 

Since the Kelvin transform is linear and preserves the $D^{1,2}(\mathbb R^N )$ norm, 
it follows from (5.38), (5.41), Lemma 5.15 and Lemma 3.12 (i) that for any $r>1$, there exists $\lambda_r>0$ such that 
$$
\sup_{\lambda\in (0,\lambda_r)}\|K[\tilde w_\lambda]^r\|_{H^1(B_1)}\lesssim \sup_{\lambda\in (0,\lambda_r)}\|K[\tilde w_\lambda]\|_{L^{2^*}(\mathbb R^N)}\le C_r.
\eqno(5.42)
$$

Since $\alpha>N-4$, we have $\frac{2N}{N-\alpha}>\frac{N}{2}$.  Firstly,  we  show that  for some $r_0\in (\frac{N}{2}, \frac{2N}{N-\alpha})$, there holds
$$
\lim_{\lambda\to 0}I_\lambda(r_0)=0, 
\eqno(5.43)
$$
where
$$
I_\lambda(r_0):=\int_{|x|\le 4}\left|\frac{1}{|x|^4}\left [(I_\alpha\ast |\tilde w_\lambda|^{p})(\frac{x}{|x|^2})
-(I_\alpha\ast |V_1|^p)(\frac{x}{|x|^2})\right ]\tilde w_\lambda^{p-2}(\frac{x}{|x|^2})\right|^{r_0}dx.
$$

Since $\alpha>N-4$, for any  $r_0\in (\frac{N}{2}, \frac{2N}{N-\alpha})$, we have 
$$
 s_1=\frac{2N}{2N-(N-\alpha)r_0}>1, \quad  s_2=\frac{2N}{(N-\alpha)r_0}>1,
 $$
and 
 $$
 \frac{1}{s_1}+\frac{1}{s_2}=1.
 $$
Note that $(N-\alpha)r_0s_2=2N$, by the Hardy-Littlewood-Sobolev inequality, we get 
$$
\begin{array}{rl}
&\int_{\mathbb R^N}\frac{1}{|x|^{(N-\alpha)r_0s_2}}|(I_\alpha\ast (|\tilde w_\mu|^{p}-|V_1|^p))(\frac{x}{|x|^2})|^{r_0s_2}dx\\
&=\int_{\mathbb R^N}|(I_\alpha\ast (|\tilde w_\mu|^{p}-|V_1|^p))(z)|^{r_0s_2}dz\\
&\le C\left(\int_{\mathbb R^N}|\tilde w_\mu^p-V_1^p|^{\frac{2N}{(N-2)p}}\right)^{\frac{(N-2)p}{(N-2)p-2\alpha}}\to 0,  \   \   {\rm as}  \ \mu\to 0.
 \end{array}
\eqno(5.44)
$$
 By (6.42), (6.47), (6.50), Lemma  6.13 and the H\"older inequality,  we have 
$$
\begin{array}{rcl}
I_\lambda(r_0)&=&\int_{|x|\le 4}\left|\frac{1}{|x|^{N-\alpha}}(I_\alpha\ast (|\tilde w_\lambda|^{p}-|V_1|^p))(\frac{x}{|x|^2})K[\tilde w_\lambda]^{p-2}(x)\right|^{r_0}dx\\
&=&\left(\int_{|x|\le 4} K[\tilde w_\lambda]^{(p-2)r_0s_1}\right)^{\frac{1}{s_1}}\\
&\mbox{}& \cdot\left(\int_{\R^N}|x|^{-(N-\alpha)r_0s_2}\left|(I_\alpha\ast (|\tilde w_\lambda|^{p}-|V_1|^p))(\frac{x}{|x|^2})\right|^{r_0s_2}dx\right)^{\frac{1}{s_2}},
\end{array}
$$
which together with (5.42) and (5.44) yields (5.43).

Next, we consider the function 
$$
J_\lambda(r_0):=\int_{|x|\le 4}\left|\frac{1}{|x|^{4}}(I_\alpha\ast | V_1|^{p})(\frac{x}{|x|^2})\tilde w_\lambda^{p-2}(\frac{x}{|x|^2})\right|^{r_0}dx.
$$
Then it follows from the H\"older inequality that 
$$
\begin{array}{rcl}
J_\lambda(r_0)&=&\int_{\frac{1}{4}\le |z|< \infty} |z|^{4r_0-2N}\left |(I_\alpha\ast | V_1|^{p})(z)\tilde w_\lambda^{p-2}(z)\right|^{r_0}dz\\
&\le &\left(\int_{\frac{1}{4}\le |z|< \infty} |z|^{(4r_0-2N)s_1} |(I_\alpha\ast | V_1|^{p})(z)|^{r_0s_1}dz\right)^{\frac{1}{s_1}}\\
&\mbox{}& \qquad \cdot \left(\int_{\mathbb R^N}|\tilde w_\lambda^{p-2}(z)|^{r_0s_2}dz\right)^{\frac{1}{s_2}},
\end{array}
\eqno(5.45)
$$
where 
$$
s_1=\frac{2N}{2N-(\alpha-N+2)r_0}, \quad s_2=\frac{2N}{(\alpha-N+4)r_0}.
$$
Since  $p=\frac{N+\alpha}{N-2}>\frac{N}{N-2}$, we have
$$
\frac{N}{2}<\frac{2N}{\alpha-N+4}=\frac{2N}{(N-2)(p-2)}.
$$
Consider the function 
$$
h(r_0):=(4r_0-2N)s_1-(N-\alpha)r_0s_1+N.
$$ 
It is easy to check that $h(\frac{N}{2})=-N<0$,  and hence   $h(r_0)<0$ for $r_0>\frac{N}{2}$ which is close to $\frac{N}{2}$. 

Since $p>\frac{N}{N-2}$, we have $\int_{\mathbb R^N}|V_1|^p<\infty$. Notice that 
$$
|V_1|^p|x|^N\le |x|^{-(N-2)p+N}\to 0
$$ 
as $|x|\to \infty$, by Lemma 3.10, we have 
$$
\begin{array}{rcl}
&\mbox{}&\int_{\frac{1}{4}\le |z|< \infty} |z|^{(4r_0-2N)s_1} |(I_\alpha\ast | V_1|^{p})(z)|^{r_0s_1}dz\\
&\mbox{}&\lesssim 
\int_{\frac{1}{4}\le |z|< \infty} |z|^{(4r_0-2N)s_1} |z|^{-(N-\alpha)r_0s_1}dz\\
&\mbox{}&=\int_{\frac{1}{4}}^\infty r^{(4r_0-2N)s_1-(N-\alpha)r_0s_1+N-1}dr<\infty,
\end{array}
\eqno(5.46)
$$
On the other hand,  noting that $r_0s_2=\frac{2N}{(N-2)(q-2)}$, we have 
$$
\int_{\mathbb R^N}|\tilde w_\lambda^{p-2}(z)|^{r_0s_2}dz= \int_{\mathbb R^N}|\tilde w_\lambda|^{2^*}=\int_{\mathbb R^N}| w_\lambda|^{2^*}<C<\infty.
\eqno(5.47)
$$
Thus, from (5.45), (5.46) and (5.47), it follows that there is $\lambda_0>0$ such that
$$
\sup_{\lambda\in (0,\lambda_0)}J_\lambda(r_0)<+\infty,
$$
which together with (5.43) implies that for some  $\lambda_0>0$, there holds
$$
\sup_{\lambda\in (0,\lambda_0)}\int_{|x|\le 4}\left|\frac{1}{|x|^4}(I_\alpha\ast |\tilde w_\lambda|^{p})(\frac{x}{|x|^2})\right|^{r_0}dx<+\infty.
$$
It remains to prove that there exists $r_0>\frac{N}{2}$ and $\lambda_0>0$ such that 
$$
\sup_{\lambda\in (0,\lambda_0)}\int_{|x|\le 4}\left|\frac{\lambda^\sigma\xi_\lambda^{\gamma/2}}{|x|^\gamma}K[\tilde w_\lambda]^{q-2}(x)\right|^{r_0}dx\le C_{r_0}.
\eqno(5.48)
$$

It is easy to see that $0<\gamma<4$.  Therefore, we find
$$
\frac{N}{2}<\frac{2N}{\gamma}.
\eqno(5.49)
$$
If $\gamma>2$, we have 
$$
 \frac{N}{2}<\frac{2N}{\gamma}<\frac{N}{\gamma-2}.
\eqno(5.50)
$$
If $\gamma>8/3$, then
$$
\frac{N}{2}<\frac{2N}{\gamma}<\frac{N}{\gamma-2}<\frac{2N}{3\gamma-8}.
\eqno(5.51)
$$

Choose $s_0>1$ and $r_0>0$ 
such that
$$
N(q-2)s_0>2(s_0-1), \qquad \frac{N}{2}<r_0<\frac{2N}{\gamma s_0}. 
\eqno(5.52)
$$
Put $\eta_\lambda=\lambda^{\sigma/2}\xi_\lambda/L_0$. Then by (5.37) we have 
$$
\begin{array}{rcl}
I_\lambda^{(1)}(r_0):&=&\lambda^{\sigma r_0}\xi_\lambda^{\gamma r_0/2}\int_{|x|\le \eta_\lambda}\left|\frac{1}{|x|^\gamma}K[\tilde w_\lambda]^{(q-2)}(x)\right|^{r_0}dx\\
&\lesssim & \lambda^{\sigma r_0+\frac{\sigma(N-2)(q-2)}{4}r_0}\xi_{\lambda}^{\frac{\gamma}{2}r_0+\frac{(N-2)(q-2)}{2}r_0}
\int_{|x|\le \eta_\lambda}|x|^{-\gamma r_0}e^{-\frac{1}{2}(q-2)r_0\lambda^{\sigma/2}\xi_\lambda |x|^{-1}}dx\\
&\lesssim &\lambda^{\frac{\sigma}{4}[2N-(3\gamma-8)r_0]}\xi_\lambda^{\frac{1}{2}[2N-2(\gamma-2)r_0]}\int_{L_0}^{+\infty}t^{\gamma r_0-N-1}e^{-\frac{1}{2}(q-2)r_0t}dt.
\end{array}
$$
Therefore, we have 
$$
\lim_{\lambda\to 0}I_\lambda^{(1)}(r_0)=0.
\eqno(5.53)
$$
By (5.42) and the H\"older inequality, we also have
$$
\begin{array}{rcl}
I_\lambda^{(2)}(r_0):&=&\lambda^{\sigma r_0}\xi_\lambda^{\gamma r_0/2}\int_{\eta_\lambda\le |x|\le 4}\left|\frac{1}{|x|^\gamma}K[\tilde w_\lambda]^{(q-2)}(x)\right|^{r_0}dx\\
&\le & \lambda^{\sigma r_0}\xi_\lambda^{\gamma r_0/2}\left(\int_{|x|\le 4}K[\tilde w_\lambda]^{\frac{(q-2)r_0s_0}{s_0-1}}\right)^{1-\frac{1}{s_0}}\left(\int_{\eta_\lambda\le |x|\le 4}|x|^{-\gamma r_0s_0}dx\right)^{\frac{1}{s_0}}\\
&\lesssim &\lambda^{\sigma r_0}\xi_\lambda^{\gamma r_0/2}\left(\int_{\eta_\lambda}^4r^{-\gamma r_0s_0+N-1}dr\right)^{\frac{1}{s_0}}\\
&\lesssim &\lambda^{\sigma r_0}\xi_\lambda^{\gamma r_0/2}\eta_\lambda^{-\frac{\gamma r_0s_0-N}{s_0}}\\
&=& \lambda^{\frac{\sigma}{2s_0}[N-(\gamma-2)r_0s_0]}\xi_\lambda^{\frac{1}{2s_0}[2N-\gamma r_0s_0]}.
\end{array}
$$
Therefore, by the choices of $r_0$ and $s_0$, we obtain 
$$
\lim_{\lambda\to 0}I_\lambda^{(2)}(r_0)=0.
\eqno(5.54)
$$
Finally, (5.48) follows from (5.53) and (5.54),  the proof of Proposition 5.14 is complete. \end{proof}

{\bf Lemma 5.16.} { If $q>\frac{N}{N-2}$, then $\|\tilde w_\lambda\|_q^q\sim 1$ as $\lambda\to 0$. Furthermore, $\tilde w_\lambda\to V_1$ in $L^q(\mathbb R^N)$ as $\lambda\to 0$. }

\begin{proof}   Since $\tilde w_\lambda\to V_1$ in $L^{2^*}(\mathbb R^N)$, as in \cite[Lemma 4.6]{Moroz-1}, using the embeddings $L^{2^*}(B_1)\hookrightarrow L^q(B_1)$ we prove that
$\liminf_{\lambda\to 0}\|w_\lambda\|_q^q>0$.  

  On the other hand,  by virtue of  Proposition 5.14,  there exists a constant $C>0$ such that for all small $\lambda>0$, 
$$
\tilde w_\lambda(x)\le \frac{C}{(1+|x|)^{N-2}}, \qquad \forall x\in \mathbb R^N,
$$
which together with the fact that $q>\frac{N}{N-2}$ implies that $w_\lambda$ is bounded in $L^q(\mathbb R^N)$  uniformly for small $\lambda>0$,  and 
 by the dominated convergence theorem $w_\lambda\to U_1$ in $L^q(\mathbb R^N)$ as $\lambda\to 0$. 
\end{proof}

\begin{proof}[Proof of Theorem 2.2]
For $N\ge 5$, the conclusion follows directly from Lemmas 5.5, 5.6 and 5.8.  We only consider the cases $N=4$ and $N=3$.

We first note that a result similar to Lemma 3.2 holds for  $\tilde w_\lambda$ and $\tilde J_\lambda$. By  Lemma 5.9, we also have  $\tau_2(\tilde w_\lambda)=\tau_2(w_\lambda)$. Therefore, by (5.32),  we get
$$
\begin{array}{rcl}
m_0&\le& \sup_{t\ge 0} \tilde J_\lambda((\tilde w_\lambda)_t)+\lambda^\sigma\tau_2( \tilde w_\lambda)^{\frac{N}{2}}\left\{\frac{1}{q}\xi_\lambda^{N-\frac{N-2}{2}q}\int_{\mathbb R^N}|\tilde w_\lambda|^q-\frac{1}{2}\xi_\lambda^{2}\int_{\mathbb R^N}|\tilde w_\lambda|^2\right\}\\
&=&m_\lambda+\lambda^\sigma\tau_2(\tilde w_\lambda)^{\frac{N}{2}}\frac{q-2}{q(2^*-2)}\xi_\lambda^{N-\frac{N-2}{2}q}\int_{\mathbb R^N}|\tilde w_\lambda|^q,
\end{array}
\eqno(5.55)
$$
which implies that
$$
\xi_\lambda^{N-\frac{N-2}{2}q}\int_{\mathbb R^N}|\tilde w_\lambda|^q\ge \lambda^{-\sigma}\frac{q(2^*-2)}{(q-2)\tau_2(w_\lambda)^{\frac{N}{2}}}\delta_\lambda.
$$
Hence, by Lemma  5.5,  we obtain
$$
\xi_\lambda^{N-\frac{N-2}{2}q}\int_{\mathbb R^N}|\tilde w_\lambda|^q\gtrsim \lambda^{-\sigma}\delta_\lambda\gtrsim \left\{\begin{array}{rcl} 
 (\ln\frac{1}{\lambda})^{-\frac{4-q}{q-2}},  \quad  \quad &{\rm if}& \   \ N=4,\\
 \lambda^{\frac{2(6-q)}{(q-2)(q-4)}},   \qquad &{\rm if}& \    \  N=3.
 \end{array}\right.
\eqno(5.56)
 $$  
 Therefore, by Lemma 5.16, we have 
$$
\xi_\lambda\gtrsim  \left\{\begin{array}{rcl} 
 (\ln\frac{1}{\lambda})^{-\frac{1}{q-2}},  \quad  \quad &{\rm if}& \   \ N=4,\\
 \lambda^{\frac{4}{(q-2)(q-4)}},   \qquad &{\rm if}& \    \  N=3.
 \end{array}\right.
\eqno(5.57)
$$
On the other hand, if  $ N=3$, then by  (5.32) and Lemma 5.11 and Lemma 5.16, we have 
$$
\xi_\lambda^{\frac{q-2}{2}}\lesssim \frac{1}{\|\tilde w_\lambda\|_2^2}\lesssim\lambda^{\frac{\sigma}{2}}\xi_\lambda.
$$
Then 
$$
\xi_\lambda^{\frac{q-4}{2}}\lesssim \lambda^{\frac{\sigma}{2}}.
$$
Hence, noting that  $\sigma=\frac{2^*-2}{q-2}=\frac{4}{q-2}$,  for $q\in (4,6)$, we have 
$$
\xi_\lambda\lesssim \lambda^{\frac{4}{(q-2)(q-4)}}.
\eqno(5.58)
$$
 If $N=4$, then by (5.32) and Lemma 5.12 and Lemma 5.16,  we have 
$$
\xi_\lambda^{q-2}\lesssim \frac{1}{\|\tilde w_\lambda\|_2^2}\lesssim \frac{1}{-\ln(\lambda^\sigma\xi_\lambda^{2})}.
$$
Note that 
$$
-\ln(\lambda^\sigma\xi_\lambda^{2})=\sigma\ln\frac{1}{\lambda}+2\ln\frac{1}{\xi_\lambda}\ge \sigma\ln\frac{1}{\lambda},
$$
it follows that 
$$
\xi_\lambda^{q-2}\lesssim  \frac{1}{\|\tilde w_\lambda\|_2^2}\lesssim \left(\ln\frac{1}{\lambda}\right)^{-1}.
$$
Hence,  we  obtain 
$$
\xi_\lambda\lesssim   \left(\ln\frac{1}{\lambda}\right)^{-\frac{1}{q-2}}.
\eqno(5.59)
$$
Thus, it follows from (5.55), (5.58), (5.59)  and Lemma 5.16 that 
$$
\delta_\lambda=m_0-m_\lambda\lesssim \lambda^\sigma\xi_\lambda^{N-\frac{N-2}{2}q}\lesssim \left\{\begin{array}{rcl} \lambda^{\frac{2}{q-2}}(\ln\frac{1}{\lambda})^{-\frac{4-q}{q-2}},   \quad &{\rm if}& \   \ N=4,\\
 \lambda^{\frac{2}{q-4}}, \qquad  \qquad  \qquad &{\rm if}& \    \  N=3,
 \end{array}\right.
$$
which together with Lemma 5.5 implies that 
$$
\delta_\lambda\sim \lambda^\sigma\xi_\lambda^{N-\frac{N-2}{2}q}\sim \left\{\begin{array}{rcl} \lambda^{\frac{2}{q-2}}(\ln\frac{1}{\lambda})^{-\frac{4-q}{q-2}},   \quad &{\rm if}& \   \ N=4,\\
 \lambda^{\frac{2}{q-4}}, \qquad  \qquad  \qquad &{\rm if}& \    \  N=3.
 \end{array}\right.
\eqno(5.60)
$$
By (5.28) and  (5.32),  we  get
$$
m_\lambda=\frac{1}{2}\int_{\mathbb R^N}|\nabla \tilde w_\lambda|^2-\frac{1}{2p}\int_{\mathbb R^N}(I_\alpha\ast |\tilde w_\lambda|^{p})|\tilde w_\lambda^{p}-\frac{q-2}{q(2^*-2)}\lambda^\sigma\xi_\lambda^{N-\frac{N-2}{2}q}\int_{\mathbb R^N}|\tilde w_\lambda|^q.
$$
By (5.30) and (5.31), we get
$$
\int_{\mathbb R^N}|\nabla \tilde w|^2=\int_{\mathbb R^N}(I_\alpha\ast |\tilde w_\lambda|^p)|\tilde w_\lambda|^p+\frac{N(q-2)}{2q}\lambda^\sigma\xi_\lambda^{N-\frac{N-2}{2}q}\int_{\mathbb R^N}|\tilde w_\lambda |^q.
$$
Therefore, we have 
$$
m_\lambda=\frac{2+\alpha}{2(N+\alpha)}\int_{\mathbb R^N}|\nabla \tilde w_\lambda|^2-\frac{\alpha(N-2)(q-2)}{4q(N+\alpha)}\lambda^\sigma\xi_\lambda^{N-\frac{N-2}{2}q}\int_{\mathbb R^N}|\tilde w_\lambda|^q.
$$
Similarly, we have
$$
m_0=\frac{2+\alpha}{2(N+\alpha)}\int_{\mathbb R^N}|\nabla V_1|^2.
$$
Thus, by virtue of (5.60), we obtain
$$
\begin{array}{rcl}
\|\nabla V_1\|_2^2-\|\nabla \tilde w_\lambda\|_2^2
&=&\frac{2(N+\alpha)}{2+\alpha}\delta_\lambda-\frac{\alpha(N-2)(q-2)}{2q(2+\alpha)}\lambda^\sigma\xi_\lambda^{N-\frac{N-2}{2}q}\int_{\mathbb R^N}|\tilde w_\lambda|^q\\
&=& \left\{\begin{array}{rcl} O(\lambda^{\frac{2}{q-2}}(\ln\frac{1}{\lambda})^{-\frac{4-q}{q-2}}),   \quad &{\rm if}& \   \ N=4,\\
O( \lambda^{\frac{2}{q-4}}), \qquad  \quad  \qquad &{\rm if}& \    \  N=3.
 \end{array}\right.
 \end{array}
\eqno(5.61)
$$
By (5.28) and (5.31), we have 
$$
\begin{array}{lcl}
m_\lambda &=&(\frac{1}{2}-\frac{1}{2^*})\int_{\R^N}|\nabla \tilde w_\lambda|^2+(\frac{1}{2^*}-\frac{1}{2p})\int_{\R^N}(I_\alpha\ast |\tilde w_\lambda|^p)|\tilde w_\lambda|^p\\
&=&\frac{1}{N}\int_{\R^N}|\nabla \tilde w_\lambda|^2+\frac{\alpha(N-2)}{2N(N+\alpha)}\int_{\R^N}(I_\alpha\ast |\tilde w_\lambda|^p)|\tilde w_\lambda|^p.
\end{array}
$$
Similarly, we also have 
$$
m_0=\frac{1}{N}\int_{\R^N}|\nabla V_1|^2+\frac{\alpha(N-2)}{2N(N+\alpha)}\int_{\R^N}(I_\alpha\ast |V_1|^p)|V_1|^p.
$$
Then it follows from (5.60) and (5.61)  that 
$$
\begin{array}{rl}
&\int_{\R^N}(I_\alpha\ast |V_1|^{p})|V_1|^{p}-\int_{\R^N}(I_\alpha\ast |\tilde w_\lambda|^{p})|\tilde w_\lambda|^p\\
&=\frac{2N(N+\alpha)}{\alpha(N-2)}\left[(m_0-m_\lambda)-\frac{1}{N}(\int_{\R^N}|\nabla V_1|^2-\int_{\R^N}|\nabla \tilde w_\lambda|^2)\right]\\
&=\left\{\begin{array}{rcl} O(\lambda^{\frac{2}{q-2}}(\ln\frac{1}{\lambda})^{-\frac{4-q}{q-2}}),   \quad &{\rm if}& \   \ N=4,\\
 O(\lambda^{\frac{2}{q-4}}), \qquad   \qquad &{\rm if}& \    \  N=3.
 \end{array}\right.
\end{array}
\eqno(5.62)
$$
Since 
$
\|\nabla V_1\|_2^2=\int_{\R^N}(I_\alpha\ast |V_1|^p)|V_1|^p=S_\alpha^{\frac{N+\alpha}{2+\alpha}},
$
it follows from (5.61) and (5.62) that 
$$
\|\nabla \tilde w_\lambda\|_2^2=S_\alpha^{\frac{N+\alpha}{2+\alpha}}+\left\{\begin{array}{rcl} O(\lambda^{\frac{2}{q-2}}(\ln\frac{1}{\lambda})^{-\frac{4-q}{q-2}}),   \quad &{\rm if}& \   \ N=4,\\
 O(\lambda^{\frac{2}{q-4}}), \qquad   \qquad &{\rm if}& \    \  N=3,
 \end{array}\right.
$$
and 
$$
\int_{\R^N}(I_\alpha\ast |\tilde w_\lambda|^{p})|\tilde w_\lambda|^p=S_\alpha^{\frac{N+\alpha}{2+\alpha}}+\left\{\begin{array}{rcl} O(\lambda^{\frac{2}{q-2}}(\ln\frac{1}{\lambda})^{-\frac{4-q}{q-2}}),   \quad &{\rm if}& \   \ N=4,\\
 O(\lambda^{\frac{2}{q-4}}), \qquad   \qquad &{\rm if}& \    \  N=3.
 \end{array}\right.
$$
Finally, by (5.32), Lemma 5.11 and Lemma 5.12,  we obtain
$$
\|\tilde w_\lambda\|_2^2\sim \left\{\begin{array}{rcl}
\ln\frac{1}{\lambda},   \quad  \quad if \   \ N=4,\\
\lambda^{-\frac{2}{q-4}},  \quad if  \    \  N=3.
\end{array}\right.
$$
The statements on $v_\lambda$ follow from the corresponding results on $w_\lambda$ and $\tilde w_\lambda$. This completes the proof of Theorem 2.2.
\end{proof}

\vskip 5mm 

\section*{6. Proof of Theorem 2.3}

In this section, we always assume that  $q=2^*$ and $ p\in (1+\frac{\alpha}{N-2}, \frac{N+\alpha}{N-2})$.  It is easy to see that under the rescaling 
$$
w(x)=\mu^{\frac{N-2}{2(N-2)(p-1)-2\alpha}}v(\mu^{\frac{1}{(N-2)(p-1)-\alpha}}x), 
\eqno(6.1)
$$
the equation $(Q_\mu)$ is reduced to 
$$
-\Delta w+\mu^\sigma w=\mu^\sigma(I_\alpha\ast |w|^p)|w|^{p-2}w+|w|^{2^*-2}w, 
$$
where $\sigma:=\frac{2}{(N-2)(p-1)-\alpha}>1$.

The associated energy functional defined by
$$
J_\mu(w):=\frac{1}{2}\int_{\mathbb R^N}|\nabla w|^2+\mu^\sigma |w|^2-\frac{1}{2p}\mu^\sigma\int_{\mathbb R^N}(I_\alpha\ast |w|^p)|w|^p-\frac{1}{2^*}\int_{\mathbb R^N}|w|^{2^*}=I_\mu(v).
$$

{\bf Lemma 6.1.}  {\it Let $\mu>0$, $v\in H^1(\mathbb R^N)$ and $w$ is the rescaling of $v$ defined in (6.1). Then

(1) $ \ \|\nabla w\|_2^2=\|\nabla v\|_2^2, \quad  \|w\|_{2^*}^{2^*}=\|v\|_{2^*}^{2^*},$

(2)  \  $\mu^{\sigma}\|w\|_2^2=\|v\|_2^2, \quad  \mu^\sigma \int_{\mathbb R^N}(I_\alpha\ast |w|^p)|w|^p= \mu \int_{\mathbb R^N}(I_\alpha\ast |v|^p)|v|^p$,

(3)  \  $J_\mu(w)=I_\mu(v)$.
}

We define the Nehari manifolds as follows
$$
\mathcal{N}_\mu=
\left\{w\in H^1(\mathbb R^N)\setminus\{0\} \ \left | \ \int_{\mathbb R^N}|\nabla w|^2+\mu^\sigma \int_{\mathbb R^N}|w|^2=\mu^\sigma\int_{\mathbb R^N}(I_\alpha\ast |w|^p)|w|^p+
\int_{\mathbb R^N}|w|^{2^*}\  \right. \right\}
$$
and 
$$
\mathcal{N}_0=
\left\{w\in D^{1,2}(\mathbb R^N)\setminus\{0\} \ \left | \ \int_{\mathbb R^N}|\nabla w|^2=\int_{\mathbb R^N}|w|^{2^*}\  \right. \right\}. 
$$
Then 
$$
m_\mu:=\inf_{w\in \mathcal {N}_\mu}J_\mu(w) \quad {\rm and} \quad 
m_0:=\inf_{u\in \mathcal {N}_0}J_0(u)
$$
are well-defined and positive. Moreover, $J_0$ is attained on $\mathcal N_0$ and 
$m_0:=\inf_{w\in \mathcal {N}_0}J_0(w)$.

 For $w\in H^1(\mathbb R^N)\setminus \{0\}$, we set
$$
\tau_3(w)=\frac{\int_{\mathbb R^N}|\nabla w|^2}{\int_{\mathbb R^N}|w|^{2^*}}.
\eqno(6.2)
$$
Then $(\tau_3(w))^{\frac{N-2}{4}}w\in \mathcal N_0$ for any $w\in H^1(\mathbb R^N)\setminus\{0\}$,  and $w\in \mathcal N_0$  if and only if $\tau_3(w)=1$.

Define the Pohozaev manifold as follows
$$
\mathcal P_\mu:=\{w\in H^1(\mathbb R^N)\setminus\{0\} \   | \  P_\mu(w)=0 \},
$$
where 
$$
\begin{array}{rcl}
P_\mu(w):&=&\frac{N-2}{2}\int_{\mathbb R^N}|\nabla w|^2+\frac{ N}{2}\mu^\sigma \int_{\mathbb R^N}|w|^2
\\ \\
&\quad &-\frac{N+\alpha}{2p}\mu^\sigma\int_{\mathbb R^N}(I_\alpha\ast |w|^p)|w|^p-\frac{N}{2^*}\int_{\mathbb R^N}|w|^{2^*}.
\end{array}
\eqno(6.3)
$$
Then by Lemma 3.1, $w_\mu\in \mathcal P_\mu$. Moreover,  we have a minimax characterizations for the least energy level $m_\mu$ similar to

{\bf Lemma 6.2.} {\it  Assume $p\in (1+\frac{\alpha}{N-2}, \frac{N+\alpha}{N-2})$. Then the rescaled family of solutions $\{w_\mu\}$ is bounded in $H^1(\mathbb R^N)$.}
\begin{proof}   Since $\{w_\mu\}$ is bounded in $D^{1,2}(\mathbb R^N)$, it suffices to show that it is also bounded 
in $L^2(\mathbb R^N)$.  Since  $w_\mu\in \mathcal N_\mu\cap \mathcal P_\mu$, we obtain
$$
\int_{\mathbb R^N}|\nabla w_\mu|^2+\mu^\sigma\int_{\mathbb R^N}|w_\mu|^2=\mu^\sigma\int_{\mathbb R^N}(I_\alpha\ast |w_\mu|^p)|w_\mu|^p+\int_{\mathbb R^N}|w_\mu|^{2^*},
$$
$$
\frac{N-2}{2}\int_{\mathbb R^N}|\nabla w_\mu|^2+\frac{N}{2}\mu^\sigma\int_{\mathbb R^N}|w_\mu|^2=\frac{N+\alpha}{2p}\mu^\sigma\int_{\mathbb R^N}(I_\alpha\ast |w_\mu|^p)|w_\mu|^p+\frac{N}{2^*}\int_{\mathbb R^N}|w_\mu|^{2^*}.
$$
Therefore, we have 
$$
\|w_\mu\|_2^2=\frac{(N+\alpha)-p(N-2)}{2p}\int_{\mathbb R^N}(I_\alpha\ast |w_\mu|^p)|w_\mu|^p.
\eqno(6.4)
$$
By the Hardy-Littlewood-Sobolev, the Sobolev  and the interpolation inequalities,  we obtain
$$
\int_{\mathbb R^N}(I_\alpha\ast |w_\mu|^p)|w_\mu|^p\le C\|w_\mu\|_{\tilde p}^{2p}\le C\left(\int_{\mathbb R^N}|w_\mu|^2\right)^{\frac{2^*-\tilde p}{2^*-2}\frac{N+\alpha}{N}}
\left(\int_{\mathbb R^N}|w_\mu|^{2^*}\right)^{\frac{\tilde p-2}{2^*-2}\frac{N+\alpha}{N}}
\eqno(6.5)
$$
and 
$$
\int_{\mathbb R^N}(I_\alpha\ast |w_\mu|^p)|w_\mu|^p
\le C\|w_\mu\|_2^{\frac{2(2^*-\tilde p)}{2^*-2}\frac{N+\alpha}{N}}
\left(\frac{1}{S}\int_{\mathbb R^N}|\nabla w_\mu|^2\right)^{\frac{2^*(\tilde p-2)}{2(2^*-2)}\frac{N+\alpha}{N}},
$$
here and in what follows, $\tilde p=\frac{2Np}{N+\alpha}$. Therefore, we get
$$
\| w_\mu\|_2^{2-\frac{2(2^*-\tilde p)}{2^*-2}\frac{N+\alpha}{N}}\le C\frac{(N+\alpha)-p(N-2)}{2p}\left(\frac{1}{S}\int_{\mathbb R^N}|\nabla w_\mu|^2\right)^{\frac{2^*(\tilde p-2)}{2(2^*-2)}\frac{N+\alpha}{N}}.
$$
Since $p>1+\frac{\alpha}{N-2}$ implies that $2-\frac{2(2^*-\tilde p)}{2^*-2}\frac{N+\alpha}{N}>0$, it follows from the boundedness of $\|\nabla w_\mu\|_2$  that $w_\mu$ is bounded in $L^2(\mathbb R^N)$. 
\end{proof}

Now, we give  the following estimation on the least energy:

{\bf Lemma 6.3.} {\it  Let $N\ge 5$ and $p\in (1+\frac{\alpha}{N-2},\frac{N+\alpha}{N-2})$, then
  $$
  m_0-m_\mu\sim \mu^\sigma\qquad  as \quad \mu\to 0.
  $$
  If $N=4$ and $p\in (1+\frac{\alpha}{2}, 2+\frac{\alpha}{2})$,  or $N=3$ and $p\in (2+\alpha, 3+\alpha)$, then 
  $$
  m_0-m_\mu\lesssim \mu^\sigma \qquad  as \quad \mu\to 0.
  $$} 
 \begin{proof}   
By (6.5) and  the boundedness of  $\{w_\mu\}$ in $H^1(\mathbb R^N)$, we get  
$$
\frac{\int_{\mathbb R^N}(I_\alpha\ast |w_\mu|^p)|w_\mu|^p-\int_{\mathbb R^N}|w_\mu|^2}{\int_{\mathbb R^N}|w_\mu|^{2^*}}
\le C\eta_\mu^{\frac{2^*-\tilde p}{2^*-2}\frac{N+\alpha}{N}}\left(\int_{\mathbb R^N}|w_\mu|^{2^*}|\right)^{\frac{\alpha}{N}}-\eta_\mu
\le \tilde C\eta_\mu^{\frac{2^*-\tilde p}{2^*-2}\frac{N+\alpha}{N}}-\eta_\mu,
$$
where 
$
\eta_\mu:=\frac{\int_{\mathbb R^N}|w_\mu|^2}{\int_{\mathbb R^N}|w_\mu|^{2^*}}.
$

Since $\frac{2^*-\tilde p}{2^*-2}\frac{N+\alpha}{N}<1$, it is easy to show that the function $g: [0,+\infty)\to \mathbb R$ defined by
$$
g(x)= \tilde Cx^{\frac{2^*-\tilde p}{2^*-2}\frac{N+\alpha}{N}}-x\le \tilde C^{\frac{1}{1-\tau_0}}\tau_0^{\frac{\tau_0}{1-\tau_0}}(1-\tau_0), \qquad \tau_0=\frac{2^*-\tilde p}{2^*-2}\frac{N+\alpha}{N}.
$$
Therefore, 
by $w_\mu\in \mathcal N_\mu$ and the definition of $\tau_3(w_\mu)$, we get
$$
\tau_3(w_\mu)=\frac{\int_{\mathbb R^N}|\nabla w_\mu|^2}{\int_{\mathbb R^N}|w_\mu|^{2^*}}
=1+\mu^\sigma\frac{\int_{\mathbb R^N}(I_\alpha\ast |w_\mu|^p)|w_\mu|^p-\int_{\mathbb R^N}|w_\mu|^2}{\int_{\mathbb R^N}|w_\mu|^{2^*}}\le 1+C\mu^\sigma.
$$
If $N\ge 3$,  then by Lemma 3.2 and the boundedness of $\{w_\mu\}$, we find
$$
\begin{array}{rcl}
m_0&\le& \sup_{t\ge 0}J_0((w_\mu)_t)=J_0((w_\mu)_{t_\mu})\\
&\le &\sup_{t\ge 0} J_\mu((w_\mu)_t)+\mu^\sigma \left(\frac{t_\mu^{N+\alpha}}{2p}\int_{\mathbb R^N}(I_\alpha\ast |w_\mu|^p)|w_\mu|^p-\frac{t_\mu^{N}}{2}\int_{\mathbb R^N}|w_\mu|^{2}\right)\\
&\le & m_\mu+\mu^\sigma\left(\frac{(\tau_3(w_\mu))^{(N+\alpha)/2}}{2p}C\|w_\mu\|_{\tilde p}^{2p}-\frac{(\tau_3(w_\mu))^{N/2}}{2}\int_{\mathbb R^N}|w_\mu|^{2}\right)\\
&\le &m_\mu +C\mu^\sigma,
\end{array}
$$
where we have used the fact that $t_\mu=\left(\frac{\int_{\mathbb R^N}|\nabla w_\mu|^2}{\int_{\mathbb R^N}|w_\mu|^{2^*}}\right)^{1/2}=(\tau_3(w_\mu))^{1/2}$. 

For each $\rho>0$, the family  $W_\rho(x):=\rho^{-\frac{N-2}{2}}W_1(x/\rho)$ are radial ground states of  $-\Delta u=u^{2^*-1}$, and verify that
 $$
 \|W_\rho\|_2^2=\rho^2\|W_1\|_2^2, \quad \int_{\mathbb R^N}(I_\alpha\ast |W_\rho|^p)|W_\rho|^p=\rho^{(N+\alpha)-p(N-2)}\int_{\mathbb R^N}(I_\alpha\ast |W_1|^p)|W_1|^p.
 \eqno(6.6)
 $$ 
 Let $g_0(\rho)=\frac{1}{2p}\int_{\mathbb R^N}(I_\alpha\ast |W_\rho|^p)|W_\rho|^p-\frac{1}{2}\int_{\mathbb R^N}|W_\rho|^2$. Then there exists an unique $\rho_0=\rho_0(p)\in (0,+\infty)$ given by 
 $$
 \rho_0=\left(\frac{[(N+\alpha)-p(N-2)]\int_{\mathbb R^N}(I_\alpha\ast |W_1|^p)|W_1|^p}{2p\int_{\mathbb R^N}|W_1|^2}\right)^ {\frac{1}{(N-2)(p-1)-\alpha}}, 
 $$
 such that 
 $$
 g_0(\rho_0)=\sup_{\rho>0}g_0(\rho)=\frac{(N-2)(p-1)-\alpha}{2[(N+\alpha)-p(N-2)]}\rho_0^2\int_{\mathbb R^N}|W_1|^2.
 $$ 
 Let $W_0=W_{\rho_0}$, then there exists $t_\mu\in (0,+\infty)$ such that
 $$
\begin{array}{rcl}
m_\mu&\le &\sup_{t\ge 0}J_\mu(tW_0)=J_\mu(t_\mu W_0)\\
&=&\frac{t^2_\mu}{2}\int_{\mathbb R^N}|\nabla W_0|^2-\frac{t^{2^*}_\mu}{2^*}\int_{\mathbb R^N}|W_0|^{2^*}+\mu^\sigma\int_{\mathbb R^N}\frac{t^2_\mu}{2}|W_0|^2-\frac{t^{2p}_\mu}{2p}(I_\alpha\ast |W_0|^p)|W_0|^p\\
&\le &\sup_{t\ge 0}\left(\frac{t^2}{2}-\frac{t^{2^*}}{2^*}\right)\int_{\mathbb R^N}|\nabla W_0|^2 +\mu^\sigma\int_{\mathbb R^N}\frac{t^2_\mu}{2}|W_0|^2-\frac{t^{2p}_\mu}{2p}(I_\alpha\ast |W_0|^p)|W_0|^p\\
&=& m_0+\mu^\sigma\int_{\mathbb R^N}\frac{t^2_\mu}{2}|W_0|^2-\frac{t^{2p}_\mu}{2p}(I_\alpha\ast |W_0|^p)|W_0|^p.
\end{array}
\eqno(6.7)
$$
If $t_\mu\ge 1$, then
$$
\int_{\mathbb R^N}|\nabla W_0|^2+\mu^\sigma\int_{\mathbb R^N}|W_0|^2\ge t_\mu^{\min\{2p-2,2^*-2\}}\left\{\int_{\mathbb R^N}|W_0|^{2^*}+\mu^\sigma\int_{\mathbb R^N}(I_\alpha\ast |W_0|^p)|W_0|^p\right\}.
$$
Hence
$$
t_\mu\le \max\left\{1, \left(\frac{Nm_0+\mu^\sigma\int_{\mathbb R^N}|W_0|^2}{Nm_0+\mu^\sigma\int_{\mathbb R^N}(I_\alpha\ast |W_0|^p)|W_0|^p}\right)^{\frac{1}{\min\{2p-2,2^*-2\}}}\right\}.
$$
If $t_\mu\le  1$, then
$$
\int_{\mathbb R^N}|\nabla W_0|^2+\mu^\sigma\int_{\mathbb R^N}|U_0|^2\le t_\mu^{\min\{2p-2,2^*-2\}}\left\{\int_{\mathbb R^N}|W_0|^{2^*}+\mu^\sigma\int_{\mathbb R^N}(I_\alpha\ast |W_0|^p)|W_0|^p\right\}.
$$
Hence
$$
t_\mu\ge \min\left\{1, \left(\frac{Nm_0+\mu^\sigma\int_{\mathbb R^N}|W_0|^2}{Nm_0+\mu^\sigma\int_{\mathbb R^N}(I_\alpha\ast |W_0|^p)|W_0|^p}\right)^{\frac{1}{\min\{2p-2,2^*-2\}}}\right\}.
$$
Since
 $$
 \int_{\mathbb R^N}(I_\alpha\ast |W_0|^p)|W_0|^p>\int_{\mathbb R^N}|W_0|^2,
 $$
we conclude that 
$$
\left(\frac{Nm_0+\mu^\sigma\int_{\mathbb R^N}|W_0|^2}{Nm_0+\mu^\sigma\int_{\mathbb R^N}(I_\alpha\ast |W_0|^p)|W_0|^p}\right)^{\frac{1}{\min\{2p-2,2^*-2\}}}\le t_\mu
\le 1.
$$
Let $A_\mu:=\frac{\int_{\mathbb R^N}(I_\alpha\ast |U_0|^p)|U_0|^p-\int_{\mathbb R^N} |U_0|^2}{Nm_0+\mu^\sigma\int_{\mathbb R^N}(I_\alpha\ast |U_0|^p)|U_0|^p}$, then 
$$
A_\mu\le \frac{1}{Nm_0}[\int_{\mathbb R^N}(I_\alpha\ast |W_0|^p)|W_0|^p-\int_{\mathbb R^N} |W_0|^2]
$$ 
and 
$$
[1-\mu^\sigma A_\mu]^{\frac{1}{\min\{2p-2,2^*-2\}}}\le t_\mu\le 1.
$$
Let $g(t):=\frac{t^2}{2}\int_{\mathbb R^N}|W_0|^2-\frac{t^{2p}}{2p}\int_{\mathbb R^N}(I_\alpha\ast |W_0|^p)|W_0|^p$, and $h(x):=g([1-x]^{\frac{1}{\min\{2p-2,2^*-2\}}})$ for $x\in [0,1]$.
Then  $g(t)$ has an unique miximum point at $t_0=\left(\frac{\int_{\mathbb R^N}|W_0|^2}{\int_{\mathbb R^N}(I_\alpha\ast |W_0|^p)|W_0|^p}\right)^{\frac{1}{2p-2}}$ and is strictly decreasing in $(t_0,1)$, and for small $x>0$, we have 
$$
h'(x)=\frac{1}{2p-2}[1-x]^{-\frac{p-2}{p-1}}\left[-\int_{\mathbb R^N}|W_0|^2+(1-x)\int_{\mathbb R^N}(I_\alpha\ast |W_0|^p)|W_0|^p\right]>0.
$$
Therefore, for small $\mu>0$, we get
$$
\begin{array}{rcl}
g(t_\mu)&\le& g([1-\mu^\sigma A_\mu]^\frac{1}{\min\{2p-2,2^*-2\}})\\
&=&h(\mu^\sigma A_\mu)\\
&=&\frac{1}{2}\int_{\mathbb R^N}|W_0|^2-\frac{1}{2p}\int_{\mathbb R^N}(I_\alpha\ast |W_0|^p)|W_0|^p+h'(\xi)\mu^\sigma A_\mu, 
\end{array}
$$
where $\xi\in (0, \mu^\sigma A_\mu).$
Since for small $x>0$, there holds
$$
h'(x)\le \frac{2}{2p-2}\left[\int_{\mathbb R^N}(I_\alpha\ast |W_0|^p)|W_0|^p-\int_{\mathbb R^N}|W_0|^2\right], 
$$
it follows that 
$$
\begin{array}{rcl}
g(t_\mu)&\le& \frac{1}{2}\int_{\mathbb R^N}|W_0|^2-\frac{1}{2p}\int_{\mathbb R^N}(I_\alpha \ast |W_0|^p)|W_0|^p\\
&\quad &+\frac{\mu^\sigma}{Nm_0(p-1)}\left[\int_{\mathbb R^N}(I_\alpha \ast |W_0|^p)|W_0|^p-\int_{\mathbb R^N}|W_0|^2\right]^2\\
&\le &-g_0(\rho_0)+C\mu^\sigma
\end{array}
$$
from which the conclusion follows.
\end{proof}

{\bf Lemma 6.4.}  {\it There exists a constant $\varpi=\varpi(q)>0$ such that  for $\mu>0$ small, 
$$
m_\mu\le \left\{\begin{array}{rcl}
m_0-\mu^\sigma \left(\ln\frac{1}{\mu}\right)^{\frac{2p-4-\alpha}{2p-2-\alpha}}\varpi=m_0-\mu^{\frac{2}{2p-2-\alpha}}\left(\ln\frac{1}{\mu}\right)^{\frac{2p-4-\alpha}{2p-2-\alpha}}\varpi,      &if&     N=4.\\
m_0-\mu^{\sigma-\frac{p-3-\alpha}{(p-1-\alpha)(p-2-\alpha)}}\varpi=m_0-\mu^{\frac{1}{p-2-\alpha}}\varpi, \quad \quad \quad  \quad \qquad    &if&     N=3, \  p>2+\alpha.
\end{array} \right.
$$}\begin{proof}  
Let  $\rho>0$ and $R$ be a large parameter, and $\eta_R\in C_0^\infty(\mathbb R)$ is a cut-off function such that $\eta_R(r)=1$ for $|r|<R$, $0<\eta_R(r)<1$ for $R<|r|<2R$, $\eta_R(r)=0$ for $|r|>2R$ and $|\eta'_R(r)|\le 2/R$.  

For $\ell\gg 1$, a straightforward computation shows that 
$$
\int_{\mathbb R^N}|\nabla (\eta_\ell W_1)|^2=\left\{ \begin{array}{rcl} Nm_0+O(\ell^{-2}), \quad \qquad &{\rm if} & \  \ N=4,\\
Nm_0+O(\ell^{-1}),  \qquad \  \    \ &{\rm if} & \   \ N=3.
\end{array}\right.
$$ 
$$
\int_{\mathbb R^N} |\eta_\ell W_1|^{2^*}=Nm_0+O(\ell^{-N}).
$$
$$
\int_{\mathbb R^N}|\eta_\ell W_1|^2=\left\{ \begin{array}{rcl} \ln \ell(1+o(1)), \quad &{\rm if} & \  \ N=4,\\
\ell(1+o(1)), \quad  \quad &{\rm if} & \   \ N=3.
\end{array}\right.
$$ 
By Lemma 3.2, we find
$$
\begin{array}{rcl}
m_\mu&\le &\sup_{t\ge 0}J_\mu((\eta_RW_\rho)_t)=J_\mu((\eta_RW_\rho)_{t_\mu})\\
&\le &\sup_{t\ge 0}\left(\frac{t^{N-2}}{2}\int_{\mathbb R^N}|\nabla(\eta_RW_\rho)|^2-\frac{t^{N}}{2^*}\int_{\mathbb R^N}|\eta_RW_\rho|^{2^*}\right)\\
&\quad & -\mu^\sigma
t_\mu^N\left[\int_{\mathbb R^N}\frac{t_\mu^\alpha}{2p}(I_\alpha\ast |\eta_RW_\rho|^p)|\eta_RW_\rho|^p-\frac{1}{2}|\eta_RW_\rho|^2\right]\\
&= & (I)-\mu^\sigma (II),
\end{array}
\eqno(6.8)
$$
where $t_\mu\in (0, +\infty)$ is the unique critical point of the function $g(t)$ defined by 
$$
\begin{array}{rcl}
g(t)&=&\frac{t^{N-2}}{2}\int_{\mathbb R^N}|\nabla(\eta_RW_\rho)|^2+ \frac{t^N}{2}\mu^\sigma\int_{\mathbb R^N}|\eta_RW_\rho|^2\\
&\quad &-\frac{t^{N+\alpha}}{2p}\mu^\sigma\int_{\mathbb R^N}(I_\alpha\ast |\eta_RW_\rho|^p)|\eta_RW_\rho|^p-\frac{t^N}{2^*}\int_{\mathbb R^N}|\eta_RW_\rho|^{2^*}.
\end{array}
$$
That is, $t=t_\mu$ solves the following equation $\ell_1(t)=\ell_2(t)$, where
$$
\ell_1(t):=\frac{N-2}{2t^2}\int_{\mathbb R^N}|\nabla(\eta_RW_\rho)|^2
$$
and 
$$
\ell_2(t):=\left\{\frac{N+\alpha}{2p}t^\alpha\mu^\sigma\int_{\mathbb R^N}(I_\alpha\ast |\eta_RW_\rho|^p)|\eta_RW_\rho|^p
+\frac{N}{2^*}\int_{\mathbb R^N}|\eta_RW_\rho|^{2^*}-\frac{N}{2}\mu^\sigma\int_{\mathbb R^N}|\eta_RW_\rho|^2\right\}.
$$

If $t_\mu\ge 1$, then
$$
\begin{array}{rcl}\frac{N-2}{2t_\mu^{2}}\int_{\mathbb R^N}|\nabla(\eta_RW_\rho)|^2&\ge& \frac{N+\alpha}{2p}\mu^\sigma\int_{\mathbb R^N}(I_\alpha\ast |\eta_RW_\rho|^p)|\eta_RW_\rho|^p
+\frac{N}{2^*}\int_{\mathbb R^N}|\eta_RW_\rho|^{2^*}\\
&\quad &-\frac{N}{2}\mu^\sigma\int_{\mathbb R^N}|\eta_RW_\rho|^2,
\end{array}
$$
and hence 
$$
\begin{array}{rcl}
t_\mu
&\le &\left(\frac{\int_{\mathbb R^N}|\nabla(\eta_RW_\rho)|^2}{\int_{\mathbb R^N}|\eta_RW_\rho|^{2^*}+\frac{2\mu^\sigma}{N-2}\{\frac{N+\alpha}{2p}\int_{\mathbb R^N}(I_\alpha\ast |\eta_RW_\rho|^p)|\eta_RW_\rho|^p-\frac{N}{2}\int_{\mathbb R^N}|\eta_RW_\rho|^2\}}\right)^{\frac{1}{2}}\\
&\le & \left(\frac{\int_{\mathbb R^N}|\nabla(\eta_RW_\rho)|^2}{\int_{\mathbb R^N}|\eta_RW_\rho|^{2^*}}\right)^{\frac{1}{2}}.
\end{array}
\eqno(6.9)
$$
If $t_\mu\le 1$, then
$$
\begin{array}{rcl}
\frac{N-2}{2t_\mu^2}\int_{\mathbb R^N}|\nabla(\eta_RW_\rho)|^2&\le& \frac{N+\alpha}{2p}\mu^\sigma\int_{\mathbb R^N}(I_\alpha\ast |\eta_RW_\rho|^p)|\eta_RW_\rho|^p
+\frac{N}{2^*}\int_{\mathbb R^N}|\eta_RW_\rho|^{2^*}\\
&\quad &-\frac{N}{2}\mu^\sigma\int_{\mathbb R^N}|\eta_RW_\rho|^2,
\end{array}
$$
and hence 
$$
\begin{array}{rcl}
t_\mu
&\ge&  \left(\frac{\int_{\mathbb R^N}|\nabla(\eta_RW_\rho)|^2}{\int_{\mathbb R^N}|\eta_RW_\rho|^{2^*}+\frac{2\mu^\sigma}{N-2}\{\frac{N+\alpha}{2p}\int_{\mathbb R^N}(I_\alpha\ast |\eta_RW_\rho|^p)|\eta_RW_\rho|^p-\frac{N}{2}\int_{\mathbb R^N}|\eta_RW_\rho|^2\}}\right)^{\frac{1}{2}}\\
&\ge & \left(\frac{\int_{\mathbb R^N}|\nabla(\eta_RW_\rho)|^2}{\int_{\mathbb R^N}|\eta_RW_\rho|^{2^*}}\right)^{\frac{1}{2}}
\left\{1
-\mu^\sigma\frac{\frac{N+\alpha}{2p}\int_{\mathbb R^N}(I_\alpha\ast |\eta_RW_\rho|^p)|\eta_RW_\rho|^p-\frac{N}{2}\int_{\mathbb R^N}|\eta_RW_\rho|^2}
{(N-2)\int_{\mathbb R^N}|\eta_RW_\rho|^{2^*}}\right\}.
\end{array}
\eqno(6.10)
$$
Therefore, we obtain
$$
|t_\mu-1|\le \left|\left(\frac{\int_{\mathbb R^N}|\nabla(\eta_RW_\rho)|^2}{\int_{\mathbb R^N}|\eta_RW_\rho|^{2^*}}\right)^{\frac{1}{2}}-1\right|
+\mu^\sigma \left(\frac{\int_{\mathbb R^N}|\nabla(\eta_RW_\rho)|^2}{\int_{\mathbb R^N}|\eta_RW_\rho|^{2^*}}\right)^{\frac{1}{2}}
\frac{\psi(\rho)}{(N-2)\int_{\mathbb R^N}|\eta_RW_\rho|^{2^*}},
$$
where $\psi(\rho):=\frac{N+\alpha}{2p}\int_{\mathbb R^N}(I_\alpha\ast |\eta_RW_\rho|^p)|\eta_RW_\rho|^p-\frac{N}{2}\int_{\mathbb R^N}|\eta_RW_\rho|^2$.

Set $\ell=R/\rho$, then
$$
(I)=\frac{1}{N}\frac{\|\nabla(\eta_\ell W_1)\|_2^N}{\|\eta_\ell W_1\|_{2^*}^N}=\left\{\begin{array}{rcl} m_0+O(\ell^{-2}), \  \  {\rm if} \   \  N=4,\\
m_0+O(\ell^{-1}),   \   \   {\rm if}  \   \  N=3.\end{array}\right.
\eqno(6.11)
$$
Since 
$$
\begin{array}{rcl}
\varphi(\rho):&=&\frac{1}{2p}\int_{\mathbb R^N}(I_\alpha\ast |\eta_RW_\rho|^p)|\eta_RW_\rho|^p-\frac{1}{2}\int_{\mathbb R^N}|\eta_RW_\rho|^2\\
&=&\frac{1}{2p}\rho^{(N+\alpha)-p(N-2)}\int_{\mathbb R^N}(I_\alpha\ast |\eta_\ell W_1|^p)|\eta_\ell W_1|^p-\frac{1}{2}\rho^2\int_{\mathbb R^N}|\eta_\ell W_1|^2
\end{array}
$$
take its maximum value $\varphi(\rho_\ell)$ at the unique point 
$$
\begin{array}{rcl}
\rho_\ell:&=&\left(\frac{[(N+\alpha)-p(N-2)]\int_{\mathbb R^N}(I_\alpha\ast |\eta_\ell W_1|^p)|\eta_\ell W_1|^p}{2p\int_{\mathbb R^N}|\eta_\ell W_1|^2}\right)^{\frac{1}{(N-2)(p-1)-\alpha}}\\
&\sim &\left\{\begin{array}{rcl}
(\ln \ell)^{-\frac{1}{2+p(N-2)-(N+\alpha)} }\quad &{\rm if}&  \   N=4,\\
\ell^{-\frac{1}{2+p(N-2)-(N+\alpha)}} \quad &{\rm if}&  \   N=3.
\end{array}\right.
\end{array}
$$
and 
$$
\begin{array}{rcl}
\varphi(\rho_\ell)&=&\sup_{\rho\ge 0}\varphi(\rho)\\
&=&\frac{2+p(N-2)-(N+\alpha)}{2[(N+\alpha)-p(N-2)]}\rho_\ell^2\int_{\mathbb R^N}|\eta_\ell W_1|^2\\
&\le &C(\int_{\mathbb R^N}|\eta_\ell W_1|^{2^*})^{\frac{(N+\alpha)-p(N-2)}{(N-2)(p-1)-\alpha}}\\
&\le &C(\int_{\mathbb R^N}|W_1|^{2^*})^{\frac{(N+\alpha)-p(N-2)}{(N-2)(p-1)-\alpha}}.
\end{array}
$$
where we have used the interpolation inequality
$$
\int_{\mathbb R^N}(I_\alpha\ast |\eta_\ell W_1|^p)|\eta_\ell W_1|^p\le C\|\eta_\ell W_1\|_{\frac{2Np}{N+\alpha}}^{2p}\le C(\int_{\mathbb R^N}|\eta_\ell W_1|^2)^{\frac{2^-\tilde p}{2^*-2}\frac{N+\alpha}{N}}(\int_{\mathbb R^N}|\eta_\ell W_1|^{2^*})^{\frac{2^*-\tilde p}{2^*-2}\frac{N+\alpha}{N}}.
$$
Since 
$$
\int_{\mathbb R^N}(I_\alpha\ast |\eta_\ell W_1|^p)|\eta_\ell W_1|^p\to \int_{\mathbb R^N}(I_\alpha\ast | W_1|^p)|W_1|^p,
$$ 
as $\ell \to +\infty$, it follows that
$$
\begin{array}{rcl}
\varphi(\rho_\ell)&=&C\left(\frac{\int_{\mathbb R^N}(I_\alpha\ast |\eta_\ell W_1|^p)|\eta_\ell W_1|^p}{\int_{\mathbb R^N}|\eta_\ell W_1|^2}\right)^{\frac{2}{(N-2)(p-1)-\alpha}}\int_{\mathbb R^N}|\eta_\ell W_1|^2\\
&=&C\frac{\left(\int_{\mathbb R^N}(I_\alpha\ast |\eta_\ell W_1|^p\right)|\eta_\ell W_1|^p)^{\frac{2}{(N-2)(p-1)-\alpha}}}{\left(\int_{\mathbb R^N}|\eta_\ell W_1|^2\right)^{\frac{(N+\alpha)-p(N-2)}{(N-2)(p-1)-\alpha}}}\\
&=& \left\{\begin{array}{rcl} C(\ln \ell(1+o(1)))^{-\frac{(N+\alpha)-p(N-2)}{(N-2)(p-1)-\alpha}}  \qquad &{\rm if}&  \   \  N=4,\\
C(\ell(1+o(1)))^{-\frac{(N+\alpha)-p(N-2)}{(N-2)(p-1)-\alpha}} \qquad &{\rm if }& \   \ N=3.
\end{array}\right.
\end{array}
$$
Similarly, we prove that $\psi(\rho)$ is bounded and
$$
|t_\mu-1|\le \left|\left(\frac{\int_{\mathbb R^N}|\nabla(\eta_\ell W_1|^2}{\int_{\mathbb R^N}|\eta_\ell W_1|^{2^*}}\right)^{\frac{1}{2}}-1\right|
+\mu^\sigma \left(\frac{\int_{\mathbb R^N}|\nabla(\eta_\ell W_1|^2}{\int_{\mathbb R^N}|\eta_\ell W_1|^{2^*}}\right)^{\frac{1}{2}}
\frac{C}{\int_{\mathbb R^N}|\eta_\ell W_1|^{2^*}}\to \frac{C\mu^\sigma}{\int_{\mathbb R^N}|W_1|^{2^*}},
$$
as $\ell \to +\infty$. Thus, for small $\mu>0$, we have 
$$
\begin{array}{rcl}
(II)&=&\varphi(\rho_\ell)+(t_\mu^N-1)\varphi(\rho_\ell )+\frac{t_\mu^N(t_\mu^\alpha-1)}{2p}\rho_\ell^{(N+\alpha)-p(N-2)}\int_{\mathbb R^N}(I_\alpha\ast |\eta_\ell W_1|^p)|\eta_\ell W_1|^p\\
&\sim &\left\{\begin{array}{rcl}
(\ln \ell)^{-\frac{(N+\alpha)-p(N-2)}{(N-2)(p-1)-\alpha} }\quad &{\rm if}&  \   N=4,\\
\ell^{-\frac{(N+\alpha)-p(N-2)}{(N-2)(p-1)-\alpha}} \quad &{\rm if}&  \   N=3.
\end{array}\right.
\end{array}
$$
It follows that if $N=4$, then 
$$
\begin{array}{rcl}
m_\mu &\le &(I)-\mu^\sigma (II)\\
&\le &m_0+O(\ell^{-2})-C\mu^\sigma (\ln\ell)^{-\frac{(N+\alpha)-p(N-2)}{(N-2)(p-1)-\alpha} }
\end{array}
\eqno(6.12)
$$
Take $\ell=(1/\mu)^M$. Then 
$$
m_\mu\le m_0+O(\mu^{2M})-C\mu^\sigma M^{-\frac{(N+\alpha)-p(N-2)}{(N-2)(p-1)-\alpha}}(\ln\frac{1}{\mu})^{-\frac{(N+\alpha)-p(N-2)}{(N-2)(p-1)-\alpha} }
$$
If $M>\frac{1}{2p-2-\alpha}$, then $2M>\sigma$, and hence
$$
m_\mu\le m_0-\mu^\sigma(\ln\frac{1}{\mu})^{-\frac{(N+\alpha)-p(N-2)}{(N-2)(p-1)-\alpha} }\varpi=m_0-\mu^{\frac{2}{2p-2-\alpha}}
(\ln\frac{1}{\mu})^{-\frac{4+\alpha -2p}{2p-2-\alpha}}\varpi,
\eqno(6.13)
$$
for small $\mu>0$, where 
$$
\varpi=\frac{1}{2}CM^{-\frac{(N+\alpha)-p(N-2)}{(N-2)(p-1)-\alpha}}.
$$

If  $N=3$, then 
$$
\begin{array}{rcl}
m_\mu &\le &(I)-\mu^\sigma (II)\\
&\le &m_0+O(\ell^{-1})-C\mu^\sigma \ell^{-\frac{(N+\alpha)-p(N-2)}{(N-2)(p-1)-\alpha} }
\end{array}
\eqno(6.14)
$$
Take $\ell=\delta^{-1}\mu^{-\tau}$. Then 
$$
m_\mu\le m_0+\delta O(\mu^\tau)-C\mu^\sigma \delta^{\frac{(N+\alpha)-p(N-2)}{(N-2)(p-1)-\alpha}}\mu^{\tau\frac{(N+\alpha)-p(N-2)}{(N-2)(p-1)-\alpha} }
$$
If 
$$
\tau=\frac{1}{1+p(N-2)-(N+\alpha)},
$$
then
$$
m_\mu\le m_0+(\delta O(1)-C\delta^{\frac{(N+\alpha)-p(N-2)}{(N-2)(p-1)-\alpha}})\mu^{\frac{1}{1+p(N-2)-(N+\alpha)}}.
$$
Since 
$$
1>\frac{(N+\alpha)-p(N-2)}{(N-2)(p-1)-\alpha},
$$
it follows that for some small $\delta>0$, there exists $\varpi>0$ such that 
$$
m_\mu\le m_0-\mu^{\frac{1}{1+p(N-2)-(N+\alpha)}}\varpi=m_0-\mu^{\frac{1}{p-2-\alpha}}\varpi.
$$
This completes the proof. \end{proof}

Combining Lemma 6.3 and Lemma 6.4, we get the following 

{\bf Lemma 6.5.} {\it Let $\delta_\mu:=m_0-m_\mu$, then 
$$
\mu^\sigma\gtrsim \delta_\mu\gtrsim \left\{\begin{array}{rcl} 
 \mu^\sigma,   \quad  \quad \qquad  &{ if}&  \   \  N\ge 5,\\
 \mu^{\sigma}(\ln\frac{1}{\mu})^{-\frac{4+\alpha-2p}{2p-2-\alpha}},   \ &{ if}& \   \ N=4,\\
 \mu^{\frac{1}{p-2-\alpha}},  \qquad  &{ if}& \    \  N=3 \ and \ p\in (2+\alpha, 3+\alpha).
 \end{array}\right.
 $$   }

{\bf Lemma 6.6.}   {\it Let $N\ge 5$ and $p\in (1+\frac{\alpha}{N-2}, \frac{N+\alpha}{N-2})$, then
$$
 \|w_\mu\|_2^2\sim \int_{\mathbb R^N}(I_\alpha\ast |w_\mu|^p)|w_\mu|^p\sim 1.    
$$}
\begin{proof}  Since 
$$
\begin{array}{rcl}
m_0&\le& m_\mu+\mu^\sigma\left(\frac{t_\mu^{N+\alpha}}{2p}\int_{\mathbb R^N}(I_\alpha\ast |w_\mu|^p)|w_\mu|^p-\frac{t_\mu^N}{2}\int_{\mathbb R^N}|w_\mu|^2\right)\\
&=&m_\mu+\mu^\sigma\frac{t_\mu^N}{4p}[2t_\mu^\alpha+p(N-2)-(N+\alpha)]\int_{\mathbb R^N}(I_\alpha\ast |w_\mu|^p)|w_\mu|^p\\
&\le &m_\mu+\frac{\mu^\sigma}{2p}[4+p(N-2)-(N+\alpha)]\int_{\mathbb R^N}(I_\alpha\ast |w_\mu|^p)|w_\mu|^p,
\end{array}
$$
it follows that
$$
\int_{\mathbb R^N}(I_\alpha\ast |w_\mu|^p)|w_\mu|^p\ge \frac{2p(m_0-m_\mu)\mu^{-\sigma}}{4+p(N-2)-(N+\alpha)}\ge C>0.
$$

On the other hand, by the Hardy-Littlewood-Sobolev inequality and the boundedness of $\{w_\mu\}$, we also have  
$$
\int_{\mathbb R^N}(I_\alpha\ast |w_\mu|^p)|w_\mu|^p\le C\|w_\mu\|_{\frac{2N}{N+\alpha}}^{2p}\le \tilde C<\infty,
$$
Therefore, $\int_{\mathbb R^N}(I_\alpha\ast |w_\mu|^p)|w_\mu|^p\sim 1$. Finally, it follows from 
(6.4) that   $\|w_\mu\|_2^2\sim \int_{\mathbb R^N}(I_\alpha\ast |w_\mu|^p)|w_\mu|^p\sim 1$ as $\mu\to 0$. The proof is complete.
\end{proof}

In what follows, we consider that case $\mu\to 0$ and prove the theorem 2.3. To this end, some technique lemmas are needed.

{\bf Lemma 6.7.} {\it Let $N\ge 5$ and $p\in (\max\{2,1+\frac{\alpha}{N-2}\},\frac{N+\alpha}{N-2})$,  then $w_\mu\to W_{\rho_0}$ in $H^1(\mathbb R^N)$ as $\mu\to 0$, where $W_{\rho_0}$ is  a positive ground sate of the equation $-\Delta W=W^{2^*-1}$ with
$$
 \rho_0=\left(\frac{[(N+\alpha)-p(N-2)]\int_{\mathbb R^N}(I_\alpha\ast |W_1|^p)|W_1|^p}{2p\int_{\mathbb R^N}|W_1|^2}\right)^ {\frac{1}{(N-2)(p-1)-\alpha}}.
 $$
 In the lower dimension cases $N=4$ and $N=3$, there exists $\xi_\mu\in (0,+\infty)$ with $\xi_\mu\to 0$ such that 
$$
w_\mu-\xi_\mu^{-\frac{N-2}{2}}W_1(\xi^{-1}_\mu\cdot)\to 0
$$
as $\mu\to 0$ in $D^{1,2}(\mathbb R^N)$ and $L^{2^*}(\mathbb R^N)$. }
\begin{proof} Note that $w_\mu$ is a positive radially symmetric function, and by Lemma 6.2, $\{w_\mu\}$ is bounded in $H^1(\mathbb R^N)$. Then there exists $w_0\in H^1(\mathbb R^N)$ verifying $-\Delta w=w^{2^*-1}$ such that 
$$
w_\mu \rightharpoonup w_0   \quad {\rm weakly \ in} \  H^1(\mathbb R^N), \quad w_\mu\to w_0 \quad {\rm in} \ L^p(\mathbb R^N) \quad {\rm for \ any} \ p\in (2,2^*),
\eqno(6.15)
$$
and 
$$
w_\mu(x)\to w_0(x) \quad a. \ e. \  {\rm on} \ R^N,  \qquad w_\mu\to w_0 \quad {\rm in} \   L^2_{loc}(\mathbb R^N).
\eqno(6.16)
$$
Observe that
$$
J_0(w_\mu)=J_\mu(w_\mu)+\frac{\mu^\sigma}{2p}\int_{\mathbb R^N}(I_\alpha\ast |w_\mu|^p)|w_\mu|^p-\frac{\mu^\sigma}{2}\int_{\mathbb R^N}|w_\mu|^2=m_\mu+o(1)=m_0+o(1),
$$
and 
$$
J'_0(w_\mu)w=J'_\mu(w_\mu)w+\mu^\sigma\int_{\mathbb R^N}(I_\alpha\ast |w_\mu|^p)|w_\mu|^{p-2}w_\mu w-\mu^\sigma\int_{\mathbb R^N}w_\mu w=o(1).
$$
Therefore, $\{w_\mu\}$ is a PS sequence of $J_0$.

By Lemma 3.5,  it is standard to show that there exists $\zeta^{(j)}_\mu\in (0,+\infty)$, $w^{(j)}\in D^{1,2}(\mathbb R^N)$ with $j=1,2,\cdots, k$, $k$  a non-negative integer, such that
$$
w_\mu=w_0+\sum_{j=1}^k(\zeta^{(j)}_\mu)^{-\frac{N-2}{2}}w^{(j)}((\zeta^{(j)}_\mu)^{-1} x)+\tilde w_\mu,
\eqno(6.17)
$$
where $\tilde w_\mu \to 0$ in $L^{2^*}(\mathbb R^N)$ and $w^{(j)}$ are nontrivial solutions of the limit equation $-\Delta w=w^{2^*-1}$. Moreover, we have
$$
\limsup_{\mu\to 0}\|w_\mu\|^2_{D^1(\mathbb R^N)}\ge  \|w_0\|^2_{D^1(\mathbb R^N)}+\sum_{j=1}^k\|w^{(j)}\|^2_{D^1(\mathbb R^N)}
\eqno(6.18)
$$
and 
$$
m_0=J_0(w_0)+\sum_{j=1}^kJ_0(w^{(j)}).
\eqno(6.19)
$$
Moreover, $J_0(w_0)\ge 0$ and $J_0(w^{(j)})\ge m_0$ for all $j=1,2,\cdots, k.$

If $N\ge 5$, then by Lemma 6.6, we have $w_0\not=0$ and hence $J_0(w_0)=m_0$ and  $k=0$.  Thus $w_\mu\to w_0$ in $L^{2^*}(\mathbb R^N)$. Since $J_0'(w_\mu)\to 0$, it follows that
$w_\mu\to w_0$ in $D^{1,2}(\mathbb R^N)$.

By  Lemma 3.7, we have 
$$
w_\mu(x)\le C_N|x|^{-\frac{N-1}{2}}\|w_\mu\|_{H^1(\mathbb R^N)} \qquad {\rm for} \  |x|\ge \alpha_N.
$$
Since
$$
\left(-\Delta-C|x|^{-\frac{2(N-1)}{N-2}}\right)w_\mu\le \left(-\Delta +\mu^\sigma-w_\mu^{2^*-2}-\mu^\sigma (I_\alpha\ast |w_\mu|^p)w_\mu^{p-2}\right)w_\mu=0,
$$
for some constant $C>0$.  We also have
$$
\left(-\Delta-C|x|^{-\frac{2(N-1)}{N-2}}\right)\frac{1}{|x|^{N-2-\varepsilon_0}}=\left(\varepsilon_0(N-2-\varepsilon_0)-C|x|^{-\frac{2}{N-2}}\right)\frac{1}{|x|^{N-\varepsilon_0}},
$$
which is positive for $|x|$ large enough. By the maximum principle on $\mathbb R^N\setminus B_R$, we deduce that 
$$
w_\mu(x)\le \frac{w_\mu(R)R^{N-2-\varepsilon_0}}{|x|^{N-2-\varepsilon_0}}\le \frac{CR^{N/2-1-\varepsilon_0}}{|x|^{N-2-\varepsilon_0}},  \qquad {\rm for} \  \ |x|\ge R.
\eqno(6.20)
$$
When $\varepsilon_0>0$ is small enough, the domination is in $L^2(\mathbb R^N)$ for $N\ge 5$, and this shows, by the dominated convergence theorem, that 
$w_\mu\to w_0$ in $L^2(\mathbb R^N)$.  Thus, we conclude that $w_\mu\to w_0$ in  $H^1(\mathbb R^N)$. Moreover, by (6.4), we obtain
$$
\|w_0\|_2^2=\frac{(N+\alpha)-p(N-2)}{2p}\int_{\mathbb R^N}(I_\alpha\ast |w_0|^p)|w_0|^p.
$$
from which it follows that $w_0=W_{\rho_0}$ with  
$$
 \rho_0=\left(\frac{[(N+\alpha)-p(N-2)]\int_{\mathbb R^N}(I_\alpha\ast |W_1|^p)|W_1|^p}{2p\int_{\mathbb R^N}|W_1|^2}\right)^ {\frac{1}{(N-2)(p-1)-\alpha}}.
 $$

If $N=4$ or $3$. By Fatou's lemma, we have $\|w_0\|^2_2\le \liminf_{\mu\to 0}\|w_\mu\|_2^2<\infty$, therefore,  $w_0=0$ and hence $k=1$. Thus, we obtain  
$J_0(w^{(1)})=m_0$ and hence $w^{(1)}=W_\rho$ for some $\rho\in (0,+\infty)$. Therefore, we conclude that 
$$w_\mu-\xi_\mu^{-\frac{N-2}{2}}W_1(\xi_\mu^{-1}\cdot )\to 0
$$ in $L^{2^*}(\mathbb R^N)$ as $\mu\to 0$, where 
$\xi_\mu:=\rho\zeta_\mu^{(1)}\in (0,+\infty)$ satisfying $\xi_\mu\to 0$ as $\mu\to 0$.
Since $$J_0'(w_\mu-\xi_\mu^{-\frac{N-2}{2}}W_1(\xi_\mu^{-1}\cdot ))=J'_0(w_\mu)+J'_0(W_1)+o(1)=o(1)$$ as $\mu\to 0$, it follows that $w_\mu-\xi_\mu^{-\frac{N-2}{2}}W_1(\xi_\mu^{-1}\cdot )\to 0$ in $D^{1,2}(\mathbb R^N)$. 
\end{proof}

{\bf Lemma 6.8.} {\it Let $N\ge 5$ and $p\in (1+\frac{\alpha}{N-2},\frac{N+\alpha}{N-2})$,  then   there exists a $\zeta_\mu\in (0,\infty)$ verifying 
$$
\zeta_\mu\sim  \mu^{\frac{1}{(N-2)(p-1)-\alpha}}
$$
such that the rescaled ground states 
$$
w_\mu(x)=\zeta_\mu^{\frac{N-2}{2}}v_\mu(\zeta_\mu x)
$$
converges to  $W_{\rho_0}$ in $H^1(\mathbb R^N)$ as $\mu\to 0$, where $W_{\rho_0}$ is  given in Lemma 6.7.}

\begin{proof}  The proof is similar to that of Lemma 5.8 and is omitted. \end{proof}

In the lower dimension cases $N=4$ and $N=3$, we further perform a  scaling
$$
\tilde w(x)=\xi_\mu^{\frac{N-2}{2}} w(\xi_\mu x), 
\eqno(6.21)
$$
where $\xi_\mu\in (0,+\infty)$ is given in Lemma 5.7.  Then the rescaled equation is as follows
$$
-\Delta\tilde w+\mu^\sigma\xi_\mu^2\tilde w=\mu^\sigma \xi_\mu^{(N+\alpha)-p(N-2)}(I_\alpha \ast |\tilde w|^p) \tilde w^{p-1}+\tilde w^{2^*-1}.
\eqno(R_\mu)
$$
The  corresponding energy functional is given by
$$
\begin{array}{rcl}
\tilde J_\mu(\tilde w):&=&\frac{1}{2}\int_{\mathbb R^N}|\nabla \tilde w|^2+\mu^\sigma\xi_\mu^2|\tilde w|^2-\frac{1}{2p}\mu^\sigma \xi_\mu^{(N+\alpha)-p(N-2)}\int_{\mathbb R^N}(I_\alpha\ast |\tilde w|^p)|\tilde w|^p\\
&\quad &-\frac{1}{2^*}\int_{\mathbb R^N}|\tilde w|^{2^*}.
\end{array}
\eqno(6.22)
$$
Clearly, we have $\tilde J_\mu(\tilde w)=J_\mu(w)=I_\mu(v)$. 

Furthermore, we  have the following lemma.

{\bf Lemma 6.9.}  {\it  Let $v,w,\tilde w\in H^1(\mathbb R^N)$ satisfy  (6.1) and (6.21), then the following statements hold true

(1) \ $ \ \|\nabla \tilde w\|_2^2= \|\nabla w\|_{2}^{2}=\|\nabla v\|_{2}^{2}, \  \|\tilde w\|^{2^*}_{2^*}=\|w\|_{2^*}^{2^*}=\|v\|_{2^*}^{2^*},$

(2) \ $\xi_\mu^2\|\tilde w\|^2_2=\|w\|_2^2=\mu^{-\sigma}\| v\|_2^2,$

(3)  \ $\xi_\mu^{(N+\alpha)-p(N-2)}\int_{\mathbb R^N}(I_\alpha\ast |\tilde w|^p)|\tilde w|^p=\int_{\mathbb R^N}(I_\alpha\ast |w|^p)|w|^p=\mu^{1-\sigma} \int_{\mathbb R^N}(I_\alpha\ast |v|^p)|v|^p$.
}

In what follows, we set  
$$
\tilde w_\mu(x)=\xi_\mu^{\frac{N-2}{2}} w_\mu(\xi_\mu x),
$$
 then by Lemma 6.7, we have 
$$
\|\nabla(\tilde w_\mu-W_1)\|_2\to 0, \qquad \|\tilde w_\mu-W_1\|_{2^*}\to 0,   \qquad {\rm as} \  \ \mu\to 0.
\eqno(6.23)
$$

Note that the corresponding Nehari and Poho\v zaev's identities are as follows
$$
\int_{\mathbb R^N}|\nabla \tilde w_\mu|^2+\mu^\sigma\xi_\mu^2\int_{\mathbb R^N}|\tilde w_\mu|^2=\mu^\sigma \xi_\mu^{(N+\alpha)-p(N-2)}\int_{\mathbb R^N}(I_\alpha\ast |\tilde w_\mu|^p)|\tilde w_\mu|^p+\int_{\mathbb R^N}|\tilde w_\mu|^{2^*}
$$
and 
$$
\frac{1}{2^*}\int_{\mathbb R^N}|\nabla \tilde w_\mu|^2+\frac{1}{2}
\mu^\sigma\xi_\mu^2\int_{\mathbb R^N}|\tilde w_\mu|^2=\frac{N+\alpha}{2Np}\mu^\sigma\xi_\mu^{(N+\alpha)-p(N-2)}\int_{\mathbb R^N}(I_\alpha\ast |\tilde w_\mu|^p)|\tilde w_\mu|^p+\frac{1}{2^*}\int_{\mathbb R^N}|\tilde w_\mu|^{2^*},
$$
it follows that 
$$
\left(\frac{1}{2}-\frac{1}{2^*}\right)\mu^\sigma\xi_\mu^2\int_{\mathbb R^N}|\tilde w_\mu|^2=\left(\frac{N+\alpha}{2Np}-\frac{1}{2^*}\right)\mu^\sigma\xi_\mu^{(N+\alpha)-p(N-2)}\int_{\mathbb R^N}(I_\alpha\ast |\tilde w_\mu|^p)|\tilde w_\mu|^p.
$$
Thus, we obtain
$$
\xi_\mu^{(N-2)(p-1)-\alpha}\int_{\mathbb R^N}|\tilde w_\mu|^2=\frac{2(2^*-\tilde p)}{\tilde p(2^*-2)}\int_{\mathbb R^N}(I_\alpha\ast |\tilde w_\mu|^p)|\tilde w_\mu|^p,
\eqno(6.24)
$$
where $\tilde p=\frac{2Np}{N+\alpha}\in (2,2^*)$.

To control the norm $\|\tilde w_\mu\|_2$, we note that  for any $\mu>0$, $\tilde w_\mu>0$ satisfies the linear inequality
$$
-\Delta \tilde w_\mu+\mu^\sigma\xi_\mu^2\tilde w_\mu=\mu^\sigma \xi_\mu^{(N+\alpha)-p(N-2)}(I_\alpha\ast  |\tilde w_\mu|^p)|\tilde w_\mu|^{p-1}+\tilde w_\mu^{2^*-1}>0,  \qquad x\in \mathbb R^N.
\eqno(6.25)
$$

{\bf Lemma 6.10.}  {\it There exists a constant $c>0$ such that 
$$
\tilde w_\mu(x)\ge c|x|^{-(N-2)}\exp(-\mu^{\frac{\sigma}{2}}\xi_\mu |x|),  \quad |x|\ge 1.
\eqno(6.26)
$$}
The proof of the above lemma  is similar to that  in \cite[Lemma 4.8]{Moroz-1}.  As consequences, we have the following two lemmas.

{\bf Lemma 6.11.}  {\it If $N=3$, then $\|\tilde w_\mu\|_2^2\gtrsim \mu^{-\frac{\sigma}{2}}\xi_\mu^{-1}.$}

{\bf Lemma 6.12.}  {\it If $N=4$, then $\|\tilde w_\mu\|_2^2\gtrsim  - \ln(\mu^{\sigma}\xi_\mu^2).$}

We remark that $\tilde w_\mu$ is only defined for $N=4$ and $N=3$. But the following discussion also apply to the case $N\ge 5$. To prove our main result, the key point is to show the boundedness of $\int_{\mathbb R^N}(I_\alpha\ast |\tilde w_\mu|^p)|\tilde w_\mu|^p$. 
We will prove this by using the well-known Hardy-Littlewood-Sobolev inequality and the best uniform decay  estimation of $\tilde w_\mu$. 
To this end, we first give several technical lemmas.  The first one concerns the decay estimate of $\tilde w_\mu$ at infinity.

{\bf Lemma 6.13.}  {\it Assume $N\ge 3$ and $2<p<\frac{N+\alpha}{N-2}$. Then there exist constants $L_0>0$ and $C_0>0$ such that for any small $\mu>0$ and $|x|\ge L_0\mu^{-\sigma/2}\xi_\mu^{-1}$, 
$$
\tilde w_\mu(x)\le C_0\mu^{\sigma(N-2)/4}\xi_\mu^{(N-2)/2}\exp(-\frac{1}{2}\mu^{\sigma/2}\xi_\mu |x|).
\eqno(6.27)
$$}
\begin{proof}  Note that 
$$
|w_\mu(x)|^2\le \frac{1}{|B_{|x|}|}\int_{B_{|x|}}|w_\mu|^2\le C|x|^{-N},
$$ 
 and hence 
$$
| w_\mu(x)|^p|x|^N\le C|x|^{-\frac{N}{2}(p-2)}\to 0,  \  \ {\rm as} \   \mu\to 0,
$$ 
it follows from Lemma 3.10 that
$$
(I_\alpha\ast |w_\mu|^p)(x)\le C|x|^{-(N-\alpha)},
$$
for all large $|x|$. If $|x|\ge L_0\mu^{-\sigma/2}\xi_\mu^{-1}$, then for large $L_0>0$, 
$$
|\tilde w_\mu|^{2^*-2}(x)=\xi_\mu^{\frac{N-2}{2}(2^*-2)}w_\mu^{2^*-2}(\xi_\mu x)\le CL_0^{-\frac{N(2^*-2)}{2}}\mu^{\frac{\sigma N(2^*-2)}{4}}\xi_\mu^2\le \frac{1}{4}\mu^{\frac{\sigma N}{N-2}}\xi_\mu^2\le \frac{1}{4}\mu^{\sigma}\xi_\mu^2
$$
and 
$$
\begin{array}{rcl}
\mu^\sigma\xi_\mu^{(N+\alpha)-p(N-2)}(I_\alpha\ast |\tilde w_\mu|^p)(x)\tilde w_\mu^{p-2}(x)
&=&\mu^{\sigma}\xi_\mu^2(I_\alpha\ast |w_\mu|^p)(\xi_\mu x)w_\mu^{p-2}(\xi_\mu x)\\
&\le& C\mu^\sigma\xi_\mu^2|\xi_\mu x|^{-N+\alpha-\frac{N}{2}(p-2)}\\
&\le& C\mu^\sigma\xi_\mu^2\cdot L_0^{-[N-\alpha+\frac{N}{2}(p-2)]}\mu^{\frac{\sigma}{2}[N-\alpha+\frac{N}{2}(p-2)]}\\
&\le& \frac{1}{4}\mu^\sigma\xi_\mu^2.
\end{array}
$$
Therefore, we obtain 
$$
-\Delta\tilde w_\mu(x)+\frac{1}{2}\mu^{\sigma}\xi_\mu^2\tilde w_\mu(x)\le 0,\quad |x|\ge L_0\mu^{-\sigma/2}\xi_\mu^{-1}.
$$ 
The rest of the proof is similar to that of Lemma  5.13 and is omitted.  
\end{proof}

To proceed further, we need the following fact \cite[Lemma 3.5]{Li-1}.

{\bf Lemma 6.14.}  {\it For any $N \ge  3$, $q \in [2, 2^*]$ and $u \in  H^1(\R^N)$,  there exists
a positive constant $C_0$ independent of $q$ and $u$ such that 
 $$
 \| u\|_q \le C_0\|u\|_{H^1(\R^N)}.
\eqno(6.28)
 $$}

In the following three lemmas, we establish the decay estimate of $(I_\alpha\ast |\tilde w_\mu|^p)(x)$ when $(\alpha,p)$ belongs to some regions in the $(\alpha,p)$ plane.

{\bf Lemma 6.15.}  {\it Assume that 
$$
\frac{N+\alpha}{N-2}>p> \frac{2(N^2-2\alpha)}{N(N-2)}=\left\{\begin{array}{rcl}
4-\frac{1}{2}\alpha, \quad  &if & \ N=4,\\
6-\frac{4}{3}\alpha, \quad &if & \ N=3.
\end{array}
\right.
\eqno(6.29)
$$
 Then for any small $\epsilon>0$, there exists a constants $C=C(\epsilon)>0$ such that
$$
(I_\alpha\ast |\tilde w_\mu|^p)(x)\le C\xi_\mu^{-(\frac{2N-(N-2)p}{2}+\epsilon)}|x|^{-(N-\alpha)},  \quad |x|\ge 2.
\eqno(6.30)
$$}
\begin{proof}   Since $\frac{N+\alpha}{N-2}>p>\frac{2(N^2-2\alpha)}{N(N-2)}>\frac{2(N-\alpha)}{N-2}$, it follows that
$$
\frac{2(N-\alpha)}{(N-2)p}<1.
$$
Let $\theta_0=1-\frac{2(N-\alpha)}{(N-2)p}\in (0,1)$, then 
$$
\frac{\theta_0Np}{\alpha}=\frac{Np}{\alpha}\left(1-\frac{2(N-\alpha)}{(N-2)p}\right)\le \frac{2N}{N-2}.
$$
Since $p>\frac{2(N^2-2\alpha)}{N(N-2)}$, we also have
$$
\frac{\theta_0Np}{\alpha}=\frac{Np}{\alpha}\left(1-\frac{2(N-\alpha)}{(N-2)p}\right)>2.
$$
Therefore,  there exists $\tau_j>0$ satisfying
$$
\max\{2, \theta_0p\}<\tau_j<\frac{\theta_0Np}{\alpha}\le \frac{2N}{N-2}
$$
and 
$$
\lim_{j\to \infty}\tau_j=\frac{\theta_0Np}{\alpha}=\frac{Np}{\alpha}\left(1-\frac{2(N-\alpha)}{(N-2)p}\right).
$$
Since $\tau_j>\theta_0p$ and $\theta_0>\frac{\alpha \tau_j}{Np}$, it follows that
$$
m_j:=\frac{\tau_j(N-\alpha)}{2(\tau_j-\theta_0 p)}>\frac{N}{2}.
\eqno(6.33)
$$
By the definition of $\theta_0$, we see that
$$
N-\alpha=\frac{N-2}{2}(1-\theta_0)p,
$$
which yields 
$$
\frac{\tau_j(N-\alpha)}{2(\tau_j-\theta_0p)}=m'_j:=\frac{(N-2)(1-\theta_0)p\tau_j}{4(\tau_j-\theta_0p)}.
\eqno(6.34)
$$
Since $p> \frac{2(N^2-2\alpha)}{N(N-2)}>\frac{2(N-\alpha)}{N-2}$, for $|x|\ge 2$, we get
$$
\int_{B_{1}(x)}\frac{|\tilde w_\mu|^p(y)}{|x-y|^{N-\alpha}}dy\lesssim \tilde w_\mu^p(\frac{1}{2}|x|)\int_{B_{1}(0)}\frac{1}{|y|^{N-\alpha}}dy\lesssim 
\frac{1}{|x|^{\frac{N-2}{2}p}}\le \frac{1}{|x|^{N-\alpha}}.
\eqno(6.33)
$$
By the H\"older inequality, for $|x|\ge 2$,  we get
$$
\int_{B_{1}(0)}\frac{|\tilde w_\mu|^p(y)}{|x-y|^{N-\alpha}}dy\le \left(\int_{\mathbb R^N}|\tilde w_\mu|^{2^*}\right)^{\frac{p}{2^*}}\left(\int_{B_{1}(0)}\frac{1}{|x-y|^{\frac{2^*(N-\alpha)}{2^*-p}}}dy\right)^{\frac{2^*-p}{2^*}}\lesssim \frac{1}{|x|^{N-\alpha}}.
\eqno(6.34)
$$
By the H\"older inequality, Lemma 3.11 and the fact that $\lim_{j\to\infty}\frac{\tau_j-\theta_0p}{\tau_j}=\frac{N-\alpha}{N}$, for $|x|\ge 2$, we get
$$
\begin{array}{rl}
&\int_{\mathbb R^N\setminus B_1(0)}\frac{\tilde w_\mu^p(y)}{(1+|x-y|^2)^{\frac{N-\alpha}{2}}}dy=\int_{\mathbb R^N\setminus B_1(0)}\frac{\tilde w_\mu^{\theta_0 p}(y)\tilde w_\mu^{(1-\theta_0)p}(y)}{(1+|x-y|^2)^{\frac{N-\alpha}{2}}}dy\\
&\lesssim    (\int_{\mathbb R^N}\tilde w_\mu^{\tau_j}(y)dy)^{\frac{\theta_0 p}{\tau_j}}\left(\int_{\mathbb R^N\setminus B_1(0)}\frac{(\tilde w_\mu(y))^{(1-\theta_0)p\tau_j/(\tau_j-\theta_0 p)}}
{(1+|x-y|^2)^{\frac{1}{2}\tau_j(N-\alpha)/(\tau_j-\theta_0 p)}}dy\right)^{\frac{\tau_j-\theta_0 p}{\tau_j}}\\
&\lesssim    (\int_{\mathbb R^N}\tilde w_\mu^{\tau_j}(y)dy)^{\frac{\theta_0 p}{\tau_j}}\left(\int_{\mathbb R^N}\frac{1}
{(1+|x-y|^2)^{m_j}(1+|y|^2)^{m'_j}}dy\right)^{\frac{\tau_j-\theta_0 p}{\tau_j}}\\
&\le \frac{1}{|x|^{N-\alpha}} C_j(\int_{\mathbb R^N}\tilde w_\mu^{\tau_j}(y)dy)^{\frac{\theta_0 p}{\tau_j}}.
\end{array}
\eqno(6.35)
$$
Therefore, by (6.33)-(6.35), Lemma 6.14 and the boundedness of  $w_\mu$ in $H^1(\mathbb R^N)$ , it follows that there exists a constant $C_j>0$ dependent of $j$ such that for any $|x|\ge 2$, 
$$
\begin{array}{rcl}
(I_\alpha\ast |\tilde w_\mu|^p)(x)&\lesssim& \int_{B_{1}(x)}\frac{|\tilde w_\mu|^p(y)}{|x-y|^{N-\alpha}}dy+\int_{B_{1}(0)}\frac{|\tilde w_\mu|^p(y)}{|x-y|^{N-\alpha}}dy
+\int_{\mathbb R^N\setminus  B_{1}(0)}\frac{|\tilde w_\mu|^p(y)}{(1+|x-y|^2)^{\frac{N-\alpha}{2}}}dy\\
&\le& \frac{C}{|x|^{N-\alpha}}+\frac{C_j}{|x|^{N-\alpha}}\xi_\mu^{(\frac{N-2}{2}-\frac{N}{\tau_j})\theta_0p}\left(\int_{\mathbb R^N}|w_\mu|^{\tau_j}\right)^{\frac{\theta_0p}{\tau_j}}\\
&\le &\frac{\tilde C_j}{|x|^{N-\alpha}}\xi_\mu^{(\frac{N-2}{2}-\frac{N}{\tau_j})\theta_0p}.
\end{array}
$$
Since $\lim_{j\to\infty}(\frac{N-2}{2}-\frac{N}{\tau_j})\theta_0p=-\frac{2N-p(N-2)}{2}$, for any $\epsilon>0$, there exists $j_0$ sufficiently large such that
$$
(I_\alpha\ast |\tilde w_\mu|^p)(x)\le \frac{\tilde C_{j_0}}{|x|^{N-\alpha}}\xi_\mu^{(\frac{N-2}{2}-\frac{N}{\tau_{j_0}})\theta_0p}\le C\xi_\mu^{-(\frac{2N-p(N-2)}{2}+\epsilon)}|x|^{-(N-\alpha)}.
$$
The proof is complete.
\end{proof}

{\bf Lemma 6.16.}  {\it Assume that 
$$
2<p\le \frac{2(N^2-2\alpha)}{N(N-2)}
=\left\{\begin{array}{rcl}
4-\frac{1}{2}\alpha, \quad  &if & \ N=4,\\
6-\frac{4}{3}\alpha, \quad &if & \ N=3.
\end{array}
\right.
\eqno(6.36)
$$
Then for any small $\epsilon>0$, there exists a constants $C=C(\epsilon)>0$ such that
 $$
(I_\alpha\ast |\tilde w_\mu|^p)(x)\le C\xi_\mu^{-(\frac{2\alpha}{N}+\epsilon)}|x|^{-\frac{1}{2}(N-2)(p-\frac{2\alpha}{N})}, \quad |x|\ge 2.
\eqno(6.37)
$$}
\begin{proof}   
By virtue of Lemma 3.10 and the boundedness of $w_\mu$ in $H^1(\mathbb R^N)$,  for all $|x|\ge \xi_\mu^{-1}$, we have 
$$
\begin{array}{rcl}
(I_\alpha \ast |\tilde w_\mu|^p)(x)&=&\xi_\mu^{\frac{N-2}{2}p-\alpha}(I_\alpha\ast |w_\mu|^p)(\xi_\mu z)\\
&\le& C\xi_\mu^{\frac{N-2}{2}p-\alpha}|\xi_\mu x|^{-(N-\alpha)}\\
&=&C\xi_\mu^{-\frac{2\alpha}{N}}\xi_\mu^{\frac{N-2}{2}p-N+\frac{2\alpha}{N}}|x|^{-(N-\alpha)}\\
&\le &C\xi_\mu^{-\frac{2\alpha}{N}}|x|^{-\frac{1}{2}(N-2)(p-\frac{2\alpha}{N})}.
\end{array}
\eqno(6.38)
$$  Without loss of generality, in what follows, we assume $|x|\le \xi_\mu^{-1}$.

By Lemma 3.8,  for $|x|\ge 2$,  we get
$$
\int_{B_{1}(x)}\frac{|\tilde w_\mu|^p(y)}{|x-y|^{N-\alpha}}dy\lesssim \tilde w_\mu^p(\frac{1}{2}|x|)\int_{B_{1}(0)}\frac{1}{|y|^{N-\alpha}}dy\lesssim 
\frac{1}{|x|^{\frac{N-2}{2}p}}\le |x|^{-\frac{1}{2}(N-2)(p-\frac{2\alpha}{N})}.
\eqno(6.39)
$$
Since $p\le \frac{2(N^2-2\alpha)}{N(N-2)}$, by the H\"older inequality and the boundedness of $\tilde w_\mu$ in $D^{1,2}(\mathbb R^N)$, for $|x|\ge 2$,  we get  
$$
\begin{array}{rcl}
\int_{B_{1}(0)}\frac{|\tilde w_\mu|^p(y)}{|x-y|^{N-\alpha}}dy&\le& \left(\int_{\mathbb R^N}|\tilde w_\mu|^{2^*}\right)^{{p/2^*}}\left(\int_{B_{1}(0)}\frac{1}{|x-y|^{{2^*(N-\alpha)}/{(2^*-p)}}}dy\right)^{(2^*-p)/2^*}\\
&\lesssim& \frac{1}{|x|^{N-\alpha}}\lesssim |x|^{-\frac{1}{2}(N-2)(p-\frac{2\alpha}{N})}.
\end{array}
\eqno(6.40)
$$
 For any $\theta_0\in (\frac{2\alpha}{Np}, \frac{2}{p})$, we have
$$
\theta_0p<2<\frac{\theta_0Np}{\alpha}.
$$
As before,  we can easily check that 
$$
m':=\frac{N-\alpha}{2-\theta_0p}>\frac{N}{2}.
\eqno(6.41)
$$ 
Since $p\le \frac{2(N^2-2\alpha)}{N(N-2)}$, it follows that
$$
m'=\frac{N-\alpha}{2-\theta_0p}\ge m:=\frac{(N-2)(1-\theta_0)p}{2(2-\theta_0p)}.
\eqno(6.42)
$$
Again, by the H\"older inequality and the boundedness of $w_\mu$ in $H^1(\mathbb R^N)$, for $|x|\ge 2$,  we get
$$
\begin{array}{rl}
&\int_{\mathbb R^N\setminus B_1(0)}\frac{\tilde w_\mu^p(y)}{(1+|x-y|^2)^{\frac{N-\alpha}{2}}}dy=\int_{\mathbb R^N\setminus B_1(0)}\frac{\tilde w_\mu^{\theta_0p}(y)\tilde w_\mu^{(1-\theta_0)p}(y)}{(1+|x-y|^2)^{\frac{N-\alpha}{2}}}dy\\
&\lesssim    (\int_{\mathbb R^N}\tilde w_\mu^{2}(y)dy)^{\frac{\theta_0 p}{2}}\left(\int_{\mathbb R^N\setminus B_1(0)}\frac{(\tilde w_\mu(y))^{2(1-\theta_0)p/(2-\theta_0 p)}}
{(1+|x-y|^2)^{(N-\alpha)/(2-\theta_0 p)}}dy\right)^{\frac{2-\theta_0 p}{2}}\\
&\lesssim    \xi_\mu^{-\theta_0p}\left(\int_{\mathbb R^N}\frac{1}
{(1+|x-y|^2)^{m'}(1+|y|^2)^{m}}dy\right)^{\frac{2-\theta_0p}{2}}.
\end{array}
$$
Finally, by  (6.41), (6.42), (6.43) and Lemma 3.11, 
 there exists a constants $C>0$ such that for any $|x|\ge 2$, 
 $$
\int_{\mathbb R^N\setminus B_1(0)}\frac{\tilde w_\mu^p(y)}{(1+|x-y|^2)^{\frac{N-\alpha}{2}}}dy\le C\xi_\mu^{-\theta_0p}|x|^{-\frac{1}{2}(N-2)(1-\theta_0)p}.
$$
For any small $\epsilon>0$, choose $\theta_0\in (\frac{2\alpha}{Np}, \frac{2}{p})$ so that 
$$
\theta_0p\le \frac{2\alpha}{N}+\frac{\epsilon}{2}, \quad \frac{1}{2}(N-2)(\theta_0p-\frac{2\alpha}{N})\le \frac{\epsilon}{2}.
$$
Recall that $|x|\le \xi_\mu^{-1}$, it follows that
$$
\begin{array}{lcl}
\int_{\mathbb R^N\setminus B_1(0)}\frac{\tilde w_\mu^p(y)}{(1+|x-y|^2)^{\frac{N-\alpha}{2}}}dy&\le& C\xi_\mu^{-(\frac{2\alpha}{N}+\frac{\epsilon}{2})}|x|^{-\frac{1}{2}(N-2)(p-\frac{2\alpha}{N})}|x|^{\frac{1}{2}(N-2)(\theta_0p-\frac{2\alpha}{N})}\\
&\le&C\xi_\mu^{-(\frac{2\alpha}{N}+\frac{\epsilon}{2})}|x|^{-\frac{1}{2}(N-2)(p-\frac{2\alpha}{N})}\xi_\mu^{-\frac{1}{2}(N-2)(\theta_0p-\frac{2\alpha}{N})}\\
&\le &C\xi_\mu^{-(\frac{2\alpha}{N}+\epsilon)}|x|^{-\frac{1}{2}(N-2)(p-\frac{2\alpha}{N})} .
\end{array}
\eqno(6.43)
$$
Finally, by (6.38), (6.39), (6.40) and (6.43), we conclude the proof of Lemma 6.16. \end{proof}

We are in the position to state our result concerning the best decay estimation of $\tilde w_\mu$ when $(\alpha,p)$ belongs to some regions in the $(\alpha,p)$ plane. 

{\bf Proposition  6.17.}  {\it Assume $N\ge 3$, $\alpha>N-4$ and 
$$
2<p<\frac{N+\alpha}{N-2}.
\eqno(6.46)
$$
Then there exists a constant $C>0$ such that for small $\mu>0$, there holds 
$$
\tilde w_\mu(x)\le C(1+|x|)^{-(N-2)}, \qquad x\in \mathbb R^N.
$$}

{\bf Remark 6.1.} 
In our proof of Theorem 2.3, we only need  Proposition 6.17 in lower dimension cases $N=3, 4$. So, we only give the proof of Proposition 6.17 for $N=3$ and $N=4$. The proof for  $N\ge 5$  is similar to the case $N=4$ and is omitted. 

\vskip 5mm

We consider the Kelvin transform of $\tilde w_\mu$ given as follows
$$
K[\tilde w_\mu](x):=|x|^{-(N-2)}\tilde w_\mu\left(\frac{x}{|x|^2}\right). 
$$
It is easy to see that $\|K[\tilde w_\mu]\|_{L^\infty(B_1)}\lesssim 1$ implies that 
$$
\sup_{|x|\ge 1} \tilde w_\mu(x)\lesssim |x|^{-(N-2)}.
$$
Thus, to prove Proposition 6.17, it suffices to show that there exists $\mu_0>0$ such that 
$$
\sup_{\mu\in (0,\mu_0)}\|K[\tilde w_\mu]\|_{L^{\infty}(B_1)}<\infty.
\eqno(6.47)
$$

It is easy to verify that $K[\tilde w_\mu]$ satisfies 
$$
-\Delta K[\tilde w_\mu]+\frac{\mu^\sigma\xi_\mu^{2}}{|x|^4}K[\tilde w_\mu]=\frac{\mu^\sigma\xi_\mu^{\eta}}{|x|^{4}}(I_\alpha\ast |\tilde w_\mu|^{p})(\frac{x}{|x|^2})\tilde w_\mu^{p-2}(\frac{x}{|x|^2})K[\tilde w_\mu]+K[\tilde w_\mu]^{2^*-1},
\eqno(6.48)
$$
here and in what follows, we set
$$
\eta:=(N+\alpha)-p(N-2)>0.
\eqno(6.49)
$$
 
We also see from Lemma 6.13 that if $|x|\le \mu^{\sigma/2}\xi_\mu/L_0$, then 
$$
K[\tilde w_\mu](x)\lesssim \frac{1}{|x|^{N-2}}\mu^{\frac{\sigma(N-2)}{4}}\xi_\mu^{\frac{N-2}{2}}e^{-\frac{1}{2}\mu^{\sigma/2}\xi_\mu |x|^{-1}}.
\eqno(6.50)
$$

   Let 
$$
a(x)=\frac{\mu^\sigma\xi_\mu^2}{|x|^4}, 
$$
and 
$$
\begin{array}{lcl}
 b(x)&=&\frac{\mu^\sigma\xi_\mu^\eta}{|x|^{4}}(I_\alpha\ast |\tilde w_\mu|^p)(\frac{x}{|x|^2})\tilde w_\mu^{p-2}(\frac{x}{|x|^2})
+K[\tilde w_\mu]^{2^*-2}(x)\\
&=&\frac{\mu^\sigma\xi_\mu^\eta}{|x|^{\gamma}}(I_\alpha\ast |\tilde w_\mu|^p)(\frac{x}{|x|^2})K[\tilde w_\mu]^{p-2}(x)
+K[\tilde w_\mu]^{2^*-2}(x),
\end{array}
$$
where and in what follows, we set
$$
\gamma:=2N-p(N-2).
\eqno(6.51)
$$
Then (6.48) reads as 
$$
-\Delta K[\tilde w_\mu]+a(x)K[\tilde w_\mu]=b(x)K[\tilde w_\mu].
$$
We shall apply the Moser iteration to prove  (6.47).
We note that it follows from (6.50) that for any $v\in H^1_0(B_4)$, 
$$
\int_{B_4}\frac{\mu^\sigma\xi_\mu^2}{|x|^4}K[\tilde w_\mu](x)|v(x)|dx<\infty.
\eqno(6.52)
$$
Since $\tilde w_\mu\to W_1$  in $L^{2^*}(\mathbb R^N)$   as $\mu\to 0$, we have 
$$
\lim_{\mu\to 0}\|K[\tilde w_\mu]-K[W_1]\|_{L^{2^*}(\mathbb R^N)}=0,
$$
which together with fact that for any $s>1$, the Lebesgue space $L^{s}(\mathbb R^N)$ has the Kadets-Klee property,  implies that  
$$
\lim_{\mu\to 0}\|K[\tilde w_\mu]^{\frac{4}{N-2}}-K[W_1]^{\frac{4}{N-2}}\|_{L^{\frac{N}{2}}(\mathbb R^N)}=0.
$$

{\bf Lemma 6.18.}  {\it Assume $N\ge 3$ and $2<p<\frac{N+\alpha}{N-2}$. Then 
$$
\lim_{\mu\to 0}\int_{|x|\le 4}\left|\frac{\mu^\sigma\xi_\mu^{\eta}}{|x|^{4}}(I_\alpha\ast |\tilde w_\mu|^{p})(\frac{x}{|x|^2})\tilde w_\mu^{p-2}(\frac{x}{|x|^2})\right|^{N/2}dx=0.
\eqno(6.53)
$$
}
\begin{proof}  We divide the integral into two parts:
$$
I_\mu^{(1)}(\frac{N}{2}):=\int_{|x|\le \mu^{\sigma/2}\xi_\mu/L_0}\left |\frac{\mu^\sigma\xi_\mu^{\eta}}{|x|^{4}}(I_\alpha\ast |\tilde w_\mu|^{p})(\frac{x}{|x|^2})\tilde w_\mu^{p-2}(\frac{x}{|x|^2})\right |^{N/2}dx,
$$
$$
I_\mu^{(2)}(\frac{N}{2}):=\int_{\mu^{\sigma/2}\xi_\mu/L_0\le |x|\le 4}\left |\frac{\mu^\sigma\xi_\mu^{\eta}}{|x|^{4}}(I_\alpha\ast |\tilde w_\mu|^{p})(\frac{x}{|x|^2})\tilde w_\mu^{p-2}(\frac{x}{|x|^2})\right |^{N/2}dx.
$$
By virtue of Lemma 3.10,  for all $|z|\ge L_0\mu^{-\sigma/2}\xi_\mu^{-1}$, we have 
$$
\begin{array}{rcl}
(I_\alpha \ast |\tilde w_\mu|^p)(z)&=&\xi_\mu^{\frac{N-2}{2}p-\alpha}(I_\alpha\ast |w_\mu|^p)(\xi_\mu z)\\
&\le& C\xi_\mu^{\frac{N-2}{2}p-\alpha}|\xi_\mu z|^{-(N-\alpha)}\\
&\le& C\xi_\mu^{\frac{N-2}{2}p-N}|z|^{-(N-\alpha)}.
\end{array}
\eqno(6.54)
$$
Therefore, by Lemma 6.13, it follows that
$$
\begin{array}{rcl}
I_\mu^{(1)}(\frac{N}{2})&=&\int_{|z|\ge L_0\mu^{-\sigma/2}\xi_\mu^{-1}} |\mu^\sigma\xi_\mu^{\eta}(I_\alpha\ast |\tilde w_\mu|^{p})(z)\tilde w_\mu^{p-2}(z)|^{N/2}dz\\
&\lesssim &\mu^{\frac{\sigma N}{2}+\frac{N(N-2)(p-2)\sigma}{8}+\frac{\sigma}{2}[\frac{N}{2}(N-\alpha)-N]}\xi_\mu^{\frac{\eta N}{2}+\frac{N(N-2)(p-2)}{4}+\frac{N}{2}(\frac{N-2}{2}p-\alpha)-N}\\
&\mbox{} &\qquad \qquad \cdot \int_{L_0}^{+\infty}s^{-\frac{N(N-\alpha)}{2}+N-1}e^{-\frac{N(p-2)}{4}s}ds\\
&\lesssim &\mu^{\frac{\sigma N}{8}[(N-2)(p-2)+2(N-\alpha)]},
\end{array}
$$
here we have used the fact that 
$$
\frac{\eta N}{2}+\frac{N(N-2)(p-2)}{4}+\frac{N}{2}(\frac{N-2}{2}p-\alpha)-N=0
$$
If $p>\frac{2\alpha}{N-2}$, it follows from the Hardy-Littlewood-Sobolev inequality  and the boundedness of $w_\mu$ in $H^1(\mathbb R^N)$ that 
$$
\begin{array}{rcl}
I_\mu^{(2)}(\frac{N}{2})&=&\int_{\frac{1}{4}\le |z|\le L_0\mu^{-\sigma/2}\xi_\mu^{-1}} |\mu^\sigma\xi_\mu^{\eta}(I_\alpha\ast |\tilde w_\mu|^{p})(z)\tilde w_\mu^{p-2}(z)|^{N/2}dz\\
&\le &\mu^{\frac{\sigma N}{2}}\xi_\mu^{\frac{\eta N}{2}}\left(\int_{\frac{1}{4}\le |z|\le L_0\mu^{-\sigma/2}\xi_\mu^{-1}}|\tilde w_\mu|^{2^*}dz\right)^{\frac{(N-2)(p-2)}{4}}\\
&\quad & \cdot \left(\int_{\frac{1}{4}\le |z|\le L_0\mu^{-\sigma/2}\xi_\mu^{-1}} |(I_\alpha\ast |\tilde w_\mu|^p)(z)|^{\frac{2N}{2N-(N-2)p}}dz\right)^{\frac{2N-(N-2)p}{4}}\\
&\lesssim & \mu^{\sigma N/2}\xi_\mu^{\eta N/2}\left(\int_{\mathbb R^N}|\tilde w_\mu|^{\frac{2Np}{2N-(N-2)p+2\alpha}}\right)^{\frac{2N-(N-2)p+2\alpha}{4}}\\
&\lesssim & \mu^{\sigma N/2}\xi_\mu^{\frac{\eta N}{2}+\frac{N(N-2)p}{4}-\frac{N[2N-(N-2)p+2\alpha]}{4}}\cdot \left(\int_{\mathbb R^N}| w_\mu|^{\frac{2Np}{2N-(N-2)p+2\alpha}}\right)^{\frac{2N-(N-2)p+2\alpha}{4}}\\
&\lesssim &\mu^{\sigma N/2}\to 0,
\end{array}
$$
as $\mu\to 0$, here we have used the facts that 
$$
2\le \frac{2Np}{2N-(N-2)p+2\alpha}\le 2^*,
$$
which follows from $\frac{N+\alpha}{N-2}>p\ge\max\{2, \frac{2\alpha}{N-2}\}\ge \frac{N+\alpha}{N-1}$ and 
$$
\frac{\eta N}{2}+\frac{N(N-2)p}{4}-\frac{N[2N-(N-2)p+2\alpha]}{4}=0.
$$
Next, we consider the case of $2<p\le \frac{2\alpha}{N-2}$.
If $q>1$ satisfies 
$$
\frac{2N-(N-2)p}{2(N-\alpha)}<q<+\infty.
\eqno(6.55)
$$
Then 
$$
1<s:=\frac{2Nq}{2N+2q\alpha-(N-2)p}<\frac{N}{\alpha}
$$
and 
 $$
\frac{2Nq}{2N-(N-2)p}=\frac{Ns}{N-\alpha s}.
$$
By the H\"older inequality and the Hardy-Littlewood-Sobolev inequality, we obtain 
$$
\begin{array}{lcl}
Q_\mu:&=&\int_{\frac{1}{4}\le |z|\le L_0\mu^{-\sigma/2}\xi_\mu^{-1}}|(I_\alpha\ast |\tilde w_\mu|^p)(z)|^{\frac{2N}{2N-(N-2)p}}dz\\
&\lesssim & \left(\int_{\frac{1}{4}\le |z|\le L_0\mu^{-\sigma/2}\xi_\mu^{-1}}dz\right)^{\frac{q-1}{q}}
\left(\int_{\frac{1}{4}\le |z|\le L_0\mu^{-\sigma/2}\xi_\mu^{-1}}|(I_\alpha\ast |\tilde w_\mu|^p|(z)|^{\frac{2Nq}{2N-(N-2)p}}dz\right)^{\frac{1}{q}}\\
&\lesssim &\mu^{-\frac{\sigma N(q-1)}{2q}}\xi_\mu^{-\frac{N(q-1)}{q}}\|\tilde w_\mu\|_{ps}^{\frac{2Np}{2N-(N-2)p}}\\
&=&\mu^{-\frac{\sigma N(q-1)}{2q}}\xi_\mu^{-\frac{N(q-1)}{q}+[\frac{N-2}{2}ps-N]\frac{2N}{s[2N-(N-2)p]}}\|w_\mu\|_{ps}^{\frac{2Np}{2N-(N-2)p}}.
\end{array}
$$
Since $2<p<ps<\frac{Np}{\alpha}\le 2^*$, it follows from the boundedness of $w_\mu$ in $H^1(\mathbb R^N)$ that
$$
Q_\mu\lesssim \mu^{-\frac{\sigma N(q-1)}{2q}}\xi_\mu^{-\frac{N(q-1)}{q}+[\frac{N-2}{2}ps-N]\frac{2N}{s[2N-(N-2)p]}}.
\eqno(6.56)
$$
Recall that $\eta=N+\alpha-(N-2)p$, it follows that
$$
\begin{array}{lcl}
\tau(\alpha,p):&=&\frac{\eta N}{2}-\left\{\frac{N(q-1)}{q}-[\frac{N-2}{2}ps-N]\frac{2N}{s[2N-(N-2)p]}\right\}\frac{2N-(N-2)p}{4}\\
&=&\frac{\eta N}{2}-\frac{N(q-1)}{4q}(2N-(N-2)p)+(\frac{N-2}{2}ps-N)\frac{N}{2s}\\
&=&\frac{\eta N}{2}-\frac{N(q-1)}{4q}(2N-(N-2)p)+\frac{Np(N-2)}{4}-\frac{N^2}{2s}\\
&=&\frac{\eta N}{2}-\frac{N(q-1)}{4q}(2N-(N-2)p)+\frac{Np(N-2)}{4}-\frac{N}{4q}(2N+2q\alpha-(N-2)p)\\
&=&0.
\end{array}
$$
On the other hand, it follows from $p>2$ that
$$
\frac{(q-1)[2N-(N-2)p]}{4q}<1.
$$
Therefore, from (6.56) we get
$$
\begin{array}{lcl}
I_\mu^{(2)}(\frac{N}{2})&\lesssim& \mu^{\frac{\sigma N}{2}-\frac{\sigma N(q-1)}{2q}\cdot\frac{2N-(N-2)p}{4}}\xi_\mu^{\frac{\eta N}{2}-\left\{\frac{N(q-1)}{q}-[\frac{N-2}{2}ps-N]\frac{2N}{s[2N-(N-2)p]}\right\}\frac{2N-(N-2)p}{4}}\\
&\lesssim& \mu^{\frac{\sigma N}{2}\left[1-\frac{(q-1)[2N-(N-2)p]}{4q}\right]}\to 0, \quad {\rm as} \ \mu\to 0.
\end{array}
$$
Thus the conclusion follows since $\lim_{\mu\to 0}(I_\mu^{(1)}(\frac{N}{2})+I_\mu^{(2)}(\frac{N}{2}))=0$. 
The proof is complete. \end{proof}

 It follows from (6.50), (6.52), (6.53) and Lemma 3.12 (i) that for any $r>1$, there exists $\mu_r>0$ such that 
$$
\sup_{\mu\in (0,\mu_r)}\|K[\tilde w_\mu]^r\|_{H^1(B_1)}\le C_r.
\eqno(6.57)
$$

To verify the condition in Lemma 3.12 (ii), we show that there exists $r_0>\frac{N}{2}$ such that 
$$
\lim_{\mu\to 0}\int_{|x|\le 4}\left|\frac{\mu^\sigma\xi_\mu^{\eta}}{|x|^{4}}(I_\alpha\ast |\tilde w_\mu|^{p})(\frac{x}{|x|^2})\tilde w_\mu^{p-2}(\frac{x}{|x|^2})\right|^{r_0}dx=0.
\eqno(6.58)
$$
We divide the integral in (6.58) into two parts:
$$
I_\mu^{(1)}(r_0):=\int_{|x|\le \mu^{\sigma/2}\xi_\mu/L_0}\left|\frac{\mu^\sigma\xi_\mu^{\eta}}{|x|^{4}}(I_\alpha\ast |\tilde w_\mu|^{p})(\frac{x}{|x|^2})\tilde w_\mu^{p-2}(\frac{x}{|x|^2})\right|^{r_0}dx
$$ 
and 
$$
I_\mu^{(2)}(r_0):=\int_{\mu^{\sigma/2}\xi_\mu/L_0\le |x|\le 4}\left|\frac{\mu^\sigma\xi_\mu^{\eta}}{|x|^{4}}(I_\alpha\ast |\tilde w_\mu|^{p})(\frac{x}{|x|^2})\tilde w_\mu^{p-2}(\frac{x}{|x|^2})\right|^{r_0}dx,
$$
and show that $\lim_{\mu\to 0}(I_\mu^{(1)}(r_0)+I_\mu^{(2)}(r_0))=0.$  

In what follows, we denote $\eta_\mu:=\mu^{\sigma/2}\xi_\mu/L_0$. It is easy to show that 
$$
I_\mu^{(1)}(r_0)=\mu^{\sigma r_0}\xi_\mu^{\eta r_0}\int_{|x|\le \eta_\mu}\left|\frac{1}{|x|^{\gamma}}(I_\alpha\ast |\tilde w_\mu|^{p})(\frac{x}{|x|^2})K[\tilde w_\mu]^{p-2}(x)\right|^{r_0}dx
$$
and 
$$
I_\mu^{(2)}(r_0)=\mu^{\sigma r_0}\xi_\mu^{\eta r_0}\int_{\eta_\mu\le |x|\le 4}\left|\frac{1}{|x|^{\gamma}}(I_\alpha\ast |\tilde w_\mu|^{p})(\frac{x}{|x|^2})K[\tilde w_\mu]^{p-2}(x)\right|^{r_0}dx.
$$

To proceed, we  divide the associated $(\alpha,p)$ regions into three parts: 
$$\mathcal A_N:=\left\{(\alpha,p):  \  p>2,  \   \frac{2\alpha}{N-2}\le p<\frac{N+\alpha}{N-2}\right\},$$
$$\mathcal B_N:=\left\{(\alpha,p):  \ 
\frac{2(N^2-2\alpha)}{N(N-2)}<p< \frac{2\alpha}{N-2} \right\},$$
and 
$$\mathcal C_N:=\left\{(\alpha,p): \    2<p<\frac{2\alpha}{N-2}, \    \  p\le \frac{2(N^2-2\alpha)}{N(N-2)}\right\},
$$ 
which are depicted in the following pictures. We shall verify the conditions needed in Lemma 3.12 (ii) under the additional assumptions $(\alpha,p)\in \mathcal A_N$, $(\alpha,p)\in \mathcal B_N$ and $\mathcal C_N$, respectively. This will be done in Lemmas 6.19-6.22.

\begin{figure}
	\centering
	\includegraphics[width=0.5\linewidth]{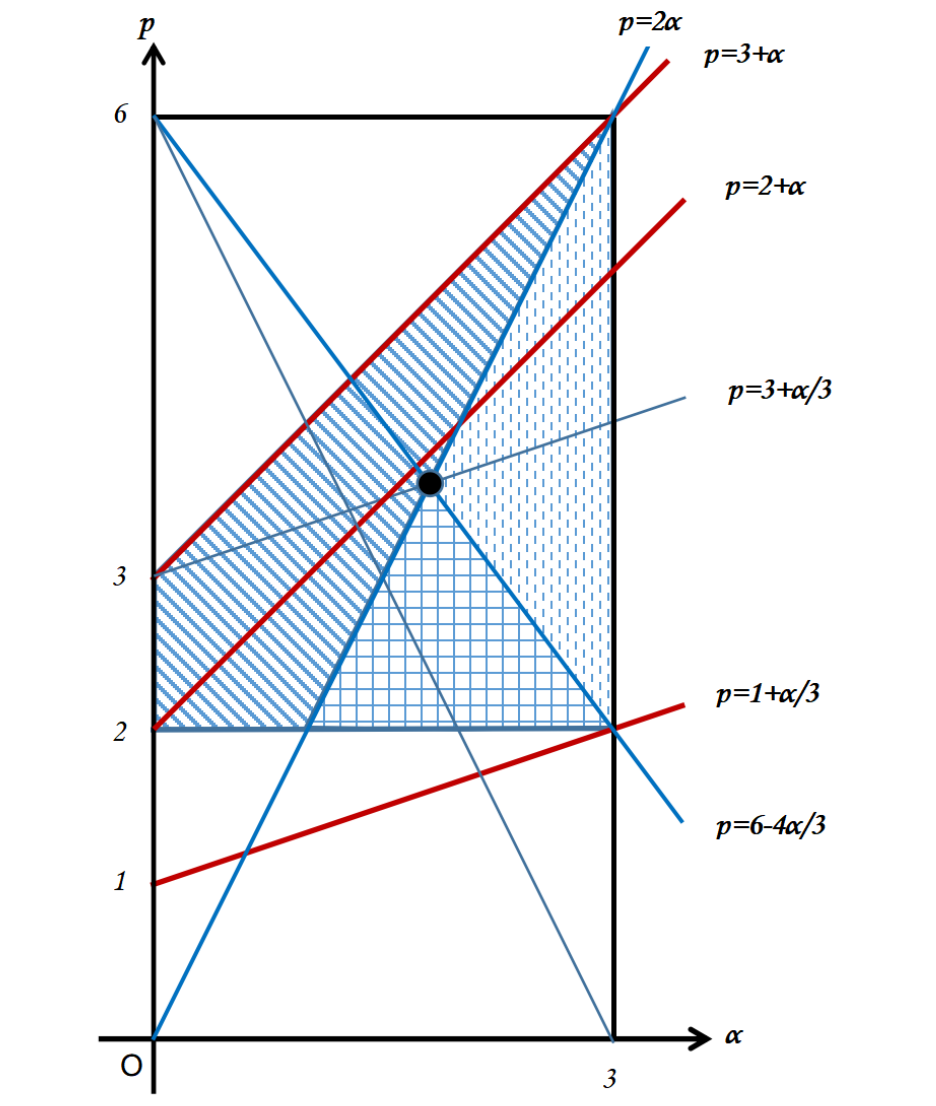}
	\caption{Three parts of the associated $(\alpha,p)$ regions in $N=3$.}
	\label{fig:figure7}
\end{figure}


\begin{figure}
	\centering
	\includegraphics[width=0.5\linewidth]{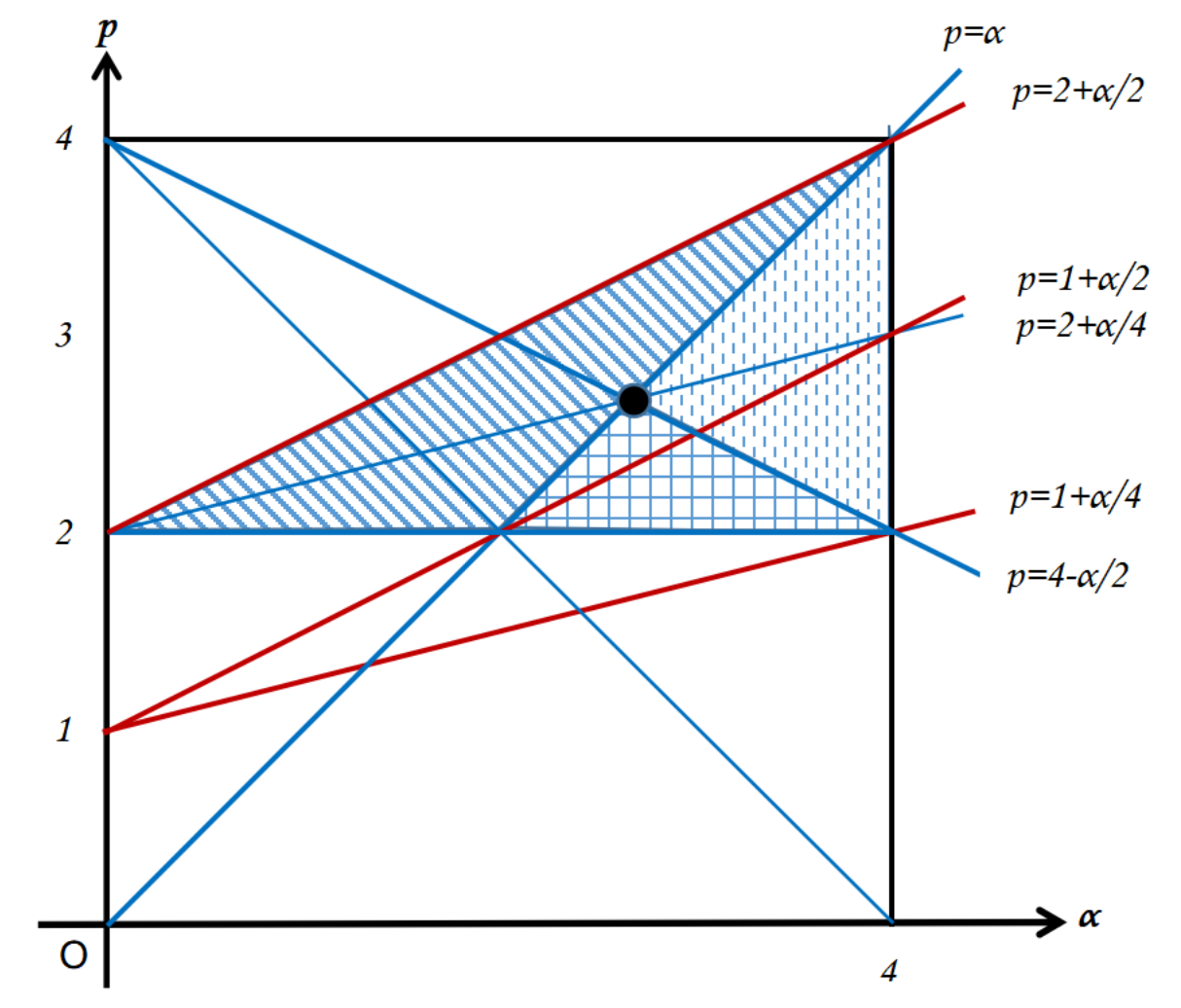}
	\caption{Three parts of the associated $(\alpha,p)$ regions in $N=4$.}
	\label{fig:figure8}
\end{figure}

%
%

\vskip 5mm

 As the first step, we check the conditions of Lemma 3.12 (ii) in an additional condition $p\ge \frac{2\alpha}{N-2}$, in this case, the well-known Hardy-Littlewood-Sobolev inequality apply. The cases $p>\frac{2\alpha}{N-2}$ and $p=\frac{2\alpha}{N-2}$  will be treated  in Lemma 6.19 and Lemma 6.20, respectively.

{\bf Lemma 6.19.}  {\it Assume $N\ge 3$, $\alpha>N-4$ and $\max\{2,\frac{2\alpha}{N-2}\}<p<\frac{N+\alpha}{N-2}$. Then 
 there exists $r_0>\frac{N}{2}$ such that 
$$
\lim_{\mu\to 0}\int_{|x|\le 4}\left|\frac{\mu^\sigma\xi_\mu^{\eta}}{|x|^{4}}(I_\alpha\ast |\tilde w_\mu|^{p})(\frac{x}{|x|^2})\tilde w_\mu^{p-2}(\frac{x}{|x|^2})\right|^{r_0}dx=0.
$$
}
\begin{proof} 
Since $\alpha>N-4$, we have $\frac{2N}{N-\alpha}>\frac{N}{2}$. We  show that  for some $r_0\in (\frac{N}{2}, \frac{2N}{N-\alpha})$, there holds
$$
\lim_{\mu\to 0}(I_\mu^{(1)}(r_0)+I_\mu^{(2)}(r_0))=0.
\eqno(6.59)
$$
Let 
$$
\eta_1=2[N+\alpha-(N-2)p], \quad \eta_2=(N-2)p-2\alpha.
$$
Then 
$$
\eta_1+\eta_2=\gamma:=2N-(N-2)p.
$$
Let $\epsilon>0$ be a small number to be specified later. Put
$$
 s_1=\frac{2N}{2N-[(N-2)p-2\alpha](r_0+\epsilon)}, \quad  s_2=\frac{2N}{[(N-2)p-2\alpha]r_0}, \quad s_3=\frac{2N}{[(N-2)p-2\alpha]\epsilon},
 $$
where $r_0>\frac{N}{2}$ and $\epsilon>0$ is such that $s_1, s_2, s_3>0$.  Then 
 $$
 \frac{1}{s_1}+\frac{1}{s_2}+\frac{1}{s_3}=1.
 $$
 By  the H\"older inequality,  we have 
$$
\begin{array}{rcl}
I_\mu^{(1)}(r_0)&=&\mu^{\sigma r_0}\xi_\mu^{\eta r_0}\int_{|x|\le \eta_\mu}\left|\frac{1}{|x|^{\gamma}}(I_\alpha\ast |\tilde w_\mu|^{p})(\frac{x}{|x|^2})K[\tilde w_\mu]^{p-2}(x)\right|^{r_0}dx\\
&\le &\mu^{\sigma r_0}\xi_\mu^{\eta r_0}\left(\int_{|x|\le \eta_\mu} |x|^{-\eta_1r_0s_1}K[\tilde w_\mu]^{(p-2)\theta r_0s_1}(x)\right)^{\frac{1}{s_1}}\\
&\mbox{}& \cdot\left(\int_{\R^N}|x|^{-\eta_2r_0s_2}\left|(I_\alpha\ast |\tilde w_\mu|^{p})(\frac{x}{|x|^2})\right|^{r_0s_2}dx\right)^{\frac{1}{s_2}}\\
&\mbox{}& \cdot \left(\int_{|x|\le 4}K[\tilde  w_\mu]^{(p-2)(1-\theta)r_0s_3}(x)\right)^{\frac{1}{s_3}}.
\end{array}
\eqno(6.60)
$$
Note that $\eta_2r_0s_2=2N$, by the Hardy-Littlewood-Sobolev inequality, we get 
$$
\begin{array}{rl}
&\int_{\mathbb R^N}\frac{1}{|x|^{\eta_2r_0s_2}}|(I_\alpha\ast |\tilde w_\mu|^{p})(\frac{x}{|x|^2})|^{r_0s_2}dx\\
&=\int_{\mathbb R^N}|(I_\alpha\ast |\tilde w_\mu|^{p})(z)|^{r_0s_2}dz\\
&\le C\left(\int_{\mathbb R^N}|\tilde w_\mu^p|^{\frac{2N}{(N-2)p}}\right)^{\frac{(N-2)p}{(N-2)p-2\alpha}}\le C<+\infty.
 \end{array}
\eqno(6.61)
$$
 Let $\theta\in (0,1)$ be a fixed number. Then $(p-2)(1-\theta)r_0s_3>2^*$ for small $\epsilon>0$,  it follows from (6.57) that
$$
\int_{|x|\le 4}|K[\tilde w_\mu]|^{(p-2)(1-\theta)r_0s_3}dx\le C<\infty.
\eqno(6.62)
$$
 By (6.50), (6.60), (6.61), (6.62),  we have 
$$
\begin{array}{rcl}
I_\mu^{(1)}(r_0)&\lesssim&\mu^{\sigma r_0}\xi_\mu^{\eta r_0}\left(\int_{|x|\le\eta_\mu} |x|^{-\eta_1r_0s_1}K[\tilde w_\mu]^{(p-2)\theta r_0s_1}(x)\right)^{\frac{1}{s_1}}\\
&\lesssim&\mu^{\sigma r_0+\frac{\sigma(N-2)(p-2)\theta}{4}r_0}\xi_\mu^{\eta r_0+\frac{(N-2)(p-2)\theta}{2}r_0}\\
&\mbox{}& \cdot \left(\int_{|x|\le \eta_\mu} |x|^{-\eta_1r_0s_1-(N-2)(p-2)\theta r_0s_1}\exp(-\frac{1}{2}(p-2)\theta r_0s_1\mu^{\frac{\sigma}{2}}\xi_\mu|x|^{-1})dx\right)^{\frac{1}{s_1}}\\
&=&\mu^{\sigma r_0+\frac{\sigma(N-2)(p-2)\theta}{4}r_0}\xi_\mu^{\eta r_0+\frac{(N-2)(p-2)\theta}{2}r_0}\\
&\mbox{}& \cdot \left(\int_0^{\eta_\mu} r^{-\eta_1r_0s_1-(N-2)(p-2)\theta r_0s_1+N-1}\exp(-\frac{1}{2}(p-2)\theta r_0s_1\mu^{\frac{\sigma}{2}}\xi_\mu r^{-1})dr\right)^{\frac{1}{s_1}}\\
&=&\mu^{\sigma r_0-\frac{\sigma(N-2)(p-2)\theta}{4}r_0+\frac{\sigma}{2}[-\eta_1r_0+\frac{N}{s_1}]}\xi_\mu^{\eta r_0-\frac{(N-2)(p-2)\theta}{2}r_0-\eta_1 r_0 +\frac{N}{s_1}}\\
&\mbox{}& \cdot \left(\int_{L_0}^\infty s^{\eta_1r_0s_1+(N-2)(p-2)\theta r_0s_1-N-1}\exp(-\frac{1}{2}(p-2)\theta r_0s_1s)ds\right)^{\frac{1}{s_1}}.
\end{array}
$$
Let 
$$
\tau_1(r_0,s_1):=\sigma r_0-\frac{\sigma(N-2)(p-2)\theta}{4}r_0+\frac{\sigma}{2}[-\eta_1r_0+\frac{N}{s_1}]
$$
and 
$$
\tau_2(r_0,s_1):=\eta r_0-\frac{(N-2)(p-2)\theta}{2}r_0-\eta_1 r_0 +\frac{N}{s_1}.
$$
Since 
$$
p>\max\{2, \frac{2\alpha}{N-2}\}\ge \frac{4-2\theta}{3-\theta}+\frac{2\alpha}{(N-2)(3-\theta)}, 
$$
it is easy to check that 
$$
\begin{array}{lcl}
\tau_1(\frac{N}{2}, \frac{4}{4-p(N-2)+2\alpha})&=&\frac{\sigma N}{2}[1-\frac{1}{2}\eta_1-\frac{1}{4}(N-2)(p-2)\theta+\frac{4-p(N-2)+2\alpha}{4}]\\
&=&\frac{\sigma N}{2}[(N-2)(\frac{1}{4}p(3-\theta)+\frac{1}{2}\theta-1)-\frac{\alpha}{2}]>0,
\end{array}
$$
and
$$
\begin{array}{lcl}
\tau_2(\frac{N}{2},  \frac{4}{4-p(N-2)+2\alpha})&=&\frac{ N}{2}[\eta-\eta_1-\frac{1}{2}(N-2)(p-2)\theta +\frac{4-p(N-2)+2\alpha}{2}]\\
&=&\frac{N}{4}(N-2)(p-2)(1-\theta)>0.
\end{array}
$$
Therefore, we can choose $r_0>\frac{N}{2}$ with $r_0-\frac{N}{2}$ being small, and $\epsilon>0$ small, 
so that 
$$
\tau_1(r_0,s_1)>0 \quad {\rm  and} \quad \tau_2(r_0,s_1)>0.
$$
 Thus, we conclude that 
$$
I_\mu^{(1)}(r_0)\lesssim \mu^{\tau_1(r_0,s_1)}\xi_\mu^{\tau_2(r_0,s_1)}\to 0, \  \ {\rm as} \  \mu\to 0.
\eqno(6.63)
$$

Let $\epsilon>0$ be a small number, and $\tau(r_0)$ be given as follows
$$
\tau(r_0):=\frac{2N-2[(N-2)(p-1)-\alpha]r_0}{(N-2)(p-2)r_0}>0, \quad r_0\in \left(\frac{N}{2},\frac{N}{(N-2)(p-1)-\alpha}\right).
$$
Put
$$ 
\tilde s_1=\frac{N(1+\epsilon)^{-1}}{N-[(N-2)(p-1)-\alpha]r_0}, \quad  s_2=\frac{2N}{[(N-2)p-2\alpha]r_0}, \quad \tilde s_3=\frac{2N(1-\epsilon\tau(r_0))^{-1}}{(N-2)(p-2)r_0}.
$$
Then it is easy to see that
$$
\frac{1}{\tilde s_1}+\frac{1}{s_2}+\frac{1}{\tilde s_3}=1.
$$
By  the H\"older inequality, we  have
$$
\begin{array}{rcl}
I_\mu^{(2)}(r_0)&=&\mu^{\sigma r_0}\xi_\mu^{\eta r_0}\int_{\eta_\mu\le |x|\le 4}\left|\frac{1}{|x|^{\gamma}}(I_\alpha\ast |\tilde w_\mu|^{p})(\frac{x}{|x|^2})K[\tilde w_\mu]^{p-2}(x)\right|^{r_0}dx\\
&\lesssim &\mu^{\sigma r_0}\xi_\mu^{\eta r_0}\left(\int_{\eta_\mu\le |x|\le 4} \frac{1}{|x|^{\eta_1r_0\tilde s_1}}dx\right)^{\frac{1}{\tilde s_1}}\left(\int_{|x|\le 4}|K[\tilde w_\mu]|^{(p-2)r_0\tilde s_3}dx\right)^{\frac{1}{\tilde s_3}}\\
& \mbox{}& \   \cdot
\left(\int_{\mathbb R^N}\frac{1}{|x|^{\eta_2r_0s_2}}|(I_\alpha\ast |\tilde w_\mu|^{p})(\frac{x}{|x|^2})|^{r_0s_2}dx\right)^{\frac{1}{s_2}}.
\end{array}
\eqno(6.64)
$$
We claim that 
$$
\eta_1r_0\tilde s_1-N>\frac{N}{2}, \quad {\rm for \ all}  \ r_0>\frac{N}{2}.
\eqno(6.65)
$$
In fact, (6.65) is equivalent to 
$$
h_3(r_0):=\frac{2[(N+\alpha)-(N-2)p]r_0}{N-[(N-2)(p-1)-\alpha]r_0}>\frac{3}{2}(1+\epsilon),
$$
 for  all $\ r_0\in (\frac{N}{2},\frac{N}{(N-2)(p-1)-\alpha})$.
The latter follows from the fact that $h_3(\frac{N}{2})=2$ and $h_3(r_0)$ is increasing.  From this fact, the claim (6.65) follows.
From (5.65), we find
$$
\begin{array}{rl}
&\int_{ \eta_\mu\le |x|\le 4} |x|^{-\eta_1r_0\tilde s_1}dx
=\int_{\eta_\mu}^4r^{-\eta_1r_0\tilde s_1+N-1}dr\\
&= \frac{1}{\eta_1r_0\tilde s_1-N}\left[L_0^{\eta_1r_0\tilde s_1-N}\mu^{\frac{\sigma}{2}[-\eta_1r_0\tilde s_1+N]}\xi_\mu^{-\eta_1r_0\tilde s_1+N}-4^{-\eta_1r_0\tilde s_1+N}\right]\\
&\lesssim  \mu^{\frac{\sigma}{2}[-\eta_1r_0\tilde s_1+N]}\xi_\mu^{-\eta_1r_0\tilde s_1+N}.
\end{array}
\eqno(6.66)
$$
Since $(p-2)r_0\tilde s_3=\frac{2N}{N-2}(1-\epsilon\tau(r_0))^{-1}>2^*$, it follows from (6.57) that
$$
\int_{|x|\le 4}|K[\tilde w_\mu]|^{(p-2)r_0\tilde s_3}dx\le C<\infty.
\eqno(6.67)
$$
Therefore, by (6.61), (6.64), (6.66) and (6.67),  we get
$$
I_\mu^{(2)}(r_0)\lesssim  \mu^{\sigma r_0+\frac{\sigma}{2}[-\eta_1r_0+\frac{N}{\tilde s_1}]}\xi_\mu^{\eta r_0-\eta_1r_0+\frac{N}{\tilde s_1}}.
$$
Let 
$$
\tau_3(r_0):=\sigma r_0+\frac{\sigma}{2}[-\eta_1r_0+\frac{N}{\tilde s_1}],
\quad 
\tau_4(r_0):=\eta r_0-\eta_1r_0+\frac{N}{\tilde s_1}.
$$
Since $p>\max\{2, \frac{2\alpha}{N-2}\}\ge 1+\frac{\alpha}{N-2}$, it is easy to check that
$$
\begin{array}{lcl}
\tau_3(\frac{N}{2})&=&\frac{\sigma N}{2}-\frac{\sigma N}{4}\eta_1+\frac{\sigma}{2}N(1+\epsilon)-\frac{\sigma N}{4}((N-2)(p-1)-\alpha)(1+\epsilon)\\
&=&\frac{\sigma N}{4}[(N-2)(p-1)-\alpha](1-\epsilon)+\frac{\sigma N}{2}\epsilon>0,
\end{array}
$$
and 
$$
\begin{array}{lcl}
\tau_4(\frac{N}{2})&=&\frac{N}{2}[\eta-\eta_1-((N-2)(p-1)-\alpha)(1+\epsilon)+2(1+\epsilon)]\\
&=&\frac{N}{2}\epsilon [N+\alpha-p(N-2)]>0.
\end{array}
$$
Therefore, we can choose $r_0>\frac{N}{2}$ such that $\tau_3(r_0)>0$ and $\tau_4(r_0)>0$. Thus,   we obtain
$$
I_\mu^{(2)}(r_0)\lesssim \mu^{\tau_3(r_0)}\xi_\mu^{\tau_4(r_0)}\to 0, \quad {\rm as} \ \mu\to 0.
\eqno(6.68)
$$
Finally, (6.59) follows from  (6.63) and (6.68), the proof is complete. 
\end{proof}

{\bf Lemma 6.20.}  {\it Assume $N\ge 3$, $\alpha>N-2$ and $p=\frac{2\alpha}{N-2}$. Then 
 there exists $r_0>\frac{N}{2}$ such that 
$$
\lim_{\mu\to 0}\int_{|x|\le 4}\left|\frac{\mu^\sigma\xi_\mu^{\eta}}{|x|^{4}}(I_\alpha\ast |\tilde w_\mu|^{p})(\frac{x}{|x|^2})\tilde w_\mu^{p-2}(\frac{x}{|x|^2})\right|^{r_0}dx=0.
$$
}
\begin{proof} 
Since $\alpha>N-2$, we have $\frac{N}{N-\alpha}>\frac{N}{2}$.
We  show that  for some $r_0\in (\frac{N}{2}, \frac{N}{N-\alpha})$, there holds
$$
\lim_{\mu\to 0}(I_\mu^{(1)}(r_0)+I_\mu^{(2)}(r_0))=0.
\eqno(6.69)
$$
Let $\epsilon>0$ and $\theta\in (0,1)$  be  small numbers satisfying
$$
\epsilon<\frac{2+\alpha-N}{2}
$$
and 
$$
0<\theta<\frac{2+\alpha-N-2\epsilon}{2+\alpha-N}.
$$
Put 
$$
\eta_1=(2N-p(N-2))\epsilon=2(N-\alpha)\epsilon, \quad \eta_2=(2N-p(N-2))(1-\epsilon)=2(N-\alpha)(1-\epsilon).
$$
Then 
$$
\eta_1+\eta_2=\gamma:=2N-(N-2)p=2(N-\alpha).
$$
Let $r_0\in (\frac{N}{2},\frac{N}{N-\alpha})$ be a number  to be specified later. Let
$$
s_1=\frac{N}{[N-(N-\alpha)r_0](1-\epsilon)}, \quad s_2=\frac{N}{(N-\alpha)r_0(1-\epsilon)},\quad s_3=\frac{1}{\epsilon}.
$$
Then
 $$
 \frac{1}{s_1}+\frac{1}{s_2}+\frac{1}{s_3}=1.
 $$
 By  the H\"older inequality,  we have 
$$
\begin{array}{rcl}
I_\mu^{(1)}(r_0)&=&\mu^{\sigma r_0}\xi_\mu^{\eta r_0}\int_{|x|\le \eta_\mu}\left|\frac{1}{|x|^{\gamma}}(I_\alpha\ast |\tilde w_\mu|^{p})(\frac{x}{|x|^2})K[\tilde w_\mu]^{p-2}(x)\right|^{r_0}dx\\
&\le &\mu^{\sigma r_0}\xi_\mu^{\eta r_0}\left(\int_{|x|\le \eta_\mu} |x|^{-\eta_1r_0s_1}K[\tilde w_\mu]^{(p-2)\theta r_0s_1}(x)\right)^{\frac{1}{s_1}}\\
&\mbox{}& \cdot\left(\int_{\R^N}|x|^{-\eta_2r_0s_2}\left|(I_\alpha\ast |\tilde w_\mu|^{p})(\frac{x}{|x|^2})\right|^{r_0s_2}dx\right)^{\frac{1}{s_2}}\\
&\mbox{}& \cdot \left(\int_{|x|\le 4}K[\tilde  w_\mu]^{(p-2)(1-\theta)r_0s_3}(x)\right)^{\frac{1}{s_3}}.
\end{array}
\eqno(6.70)
$$
Note that $\eta_2r_0s_2=2N$, by the Hardy-Littlewood-Sobolev inequality, we get 
$$
\begin{array}{rl}
&\int_{\mathbb R^N}\frac{1}{|x|^{\eta_2r_0s_2}}|(I_\alpha\ast |\tilde w_\mu|^{p})(\frac{x}{|x|^2})|^{r_0s_2}dx\\
&=\int_{\mathbb R^N}|(I_\alpha\ast |\tilde w_\mu|^{p})(z)|^{r_0s_2}dz\\
&\le C\left(\int_{\mathbb R^N}|\tilde w_\mu|^{ps_0}\right)^{\frac{r_0s_2}{s_0}}\lesssim \xi_\mu^{[\frac{N-2}{2}p-\frac{N}{s_0}]r_0s_2},
 \end{array}
\eqno(6.71)
$$
where 
$$
1<s_0:=\frac{Nr_0s_2}{r_0s_2\alpha+N}=\frac{N}{N-(N-\alpha)\epsilon}<\frac{N}{\alpha}.
$$
 Fix $\theta\in (0,1)$ and take $\epsilon>0$ small if necessary,   we then have  $(p-2)(1-\theta)r_0s_3>2^*$.  Therefore, it follows from (6.56) that
$$
\int_{|x|\le 4}|K[\tilde w_\mu]|^{(p-2)(1-\theta)r_0s_3}dx\le C<\infty.
\eqno(6.72)
$$
 By (6.50), (6.70), (6.71), (6.72),  we have 
$$
\begin{array}{rcl}
I_\mu^{(1)}(r_0)&\lesssim&\mu^{\sigma r_0}\xi_\mu^{\eta r_0}\left(\int_{|x|\le \eta_\mu} |x|^{-\eta_1r_0s_1}K[\tilde w_\mu]^{(p-2)\theta r_0s_1}(x)\right)^{\frac{1}{s_1}}\\
&\lesssim&\mu^{\sigma r_0+\frac{\sigma(N-2)(p-2)\theta}{4}r_0}\xi_\mu^{\eta r_0+\frac{(N-2)(p-2)\theta}{2}r_0+[\frac{N-2}{2}p-\frac{N}{s_0}]r_0}\\
&\mbox{}& \cdot \left(\int_{|x|\le \eta_\mu} |x|^{-\eta_1r_0s_1-(N-2)(p-2)\theta r_0s_1}\exp(-\frac{1}{2}(p-2)\theta r_0s_1\mu^{\frac{\sigma}{2}}\xi_\mu|x|^{-1})dx\right)^{\frac{1}{s_1}}\\
&=&\mu^{\sigma r_0+\frac{\sigma(N-2)(p-2)\theta}{4}r_0}\xi_\mu^{\eta r_0+\frac{(N-2)(p-2)\theta}{2}r_0+[\frac{N-2}{2}p-\frac{N}{s_0}]r_0}\\
&\mbox{}& \cdot \left(\int_0^{\eta_\mu} r^{-\eta_1r_0s_1-(N-2)(p-2)\theta r_0s_1+N-1}\exp(-\frac{1}{2}(p-2)\theta r_0s_1\mu^{\frac{\sigma}{2}}\xi_\mu r^{-1})dr\right)^{\frac{1}{s_1}}\\
&=&\mu^{\sigma r_0-\frac{\sigma(N-2)(p-2)\theta}{4}r_0+\frac{\sigma}{2}[-\eta_1r_0+\frac{N}{s_1}]}\xi_\mu^{\eta r_0-\frac{(N-2)(p-2)\theta}{2}r_0-\eta_1 r_0 +\frac{N}{s_1}+[\frac{N-2}{2}p-\frac{N}{s_0}]r_0}\\
&\mbox{}& \cdot \left(\int_{L_0}^\infty s^{\eta_1r_0s_1+(N-2)(p-2)\theta r_0s_1-N-1}\exp(-\frac{1}{2}(p-2)\theta r_0s_1s)ds\right)^{\frac{1}{s_1}}.
\end{array}
$$
Let 
$$
\tau_5(r_0,s_1):=\sigma r_0-\frac{\sigma(N-2)(p-2)\theta}{4}r_0+\frac{\sigma}{2}[-\eta_1r_0+\frac{N}{s_1}]
$$
and 
$$
\tau_6(r_0,s_1):=\eta r_0-\frac{(N-2)(p-2)\theta}{2}r_0-\eta_1 r_0 +\frac{N}{s_1}+[\frac{N-2}{2}p-\frac{N}{s_0}]r_0.
$$
Recall that $\eta_2r_0s_2=2N$ and 
$p=\frac{2\alpha}{N-2},$
it is easy to check that 
$$
\tau_5(r_0, s_1)=\sigma\left\{[1 -\frac{(N-2)(p-2)\theta}{4}-\frac{1}{2}\eta_1]r_0+\frac{N}{2s_1}\right\}
$$
and 
$$
\begin{array}{lcl}
\tau_6(r_0,s_1)&=&\eta r_0-\frac{(N-2)(p-2)\theta}{2}r_0-\eta_1r_0+\frac{N}{s_1}-\frac{N}{s_2}\\
&=&\eta r_0-\frac{(N-2)(p-2)\theta}{2}r_0-\eta_1r_0+\frac{N}{s_1}-\frac{1}{2}\eta_2 r_0\\
&=&[-\frac{1}{2}\eta_1-\frac{1}{2}(N-2)(p-2)\theta]r_0+\frac{N}{s_1}.
\end{array}
$$
Hence, we get
$$
\begin{array}{lcl}
\tau_5(\frac{N}{2}, \frac{2}{(2+\alpha-N)(1-\epsilon)})&=&\frac{\sigma N}{2}[3+\alpha-N-2\epsilon-\frac{1}{2}\eta_1-\frac{1}{4}(N-2)(p-2)\theta]\\
&=&\frac{\sigma N}{4}[2-(N-\alpha)\epsilon+2+\alpha-N-2\epsilon-\frac{1}{2}(N-2)(p-2)\theta]
>0
\end{array}
$$
and
$$
\begin{array}{lcl}
\tau_6(\frac{N}{2},  \frac{2}{(2+\alpha-N)(1-\epsilon)})&=&\frac{ N}{2}[2+\alpha-N-2\epsilon-\frac{1}{2}(N-2)(p-2)\theta]>0.
\end{array}
$$
Therefore, we can choose $r_0>\frac{N}{2}$ with $r_0-\frac{N}{2}$ being small 
so that 
$$
\tau_5(r_0,s_1)>0 \quad {\rm  and} \quad \tau_6(r_0,s_1)>0.
$$
 Thus, we conclude that 
$$
I_\mu^{(1)}(r_0)\lesssim \mu^{\tau_5(r_0,s_1)}\xi_\mu^{\tau_6(r_0,s_1)}\to 0, \  \ {\rm as} \  \mu\to 0.
\eqno(6.73)
$$

Let 
$$
\tilde s_1=\frac{N}{(N-\alpha)r_0\epsilon}, \quad s_2=\frac{N}{(N-\alpha)r_0(1-\epsilon)},\quad \tilde  s_3=\frac{N}{N-(N-\alpha)r_0}.
$$
Then it is easy to see that
$$
\frac{1}{\tilde s_1}+\frac{1}{s_2}+\frac{1}{\tilde s_3}=1.
$$
By  the H\"older inequality, we  have
$$
\begin{array}{rcl}
I_\mu^{(2)}(r_0)&=&\mu^{\sigma r_0}\xi_\mu^{\eta r_0}\int_{\eta_\mu\le |x|\le 4}\left|\frac{1}{|x|^{\gamma}}(I_\alpha\ast |\tilde w_\mu|^{p})(\frac{x}{|x|^2})K[\tilde w_\mu]^{p-2}(x)\right|^{r_0}dx\\
&\lesssim &\mu^{\sigma r_0}\xi_\mu^{\eta r_0}\left(\int_{\eta_\mu\le |x|\le 4} \frac{1}{|x|^{\eta_1r_0\tilde s_1}}dx\right)^{\frac{1}{\tilde s_1}}\left(\int_{|x|\le 4}|K[\tilde w_\mu]|^{(p-2)r_0\tilde s_3}dx\right)^{\frac{1}{\tilde s_3}}\\
& \mbox{}& \   \cdot
\left(\int_{\mathbb R^N}\frac{1}{|x|^{\eta_2r_0s_2}}|(I_\alpha\ast |\tilde w_\mu|^{p})(\frac{x}{|x|^2})|^{r_0s_2}dx\right)^{\frac{1}{s_2}}.
\end{array}
\eqno(6.74)
$$

Note that  $\eta_1r_0\tilde s_1=2N$, it follows  that
$$
\int_{ \eta_\mu\le |x|\le 4} |x|^{-\eta_1r_0\tilde s_1}dx
=\int_{\eta_\mu}^4r^{-\eta_1r_0\tilde s_1+N-1}dr
\lesssim  \mu^{-\frac{\sigma N}{2}}\xi_\mu^{-N}.
\eqno(6.75)
$$
Since $(p-2)r_0\tilde s_3\ge  (p-2)\frac{N}{2}\tilde s_3\ge 2^*$, it follows from (6.57) that
$$
\int_{|x|\le 4}|K[\tilde w_\mu]|^{(p-2)r_0\tilde s_3}dx\le C<\infty.
\eqno(6.76)
$$
Therefore, by (6.71), (6.74), (6.75) and (6.76),  we get
$$
I_\mu^{(2)}(r_0)\lesssim  \mu^{\sigma r_0-\frac{\sigma N}{2\tilde s_1}}\xi_\mu^{\eta r_0-\frac{N}{\tilde s_1}+[\frac{N-2}{2}p-\frac{N}{s_0}]r_0}.
$$
It is easy to check that
$$
\tau_7(r_0):=\sigma r_0-\frac{\sigma N}{2\tilde s_1}=\frac{\sigma r_0}{2}[2-(N-\alpha)\epsilon]>0
$$
and 
$$
\tau_8(r_0):=\eta r_0-\frac{N}{\tilde s_1}+[\frac{N-2}{2}p-\frac{N}{s_0}]r_0=0.
$$
Thus,   we obtain
$$
I_\mu^{(2)}(r_0)\lesssim \mu^{\tau_7(r_0)}\to 0, \quad {\rm as} \ \mu\to 0.
\eqno(6.77)
$$
Finally, (6.69) follows from (6.73) and (6.77), the proof is complete. 
\end{proof}

Next, we consider the case $p< \frac{2\alpha}{N-2}$. In this case, it seems that  the Hardy-Littlewood-Sobolev inequality does not apply directly. So we apply the decay estimates of the Riesz potential  established in Lemma 6.15 and Lemma 6.16.

{\bf Lemma 6.21.}  {\it Assume $N\ge 3$ and 
$$
\frac{2(N^2-2\alpha)}{N(N-2)}<p<\frac{2\alpha}{N-2}.
$$ 
Then there exists $r_0>\frac{N}{2}$ such that
$$
\lim_{\mu\to 0}\int_{|x|\le 4}\left|\frac{\mu^\sigma\xi_\mu^{\eta}}{|x|^{4}}(I_\alpha\ast |\tilde w_\mu|^{p})(\frac{x}{|x|^2})\tilde w_\mu^{p-2}(\frac{x}{|x|^2})\right|^{r_0}dx=0.
$$}
\begin{proof}  
Let $\theta\in (0,1)$ be a small constant and $r_0>\frac{N}{2}$,  for any $s_1,s_2>1$ verifying $\frac{1}{s_1}+\frac{1}{s_2}=1$,  by the H\"older inequality and Lemma 6.15, we have 
$$
\begin{array}{rcl}
I_\mu^{(1)}(r_0)&=& \int_{|x|\le \eta_\mu}\left|\frac{\mu^\sigma\xi_\mu^{\eta}}{|x|^{\gamma}}(I_\alpha\ast |\tilde w_\mu|^{p})(\frac{x}{|x|^2})K[\tilde w_\mu]^{p-2}(x)\right|^{r_0}dx\\
&\le &\left(\int_{|x|\le \eta_\mu}\left|\frac{\mu^\sigma\xi_\mu^{\eta}}{|x|^{\gamma }}(I_\alpha\ast |\tilde w_\mu|^{p})(\frac{x}{|x|^2})K[\tilde w_\mu]^{(p-2)\theta}(x)\right|^{r_0s_1}dx\right)^{\frac{1}{s_1}}\\
&\mbox{}& \quad\cdot \left(\int_{|x|\le \eta_\mu}|K[\tilde w_\mu]|^{(p-2)(1-\theta)r_0s_2}dx\right)^{\frac{1}{s_2}}\\
&\lesssim & \mu^{\sigma r_0}\xi_\mu^{\eta r_0}\left(\int_{ |x|\le\eta_\mu}|x|^{-\gamma r_0s_1}\left|(I_\alpha\ast |\tilde w_\mu|^{p})(\frac{x}{|x|^2})K[\tilde w_\mu|^{(p-2)\theta}(x)\right|^{r_0s_1}dx\right)^{\frac{1}{s_1}}\\
&\le &\mu^{\sigma r_0}\xi_\mu^{\eta r_0}\cdot \xi_\mu^{-(\frac{2N-(N-2)p}{2}+\epsilon)r_0}\cdot \mu^{\frac{\sigma(N-2)(p-2)}{4}r_0\theta}\xi_\mu^{\frac{(N-2)(p-2)}{2}r_0\theta}\\
&\mbox{}&\quad \cdot \left(\int_{|x|\le\eta_\mu}|x|^{[-\gamma +N-\alpha-(N-2)(p-2)\theta] r_0s_1}
e^{-\frac{(p-2)\theta}{2}r_0s_1\mu^{\sigma/2}\xi_\mu |x|^{-1}}dx\right)^{\frac{1}{s_1}}\\
&\le &\mu^{\sigma r_0}\xi_\mu^{\eta r_0}\cdot \xi_\mu^{-(\frac{2N-(N-2)p}{2}+\epsilon)r_0}\cdot \mu^{\frac{\sigma(N-2)(p-2)}{4}r_0\theta}\xi_\mu^{\frac{(N-2)(p-2)}{2}r_0\theta}\\
&\mbox{}&\quad \cdot \left(\int_0^{\eta_\mu}r^{[-\gamma +N-\alpha-(N-2)(p-2)\theta] r_0s_1+N-1}
e^{-\frac{(p-2)\theta}{2}r_0s_1\mu^{\sigma/2}\xi_\mu r^{-1}}dr\right)^{\frac{1}{s_1}}\\
&\le &\mu^{\tau_9(r_0,s_1)} \xi_\mu^{\tau_{10}(r_0,s_1)}\\
&\mbox{}&\quad \cdot \left(\int_{L_0}^\infty s^{[\gamma +(N-2)(p-2)\theta -\frac{1}{2}(N-2)(p-\frac{2\alpha}{N})]r_0s_1-N-1}
e^{-\frac{(p-2)\theta}{2}r_0s_1s}ds\right)^{\frac{1}{s_1}},
\end{array}
\eqno(6.78)
$$
where
$$
\tau_9(r_0,s_1):=\sigma r_0+\frac{\sigma (N-2)(p-2)\theta}{4}r_0+\frac{\sigma}{2}[-\gamma+N-\alpha-(N-2)(p-2)\theta]r_0+\frac{\sigma N}{2s_1},
$$
$$
\tau_{10}(r_0,s_1):=\eta r_0-(\frac{2N-(N-2)p}{2}+\epsilon)r_0-\gamma r_0-\frac{(N-2)(p-2)\theta}{2}r_0+(N-\alpha)r_0+N/s_1.
$$
A direct computation shows that
$$
\tau_9(\frac{N}{2},1)=\frac{\sigma N}{4}[\frac{1}{2}(N-2)(p-2)(2-\theta)+N-\alpha]>0,
$$
$$
\tau_{10}(\frac{N}{2},1)=\frac{1}{4}N(N-2)(p-2)(1-\theta)-\frac{1}{2}N\epsilon>0.
$$
Therefore, we can choose $r_0>\frac{N}{2}$ and $s_1>1$ with $r_0-\frac{N}{2}$ and $s_1-1$ being small, such that  $\tau_9(r_0,s_1)>0$ and $\tau_{10}(r_0,s_1)>0$. Thus, we obtain  
$$
I_\mu^{(1)}(r_0)\lesssim \mu^{\tau_9(r_0,s_1)} \xi_\mu^{\tau_{10}(r_0,s_1)}\to 0, \   \  {\rm as}  \   \mu\to 0.
$$
By the H\"older inequality and Lemma 6.15, we have 
$$
\begin{array}{rcl}
I_\mu^{(2)}(r_0)&=& \int_{\eta_\mu\le |x|\le 4}\left|\frac{\mu^\sigma\xi_\mu^{\eta}}{|x|^{\gamma}}(I_\alpha\ast |\tilde w_\mu|^{p})(\frac{x}{|x|^2})K[\tilde w_\mu]^{p-2}(x)\right|^{r_0}dx\\
&\le& \mu^{\sigma r_0}\xi_\mu^{\eta r_0}\left(\int_{\eta_\mu\le |x|\le 4}|x|^{-\gamma r_0s_1}|(I_\alpha\ast |\tilde w_\mu|^{p})(\frac{x}{|x|^2})|^{r_0s_1}\right)^{\frac{1}{s_1}}\\
&\qquad &\cdot \left(\int_{|x|\le 4}|K[\tilde w_\mu]|^{(p-2)r_0s_2}\right)^{\frac{1}{s_2}}\\
&\lesssim &\mu^{\sigma r_0}\xi_\mu^{\eta r_0-(\frac{2N-(N-2)p}{2}+\epsilon)r_0} 
\left(\int_{\eta_\mu\le |x|\le 4}|x|^{-\gamma r_0s_1+(N-\alpha)r_0s_1}dx\right)^{\frac{1}{s_1}}\\
&=&\mu^{\sigma r_0}\xi_\mu^{\eta r_0-(\frac{2N-(N-2)p}{2}+\epsilon)r_0}
\cdot \left(\int_{\eta_\mu}^4r^{-\gamma r_0s_1+(N-\alpha)r_0s_1+N-1}dr\right)^{\frac{1}{s_1}}.
\end{array}
\eqno(6.79)
$$
If 
$
-\gamma+N-\alpha+2>0,
$
then  we can choose $r_0>\frac{N}{2}$ and $s_1>1$ with $r_0-\frac{N}{2}>0$ and $s_1-1>0$ being small, such that
$$
-\gamma r_0s_1+(N-\alpha)r_0s_1+N>0.
$$
Therefore, we obtain 
$$
\begin{array}{rcl}
I_\mu^{(2)}(r_0)&\lesssim&\mu^{\sigma r_0}\xi_\mu^{\eta r_0-(\frac{2N-(N-2)p}{2}+\epsilon)r_0}
\cdot \left(\int_{0}^4r^{-\gamma r_0s_1+(N-\alpha)r_0s_1+N-1}dx\right)^{\frac{1}{s_1}}\\
&\lesssim& \mu^{\sigma r_0}\xi_\mu^{\eta r_0-(\frac{2N-(N-2)p}{2}+\epsilon)r_0}\to 0,
\end{array}
\eqno(6.80)
$$ 
as $\mu\to 0$, here we have used the fact that 
$$
\eta r_0-(\frac{2N-(N-2)p}{2}+\epsilon)r_0=[\alpha-\frac{1}{2}(N-2)p-\epsilon]r_0\ge 0.
$$ 
If 
$
-\gamma+N-\alpha+2\le 0,
$
then for $r_0-\frac{N}{2}>0$ and $s_1-1>0$ small, we have 
$$
\tau_{11}(r_0,s_1):=-\gamma r_0s_1+(N-\alpha)r_0s_1+N<0.
$$
Therefore, we obtain
$$
\begin{array}{rcl}
I_\mu^{(2)}(r_0)&\lesssim&\mu^{\sigma r_0}\xi_\mu^{\eta r_0-(\frac{2N-(N-2)p}{2}+\epsilon)r_0}\\
&\mbox{}&\quad \cdot \left(\frac{1}{\tau_{11}(r_0,s_1)}[4^{\tau_{11}(r_0,s_1)}-L_0^{-\tau_{11}(r_0,s_1)}\mu^{\frac{\sigma}{2}\tau_{11}(r_0s_1)}\xi_\mu^{\tau_{11}(r_0s_1)}]\right)^{\frac{1}{s_1}}\\
&\lesssim& \mu^{\frac{\sigma}{2}[2-\gamma+N-\alpha]r_0+\frac{\sigma N}{2s_1}}\xi_\mu^{[\eta-(\frac{2N-(N-2)p}{2}+\epsilon)-\gamma+N-\alpha]r_0+\frac{N}{s_1}}.
\end{array}
\eqno(6.81)
$$
Since
$$
[2-\gamma+N-\alpha]\frac{N}{2}+N=[(N-2)(p-2)+N-\alpha]\frac{N}{2}>0,
$$
$$
[\eta-(\frac{2N-(N-2)p}{2}+\epsilon)-\gamma+N-\alpha]\frac{N}{2}+N=\frac{1}{4}N(N-2)(p-2)-\frac{1}{2}N\epsilon>0,
$$
 we can choose $r_0>\frac{N}{2}$ and $s_1>1$ with $r_0-\frac{N}{2}>0$ and $s_1-1>0$ being small, such that
$$
\frac{\sigma}{2}[2-\gamma+N-\alpha]r_0+\frac{\sigma N}{2s_1}>0,
$$
and 
$$
[\eta-(\frac{2N-(N-2)p}{2}+\epsilon)-\gamma+N-\alpha]r_0+\frac{N}{s_1}>0.
$$
It follows from (6.81) that $\lim_{\mu\to 0}I_\mu^{(2)}(r_0)=0$.
Thus, $\lim_{\mu\to 0}(I_\mu^{(1)}(r_0)+I_\mu^{(2)}(r_0))=0$ and the proof is complete. 
\end{proof}

{\bf Lemma 6.22.}  {\it Assume $N\ge 3$ and 
$$
2<p< \frac{2\alpha}{N-2} \quad {and} \quad 
p\le \frac{2(N^2-2\alpha)}{N(N-2)}.
 $$ 
Then there exists $r_0>\frac{N}{2}$ such that
$$
\lim_{\mu\to 0}\int_{|x|\le 4}\left|\frac{\mu^\sigma\xi_\mu^{\eta}}{|x|^{4}}(I_\alpha\ast |\tilde w_\mu|^{p})(\frac{x}{|x|^2})\tilde w_\mu^{p-2}(\frac{x}{|x|^2})\right|^{r_0}dx=0.
$$}
\begin{proof} 
Fix $\theta\in (\frac{\alpha}{N},1)$, for $r_0>\frac{N}{2}$ and any $s_1,s_2>1$ verifying $\frac{1}{s_1}+\frac{1}{s_2}=1$,  by the H\"older inequality, we have 
$$
\begin{array}{rcl}
I_\mu^{(1)}(r_0)&=& \int_{|x|\le \eta_\mu}\left|\frac{\mu^\sigma\xi_\mu^{\eta}}{|x|^{\gamma}}(I_\alpha\ast |\tilde w_\mu|^{p})(\frac{x}{|x|^2})K[\tilde w_\mu]^{p-2}(x)\right|^{r_0}dx\\
&\le & \mu^{\sigma r_0}\xi_\mu^{\eta r_0}\left(\int_{|x|\le \eta_\mu}\left|\frac{1}{|x|^{\gamma }}(I_\alpha\ast |\tilde w_\mu|^{p})(\frac{x}{|x|^2})K[\tilde w_\mu]^{(p-2)\theta}(x)\right|^{r_0s_1}dx\right)^{\frac{1}{s_1}}\\
&\mbox{}& \cdot \left(\int_{|x|\le \eta_\mu}|K[\tilde w_\mu]^{(p-2)(1-\theta)}(x)|^{r_0s_2}dx\right)^{\frac{1}{s_2}}.
\end{array}
\eqno(6.82)
$$
By Lemma 6.16 and (6.50), we find that
$$
\begin{array}{rl}
&\left(\int_{|x|\le \eta_\mu}\left|\frac{1}{|x|^{\gamma }}(I_\alpha\ast |\tilde w_\mu|^{p})(\frac{x}{|x|^2})K[\tilde w_\mu]^{(p-2)\theta}(x)\right|^{r_0s_1}dx\right)^{\frac{1}{s_1}}\\
&= \xi_\mu^{-(\frac{2\alpha}{N}+\epsilon)r_0}\cdot \mu^{\frac{\sigma(N-2)(p-2)}{4}r_0\theta}\xi_\mu^{\frac{(N-2)(p-2)}{2}r_0\theta}\\
&\mbox{} \cdot \left(\int_{|x|\le\eta_\mu}|x|^{[-\gamma +\frac{1}{2}(N-2)(p-\frac{2\alpha}{N})-(N-2)(p-2)\theta] r_0s_1}
e^{-\frac{(p-2)\theta}{2}r_0s_1\mu^{\sigma/2}\xi_\mu |x|^{-1}}dx\right)^{\frac{1}{s_1}}\\
&\le \xi_\mu^{-(\frac{2\alpha}{N}+\epsilon)r_0}\cdot \mu^{\frac{\sigma(N-2)(p-2)}{4}r_0\theta}\xi_\mu^{\frac{(N-2)(p-2)}{2}r_0\theta}\\
&\mbox{} \cdot \left(\int_{L_0}^\infty s^{[\gamma +(N-2)(p-2)\theta -\frac{1}{2}(N-2)(p-\frac{2\alpha}{N})]r_0s_1-N-1}
e^{-\frac{(p-2)\theta}{2}r_0s_1s}ds\right)^{\frac{1}{s_1}},
\end{array}
$$
which together with (6.82)  and (6.57) implies that
$$
I_\mu^{(1)}(r_0) \lesssim \mu^{\tau_{12}(r_0,s_1)}\xi_\mu^{\tau_{13}(\epsilon,r_0,s_1)},
\eqno(6.83)
$$
where
$$
\tau_{12}(r_0,s_1):=\frac{\sigma}{2}[2-\gamma+\frac{1}{2}(N-2)(p-\frac{2\alpha}{N})-\frac{1}{2}(N-2)(p-2)\theta]r_0+\frac{\sigma N}{2s_1},
$$
$$
\tau_{13}(\epsilon, r_0,s_1):=[\eta-(\frac{2\alpha}{N}+\epsilon)-\gamma -\frac{1}{2}(N-2)(p-2)\theta+\frac{1}{2}(N-2)(p-\frac{2\alpha}{N})]r_0+N/s_1.
$$
Since 
$$
\tau_{12}(\frac{N}{2},1)=\frac{\sigma}{8}N(N-2)[N(p-2)(3-\theta)+2(N-\alpha)]>0,
$$
$$
\tau_{13}(0, \frac{N}{2},1)=\frac{1}{4}N(N-2)(p-2)(1-\theta)>0.
$$
Therefore, we can choose $\epsilon>0$ sufficiently small and choose $r_0>\frac{N}{2}$ and $s_1>1$ with $r_0-\frac{N}{2}$ and $s_1-1$ being small, such that  $\tau_{12}(r_0,s_1)>0$ and $\tau_{13}(\epsilon,r_0,s_1)>0$. Thus, by (6.83), we obtain  
$$
I_\mu^{(1)}(r_0)\lesssim \mu^{\tau_{12}(r_0,s_1)} \xi_\mu^{\tau_{13}(\epsilon, r_0,s_1)}\to 0, \   \  {\rm as}  \   \mu\to 0.
$$
By the H\"older inequality and Lemma 6.16, we have 
$$
\begin{array}{rcl}
I_\mu^{(2)}(r_0)&=& \int_{\eta_\mu\le |x|\le 4}\left|\frac{\mu^\sigma\xi_\mu^{\eta}}{|x|^{\gamma}}(I_\alpha\ast |\tilde w_\mu|^{p})(\frac{x}{|x|^2})K[\tilde w_\mu]^{p-2}(x)\right|^{r_0}dx\\
&\le& \mu^{\sigma r_0}\xi_\mu^{\eta r_0}\left(\int_{\eta_\mu\le |x|\le 4}|x|^{-\gamma r_0s_1}|(I_\alpha\ast |\tilde w_\mu|^{p})(\frac{x}{|x|^2})|^{r_0s_1}\right)^{\frac{1}{s_1}}\\
&\qquad &\cdot \left(\int_{|x|\le 4}|K[\tilde w_\mu]|^{(p-2)r_0s_2}\right)^{\frac{1}{s_2}}\\
&\lesssim &\mu^{\sigma r_0}\xi_\mu^{\eta r_0-(\frac{2\alpha}{N}+\epsilon)r_0} 
\left(\int_{\eta_\mu\le |x|\le 4}|x|^{-\gamma r_0s_1+\frac{1}{2}(N-2)(p-\frac{2\alpha}{N})r_0s_1}dx\right)^{\frac{1}{s_1}}.
\end{array}
\eqno(6.84)
$$
If 
$
-\gamma+\frac{1}{2}(N-2)(p-\frac{2\alpha}{N})+2>0,
$
then   we can choose $r_0>\frac{N}{2}$ and $s_1>1$ with $r_0-\frac{N}{2}>0$ and $s_1-1>0$ being small, such that
$$
-\gamma r_0s_1+\frac{1}{2}(N-2)(p-\frac{2\alpha}{N})r_0s_1+N>0.
$$
Therefore, from (6.84) we obtain 
$$
\begin{array}{rcl}
I_\mu^{(2)}(r_0)&\lesssim&\mu^{\sigma r_0}\xi_\mu^{\eta r_0-(\frac{2\alpha}{N}+\epsilon)r_0}
\cdot \left(\int_{0}^4r^{-\gamma r_0s_1+\frac{1}{2}(N-2)(p-\frac{2\alpha}{N})r_0s_1+N-1}dx\right)^{\frac{1}{s_1}}\\
&\lesssim& \mu^{\sigma r_0}\xi_\mu^{\eta r_0-(\frac{2\alpha}{N}+\epsilon)r_0}.
\end{array}
\eqno(6.85)
$$ 
Since $p\le\frac{2(N^2-2\alpha)}{N(N-2)}$ and $p<\frac{2\alpha}{N-2}$, it is easy to see that
$$
p<\frac{N^2+(N-2)\alpha}{N(N-2)}=\left\{\begin{array}{rcl}
2+\frac{1}{4}\alpha, \quad  &if & \ N=4,\\
3+\frac{1}{3}\alpha, \quad &if & \ N=3,
\end{array}
\right. 
$$
which implies that
$$
\eta-\frac{2\alpha}{N}=N+\alpha-(N-2)p-\frac{2\alpha}{N}=(N-2)\left(\frac{N^2+(N-2)\alpha}{N(N-2)}-p\right)>0.
$$
Therefore,  we can choose $\epsilon>0$ sufficiently small such that $\eta-(\frac{2\alpha}{N}+\epsilon)> 0$, and hence
$\lim_{\mu\to 0}I_\mu^{(2)}(r_0)=0$.

If 
$
-\gamma+\frac{1}{2}(N-2)(p-\frac{2\alpha}{N})+2\le 0,
$
then we can choose  $r_0-\frac{N}{2}>0$ and $s_1-1>0$ small, such that 
$$
\tau_{14}(r_0,s_1):=-\gamma r_0s_1+\frac{1}{2}(N-2)(p-\frac{2\alpha}{N})r_0s_1+N<0.
$$
Therefore, from (6.84) we obtain
$$
\begin{array}{rcl}
I_\mu^{(2)}(r_0)&\lesssim&\mu^{\sigma r_0}\xi_\mu^{\eta r_0-(\frac{2\alpha}{N}+\epsilon)r_0}
\mu^{\frac{\sigma}{2s_1}\tau_{14}(r_0,s_1)}\xi_\mu^{\tau_{14}(r_0,s_1)/s_1}\\
&\lesssim& \mu^{\frac{\sigma}{2}[2-\gamma+\frac{1}{2}(N-2)(p-\frac{2\alpha}{N})]r_0+\frac{\sigma N}{2s_1}}\xi_\mu^{[\eta-(\frac{2\alpha}{N}+\epsilon)-\gamma+\frac{1}{2}(N-2)(p-\frac{2\alpha}{N})]r_0+\frac{N}{s_1}}.
\end{array}
\eqno(6.86)
$$
Since
$$
4-\gamma+\frac{1}{2}(N-2)(p-\frac{2\alpha}{N})=(N-2)\left(\frac{3}{2}p-2-\frac{\alpha}{N}\right)>0,
$$
$$
\eta-\frac{2\alpha}{N}-\gamma+\frac{1}{2}(N-2)(p-\frac{2\alpha}{N})+2
=\frac{1}{4}N(N-2)(p-2)>0,
$$
we can choose $\epsilon>0$ sufficiently small,  and  choose $r_0>\frac{N}{2}$ and $s_1>1$ with $r_0-\frac{N}{2}>0$ and $s_1-1>0$ being small, such that
$$
\frac{\sigma}{2}[2-\gamma+\frac{1}{2}(N-2)(p-\frac{2\alpha}{N})]r_0+\frac{\sigma N}{2s_1}>0,
$$
and 
$$
[\eta-(\frac{2\alpha}{N}+\epsilon)-\gamma+\frac{1}{2}(N-2)(p-\frac{2\alpha}{N})]r_0+\frac{N}{s_1}>0.
$$
Therefore, by (6.86), we obtain
$$
I_\mu^{(2)}(r_0)\lesssim\mu^{\frac{\sigma}{2}[2-\gamma+\frac{1}{2}(N-2)(p-\frac{2\alpha}{N})]r_0+\frac{\sigma N}{2s_1}}\xi_\mu^{[\eta-(\frac{2\alpha}{N}+\epsilon)-\gamma+\frac{1}{2}(N-2)(p-\frac{2\alpha}{N})]r_0+\frac{N}{s_1}}\to 0, 
$$
as $\mu\to 0$. The proof is complete.
\end{proof}

{\bf Lemma 6.23.} { If $\max\{2,\frac{N+\alpha}{2(N-2)}\}<p<\frac{N+\alpha}{N-2}$, then $\int_{\mathbb R^N}(I_\alpha\ast |\tilde w_\mu|^p)|\tilde w_\mu|^p\sim 1$ as $\mu\to 0$. Furthermore, $\tilde w_\mu\to W_1$ in $L^{\frac{2Np}{N+\alpha}}(\mathbb R^N)$ as $\mu\to 0$. }
\begin{proof}   
By virtue of  Proposition 6.18,  there exists a constant $C>0$ such that for all small $\mu>0$, 
$$
\tilde w_\mu(x)\le \frac{C}{(1+|x|)^{N-2}}, \qquad \forall x\in \mathbb R^N,
$$
which together with the fact that $\frac{2Np}{N+\alpha}>\frac{N}{N-2}$ implies that $\tilde w_\mu$ is bounded in $L^{\frac{2Np}{N+\alpha}}(\mathbb R^N)$  uniformly for small $\mu>0$,  and 
 by the dominated convergence theorem $\tilde w_\mu\to W_1$ in $L^{\frac{2Np}{N+\alpha}}(\mathbb R^N)$ as $\mu\to 0$. 
 Therefore, by the Hardy-Littlewood-Sobolev inequality, we find that 
 $$
 \int_{\mathbb R^N}(I_\alpha\ast |\tilde w_\mu-W_1|^p)|\tilde w_\mu-W_1|^p \le C\|\tilde w_\mu-W_1\|_{\frac{2Np}{N+\alpha}}^{2p}\to 0,
 $$
 as $\mu\to 0$.  Since 
 $$
 \lim_{\mu\to 0}\int_{\mathbb R^N}(I_\alpha\ast |\tilde w_\mu|^p)|\tilde w_\mu|^p-\int_{\mathbb R^N}(I_\alpha\ast |\tilde w_\mu-W_1|^p)|\tilde w_\mu-W_1|^p=
 \int_{\mathbb R^N}(I_\alpha\ast |W_1|^p)|W_1|^p,
 $$ 
 it follows that
$$ 
\lim_{\mu\to 0}\int_{\mathbb R^N}(I_\alpha\ast |\tilde w_\mu|^p)|\tilde w_\mu|^p=
 \int_{\mathbb R^N}(I_\alpha\ast |W_1|^p)|W_1|^p.
 $$
Particularly, we have  that $\int_{\mathbb R^N}(I_\alpha\ast |\tilde w_\mu|^p)|\tilde w_\mu|^p\sim 1$ as $\mu\to 0$. 
\end{proof}

\begin{proof}[Proof of Theorem 2.3.]   For $N\ge 5$, the conclusion follows directly from Lemmas 6.5, 6.6 and 6.8.  We only consider the cases $N=4$ and $N=3$.

We first note that for a result similar to Lemma 3.2 holds for  $\tilde w_\mu$ and $\tilde J_\mu$. By  (6.2) and Lemma 6.9, we also have  $\tau_3(\tilde w_\mu)=\tau_3(w_\mu)$. Therefore, by (6.24),  we get
$$
\begin{array}{rcl}
m_0&\le& \sup_{t\ge 0} \tilde J_\lambda((\tilde w_\mu)_t)\\
&\quad &+\mu^\sigma\left\{\frac{\tau_3(\tilde w_\mu)^{\frac{N+\alpha}{2}}}{2p}\xi_\mu^{(N+\alpha)-p(N-2)}\int_{\mathbb R^N}(I_\alpha\ast |\tilde w_\mu|^p)|\tilde w_\mu|^p-\frac{\tau_3(\tilde w_\mu)^{\frac{N}{2}}}{2}\xi_\mu^{2}\int_{\mathbb R^N}|\tilde w_\mu|^2\right\}\\
&=&m_\mu+\mu^\sigma\tau_3(w_\mu)^{\frac{N}{2}}\frac{4+p(N-2)-(N+\alpha)}{4p}\xi_\mu^{(N+\alpha)-p(N-2)}\int_{\mathbb R^N}(I_\alpha\ast |\tilde w_\mu|^p)|\tilde w_\mu|^p,
\end{array}
\eqno(6.87)
$$
which implies that
$$
\xi_\mu^{(N+\alpha)-p(N-2)}\int_{\mathbb R^N}(I_\alpha\ast |\tilde w_\mu|^p)|\tilde w_\mu|^p\ge \mu^{-\sigma}\frac{4p}{[4+p(N-2)-(N+\alpha)]\tau_3(w_\mu)^{\frac{N}{2}}}\delta_\mu.
$$
Hence, by Lemma 6.5, we obtain
$$
\begin{array}{rl}
&\xi_\mu^{(N+\alpha)-p(N-2)}\int_{\mathbb R^N}(I_\alpha\ast |\tilde w_\mu|^p)|\tilde w_\mu|^p\gtrsim \mu^{-\sigma}\delta_\mu\\
&\qquad \gtrsim  \left\{\begin{array}{rcl} 
(\ln\frac{1}{\mu})^{-\frac{4+\alpha-2p}{2p-2-\alpha}}, \    \qquad &{ if}& \   \ N=4,\\
 \mu^{\frac{3+\alpha-p}{(p-2-\alpha)(p-1-\alpha)}},    \quad &{ if}& \    \  N=3 \ and \ p\in (2+\alpha, 3+\alpha).
 \end{array}\right.
 \end{array}
\eqno(6.88)
 $$  
 Therefore, by Lemma 6.23, we have 
$$
\xi_\mu\gtrsim    \left\{\begin{array}{rcl} 
(\ln\frac{1}{\mu})^{-\frac{1}{2p-2-\alpha}}, \    \qquad &{ if}& \   \ N=4,\\
 \mu^{\frac{1}{(p-2-\alpha)(p-1-\alpha)}},    \quad &{ if}& \    \  N=3 \ and \ p\in (2+\alpha, 3+\alpha).
 \end{array}\right.
\eqno(6.89)
$$
On the other hand, if  $ N=3$, then by  (6.24) and Lemma 6.11 and Lemma 6.23, we have 
$$
\xi_\mu^{p-1-\alpha}\lesssim \frac{1}{\|\tilde w_\mu\|_2^2}\lesssim\mu^{\frac{\sigma}{2}}\xi_\mu.
$$
Then 
$$
\xi_\mu^{p-2-\alpha}\lesssim \mu^{\frac{\sigma}{2}}.
$$
Hence, noting that $\sigma=\frac{2}{p-1-\alpha}$, we have 
$$
\xi_\mu\lesssim \mu^{\frac{1}{(p-1-\alpha)(p-2-\alpha)}}.
\eqno(6.90)
$$
 If $N=4$, then by (6.24) and Lemma 6.12 and Lemma 6.23,  we have 
$$
\xi_\mu^{2p-2-\alpha}\lesssim \frac{1}{\|\tilde w_\mu\|_2^2}\lesssim \frac{1}{-\ln(\mu^\sigma\xi_\mu^{2})}.
$$
Note that 
$$
-\ln(\mu^\sigma\xi_\mu^{2})=\sigma\ln\frac{1}{\mu}+2\ln\frac{1}{\xi_\mu}\ge \sigma\ln\frac{1}{\mu},
$$
it follows that 
$$
\xi_\mu^{2p-2-\alpha}\lesssim  \frac{1}{\|\tilde w_\mu\|_2^2}\lesssim \left(\ln\frac{1}{\mu}\right)^{-1},
$$
we then obtain 
$$
\xi_\mu\lesssim   \left(\ln\frac{1}{\mu}\right)^{-\frac{1}{2p-2-\alpha}}.
\eqno(6.91)
$$
Thus,  it follows from (6.87), (6.90), (6.91)  and Lemma 6.23 that 
$$
\delta_\mu=m_0-m_\mu\lesssim \mu^\sigma\xi_\mu^{(N+\alpha)-p(N-2)}\lesssim \left\{\begin{array}{rcl} \mu^{\frac{2}{2p-2-\alpha}}(\ln\frac{1}{\mu})^{-\frac{4+\alpha-2p}{2p-2-\alpha}},   \  &{\rm if}& \   N=4,\\
 \mu^{\frac{1}{p-2-\alpha}}, \qquad  \qquad  \qquad    &{\rm if}& \     N=3,
 \end{array}\right.
$$
which together with Lemma 6.5 implies that 
$$
\delta_\mu\sim \mu^\sigma\xi_\mu^{(N+\alpha)-p(N-2)}\sim\left\{\begin{array}{rcl} \mu^{\frac{2}{2p-2-\alpha}}(\ln\frac{1}{\mu})^{-\frac{4+\alpha-2p}{2p-2-\alpha}},   \  &{\rm if}& \   N=4,\\
 \mu^{\frac{1}{p-2-\alpha}}, \qquad  \qquad  \qquad    &{\rm if}& \     N=3.
 \end{array}\right.
\eqno(6.92)
$$
Arguing as in the proof of Theorem 2.2,  we also have 
$$
\|\nabla \tilde w_\mu\|_2^2=S^{\frac{N}{2}}+\left\{\begin{array}{rcl} O(\mu^{\frac{2}{2p-2-\alpha}}(\ln\frac{1}{\mu})^{-\frac{4+\alpha-2p}{2p-2-\alpha}}),   \quad &{\rm if}& \   \ N=4,\\
  O(\mu^{\frac{1}{p-2-\alpha}}), \qquad  \qquad  \qquad &{\rm if}& \    \  N=3,
 \end{array}\right.
\eqno(6.93)
$$
and 
$$
\|\tilde w_\mu\|_{2^*}^{2^*}=S^{\frac{N}{2}}+\left\{\begin{array}{rcl} O(\mu^{\frac{2}{2p-2-\alpha}}(\ln\frac{1}{\mu})^{-\frac{4+\alpha-2p}{2p-2-\alpha}}),   \quad &{\rm if}& \   \ N=4,\\
 O(\mu^{\frac{1}{p-2-\alpha}}), \qquad  \qquad  \qquad &{\rm if}& \    \  N=3.
 \end{array}\right.
$$
Finally, by (6.24), Lemma 6.11 and Lemma 6.12, we obtain
$$
\|\tilde w_\mu\|_2^2\sim\left\{\begin{array}{rcl}
\ln\frac{1}{\mu},   \quad  \quad if \   \ N=4,\\
\mu^{-\frac{1}{p-2-\alpha}},  \quad if  \    \  N=3.
\end{array}\right.
$$
The statements on $v_\mu$ follow from the corresponding results on $w_\mu$ and $\tilde w_\mu$. This completes the proof of Theorem 2.3.
\end{proof}

\vskip 5mm

\section*{7. Proofs of Other Results and Final Remarks}

At first, we prove Theorem 2.4  and Theorem 2.5. 
We consider $(Q_\lambda)$ and its limit equation
$$
-\Delta v+v=(I_\alpha\ast |v|^p)|v|^{p-2}v.
\eqno(7.1)
$$
The corresponding energies of ground states  are given by $m_\lambda=\inf_{v\in \mathcal M_\lambda}I_\lambda(v)$ and 
$$
m_0:=\inf_{v\in \mathcal M_0}I_0(v),
$$
where 
$$
I_0(v)=\frac{1}{2}\int_{\mathbb R^N}|\nabla v|^2+|v|^2-\frac{1}{2p}\int_{\mathbb R^N}(I_\alpha\ast |v|^p)|v|^p, 
\eqno(7.2)
$$
and 
$$
\mathcal{M}_0=
\left\{v\in H^1(\mathbb R^N)\setminus\{0\} \ \left | \ \int_{\mathbb R^N}|v|^2+|v|^2=\int_{\mathbb R^N}(I_\alpha\ast |v|^p)|v|^p  \  \right. \right\}.
$$
Then $m_\lambda$ and $m_0$ are well-defined and positive. Moreover, $I_0$ is attained on $\mathcal M_0$ by  positive solutions of (7.1).

{\bf Lemma 7.1.}  {\it Assume that the assumptions of Theorems 2.5 holds. Then 
for small $\lambda>0$, there holds
$$
m_0-m_\lambda\sim \lambda.
$$
}
\begin{proof}
The proof is similar to that of Lemma 4.5 and is omitted. 
\end{proof}

\begin{proof}[Proof of Theorem 2.4.] 
Since $v_\lambda$ is bounded in $H^1(\mathbb R^N)$,  there exists $v_0\in H^1(\mathbb R^N)$ verifying $-\Delta v+v=(I_\alpha\ast |v|^p)|v|^{p-2}v$ such that  up to a subsequence, we have 
$$
v_\lambda \rightharpoonup v_0   \quad {\rm weakly \ in} \  H^1(\mathbb R^N), \quad v_\lambda\to v_0 \quad {\rm in} \ L^p(\mathbb R^N) \quad {\rm for \ any} \ p\in (2,2^*),
$$
and 
$$
v_\lambda(x)\to v_0(x) \quad a. \ e. \  {\rm on} \   {\R^N},  \qquad v_\lambda\to v_0 \quad {\rm in} \   L^2_{loc}(\mathbb R^N).
$$

Being a ground state solution, $v_\lambda$ satisfies 
$$
m_\lambda=\frac{1}{2}\int_{\mathbb R^N}|\nabla v_\lambda|^2+|v_\lambda|^2-\frac{1}{2p}\int_{\mathbb R^N}(I_\alpha\ast |v_\lambda|^p)|v_\lambda|^p-\frac{1}{q}\lambda\int_{\mathbb R^N}|v_\lambda|^q,
\eqno(7.3)
$$
Since $v_\lambda\in \mathcal M_\lambda$, we also have
$$
\int_{\mathbb R^N}|\nabla v_\lambda|^2+|v_\lambda|^2=\int_{\mathbb R^N}(I_\alpha\ast |v_\lambda|^p)|v_\lambda|^p+\lambda\int_{\mathbb R^N}|v_\lambda|^q.
\eqno(7.4)
$$
Furthermore, by Lemma 3.1, the Poho\v{z}aev identity is given by 
$$
\frac{N-2}{2}\int_{\mathbb R^N}|\nabla v_\lambda|^2+\frac{N}{2}\int_{\mathbb R^N}|v_\lambda|^2=\frac{N+\alpha}{2p}\int_{\mathbb R^N}(I_\alpha\ast |v_\lambda|^p)|v_\lambda|^p+\frac{N}{q}\lambda\int_{\mathbb R^N}|v_\lambda|^q.
\eqno(7.5)
$$
Let 
$$
A_\lambda=\int_{\mathbb R^N}|\nabla v_\lambda|^2,  \  B_\lambda=\int_{\mathbb R^N}|v_\lambda|^2, \  C_\lambda=\int_{\mathbb R^N}(I_\alpha\ast |v_\lambda|^p)|v_\lambda|^p, \  D_\lambda=\int_{\mathbb R^N}|v_\lambda|^q.
$$
Then  by (7.4) and (7.5), we have 
$$
\left\{\begin{array}{rl}
&A_\lambda-C_\lambda=-B_\lambda+\lambda D_\lambda,\\
&\frac{N-2}{2}A_\lambda-\frac{N+\alpha}{2p}C_\lambda=-\frac{N}{2}B_\lambda+\frac{N}{q}\lambda D_\lambda.
\end{array}\right.
$$
Solving this system to obtain 
$$
A_\lambda=\frac{1}{\eta}\left[ (N(p-1)-\alpha)B_\lambda+(N+\alpha-\frac{2Np}{q})\lambda D_\lambda\right],
$$
$$
C_\lambda=\frac{1}{\eta}\left[ 2pB_\lambda+(p(N-2)-\frac{2Np}{q})\lambda D_\lambda\right],
$$
where $\eta=N+\alpha-p(N-2)>0$.
Therefore, we obtain
$$
\begin{array}{rcl}
m_\lambda&=&\frac{1}{2}A_\lambda+\frac{1}{2}B_\lambda-\frac{1}{2p}C_\lambda-\frac{1}{q}\lambda D_\lambda\\
&
=&\frac{p-1}{\eta}B_\lambda+\frac{q(2+\alpha)-2(2p+\alpha)}{2\eta q}\lambda D_\lambda\\
&=&\frac{p-1}{\eta}\int_{\mathbb R^N}| v_\lambda|^2+\frac{q(2+\alpha)-2(2p+\alpha)}{2\eta q}\lambda \int_{\mathbb R^N}|v_\lambda|^q.
\end{array}
$$

Moreover, the  ground state $v_0$ of (7.1) also satisfies the following identities 
$$
m_0=\frac{1}{2}\int_{\mathbb R^N}|\nabla v_0|^2+|v_0|^2-\frac{1}{2p}\int_{\mathbb R^N}(I_\alpha\ast |v_0|^p)|v_0|^p,
\eqno(7.6)
$$
$$
\int_{\mathbb R^N}|\nabla v_0|^2+|v_0|^2=\int_{\mathbb R^N}(I_\alpha\ast |v_0|^p)|v_0|^p,
\eqno(7.7)
$$
$$
\frac{N-2}{2}\int_{\mathbb R^N}|\nabla v_0|^2+\frac{N}{2}\int_{\mathbb R^N}|v_0|^2=\frac{N+\alpha}{2p}\int_{\mathbb R^N}(I_\alpha\ast |v_0|^p)|v_0|^p.
\eqno(7.8)
$$
In a similar way, we show that 
$$
m_0=\frac{p-1}{\eta}\int_{\mathbb R^N}|v_0|^2.
$$
Therefore, 
$$
\frac{p-1}{\eta}\int_{\mathbb R^N}|v_\lambda|^2-|v_0|^2=m_\lambda-m_0-\frac{q(2+\alpha)-2(2p+\alpha)}{2\eta q}\lambda \int_{\mathbb R^N}|v_\lambda|^q,
$$
which together with Lemma 7.1 and $q\ge \frac{2(2p+\alpha)}{2+\alpha}$ implies that
$$
\|v_0\|_2^2-\|v_\lambda\|_2^2\sim \lambda.
$$
Arguing in a similar way, we show that 
$$
\frac{p-1}{Np-(N+\alpha)}\int_{\mathbb R^N}|\nabla v_\lambda|^2-|\nabla v_0|^2=m_\lambda-m_0+\frac{\alpha(q-2)}{2q[Np-(N+\alpha)]}\lambda\int_{\mathbb R^N}|v_\lambda|^q.
$$
Therefore, we obtain 
$$
\|\nabla v_\lambda\|_2^2-\|\nabla v_0\|_2^2=O(\lambda).
$$
We have proved $\|\nabla v_\lambda\|_2\to \|\nabla v_0\|_2$ and $\|v_\lambda\|_2\to \|v_0\|_2$, therefore, we obtain $v_\lambda\to v_0$ in $H^1(\mathbb R^N)$.
\end{proof}

Now, we consider $(Q_\mu)$ and its limit equation
$$
-\Delta v+v=|v|^{p-2}v.
\eqno(7.9)
$$
The corresponding energies of ground states  are given by $m_\mu=\inf_{v\in \mathcal M_\mu}I_\mu(v)$ and 
$$
m_0:=\inf_{v\in \mathcal M_0}I_0(v),
$$
where 
$$
I_0(v)=\frac{1}{2}\int_{\mathbb R^N}|\nabla v|^2+|v|^2-\frac{1}{q}\int_{\mathbb R^N}|v|^q, 
\eqno(7.10)
$$
and 
$$
\mathcal{M}_0=
\left\{v\in H^1(\mathbb R^N)\setminus\{0\} \ \left | \ \int_{\mathbb R^N}|v|^2+|v|^2=\int_{\mathbb R^N}|v|^q \  \right. \right\}. 
$$
Then $m_\mu$ and $m_0$ are well-defined and positive. Moreover, $I_0$ is attained on $\mathcal M_0$ by  the unique positive solution of (7.9).

{\bf Lemma 7.2.}  {\it Assume that the assumptions of Theorems 2.5 holds. Then 
for small $\mu>0$, there holds
$$
m_0-m_\mu\sim \mu.
$$
}
\begin{proof}
The proof is similar to that of Lemma 4.5 and is omitted. 
\end{proof}

\begin{proof}[Proof of Theorem 2.5.] 
Since $v_\mu$ is bounded in $H^1(\mathbb R^N)$,  there exists $v_0\in H^1(\mathbb R^N)$ verifying $-\Delta v+v=v^{q-1}$ such that 
$$
v_\mu \rightharpoonup v_0   \quad {\rm weakly \ in} \  H^1(\mathbb R^N), \quad v_\mu\to v_0 \quad {\rm in} \ L^p(\mathbb R^N) \quad {\rm for \ any} \ p\in (2,2^*),
$$
and 
$$
v_\mu(x)\to v_0(x) \quad a. \ e. \  {\rm on} \  {\R^N},  \qquad v_\mu\to v_0 \quad {\rm in} \   L^2_{loc}(\mathbb R^N).
$$

Being a ground state solution, $v_\mu$ satisfies 
$$
m_\mu=\frac{1}{2}\int_{\mathbb R^N}|\nabla v_\mu|^2+|v_\mu|^2-\frac{1}{2p}\mu\int_{\mathbb R^N}(I_\alpha\ast |v_\mu|^p)|v_\mu|^p-\frac{1}{q}\int_{\mathbb R^N}|v_\mu|^q,
\eqno(7.11)
$$
Since $v_\mu\in \mathcal M_\mu$, we also have
$$
\int_{\mathbb R^N}|\nabla v_\mu|^2+|v_\mu|^2=\mu\int_{\mathbb R^N}(I_\alpha\ast |v_\mu|^p)|v_\mu|^p+\int_{\mathbb R^N}|v_\mu|^q.
\eqno(7.12)
$$
Furthermore, by Lemma 3.1, the Poho\v{z}aev identity is given by 
$$
\frac{N-2}{2}\int_{\mathbb R^N}|\nabla v_\mu|^2+\frac{N}{2}\int_{\mathbb R^N}|v_\mu|^2=\frac{N+\alpha}{2p}\mu\int_{\mathbb R^N}(I_\alpha\ast |v_\mu|^p)|v_\mu|^p+\frac{N}{q}\int_{\mathbb R^N}|v_\mu|^q.
\eqno(7.13)
$$
Therefore, it follows from (7.11) and (7.12) that
$$
m_\mu=\frac{p-1}{2p}\mu\int_{\mathbb R^N}(I_\alpha\ast |v_\mu|^p)|v_\mu|^p+\frac{q-2}{2q}\int_{\mathbb R^N}|v_\mu|^q,
$$
which together with Lemma 7.2 and the Hardy-Littlewood-Sobolev inequality implies that
$$
m_0=\frac{q-2}{2q}\int_{\mathbb R^N}|v_0|^q.
$$
This implies that $v_0$ is a ground state solution of the limit equation (7.9).

Moreover, the unique ground state $v_0$ of (7.9) also satisfies the following identities 
$$
m_0=\frac{1}{2}\int_{\mathbb R^N}|\nabla v_0|^2+|v_0|^2-\frac{1}{q}\int_{\mathbb R^N}|v_0|^q,
\eqno(7.14)
$$
$$
\int_{\mathbb R^N}|\nabla v_0|^2+|v_0|^2=\int_{\mathbb R^N}|v_0|^q,
\eqno(7.15)
$$
$$
\frac{N-2}{2}\int_{\mathbb R^N}|\nabla v_0|^2+\frac{N}{2}\int_{\mathbb R^N}|v_0|^2=\frac{N}{q}\int_{\mathbb R^N}|v_0|^q.
\eqno(7.16)
$$
By(7.14),  (7.15) and (7.16), we have 
$$
\int_{\mathbb R^N}|\nabla v_0|^2=\frac{N(q-2)}{2q}\int_{\mathbb R^N}|\nabla v_0|^2+|v_0|^2=Nm_0=\frac{N(q-2)}{2q}S_q^{\frac{q}{q-2}}.
$$
Therefore, we obtain
$$
\frac{N(q-2)}{2q}\int_{\mathbb R^N}|v_0|^2=(1-\frac{N(q-2)}{2q})\int_{\mathbb R^N}|\nabla v_0|^2=\frac{2N-q(N-2)}{2q}\int_{\mathbb R^N}|\nabla v_0|^2,
$$
which yields
$$
\int_{\mathbb R^N}|v_0|^2=\frac{2N-q(N-2)}{q-2}m_0=\frac{2N-q(N-2)}{2q}S_q^{\frac{q}{q-2}}.
$$

It follows from (7.11), (7.13), (7.14) and (7.16)  that
$$
\begin{array}{rcl}
m_\mu-m_0&=&\frac{1}{2}\int_{\mathbb R^N}|\nabla v_\mu|^2-|\nabla v_0|^2+\frac{1}{2}\int_{\mathbb R^N}|v_\mu|^2-|v_0|^2\\
&\mbox{} &
-\frac{1}{2p}\mu\int_{\mathbb R^N}(I_\alpha\ast |v_\mu|^p)|v_\mu|^p-\frac{1}{q}\int_{\mathbb R^N}|v_\mu|^q-|v_0|^q\\
&=&(\frac{N+\alpha}{2Np}-\frac{1}{2p})\mu\int_{\mathbb R^N}(I_\alpha\ast |v_\mu|^p)|v_\mu|^p+\frac{1}{N}\int_{\mathbb R^N}|\nabla v_\mu|^2-|\nabla v_0|^2.
\end{array}
$$
Therefore, we obtain
$$
\int_{\mathbb R^N}|\nabla v_\mu|^2-|\nabla v_0|^2=-N(m_0-m_\mu)-\frac{\alpha}{2p}\mu\int_{\mathbb R^N}(I_\alpha\ast |v_\mu|^p)|v_\mu|^p.
\eqno(7.17)
$$
Hence, Lemma 7.2 implies that
$$
\int_{\mathbb R^N}|\nabla v_0|^2-|\nabla v_\mu|^2\sim \mu.
$$
By (7.11), (7.12), (7.14) and (7.15), we obtain
$$
\begin{array}{rcl}
m_\mu-m_0&=&\frac{q-2}{2q}\int_{\mathbb R^N}|\nabla v_\mu|^2-|\nabla v_0|^2+\frac{q-2}{2q}\int_{\mathbb R^N}|v_\mu|^2-|v_0|^2\\
&\mbox{} &
+\frac{2p-q}{2pq}\mu\int_{\mathbb R^N}(I_\alpha\ast |v_\mu|^p)|v_\mu|^p,
\end{array}
$$
it follows from (7.17) that
$$
\begin{array}{rcl}
\frac{N(q-2)}{2q}\int_{\mathbb R^N}|v_\mu|^2-|v_0|^2&=&\frac{(N+\alpha)q-2Np}{2pq}\mu\int_{\mathbb R^N}(I_\alpha\ast |v_\mu|^p)|v_\mu|^p\\
&\mbox{}& 
+\frac{2N-q(N-2)}{2q}\int_{\mathbb R^N}|\nabla v_\mu|^2-|\nabla v_0|^2\\
&=&-\frac{2N-q(N-2)}{2q}N(m_0-m_\mu)\\
&\mbox{}&+\left[\frac{(N+\alpha)q-2Np}{2pq}-\frac{2N-q(N-2)}{4pq}\alpha\right]\mu \int_{\mathbb R^N}(I_\alpha\ast |v_\mu|^p)|v_\mu|^p\\
&=&-\frac{2N-q(N-2)}{2q}N(m_0-m_\mu)\\
&\mbox{}&-\frac{N[2(2p+\alpha)- q(2+\alpha)]   }{4pq}\mu \int_{\mathbb R^N}(I_\alpha\ast |v_\mu|^p)|v_\mu|^p,
\end{array}
$$
from which it follows that
$$
\int_{\mathbb R^N}|v_\mu|^2-|v_0|^2=O(\mu).
$$
If $q\le \frac{2(2p+\alpha)}{2+\alpha}$, we also have that $\int_{\mathbb R^N}|v_0|^2-|v_\mu|^2\sim \mu.$

We have proved $\|\nabla v_\mu\|_2\to \|\nabla v_0\|_2$ and $\|v_\mu\|_2\to \|v_0\|_2$, therefore, we obtain $v_\mu\to v_0$ in $H^1(\mathbb R^N)$.
\end{proof}

Finally,  we consider the case $\lambda\to \infty$ in $(Q_\lambda)$ and $\mu\to \infty$ in $(Q_\mu)$.

It is easy to see that under the rescaling 
$$
w(x)=\lambda^{\frac{1}{q-2}}v(x),
$$
the equation $(Q_\lambda)$ is reduced to 
$$
-\Delta w+w=\lambda^{-\frac{2(p-1)}{q-2}}(I_\alpha\ast |w|^p)|w|^{p-2}w+|w|^{q-2}w.
$$
Therefore,  by Theorem 2.5, we have the following 

{\bf Theorem 7.1.} {\it  Let $N\ge 3$, $p\in [\frac{N+\alpha}{N}, \frac{N+\alpha}{N-2}]$, $q\in (2,2^*)$ and  $ v_\lambda$ be  the ground state of $(Q_\lambda)$, then as $\lambda\to \infty$, 
the rescaled family of ground states $\tilde w_\lambda=\lambda^{\frac{1}{q-2}}v_\lambda$ converges in $H^1(\mathbb R^N)$ to the unique positive solution $w_0\in H^1(\mathbb R^N)$ of the equation
$$
-\Delta w+w=w^{q-1}.
$$
Moreover,  as $\lambda\to \infty$, there holds
$$
\|v_\lambda\|_2^2=\left\{\begin{array}{lcl}
\lambda^{-\frac{2}{q-2}}\left(\frac{2N-q(N-2)}{2q}S_q^{\frac{q}{q-2}}+O(\lambda^{-\frac{2(p-1)}{q-2}})\right), \   \  &if & \  q>\frac{2(2p+\alpha)}{q-2},\\
\lambda^{-\frac{2}{q-2}}\left(\frac{2N-q(N-2)}{2q}S_q^{\frac{q}{q-2}}-\Theta(\lambda^{-\frac{2(p-1)}{q-2}})\right), \   \  &if & \  q\le \frac{2(2p+\alpha)}{q-2},
\end{array}
\right.
$$
$$
\|\nabla v_\lambda\|_2^2=\lambda^{-\frac{2}{q-2}}\left(
\frac{N(q-2)}{2q}S_q^{\frac{q}{q-2}}+O(\lambda^{-\frac{2(p-1)}{q-2}})\right),
$$
and the least energy $m_\lambda$ of the ground state satisfies 
$$
\frac{q-2}{2q}S_q^{\frac{q}{q-2}}-\lambda^{\frac{2}{q-2}}m_\lambda \sim\lambda^{-\frac{2(p-1)}{q-2}},
$$
as $\lambda\to \infty$, where $S_q$ is given in (2.19).}

Under the rescaling 
$$
w(x)=\mu^{\frac{1}{2(p-1)}}v(x),
$$
the equation $(Q_\mu)$ is reduced to 
$$
-\Delta w+w=(I_\alpha\ast |w|^p)|w|^{p-2}w+\mu^{-\frac{q-2}{2(p-1)}}|w|^{q-2}w.
$$
Then by Theorem 2.4 we have the following 

{\bf Theorem 7.2.} {\it  Let $N\ge 3$, $p\in (\frac{N+\alpha}{N}, \frac{N+\alpha}{N-2})$, $q\in (2,2^*]$ and $v_\mu$ be  the ground state of $(Q_\mu)$, then as $\mu\to \infty$, 
the rescaled family of ground states $\tilde w_\mu=\mu^{\frac{1}{2(p-1)}}v_\mu$ converges up to a subsequence  in $H^1(\mathbb R^N)$ to a positive solution $w_0\in H^1(\mathbb R^N)$ of the equation
$$
-\Delta v+v=(I_\alpha \ast |v|^p)v^{p-1}.
$$
Moreover, as $\lambda\to \infty$, there holds
$$
\|v_\lambda\|_2^2=\left\{\begin{array}{lcl}
\lambda^{-\frac{1}{p-1}}\left(\frac{N+\alpha-p(N-2)}{2p}S_p^{\frac{p}{p-1}}+O(\lambda^{-\frac{q-2}{2(p-1)}})\right), \   \   
&if &  \   q<\frac{2(2p+\alpha)}{2+\alpha},\\
\lambda^{-\frac{1}{p-1}}\left(\frac{N+\alpha-p(N-2)}{2p}S_p^{\frac{p}{p-1}}-\Theta(\lambda^{-\frac{q-2}{2(p-1)}})\right), \   \   
&if &  \   q\ge \frac{2(2p+\alpha)}{2+\alpha},
\end{array}\right.
$$
$$
\|\nabla v_\lambda\|_2^2=\lambda^{-\frac{1}{p-1}}\left(\frac{N(p-1)-\alpha}{2p}S_p^{\frac{p}{p-1}}+O(\lambda^{-\frac{q-2}{2(p-1)}})\right), 
$$
and  the least energy $m_\mu$ of the ground state satisfies 
$$
\frac{p-1}{2p}S_p^{\frac{p}{p-1}}-\mu^{\frac{1}{p-1}}m_\mu\sim \mu^{-\frac{q-2}{2(p-1)}},
$$
as $\mu\to \infty$, where $S_p$ is given in (2.17).}

Our main results may be applied in other context. As an example, we consider 
$$
-\varepsilon^2\Delta u+u=(I_\alpha\ast |u|^p)|u|^{p-2}u+|u|^{q-2}u,  \   in \ \R^N,
\eqno(Q_\varepsilon)
$$
where $p, q$ and $\alpha$ are the same as before, and $\varepsilon>0$ is a parameter. 
Set 
$v(x)=u(\varepsilon x)$, then we have 
$$
-\Delta v+v=\varepsilon^\alpha(I_\alpha\ast |v|^p)|v|^{p-2}v+|v|^{q-2}v,  \   in \ \R^N,
$$
then as direct consequences of main results in this paper, we have the following 

{\bf Theorem 7.3.} {\it Let $q=2^*$,  $p\in (1+\frac{\alpha}{N-2}, \frac{N+\alpha}{N-2})$ for $N\ge 4$ and  $ p\in (2+\alpha, 3+\alpha)$ for $N=3$,  then the problem  $(Q_\varepsilon)$ admits a positive ground state $u_\varepsilon\in H^1(\mathbb R^N)$,  which is radially symmetric and radially nonincreasing. Furthermore, the following statements hold true:

 If $N\ge 5$, then for small $\varepsilon>0$
 $$
u_\varepsilon(0)\sim \varepsilon^{-\frac{\alpha(N-2)}{2(N-2)(p-1)-2\alpha}},
$$
$$
\int_{\mathbb R^N}(I_\alpha\ast |u_\varepsilon|^p)|u_\varepsilon|^p\sim \|u_\varepsilon\|_2^2\sim \varepsilon^{\frac{(N-2)(Np-N-\alpha)}{(N-2)(p-1)-\alpha}},  
$$
$$
\|\nabla u_\varepsilon\|_2^2=\varepsilon^{N-2}\left[S^{\frac{N}{2}}+O(\varepsilon^{\frac{2\alpha}{(N-2)(p-1)-\alpha}})\right], 
$$
$$
\| u_\varepsilon\|_{2^*}^{2^*}=\varepsilon^N\left[S^{\frac{N}{2}}+O(\varepsilon^{\frac{2\alpha}{(N-2)(p-1)-\alpha}})\right].
$$
Moreover,  there exists $\zeta_\varepsilon\in (0,+\infty)$
verifying  
$$
\zeta_\varepsilon \sim \varepsilon^{\frac{\alpha}{(N-2)(p-1)-\alpha}},
$$
such that for small $\varepsilon>0$, the rescaled family of ground states
$$
w_\varepsilon(x)=\zeta_\varepsilon^{\frac{N-2}{2}}u_\varepsilon(\varepsilon \zeta_\varepsilon x)
$$
satisfies 
$$
\|\nabla w_\varepsilon\|^2_2\sim \|w_\varepsilon\|_{2^*}^{2^*} \sim \int_{\mathbb R^N}(I_\alpha \ast |w_\varepsilon|^p)|w_\varepsilon|^p\sim \|w_\varepsilon\|_2^2\sim 1,
$$
and as $\varepsilon\to 0$, $w_\varepsilon$ converges in $H^1(\mathbb R^N)$ to $W_{\rho_0}$ with 
 $$
 \rho_0=\left(\frac{[(N+\alpha)-p(N-2)]\int_{\mathbb R^N}(I_\alpha\ast |W_1|^p)|W_1|^p}{2p\int_{\mathbb R^N}|W_1|^2}\right)^ {\frac{1}{(N-2)(p-1)-\alpha}}.
 $$
 In the lower dimension cases, we assume that $p\in (\max\{2, 1+\frac{\alpha}{2}\}, 2+\frac{\alpha}{2})$ if $N=4$, and $p\in (2+\alpha, 3+\alpha)$ if $N=3$, 
  then for small $\varepsilon>0$
 $$
u_\varepsilon(0)  \sim\left\{\begin{array}{rcl}
\varepsilon^{-\frac{\alpha}{2p-2-\alpha}}(\ln\frac{1}{\varepsilon})^{\frac{1}{2p-2-\alpha}},   \  \quad if \   \ N=4,\\
\varepsilon^{-\frac{\alpha}{2(p-2-\alpha)}},  \qquad \qquad \qquad if  \    \  N=3,
\end{array}\right.
$$
$$
\int_{\mathbb R^N}(I_\alpha\ast |u_\varepsilon|^p)|u_\varepsilon|^p\sim \|u_\varepsilon\|_2^2\sim \left\{\begin{array}{rcl}
\varepsilon^{\frac{2(4p-4-\alpha)}{2p-2-\alpha}}(\ln\frac{1}{\varepsilon})^{-\frac{4+\alpha-2p}{2p-2-\alpha}},   \  \quad if \   \ N=4,\\
\varepsilon^{\frac{3p-6-2\alpha}{p-2-\alpha}},  \qquad \qquad \qquad if  \    \  N=3,
\end{array}\right.
$$
$$
\|\nabla u_\varepsilon\|_2^2=\varepsilon^{N-2}S^{\frac{N}{2}}+\left\{\begin{array}{rcl} 
O(\varepsilon^{\frac{4p-4}{2p-2-\alpha}}(\ln\frac{1}{\varepsilon})^{-\frac{4+\alpha-2p}{2p-2-\alpha}}),   \quad  if \   \ N=4,\\
 O( \varepsilon^{\frac{p-2}{p-2-\alpha}}), \  \quad \quad  \qquad  \qquad   if  \    \  N=3,
 \end{array}\right. 
$$
$$
\| u_\varepsilon\|_{2^*}^{2^*}=\varepsilon^NS^{\frac{N}{2}}+\left\{\begin{array}{rcl} 
O(\varepsilon^{\frac{2(4p-4-\alpha)}{2p-2-\alpha}}(\ln\frac{1}{\varepsilon})^{-\frac{4+\alpha-2p}{2p-2-\alpha}}),   \quad if \   \ N=4,\\
 O( \varepsilon^{\frac{3p-6-2\alpha}{p-2-\alpha}}), \  \quad \quad  \qquad  \qquad  if  \    \  N=3.
 \end{array}\right. 
 $$
and  there exists $\zeta_\varepsilon\in (0,+\infty)$
verifying  
$$
\zeta_\varepsilon \sim\left\{\begin{array}{rcl}
\varepsilon^{\frac{\alpha}{2p-2-\alpha}}(\ln\frac{1}{\varepsilon})^{-\frac{1}{2p-2-\alpha}},   \  \quad if \   \ N=4,\\
\varepsilon^{\frac{\alpha}{p-2-\alpha}},  \qquad \qquad \qquad if  \    \  N=3,
\end{array}\right.
$$
such that for small $\varepsilon>0$, the rescaled family of ground states
$$
w_\varepsilon(x)=\zeta_\varepsilon^{\frac{N-2}{2}}u_\varepsilon(\varepsilon\zeta_\varepsilon x)
$$
satisfies 
$$
\|\nabla w_\varepsilon\|^2_2\sim \|w_\varepsilon\|_{2^*}^{2^*}\sim \int_{\mathbb R^N}(I_\alpha\ast |w_\varepsilon|^p)|w_\varepsilon|^p\sim 1, 
$$ 
$$
\|w_\varepsilon\|_2^2\sim \left\{\begin{array}{rcl}
\ln\frac{1}{\varepsilon},   \quad  \quad if \   \ N=4,\\
\varepsilon^{-\frac{\alpha}{p-2-\alpha}},  \quad if  \    \  N=3,
\end{array}\right.
$$
and as $\varepsilon\to 0$, $w_\varepsilon$ converges in $D^{1,2}(\mathbb R^N)$ and $L^{\frac{2Np}{N+\alpha}}(\mathbb R^N)$  to $W_1$.
Furthermore, the least energy $m_\varepsilon$ of the ground state satisfies 
$$
\frac{1}{N}S^{\frac{N}{2}}-\varepsilon^{-N}m_\varepsilon \sim \left\{\begin{array}{rcl} \varepsilon^{\frac{2\alpha}{(N-2)(p-1)-\alpha}}, \   \quad  \qquad  \quad  if  \   \ N\ge 5,\\
\varepsilon^{\frac{2\alpha}{2p-2-\alpha}}(\ln\frac{1}{\varepsilon})^{-\frac{4+\alpha-2p}{2p-2-\alpha}},   \quad  if  \   \ N=4,\\
  \varepsilon^{\frac{\alpha}{p-2-\alpha}}, \  \quad \quad  \qquad  \qquad  if \    \  N=3,
 \end{array}\right. 
$$
as $\varepsilon\to 0$, where  $S$ is the best Sobolev constant given in (2.15).}

{\bf Theorem 7.4.} {\it  If $p\in (\frac{N+\alpha}{N}, \frac{N+\alpha}{N-2})$ and $q\in (2,2^*)$.  Then the problem  $(Q_\varepsilon)$ admits a positive ground state $u_\varepsilon\in H^1(\mathbb R^N)$,  which is radially symmetric and radially nonincreasing, and the rescaled family $u_\varepsilon(\varepsilon\cdot)$ converges in $H^1(\mathbb R^N)$ to the unique positive solution $v_0\in H^1(\mathbb R^N)$ of the equation
$$
-\Delta v+v=v^{q-1}.
$$
Moreover, as $\varepsilon\to 0$, there holds
$$
\|u_\varepsilon\|_2^2=\frac{2N-q(N-2)}{2q}S_q^{\frac{q}{q-2}}\varepsilon^N+O(\varepsilon^{N+\alpha}), \quad if \  \ q> \frac{2(2p+\alpha)}{2+\alpha},
$$
$$
\|u_\varepsilon\|_2^2=\frac{2N-q(N-2)}{2q}S_q^{\frac{q}{q-2}}\varepsilon^N-\Theta(\varepsilon^{N+\alpha}), \quad if \  \ q\le \frac{2(2p+\alpha)}{2+\alpha},
$$
$$
\|\nabla u_\varepsilon\|_2^2=\frac{N(q-2)}{2q}S_q^{\frac{q}{q-2}}\varepsilon^{N-2}+O(\varepsilon^{N-2+\alpha}),
$$
and the least energy $m_\varepsilon$ of the ground state satisfies 
$$
\frac{q-2}{2q}S_q^{\frac{q}{q-2}}-\varepsilon^{-N}m_\varepsilon\sim\varepsilon^\alpha,
$$
as $\varepsilon\to 0$, where $S_q$ is given in (2.19).}

\bigskip
{\small
\noindent {\bf Acknowledgements.} 
Part of this research was carried out while S.M. was visiting Swansea University. S.M. thanks the Department of Mathematics at Swansea University for its hospitality.
S.M. was supported by National Natural Science Foundation of China
(Grant Nos.11571187, 11771182)


\vskip 10mm

\begin {thebibliography}{99}
\footnotesize

\bibitem{Akahori-2}
T. Akahori, S. Ibrahim, N. Ikoma, H. Kikuchi and H. Nawa, 
{Uniqueness and nondegeneracy of ground states to nonlinear scalar field equations involving the Sobolev critical exponent in their nonlinearities for high frequencies.} 
Calc. Var. Partial Differential Equations {\bf 58} (2019), Paper No. 120, 32 pp.

\bibitem{Akahori-3} 
T. Akahori, S. Ibrahim, H. Kikuchi and H. Nawa, 
{Global dynamics above the ground state energy for the combined power type nonlinear Schrodinger equations with energy critical growth at low frequencies.}
Mem. Amer. Math. Soc. {\bf 272} (2021), no. 1331, v+130 pp.


\bibitem{Berestycki-1}  
H. Berestycki and P.-L. Lions,  Nonlinear scalar field equations. I. Existence of a ground state, 
{Archive for Rational Mechanics and Analysis} {\bf  82} (1983), 313--345.



\bibitem{Bohmer-Harko}
C. G. B\"ohmer and T. Harko,
{Can dark matter be a Bose--Einstein condensate?}
J. Cosmol. Astropart. Phys., June 2007, pp.025--025.


\bibitem{Brezis-1}   
H. Brezis and E. Lieb, A relation between pointwise convergence of functions and convergence of functionals, 
{Proc. Amer. Math. Soc.}  {\bf 88} (1983),  486--490.


\bibitem{Cassani-1}  D. Cassani, J. Van Schaftingen and Jianjun Zhang, 
Groundstates for Choquard type equations with Hardy-Littlewood-Sobolev lower critical exponent, 
Proc. Roy. Soc. Edinburgh Sect. A {\bf 150} (2020), 1377--1400.
   
 \bibitem{Cazenave-1}  
T. Cazenave and P.-L. Lions, 
Orbital stability of standing waves
for some nonlinear Schr\"odinger equations, {Comm. Math. Phys.}
{\bf 85} (1982), 549--561.   
   

\bibitem{Chavanis-11}
P.-H. Chavanis,
{Mass-radius relation of Newtonian self-gravitating Bose-Einstein condensates with short-range interactions: I. Analytical results},
Phys. Rev. D {\bf 84} (2011), pp.043531.

\bibitem{Chavanis-15}
P.-H. Chavanis,
{Self-gravitating Bose-Einstein condensates},
in {\em Quantum aspects of black holes}, pp. 151--194.
Fundam. Theor. Phys. {\bf 178}, Springer, Cham,  2015.

\bibitem{Coles}
M. Coles and S. Gustafson,
{Solitary Waves and Dynamics for Subcritical Perturbations of Energy Critical NLS.}
Publ. Res. Inst. Math. Sci. {\bf 56} (2020), 647--699.


\bibitem{Dovetta-1}  
S. Dovetta, E. Serra and P. Tilli, 
Action versus energy ground states in nonlinear Schr\"odinger equations,  
arXiv:2107.08655v2.

\bibitem{Duo-1} J. Duoandikoetxea, 
Fractional integrals on radial functions with applications to weighted inequalities,  
Ann. Mat. Pura Appl. (4) {\bf 192} (2013), 553--568.


\bibitem{GT}   D. Gilbarg and N. S. Trudinger,  
{Elliptic Partial Differential Equations of Second Order}. 
Springer-Verlag, Berlin, 1983. xiii+513 pp.

\bibitem{Grillakis-1} M. Grillakis, J. Shatah and W. Strauss, Stability theory of solitary waves in
the presence of symmetry. I, {J. Functional Analysis} {\bf 74} (1987),  160--197.

\bibitem{Grillakis-2} M. Grillakis, J. Shatah and W. Strauss, Stability theory of solitary waves in the presence of symmetry. II, {J. Functional Analysis} {\bf  94} (1990),  308--348.


\bibitem{Ilyasov-1} 
Y. Il'yasov, On orbital stability of the physical ground states of the NLS equations,  arXiv:2103.16353v2.

\bibitem{Jeanjean-2} L. Jeanjean, Existence of solutions with prescribed norm for semilinear elliptic equations, 
{Nonlinear Anal.} {\bf 28} (1997), 1633--1659.

\bibitem{Jeanjean-5}  L. Jeanjean, J. Jendrej, T. Le and N. Visciglia, 
Orbital stability of ground states for a Sobolev critical schr\"odinger equation,
J. Math. Pures Appl. (9) {\bf 164} (2022), 158--179.


\bibitem{Jeanjean-3}   L. Jeanjean and   T. Le, 
Multiple normalized solutions for a Sobolev critical Schr\"odinger equation, 
Math. Ann. {\bf 384} (2022), 101--134.

   
\bibitem{Jeanjean-1} L. Jeanjean and Sheng-Shen Lu, 
On global minimizers for  a mass constrained problem, 
Calc. Var. Partial Differential Equations {\bf 61} (2022), Paper No. 214, 18 pp.


\bibitem{Jeanjean-4}   L. Jeanjean, J. Zhang and X. Zhong, 
A global branch approach to normalized solutions for the Schr\"odinger equation, arXiv:2112.05869v1.


  
\bibitem{Lewin-1} 
M. Lewin and S. R. Nodari, 
{The double-power nonlinear Schr\"odingger equation and its generalizations: uniqueness, non-degeneracy and applications,}
Calc. Var. Partial Differ. Equ. {\bf 59} (6), 1--49 (2020).

\bibitem{Li-2} Xinfu Li, Shiwang Ma and Guang Zhang, Existence and qualitative properties of solutions for
Choquard equations with a local term,  {Nonlinear Analysis: Real World Applications} {\bf  45}
(2019), 1--25.   

\bibitem{Li-1} Xinfu Li and Shiwang Ma, Choquard equations with critical nonlinearities, {Commun. Contemp. Math.}, {\bf 22}
(2019), 1950023.

\bibitem{Li-3} Xinfu Li, Existence and symmetry of normalized ground state to Choquard equation with
local perturbation, arXiv:2103.07026v1.

\bibitem{Li-4} Xinfu Li, Standing waves to upper critical Choquard equation with a local perturbation: multiplicity, qualitative properties and stability, arXiv:2104.09317v1.

\bibitem{Li-5} Xinfu Li, Nonexistence, existence and symmetry of normalized ground states to Choquard equations with a local perturbation, arXiv:2103.07026v2.

\bibitem{Li-6} Xinfu Li, Jianguang Bao and  Wenguang Tang, Normalized solutions to lower critical Choquard equation with a local perturbation, 
arXiv:2207.10377v2.

\bibitem{Lieb-1}   E. H. Lieb,  Existence and uniqueness of the minimizing solution of Choquard’s
nonlinear equation. {Stud. Appl. Math.},   {\bf 57} (1976/77), 93--105.

\bibitem{Lieb-Loss 2001} E.H. Lieb and  M. Loss, Analysis, American Mathematical Society, Providence, RI, 
2001.

\bibitem{Lions-1}   P. H. Lions, The concentration-compactness principle in the calculus of variations: The locally compact cases, Part I and Part II, {Ann. Inst. H. Poincar\'{e} Anal. Non Lin\'{e}aire} {\bf 1} (1984), 223--283.


\bibitem{LiuXQ} X.Q. Liu, J. Q. Liu and Z.-Q. Wang,  Quasilinear elliptic equations with critical growth via perturbation method. {J. Differential Equations}  {\bf 254} (2013), 102--124.

\bibitem{Liu-1} Zeng Liu and V. Moroz, Limit profiles for singularly perturbed Choquard equations with local repulsion, 
Calc. Var. Partial Differential Equations {\bf 61} (2022), Paper No. 160, 59 pp.


\bibitem{Ma-1} S. Ma and V. Moroz, Asymptotic profiles for a nonlinear Schr\"odinger equation with critical combined
powers nonlinearity,  arXiv:2108.01421.


\bibitem{Moroz-1}   V. Moroz and C. B. Muratov,  Asymptotic properties of ground states of scalar field equations with a vanishing
parameter, {J. Eur. Math. Soc.}  {\bf 16} (2014), 1081--1109.


\bibitem {Moroz-2}   V. Moroz  and J. Van Schaftingen, Groundstates of nonlinear Choquard equations: Existence, qualitative properties and decay asymptotics, 
{J. Functional Analysis} {\bf  265}  (2013),  153--184.

\bibitem{MvS-survey}
V. Moroz and J. Van Schaftingen,
{A guide to the Choquard equation},
J. Fixed Point Theory Appl. {\bf 19} (2017), 773--813.

\bibitem{Paredes}
A. Paredes, D. N. Olivieri and H. Michinel,
{From optics to dark matter: A review on nonlinear Schrödinger--Poisson systems},
Physica D: Nonlinear Phenomena {\bf 403} (2020), 132301.

\bibitem{Ruffini}
R. Ruffini and S. Bonazzola,
{Systems of self-gravitating particles in general relativity and the concept of an equation of state},
Physical Review {\bf 187} (1969), 1767--1783.


\bibitem{Shatah-1} J. Shatah and W. Strauss, Instability of nonlinear bound states, {Comm. Math. Phys.} {\bf 100} (1985), 173--190.

\bibitem{Sie-1} D. Siegel and E. Talvila, Pointwise growth estimates of the Riesz potential, 
Dynam. Contin. Discrete Impuls. Systems {\bf 5} (1999), 185--194.

 \bibitem{Soave-1} N. Soave, Normalized ground states for the NLS equation with combined nonlinearities: The Sobolev critical case, {J. Functional Analysis} {\bf 279} (2020), 108610.
 
\bibitem{Soave-2}  N. Soave, Normalized ground state for the NLS equations with combined nonlinearities, {J. Differential Equations} {\bf 269} (2020), 6941--6987.



\bibitem{Tao}
T. Tao, M. Visan and X. Zhang,
{The nonlinear Schr\"odinger equation with combined power-type nonlinearities.} 
Commun. Partial Differ. Equ. {\bf 32} (2007), 1281--1343.





\bibitem{Wang}
X. Z. Wang, {Cold Bose stars: Self-gravitating Bose-Einstein condensates},
Phys. Rev. D {\bf 64} (2001), pp.124009



\bibitem{Wei-1}  J. Wei and Y. Wu, Normalized solutions for Schr\"odinger equations with critical Soblev exponent and mixed nonlinearities, {J. Functional Analysis} {\bf  283} (2022) 109574.

\bibitem{Wei-2}  J. Wei and Y. Wu,     
On some nonlinear Schr\"odinger equations in $\mathbb R^N$,  
Proc. Roy. Soc. Edinburgh Sect. A,
doi:10.1017/prm.2022.56

\bibitem{Weinstein-1} M. I. Weinstein, Modulational stability of ground states of nonlinear Schr\"odinger
equations, {SIAM J. Math. Anal.} {\bf 16} (1985), 472--491.

\bibitem{Sun-1} Shuai Yao, Juntao Sun and Tsung-fang Wu,  Normalized solutions for the Schr\"odinger equation with combined Hartree type and power nonlinearities, 
arXiv:2102.10268v1.

\bibitem{Sun-2}   Shuai Yao, Haibo Chen, 
V. D. R\u adulescu and Juntao Sun,   Normalized solutions for lower critical Choquard equations with critical Sobolev pertubation,   {SIAM J. Math. Anal.}
{\bf 54} (2022),  3696--3723.











\end{thebibliography}
\end{document}